\documentclass[reqno,10pt]{amsart}
\UseRawInputEncoding
\pdfoutput=1
\usepackage{amsfonts}
\usepackage{graphicx}
\usepackage{amsfonts,amsmath, amssymb}
\usepackage{bbm,mathrsfs,mathscinet}
\usepackage[bookmarks=false,hidelinks]{hyperref}
\usepackage{marginnote}
\usepackage{color}
\usepackage{mathrsfs}
\usepackage{amsthm}
\usepackage{enumitem}
\allowdisplaybreaks

\numberwithin{equation}{section}

\newcommand{\R}{\mathbb{R}}

\newtheorem{theorem}{Theorem}[section]
\newtheorem{corollary}[theorem]{Corollary}
\newtheorem{lemma}[theorem]{Lemma}
\newtheorem{proposition}[theorem]{Proposition}
\newtheorem{remark}[theorem]{Remark}
\newtheorem{definition}[theorem]{Definition}

\def\f{\frac}

\begin{document}
	
	\title[the hydrodynamic and newtonian limits of rVMB]{hydrodynamic limit and newtonian limit from the relativistic  Vlasov-Maxwell-Boltzmann system to the classical Euler-Poisson system}
	

	\author[Y. Wang]{Yong Wang}
	\address[Y. Wang]{Academy of Mathematics and Systems Science, Chinese Academy of Sciences, Beijing 100190, China; School of Mathematical Sciences, University of Chinese Academy of Sciences, Beijing 100049, China.}
	\email{yongwang@amss.ac.cn}
	
	\author[H. Xiong]{Hang Xiong}
	\address[H. Xiong]{School of Mathematics and Computational Science, Xiangtan University, Xiangtan 411105, China.} \email{hbear0810@amss.ac.cn}

	\author[H. Y. Zhang]{Hongyao Zhang}
	\address[H.Y. Zhang]{Academy of Mathematics and Systems Science, Chinese Academy of Sciences, Beijing 100190, China; School of Mathematical Sciences, University of Chinese Academy of Sciences, Beijing 100049, China.}
	\email{zhanghongyao@amss.ac.cn}
	
    \begin{abstract}
    In this paper, around a global smooth irrotational solution to the classical isentropic compressible Euler-Poisson system, we construct classical solutions to the one-species relativistic Vlasov-Maxwell-Boltzmann system on any finite time interval $[0,T]$, and rigorously justify the combined hydrodynamic and Newtonian limits to the Euler-Poisson system. In particular, this yields a rigorous derivation of the compressible Euler-Poisson system, whose Poisson coupling induces an instantaneous electrostatic response and thus no longer preserves a strict finite-speed propagation structure, from a relativistic kinetic model with finite propagation speed. The analysis is based on a Hilbert expansion in $\varepsilon$ for the relativistic Vlasov-Maxwell-Boltzmann system, an asymptotic expansion in $\mathfrak{c}^{-1}$ for the relativistic Euler-Maxwell system, and estimates for both the expansion coefficients and the remainder terms that are uniform in $\mathfrak{c}$ and $\varepsilon$, without imposing any \emph{a priori} relation between these two parameters.
    \end{abstract}
	
	\date{\today}
	\maketitle
	
	\setcounter{tocdepth}{1}
	\tableofcontents
	
	\thispagestyle{empty}

	\section{Introduction}
	\subsection{Background and formulation} Plasma, often referred to as the fourth state of matter, consists of moving electrons and ions. It arises in a wide variety of physical settings, including stellar interiors, fusion devices, and space environments such as the solar wind and the Earth's ionosphere. Since the ion mass is much larger than the electron mass, ion dynamics can in some regimes be neglected as a first approximation. In extreme environments with high temperature, high density, or strong electromagnetic fields, plasma particles may attain relativistic velocities, and relativistic effects then become significant. In such a regime, the plasma can be modeled at the kinetic level by the one-species relativistic Vlasov-Maxwell-Boltzmann system (rVMB):
	\begin{equation} \label{e1.1}
		\left\{\begin{aligned}
			&\partial_t F+ \hat{p} \cdot \nabla_x F-e_{-}\Big(E+\frac{p}{p^0} \times B\Big) \cdot \nabla_p F=\frac{1}{\varepsilon} Q_{\mathfrak{c}}\left(F, F\right),\\
			& \partial_t E-\mathfrak{c} \nabla_x \times B=4 \pi e_{-} \int_{\mathbb{R}^3} \hat{p} F d p, \\
			& \partial_t B+\mathfrak{c} \nabla_x \times E=0, \\
			& \operatorname{div} E=4 \pi e_{-} \Big(\bar{n}-\int_{\mathbb{R}^3} F d p\Big), \\
			& \operatorname{div}  B=0,
		\end{aligned}\right.
	\end{equation}
	where $\varepsilon>0$ is the Knudsen number (the mean free path), $F=F(t, x, p)$ is the number density function for electrons at time $t \geq 0$, $E(t, x)$, $B(t, x)$ are the electromagnetic fields, position $x=\left(x_1, x_2, x_3\right) \in \mathbb{R}^3$ and momentum $p=$ $\left(p_1, p_2, p_3\right) \in \mathbb{R}^3$. Here  $\mathfrak{c}$ denotes the speed of light and $p^0$ denotes the energy of an electron with
	$$p^0=\sqrt{m^2 \mathfrak{c}^2+|p|^2}.$$
	Here, $\hat{p}$ denotes the normalized particle velocity 
	$$
    \hat{p}:=\mathfrak{c} \frac{p}{p^0}=\frac{\mathfrak{c} p}{\sqrt{m^2 \mathfrak{c}^2+|p|^2}}.
    $$
	The constants $-e_{-}$and $m$ are the magnitude of the electrons' charge and rest mass, respectively. The constant ion background charge is denoted by $e_{-}\bar{n}>0$. 
	
	The relativistic Boltzmann collision operator $Q_{\mathfrak{c}}(\cdot, \cdot)$ in \eqref{e1.1} takes the form of
	\begin{equation}\label{e1.2}
		Q_{\mathfrak{c}}(F, G)=\frac{\mathfrak{c}}{2p^0} \int_{\mathbb{R}^3} \frac{d q}{q^0} \int_{\mathbb{R}^3} \frac{d q^{\prime}}{q^{\prime 0}} \int_{\mathbb{R}^3} \frac{d p^{\prime}}{p^{\prime 0}} W\left[F\left(p^{\prime}\right) G\left(q^{\prime}\right)-F(p) G(q)\right] .
	\end{equation}
	Here the ``transition rate" $W=W\left(p, q \mid p^{\prime}, q^{\prime}\right)$ is defined as
	\begin{equation}\label{e1.3}
		W=\mathfrak{s} \sigma(g, \theta) \delta\left(p^0+q^0-p^{\prime 0}-q^{\prime 0}\right) \delta^{(3)}\left(p+q-p^{\prime}-q^{\prime}\right) .
	\end{equation}
	
	Denote the four-momenta $p^\mu=\left(p^0, p\right)$ and $q^\mu=\left(q^0, q\right)$. We use the Einstein convention that repeated up-down indices are summed, and we raise and lower indices using the
	Minkowski metric $g_{\mu \nu}=g^{\mu \nu}:=\operatorname{diag}(-1,1,1,1)$. The Lorentz inner product is then given by
	$$
	p^\mu q_\mu := p^\mu g_{\mu \nu} q_\mu= -p^0 q^0+\sum_{i=1}^3 p_i q_i .
	$$
	
	The quantity $\mathfrak{s}$ in \eqref{e1.3} is the square of the energy in the \emph{center of momentum}, $p+q=0$, and is given as
	$$
	\mathfrak{s}=\mathfrak{s}(p, q)=-\left(p^\mu+q^\mu\right)\left(p_\mu+q_\mu\right)=2\left(p^0 q^0-p \cdot q+m^2 \mathfrak{c}^2\right) \geq 4 m^2 \mathfrak{c}^2 .
	$$
	And the relativistic momentum $g$ in \eqref{e1.3} is denoted as
	$$
	g=g(p, q)=\sqrt{-\left(p^\mu-q^\mu\right)\left(p_\mu-q_\mu\right)}=\sqrt{2\left(p^0 q^0-p \cdot q-m^2 \mathfrak{c}^2\right)} \geq 0 .
	$$
	
	The post-collision momentum pair $\left(p^{\prime \mu}, q^{\prime \mu}\right)$ and the pre-collision momentum pair $\left(p^\mu, q^\mu\right)$ satisfy the relation
	\begin{equation}\label{e1.4}
		p^\mu+q^\mu=p^{\prime \mu}+q^{\prime \mu} .
	\end{equation}
	One may also write \eqref{e1.4} as
	$$
	\begin{gathered}
		p^0+q^0=p^{\prime 0}+q^{\prime 0}, \\
		p+q=p^{\prime}+q^{\prime},
	\end{gathered}
	$$
	where the former represents the principle of conservation of energy and the latter represents the conservation of momentum after a binary collision.
	
	\
	
	\textit{Hypothesis on the collision kernel.}
    The function $\sigma(g,\theta)$ in \eqref{e1.3} is called the differential cross section, or scattering kernel, and characterizes binary interactions between relativistic particles. Throughout this paper, we consider the hard-sphere case and normalize
    \[
    \sigma(g,\theta)=1
    \]
    for simplicity.
	
	\subsection{Hilbert expansion} In the present paper, we are concerned with both the hydrodynamic
	limit and Newtonian limit from rVMB to the classical Euler-Poisson system. To achieve this, we perform a Hilbert expansion for rVMB $\eqref{e1.1}$ with small Knudsen number $\varepsilon$. To emphasize the dependence on $\varepsilon$ and $\mathfrak{c}$ for the solutions of $\eqref{e1.1}$, we denote them as $F^{\varepsilon,\mathfrak{c}}$, $E^{\varepsilon,\mathfrak{c}}$ and $B^{\varepsilon,\mathfrak{c}}$. Define the Hilbert expansion by
	\begin{equation}\label{e1.5}
		\left\{\begin{aligned}
			&F^{\varepsilon,\mathfrak{c}}  =\sum_{n=0}^{2 k-1} \varepsilon^n F_n^{\mathfrak{c}}+\varepsilon^k F_R^{\varepsilon,\mathfrak{c}}, \\
			&E^{\varepsilon,\mathfrak{c}}  =\sum_{n=0}^{2 k-1} \varepsilon^n E_n^{\mathfrak{c}}+\varepsilon^k E_R^{\varepsilon,\mathfrak{c}},\quad(k \geq 5), \\
			&B^{\varepsilon,\mathfrak{c}}  =\sum_{n=0}^{2 k-1} \varepsilon^n B_n^{\mathfrak{c}}+\varepsilon^k B_R^{\varepsilon,\mathfrak{c}},
		\end{aligned}\right.
	\end{equation}
	where $F_0^{\mathfrak{c}}, F_1^{\mathfrak{c}}, \ldots, F_{2 k-1}^{\mathfrak{c}}$ in $\eqref{e1.5}$ will depend upon $\mathfrak{c}$ but be independent of $\varepsilon$. Also, $F_R^{\varepsilon, \mathfrak{c}}$, $E_R^{\varepsilon,\mathfrak{c}}$, and $B_R^{\varepsilon,\mathfrak{c}}$ are called the remainder terms which will depend upon $\varepsilon$ and $\mathfrak{c}$. For $\mathfrak{c}=1$, Guo-Xiao \cite{Guo-CMP-2021} established the Hilbert expansion for rVMB. Since we shall consider both the hydrodynamic limit $\varepsilon \to 0$ and the Newtonian limit $\mathfrak{c} \to \infty$ of rVMB, it is crucial to derive estimates for the expansion coefficients $\big(F_n^{\mathfrak{c}},E_n^{\mathfrak{c}},B_n^{\mathfrak{c}}\big)$ $(n=0,1,\cdots,2k-1)$ that are uniform in $\mathfrak{c}$, together with uniform in $\mathfrak{c}$ and $\varepsilon$ estimates for the remainder terms $F_R^{\varepsilon,\mathfrak{c}}$, $E_R^{\varepsilon,\mathfrak{c}}$, and $B_R^{\varepsilon,\mathfrak{c}}$.
	
	To determine the coefficients $F_i^{\mathfrak{c}}(t, x, p), E_i^{\mathfrak{c}}(t, x), B_i^{\mathfrak{c}}(t, x), i=0,1, \ldots, 2k-1$, we plug the formal expansion \eqref{e1.5} into the rescaled equations \eqref{e1.1}.
	Comparing the coefficients to each order of $\varepsilon$, one obtains the following equations:
	\begin{align} 
		\frac{1}{\varepsilon}: &Q_{\mathfrak{c}}\left(F_0^{\mathfrak{c}}, F_0^{\mathfrak{c}}\right)=0, \label{e1.7} \\
		\varepsilon^0: &\partial_t F_0^{\mathfrak{c}}+ \hat{p} \cdot \nabla_x F_0^{\mathfrak{c}}-e_{-}\big(E_0^{\mathfrak{c}}+\frac{p}{p^0} \times B_0^{\mathfrak{c}}\big) \cdot \nabla_p F_0^{\mathfrak{c}}=Q_{\mathfrak{c}}\left(F_1^{\mathfrak{c}}, F_0^{\mathfrak{c}}\right)+Q_{\mathfrak{c}}\left(F_0^{\mathfrak{c}}, F_1^{\mathfrak{c}}\right), \nonumber \\
		& \partial_t E_0^{\mathfrak{c}}-\mathfrak{c} \nabla_x \times B_0^{\mathfrak{c}}=4 \pi e_{-} \int_{\mathbb{R}^3} \hat{p} F_0^{\mathfrak{c}} d p, \nonumber \\
		& \partial_t B_0^{\mathfrak{c}}+\mathfrak{c} \nabla_x \times E_0^{\mathfrak{c}}=0, \label{e1.8} \\
		& \operatorname{div} E_0^{\mathfrak{c}}=4 \pi e_{-}\Big(\bar{n}-\int_{\mathbb{R}^3} F_0^{\mathfrak{c}} d p\Big), \nonumber \\
		& \operatorname{div} B_0^{\mathfrak{c}}=0, \nonumber \\
		& \cdots \cdots  \nonumber\\
		\varepsilon^n: &\partial_t F_n^{\mathfrak{c}}+ \hat{p} \cdot \nabla_x F_n^{\mathfrak{c}}-e_{-}\Big(E_n^{\mathfrak{c}}+\frac{p}{p^0} \times B_n^{\mathfrak{c}}\Big) \cdot \nabla_p F_0^{\mathfrak{c}}-e_{-}\big(E_0^{\mathfrak{c}}+\frac{p}{p^0} \times B_0^{\mathfrak{c}}\big) \cdot \nabla_p F_n^{\mathfrak{c}} \nonumber \\
		&=\sum_{\substack{i+j=n+1 \\ i, j \geq 0}} Q_{\mathfrak{c}}\left(F_i^{\mathfrak{c}}, F_j^{\mathfrak{c}}\right)+e_-\sum_{\substack{i+j=n \\ i, j \geq 1}}\Big(E_i^{\mathfrak{c}}+\frac{p}{p^0} \times B_i^{\mathfrak{c}}\Big) \cdot \nabla_p F_j^{\mathfrak{c}}, \nonumber \\
		&\partial_t E_n^{\mathfrak{c}}-\mathfrak{c} \nabla_x \times B_n^{\mathfrak{c}}=4 \pi e_{-} \int_{\mathbb{R}^3} \hat{p} F_n^{\mathfrak{c}} d p, \label{e1.9} \\
		& \partial_t B_n^{\mathfrak{c}}+\mathfrak{c} \nabla_x \times E_n^{\mathfrak{c}}=0, \nonumber \\
		& \operatorname{div} E_n^{\mathfrak{c}}=-4 \pi e_{-} \int_{\mathbb{R}^3} F_n^{\mathfrak{c}} d p, \nonumber \\
		& \operatorname{div} B_n^{\mathfrak{c}}=0. \nonumber 
	\end{align}
	The remainder terms $F_R^{\varepsilon,\mathfrak{c}}, E_R^{\varepsilon,\mathfrak{c}}$ and $B_R^{\varepsilon,\mathfrak{c}}$ satisfy the following equations:
		\begin{align}\label{e1.10}
			 & \partial_t F_R^{\varepsilon,\mathfrak{c}}+ \hat{p} \cdot \nabla_x F_R^{\varepsilon,\mathfrak{c}}-e_{-}\Big(E_R^{\varepsilon,\mathfrak{c}}+\frac{p}{p^0} \times B_R^{\varepsilon,\mathfrak{c}}\Big) \cdot \nabla_p F_0^{\mathfrak{c}}-e_{-}\big(E_0^{\mathfrak{c}}+\frac{p}{p^0} \times B_0^{\mathfrak{c}}\big) \cdot \nabla_p F_R^{\varepsilon,\mathfrak{c}} \nonumber\\
		&\quad-\frac{1}{\varepsilon}\left[Q_{\mathfrak{c}}\left(F_R^{\varepsilon,\mathfrak{c}}, F_0^{\mathfrak{c}}\right)+Q_{\mathfrak{c}}\left(F_0^{\mathfrak{c}}, F_R^{\varepsilon,\mathfrak{c}}\right)\right]-\varepsilon^{k-1} Q_{\mathfrak{c}}\left(F_R^{\varepsilon,\mathfrak{c}}, F_R^{\varepsilon,\mathfrak{c}}\right) \nonumber\\
		&=\sum_{i=1}^{2 k-1} \varepsilon^{i-1}\left[Q_{\mathfrak{c}}\left(F_i^{\mathfrak{c}}, F_R^{\varepsilon,\mathfrak{c}}\right)+Q_{\mathfrak{c}}\left(F_R^{\varepsilon,\mathfrak{c}}, F_i^{\mathfrak{c}}\right)\right]+\varepsilon^k e_{-}\Big(E_R^{\varepsilon,\mathfrak{c}}+\frac{p}{p^0} \times B_R^{\varepsilon,\mathfrak{c}}\Big) \cdot \nabla_p F_R^{\varepsilon,\mathfrak{c}} \nonumber\\
		&\quad +\sum_{i=1}^{2 k-1} \varepsilon^i e_{-}\left[\Big(E_i^{\mathfrak{c}}+\frac{p}{p^0} \times B_i^{\mathfrak{c}}\Big) \cdot \nabla_p F_R^{\varepsilon,\mathfrak{c}}+\Big(E_R^{\varepsilon,\mathfrak{c}}+\frac{p}{p^0} \times B_R^{\varepsilon,\mathfrak{c}}\Big) \cdot \nabla_p F_i^{\mathfrak{c}}\right]+\varepsilon^k A,\nonumber\\
			& \partial_t E_R^{\varepsilon,\mathfrak{c}}-\mathfrak{c} \nabla_x \times B_R^{\varepsilon,\mathfrak{c}}=4 \pi e_{-} \int_{\mathbb{R}^3} \hat{p} F_R^{\varepsilon,\mathfrak{c}} d p, \\
			&\partial_t B_R^{\varepsilon,\mathfrak{c}}+\mathfrak{c} \nabla_x \times E_R^{\varepsilon,\mathfrak{c}}=0, \nonumber\\
			&\operatorname{div} E_R^{\varepsilon,\mathfrak{c}}=-4 \pi e_{-} \int_{\mathbb{R}^3} F_R^{\varepsilon,\mathfrak{c}} d p, \nonumber\\
			&\operatorname{div} B_R^{\varepsilon,\mathfrak{c}}=0,\nonumber
		\end{align}
	where
	$$
	A=\sum_{\substack{i+j \geq 2 k+1 \\ 2 \leq i, j \leq 2 k-1}} \varepsilon^{i+j-2 k-1} Q_{\mathfrak{c}}\left(F_i^{\mathfrak{c}}, F_j^{\mathfrak{c}}\right)+\sum_{\substack{i+j \geq 2 k \\ 1 \leq i, j \leq 2 k-1}} \varepsilon^{i+j-2 k} e_{-}\left(E_i^{\mathfrak{c}}+\frac{p}{p^0} \times B_i^{\mathfrak{c}}\right) \cdot \nabla_p F_j^{\mathfrak{c}} .
	$$
	It follows from $\eqref{e1.7}$ that $F_0^{\mathfrak{c}}$ should be a relativistic local Maxwellian, i.e.,
	\begin{equation}\label{e1.11}
		\mathbf{M}_{\mathfrak{c}}(t, x, p):=F_0^{\mathfrak{c}}(t, x, p)\equiv \frac{n_{\mathfrak{e}} \gamma}{4 \pi m^3 \mathfrak{c}^3 K_2(\gamma)} \exp \left\{\frac{u_{\mathfrak{e}}^\mu p_\mu}{k_B T_\mathfrak{e}}\right\},
	\end{equation}
	where $k_B$ is Boltzmann's constant and $\gamma$ is a dimensionless variable defined as
	$$
	\gamma=\frac{m \mathfrak{c}^2}{k_B T_\mathfrak{e}}.
	$$
	 Here $K_j(\gamma)(j=0,1,2, \ldots)$ are the modified Bessel functions of the second kind defined in \eqref{2.1}. Here $n_{\mathfrak{e}}(t, x)>0, u_{\mathfrak{e}}(t, x)$ and $T_\mathfrak{e}(t, x)>0$ are the number density, velocity and temperature, which are part of the solution to the 1-fluid relativistic Euler-Maxwell system \eqref{e0.1} constructed in Theorem \ref{c3.4}.
	
	Throughout the paper, since $e_-$, $m$, and $k_B$ are fixed physical constants, we normalize them by setting
	$$
	e_-=m=k_B=1.
	$$
	\subsection{The relativistic Euler-Maxwell system and the classical Euler-Poisson system}
	Similar to \cite{Cercignani-2002,Guo-CMP-2021}, for $\alpha, \beta \in\{0,1,2,3\}$, we define the first momentum as
	$$
	I^\alpha\left[\mathbf{M}_{\mathfrak{c}}\right]:=\int_{\mathbb{R}^3} \frac{p^\alpha}{p^0} \mathbf{M}_{\mathfrak{c}} d p
	$$
	and the second momentum as
	$$
	T^{\alpha \beta}\left[\mathbf{M}_{\mathfrak{c}}\right]:=\int_{\mathbb{R}^3} \frac{p^\alpha p^\beta}{p^0} \mathbf{M}_{\mathfrak{c}} d p .
	$$
	It has been shown in \cite[Proposition 3.3]{Speck-CMP-2011} that
	\begin{equation}\label{e1.12}
		\begin{aligned}
			I^\alpha\left[\mathbf{M}_{\mathfrak{c}}\right] & =\frac{n_{\mathfrak{e}} u_{\mathfrak{e}}^\alpha}{\mathfrak{c}}, \\
			T^{\alpha \beta}\left[\mathbf{M}_{\mathfrak{c}}\right] & =\frac{e_{\mathfrak{e}}+P_{\mathfrak{e}}}{\mathfrak{c}^3} u_{\mathfrak{e}}^\alpha u_{\mathfrak{e}}^\beta+\frac{P_{\mathfrak{e}} g^{\alpha \beta}}{\mathfrak{c}},
		\end{aligned}
	\end{equation}
	where $e_{\mathfrak{e}}(t, x)>0$ is the proper energy density and $P_{\mathfrak{e}}(t, x)>0$ is the pressure.
	
	Projecting the first equation in \eqref{e1.8} onto $1$, $p$, $p^0$, which are five collision invariants for the relativistic Boltzmann collision operator $Q_{\mathfrak{c}}(\cdot,\cdot)$, and using \eqref{e1.12}, we can obtain the following relativistic Euler-Maxwell system (rEM):
	\begin{equation}\label{e0.1}
		\left\{\begin{aligned}
			&\frac{1}{\mathfrak{c}} \partial_t\left(n_{\mathfrak{e}} u_{\mathfrak{e}}^0\right)+\operatorname{div}\left(n_{\mathfrak{e}} u_{\mathfrak{e}}\right)=0, \\
			&\frac{1}{\mathfrak{c}} \partial_t\left[\left(e_{\mathfrak{e}}+P_{\mathfrak{e}}\right) u_{\mathfrak{e}}^0 u_{\mathfrak{e}}\right]+\operatorname{div}\left[\left(e_{\mathfrak{e}}+P_{\mathfrak{e}}\right) u_{\mathfrak{e}} \otimes u_{\mathfrak{e}}\right]+\mathfrak{c}^2 \nabla_x P_{\mathfrak{e}}+\mathfrak{c} n_{\mathfrak{e}}\left[u_{\mathfrak{e}}^0 E_0^{\mathfrak{c}}+u_{\mathfrak{e}} \times B_0^{\mathfrak{c}}\right]=0, \\
			&\frac{1}{\mathfrak{c}} \partial_t\left[\left(e_{\mathfrak{e}}+P_{\mathfrak{e}}\right)\left(u_{\mathfrak{e}}^0\right)^2-\mathfrak{c}^2 P_{\mathfrak{e}}\right]+\operatorname{div}\left[\left(e_{\mathfrak{e}}+P_{\mathfrak{e}}\right) u_{\mathfrak{e}}^0 u_{\mathfrak{e}}\right]+\mathfrak{c} n_{\mathfrak{e}} u_{\mathfrak{e}} \cdot E_0^{\mathfrak{c}}=0,\\
			&\partial_t E_0^{\mathfrak{c}}-\mathfrak{c} \nabla_x \times B_0^{\mathfrak{c}}=4 \pi n_{\mathfrak{e}}u_{\mathfrak{e}}, \quad \operatorname{div} E_0^{\mathfrak{c}}=4 \pi \Big(\bar{n}-\frac{1}{\mathfrak{c}} n_{\mathfrak{e}} u_{\mathfrak{e}}^0\Big), \\
			&\partial_t B_0^{\mathfrak{c}}+\mathfrak{c} \nabla_x \times E_0^{\mathfrak{c}}=0, \quad \operatorname{div} B_0^{\mathfrak{c}}=0 .
		\end{aligned}\right.
	\end{equation}
	The fluid quantities $n_{\mathfrak{e}}, T_\mathfrak{e}, S_\mathfrak{e}, P_{\mathfrak{e}}, e_{\mathfrak{e}}$ in \eqref{e0.1} satisfy the following relations
	\begin{align}
		& P_{\mathfrak{e}}=n_{\mathfrak{e}} T_\mathfrak{e}=\mathfrak{c}^2 \frac{n_{\mathfrak{e}}}{\gamma}, \label{e0.2}\\
		& e_{\mathfrak{e}}=\mathfrak{c}^2 n_{\mathfrak{e}} \frac{K_1(\gamma)}{K_2(\gamma)}+3 P_{\mathfrak{e}}=\mathfrak{c}^2 n_{\mathfrak{e}} \frac{K_3(\gamma)}{K_2(\gamma)}-P_{\mathfrak{e}}, \label{e0.3}\\
		& n_{\mathfrak{e}}=4 \pi e^4 \mathfrak{c}^3 \exp \left(-S_\mathfrak{e}\right) \frac{K_2(\gamma)}{\gamma} \exp \Big(\gamma \frac{K_1(\gamma)}{K_2(\gamma)}\Big),\label{e0.4}
	\end{align}
	where $S_\mathfrak{e}(t, x)>0$ is the entropy and $K_i(\gamma)$ is the modified Bessel function of the second kind defined in \eqref{2.1}. 
    
    Throughout this paper, we adopt the following isentropic assumption:
    \begin{equation}\label{e.1}
		e_{\mathfrak{e}}+P_{\mathfrak{e}}=n_{\mathfrak{e}} h\left(n_{\mathfrak{e}}\right)>0 \quad \text { and } \quad P_{\mathfrak{e}}^{\prime}\left(n_{\mathfrak{e}}\right)=n_{\mathfrak{e}} h^{\prime}\left(n_{\mathfrak{e}}\right),
	\end{equation}
	which yields that $S_\mathfrak{e}$ is a constant according to the following Gibbs relation (see \cite{Yvonne-OMM-2009})
    \begin{align}\label{e.18}
    T_\mathfrak{e} d S_\mathfrak{e}=d\Big(\frac{e_{\mathfrak{e}}}{n_{\mathfrak{e}}}\Big)+P_{\mathfrak{e}} d\Big(\frac{1}{n_{\mathfrak{e}}}\Big).
    \end{align}
	Consequently, under the isentropic assumption \eqref{e.1}, $T_\mathfrak{e}$  is determined by $n_{\mathfrak{e}}$ through \eqref{e0.4}. In fact, noting $\frac{K_1(\gamma)}{K_2(\gamma)}=\frac{K_3(\gamma)}{K_2(\gamma)}-\frac{4}{\gamma}$, we have from \eqref{e0.4} and \eqref{2.2}-\eqref{2.4} that
	\begin{align}\label{e.7}
	 n_{\mathfrak{e}}&=4 \pi \mathfrak{c}^3 \exp \left(-S_\mathfrak{e}\right) \frac{K_2(\gamma)}{\gamma} \exp \Big(\gamma \frac{K_3(\gamma)}{K_2(\gamma)}\Big)\nonumber\\
	 & =4 \pi \mathfrak{c}^3 e^{-S_\mathfrak{e}}\frac{1}{\gamma} \sqrt{\frac{\pi}{2 \gamma}} e^{-\gamma}\Big(1+\frac{15}{8 \gamma}+O\left(\gamma^{-2}\right)\Big) \cdot e^{\gamma\big(1+\frac{5}{2 \gamma}+\frac{15}{8 \gamma^2}+O\left(\gamma^{-3}\right)\big)} \nonumber\\
	 & =e^{\frac{5}{2}-S_\mathfrak{e}}(2 \pi)^{\frac{3}{2}} \mathfrak{c}^3 \gamma^{-\frac{3}{2}}\Big(1+\frac{15}{4} \frac{1}{\gamma}+O\left(\gamma^{-2}\right)\Big) \\
    & =e^{\frac{5}{2}-S_\mathfrak{e}}(2 \pi)^{\frac{3}{2}}\Big(T_\mathfrak{e}^{\frac{3}{2}}+\frac{15}{4} \frac{T_\mathfrak{e}^{\frac{5}{2}}}{\mathfrak{c}^2}+O\Big(\frac{T_\mathfrak{e}^{\frac{7}{2}}}{\mathfrak{c}^4}\Big)\Big). \nonumber
	\end{align}
	
	Using $\eqref{e0.1}_1$ and \eqref{e.1}, we have from $\eqref{e0.1}_{2,3}$ that
	\begin{equation}\label{e.2}
		\begin{aligned}
			& \frac{n_{\mathfrak{e}} u_{\mathfrak{e}}^0}{\mathfrak{c}} \partial_t\left(h u_{\mathfrak{e}}\right)+n_{\mathfrak{e}} u_{\mathfrak{e}} \cdot \nabla_x\left(h u_{\mathfrak{e}}\right)+\mathfrak{c}^2 \nabla_x P_{\mathfrak{e}}+\mathfrak{c}  n_{\mathfrak{e}}\left[u_{\mathfrak{e}}^0 E_0^{\mathfrak{c}}+u_{\mathfrak{e}} \times B_0^{\mathfrak{c}}\right]=0, \\
			& \frac{n_{\mathfrak{e}} u_{\mathfrak{e}}^0}{\mathfrak{c}} \partial_t\left(h u_{\mathfrak{e}}^0\right)-\mathfrak{c} \partial_t P_{\mathfrak{e}}+n_{\mathfrak{e}} u_{\mathfrak{e}} \cdot \nabla_x \left(h u_{\mathfrak{e}}^0\right)+\mathfrak{c} n_{\mathfrak{e}} u_{\mathfrak{e}} \cdot E_0^{\mathfrak{c}}=0 .
		\end{aligned}
	\end{equation}
	We further employ $u_{\mathfrak{e}}^0\eqref{e.2}_2-u_{\mathfrak{e}} \cdot\eqref{e.2}_1$ to obtain
	\begin{equation}\label{e.3}
		\mathfrak{c}u_{\mathfrak{e}}^0\big[n_{\mathfrak{e}} \partial_t h-\partial_t P_{\mathfrak{e}}\big]+\mathfrak{c}^2u_{\mathfrak{e}} \cdot\big[n_{\mathfrak{e}} \nabla_x h-\nabla_x P_{\mathfrak{e}}\big]=0,
	\end{equation}
	which, together with \eqref{e.1}, yields that $\eqref{e0.1}_3$ can be derived from $\eqref{e0.1}_{1,2}$ and \eqref{e.1}. 
	 For consistency of presentation, denote
	 $$
	 E_{\mathfrak{e}}:=E_0^{\mathfrak{c}}, \quad B_{\mathfrak{e}}:=B_0^{\mathfrak{c}}.
	 $$
	 Thus, under the assumption \eqref{e.1}, we can rewrite \eqref{e0.1} as:
	\begin{equation}\label{e.4}
		\left\{\begin{aligned}
			&\frac{1}{\mathfrak{c}} \partial_t\left(n_{\mathfrak{e}} u_{\mathfrak{e}}^0\right)+\operatorname{div}\left(n_{\mathfrak{e}} u_{\mathfrak{e}}\right)=0, \\
			&\frac{u_{\mathfrak{e}}^0}{\mathfrak{c}} \partial_t(\frac{h u_{\mathfrak{e}}}{\mathfrak{c}^2})+(u_{\mathfrak{e}}\cdot \nabla_x)(\frac{h u_{\mathfrak{e}}}{\mathfrak{c}^2}) + \nabla_x h(n_{\mathfrak{e}})+  \Big[\frac{u_{\mathfrak{e}}^0}{\mathfrak{c}} E_{\mathfrak{e}}+\frac{u_{\mathfrak{e}}}{\mathfrak{c}} \times B_{\mathfrak{e}}\Big]=0, \\
			&\partial_t E_{\mathfrak{e}}-\mathfrak{c} \nabla_x \times B_{\mathfrak{e}}=4\pi n_{\mathfrak{e}}u_{\mathfrak{e}}, \quad \operatorname{div} E_{\mathfrak{e}}= 4\pi\Big(\bar{n}-\frac{1}{\mathfrak{c}} n_{\mathfrak{e}} u_{\mathfrak{e}}^0\Big),\\
			&\partial_t B_{\mathfrak{e}}+\mathfrak{c} \nabla_x \times E_{\mathfrak{e}}=0, \quad \operatorname{div} B_{\mathfrak{e}}=0.
		\end{aligned}\right.
	\end{equation}
	We impose \eqref{e.4} with the initial condition 
	\begin{align}\label{e.4-1}
	    \left.(n_{\mathfrak{e}}, u_{\mathfrak{e}}, E_{\mathfrak{e}}, B_{\mathfrak{e}})\right|_{t=0}= (n_{\mathfrak{e}}, u_{\mathfrak{e}}, E_{\mathfrak{e}}, B_{\mathfrak{e}})(0,x)
	\end{align}
	satisfying the following compatibility conditions 
    \begin{align}\label{e.19-1}
        \operatorname{div} E_{\mathfrak{e}}(0)= 4\pi\Big(\bar{n}-\frac{1}{\mathfrak{c}} n_{\mathfrak{e}}(0) u_{\mathfrak{e}}^0(0)\Big), \quad
	    \operatorname{div} B_{\mathfrak{e}}(0)=0.
    \end{align}
    
 Formally, as $\mathfrak{c}\to\infty$, rEM \eqref{e.4} reduces to the classical isentropic compressible Euler-Poisson system (EP)
\begin{equation}\label{e0.5}
\left\{
\begin{aligned}
&\partial_t\rho+\operatorname{div}(\rho \mathfrak{u})=0, \\
&\partial_t(\rho \mathfrak{u})+\operatorname{div}(\rho \mathfrak{u}\otimes \mathfrak{u})+\nabla_x \mathcal{P}
=-\rho\nabla\phi, \\
&\Delta\phi=4\pi(\bar{\rho}-\rho),
\end{aligned}
\right.
\end{equation}
subject to the far-field condition
\[
|\phi(t,x)|\to 0 \qquad \text{as } |x|\to\infty.
\]
Here, $\rho(t,x)$ and $\mathfrak{u}(t,x)$ denote the electron density and velocity, respectively, $\bar{\rho}>0$ is the constant ion background density, $\phi(t,x)$ is the electric potential, and $\mathcal{P}(t,x)$ is the pressure. Under the isentropic assumption, the entropy is constant, and the pressure is given by the polytropic gas law
\begin{equation}\label{e0.6}
\mathcal{P}=\rho \theta, \quad  \theta:=K_{\bar{\eta}}\,\rho^{\frac{2}{3}}, 
\quad
K_{\bar{\eta}}=(2\pi)^{-1}e^{-{\frac{5}{3}}+\frac{2}{3}\bar{\eta}},
\end{equation}
where $\theta$ is the temperature and $\bar{\eta}$ is the constant entropy.

    For the classical isentropic compressible Euler-Poisson system \eqref{e0.5}, Guo~\cite{Guo-CMP-1998} constructed global smooth irrotational solutions in $\mathbb{R}^3$ near the constant equilibrium state $(\bar{\rho},0,0)$. In this setting, the neutrality condition is used to rule out non-decaying electric modes and to obtain the desired decay of the electric field, while the irrotationality assumption is essential for exploiting the dispersive structure generated by the Poisson coupling. Since then, further results on the global existence and asymptotic stability of smooth solutions for classical Euler-Poisson and Euler-Maxwell systems have also been obtained; see, for example, \cite{Guo-CMP-2011,Guo-AOM-2016}.
   
\subsection{A brief review of the Newtonian and hydrodynamic limits for rVMB}
For the relativistic Boltzmann equation, the first rigorous justification of the Newtonian limit in the periodic domain $\mathbb{T}^3$ was obtained by Calogero \cite{Calogero-JMP-2004}, who proved a local-in-time convergence result. Subsequently, Strain \cite{Strain-JMA-2010} considered the case of unique global-in-time mild solutions near vacuum in $\mathbb{R}^3$ and derived explicit convergence rates on any finite time interval $[0,T]$. More recently, Cao et al.~\cite{Cao-Ouyang-Wang-Xiao-JFA-2026} established global well-posedness and the global Newtonian limit in a periodic box with estimates uniform in $\mathfrak{c}$. These works provide a rather complete picture for the field-free relativistic Boltzmann equation.

For field-coupled models, however, the situation is substantially more delicate. In the classical Vlasov-Maxwell-Boltzmann setting, Jiang-Lei-Zhao \cite{Jiang-Lei-Zhao-JFA-2024} justified the limit toward the Vlasov-Poisson-Boltzmann system, highlighting the additional difficulty caused by the degeneration of the Maxwell system when $\mathfrak{c}\to\infty$. This indicates that, for electromagnetic kinetic models, the Newtonian limit must account not only for the particle dynamics but also for the singular transition of the field equations.

The hydrodynamic limit from Boltzmann-type equations to Euler or Navier-Stokes type systems is a classical topic in kinetic theory. Early rigorous works of Nishida \cite{Nishida-CMP-1978} and Caflisch \cite{Caflisch-CPAM-1980} laid the foundation for the derivation of compressible fluid equations from the nonlinear Boltzmann equation. For non-relativistic plasma models, Guo-Jang \cite{Guo-CMP-2010} established the global Hilbert expansion for the classical Vlasov-Poisson-Boltzmann system and derived the Euler-Poisson system.

In the relativistic setting, Speck-Strain \cite{Speck-CMP-2011} gave the first rigorous hydrodynamic limit from the relativistic Boltzmann equation to relativistic Euler equations in a local-in-time framework. Later Guo-Xiao \cite{Guo-CMP-2021} established a global Hilbert expansion for rVMB with a fixed ion background and derived rEM, working in normalized units with the speed of light $\mathfrak{c}$ set to one. Their work shows that the relativistic hydrodynamic theory can be extended from field-free models to genuinely electromagnetic coupled systems.

Compared with the separate studies of $\mathfrak{c}\to\infty$ and $\varepsilon\to0$, the joint limit is considerably more subtle, since one must derive estimates that remain uniform with respect to both singular parameters. A major step in this direction was made by Wang-Xiao \cite{Wang-Xiao-JLMS-2026}, who established, for the relativistic Boltzmann equation, both the hydrodynamic limit and the Newtonian limit without imposing an \emph{a priori} relation between $\varepsilon$ and $\mathfrak{c}$. Their work provides a unified two-parameter framework and offers a basic paradigm for the study of concurrent limits in relativistic kinetic theory.

In the present paper, we study the joint limit $\varepsilon\to0$ and $\mathfrak{c}\to\infty$ for rVMB, without imposing any \emph{a priori} coupling between the two parameters, for any given $T>0$, and show that the limiting dynamics are described by the classical isentropic EP.
\subsection{Main results} 
	For later use, define 
	\begin{equation}\label{e1.13}
		f_R^{\varepsilon,\mathfrak{c}}=\frac{F_R^{\varepsilon,\mathfrak{c}}}{\sqrt{\mathbf{M}_{\mathfrak{c}}}}.
	\end{equation}
	To use the $L^2$-$L^{\infty}$ framework \cite{Guo-KRM-2009}, we also introduce a relativistic global Maxwellian
	$$
	J_{\mathbf{M}}=\frac{n_M \gamma_M}{4 \pi \mathfrak{c}^3 K_2\left(\gamma_M\right)} \exp \left\{-\frac{\mathfrak{c}p^0}{T_M}\right\}
	$$
	with $\gamma_M=\frac{\mathfrak{c}^2}{T_M}$,
	where $n_M>0$ and $T_M>0$ are chosen so that for all $(t,x)\in [0,T]\times\mathbb{R}^3$
    \begin{align}\label{e1.21}
    n_M< n_{\mathfrak{e}}(t, x) <2 n_M, \quad 
    T_M< T_\mathfrak{e}(t, x) <2 T_M .
    \end{align}
    In what follows, we assume that there exist
    $\alpha\in\left(\frac12,1\right)$ and $C>0$ such that
    \begin{equation}\label{e1.15}
    \frac{J_{\mathbf M}}{C}\le \mathbf M_{\mathfrak{c}}(t,x,p)\le C J_{\mathbf M}^{\alpha},
    \end{equation}
	where the constant $C$ is independent of $\mathfrak{c}$. We further define
	\begin{equation}\label{e1.14}
		h_R^{\varepsilon,\mathfrak{c}} :=\frac{w_{\ell}}{\sqrt{J_{\mathbf{M}}}}F_R^{\varepsilon,\mathfrak{c}},
	\end{equation}
	with $w_{\ell} :=\langle p\rangle^{\ell}, \langle p\rangle := \sqrt{1+|p|^2}$ for some $\ell \geq 14$. 	
	
	Before stating Theorem \ref{t1.1}, we introduce the weighted quantity measuring the spatial decay of the initial remainder:
\begin{align}\label{e1.16-0}
\mathcal{I}_0 := \sup_{(x,p) \in \mathbb{R}^3 \times \mathbb{R}^3} \Bigg\{ (1+|x|) \Bigg( \sum_{|\alpha|+|\beta|\leq 1} \big|\partial_x^\alpha \partial_p^\beta h_R^{\varepsilon, \mathfrak{c}}(0,x,p)\big|+\sum_{1\leq |\alpha|\leq 2} \big( \big|\partial_x^\alpha E_R^{\varepsilon, \mathfrak{c}}(0,x)\big| + \big|\partial_x^\alpha B_R^{\varepsilon, \mathfrak{c}}(0,x)\big| \big) \Bigg) \Bigg\}. 
\end{align}
We assume that the initial remainder decays sufficiently fast at spatial infinity so that
\[
\mathcal{I}_0<\infty
\]
and that $\mathcal I_0$ is independent of the speed of light $\mathfrak{c}$.
	\begin{theorem}\label{t1.1}
		Recall $F_0^{\mathfrak{c}}=\mathbf{M}_{\mathfrak{c}}$ in \eqref{e1.11}. Let $\left(n_{\mathfrak{e}}, u_{\mathfrak{e}}, T_{\mathfrak{e}}, E_{\mathfrak{e}}, B_{\mathfrak{e}}\right)(t, x)$ be the smooth solution of rEM \eqref{e.4} given by Theorem \ref{c3.4}.
		Assume that the initial data are of the form
		\begin{align*}
		\begin{cases}
		F^{\varepsilon, \mathfrak{c}}(0, x, p)=\mathbf{M}_{\mathfrak{c}}(0, x, p)+\sum_{n=1}^{2 k-1} \varepsilon^n F_n^{\mathfrak{c}}(0, x, p)+\varepsilon^k F_R^{\varepsilon, \mathfrak{c}}(0, x, p) \geq 0,\\
		E^{\varepsilon,\mathfrak{c}}(0, x)  =E_{\mathfrak{e}}(0, x)+\sum_{n=1}^{2 k-1} \varepsilon^n E_n^{\mathfrak{c}}(0, x)+\varepsilon^k E_R^{\varepsilon,\mathfrak{c}}(0, x), \\
		B^{\varepsilon,\mathfrak{c}} (0, x) =B_{\mathfrak{e}}(0, x)+\sum_{n=1}^{2 k-1} \varepsilon^n B_n^{\mathfrak{c}}(0, x)+\varepsilon^k B_R^{\varepsilon,\mathfrak{c}}(0, x),
		\end{cases}
		\end{align*}
		for $k\geq 5$ with
		\begin{align}\label{e.16-1}
		&\varepsilon^{\frac{5}{2}} \big\|\big(E_R^{\varepsilon,\mathfrak{c}}, B_R^{\varepsilon,\mathfrak{c}}\big)(0)\big\|_{L^{\infty}}+\varepsilon^{\frac{3}{2}}\mathcal{I}_0+\left\|\nabla_{x,p}f_R^{\varepsilon,\mathfrak{c}}(0)\right\|+\left\|\big(\nabla_x E_R^{\varepsilon,\mathfrak{c}},\nabla_x B_R^{\varepsilon,\mathfrak{c}}\big)(0)\right\| \nonumber\\
		&\qquad+\left\|f_R^{\varepsilon,\mathfrak{c}}(0)\right\|+\left\|\big(E_R^{\varepsilon,\mathfrak{c}}, B_R^{\varepsilon,\mathfrak{c}}\big)(0)\right\|+1\leq C<\infty .
		\end{align}
		For any given finite time $T>0$, there exist two independent positive constants $\varepsilon_0 \in(0,1)$ and $\mathfrak{c}_0 \gg 1$ such that, for each $0<\varepsilon \leq \varepsilon_0$ and $\mathfrak{c} \geq \mathfrak{c}_0$,
     there exists a unique classical solution $(F^{\varepsilon, \mathfrak{c}},E^{\varepsilon, \mathfrak{c}},B^{\varepsilon, \mathfrak{c}})$ of rVMB \eqref{e1.1} in $t \in[0, T]$ of the form
		\begin{align}\label{e.16-2}
		\begin{cases}
		F^{\varepsilon, \mathfrak{c}}(t, x, p)=\mathbf{M}_{\mathfrak{c}}(t, x, p)+\sum_{n=1}^{2 k-1} \varepsilon^n F_n^{\mathfrak{c}}(t, x, p)+\varepsilon^k F_R^{\varepsilon, \mathfrak{c}}(t, x, p) \geq 0,\\
		E^{\varepsilon,\mathfrak{c}}(t, x)  =E_{\mathfrak{e}}(t, x)+\sum_{n=1}^{2 k-1} \varepsilon^n E_n^{\mathfrak{c}}(t, x)+\varepsilon^k E_R^{\varepsilon,\mathfrak{c}}(t, x), \\
		B^{\varepsilon,\mathfrak{c}} (t, x) =B_{\mathfrak{e}}(t, x)+\sum_{n=1}^{2 k-1} \varepsilon^n B_n^{\mathfrak{c}}(t, x)+\varepsilon^k B_R^{\varepsilon,\mathfrak{c}}(t, x),
		\end{cases}
		\end{align}
		where the functions $\big(F_n^{\mathfrak{c}},E_n^{\mathfrak{c}},B_n^{\mathfrak{c}}\big)(n=1, \cdots, 2 k-1)$ are constructed in Proposition \ref{t5.3}.
		Furthermore, for the remainder terms $\big(F_R^{\varepsilon,\mathfrak{c}}, E_R^{\varepsilon,\mathfrak{c}}, B_R^{\varepsilon,\mathfrak{c}}\big)$ in \eqref{e1.10}, there exists a constant $C_T > 0$ such that the following estimate holds:
		\begin{align}\label{e1.16}
		\sup _{t \in [0,T]}&\Big\{\varepsilon^{\frac{5}{2}}\left\|\nabla_{x,p}h_R^{\varepsilon,\mathfrak{c}}(t)\right\|_{L^{\infty}}+\varepsilon^{\frac{3}{2}}\left\|h_R^{\varepsilon,\mathfrak{c}}(t)\right\|_{L^{\infty}}+\varepsilon^{\frac{7}{2}}\big\|\big(\nabla_x E_R^{\varepsilon,\mathfrak{c}}, \nabla_x B_R^{\varepsilon,\mathfrak{c}}\big)(t)\big\|_{L^{\infty}}+\varepsilon^{\frac{5}{2}}\big\|\big(E_R^{\varepsilon,\mathfrak{c}}, B_R^{\varepsilon,\mathfrak{c}}\big)(t)\big\|_{L^{\infty}}\nonumber\\
	   &+\left\|\nabla_{x,p}f_R^{\varepsilon,\mathfrak{c}}(t)\right\|+\left\|\big(\nabla_{x}E_R^{\varepsilon,\mathfrak{c}},\nabla_{x}B_R^{\varepsilon,\mathfrak{c}}\big)(t)\right\|+\left\|f_R^{\varepsilon,\mathfrak{c}}(t)\right\|+\left\|\big(E_R^{\varepsilon,\mathfrak{c}}, B_R^{\varepsilon,\mathfrak{c}}\big)(t)\right\|\Big\} \nonumber\\
		&\leq C_T \Big[ \varepsilon^{\frac{5}{2}} \big\|\big(E_R^{\varepsilon,\mathfrak{c}}, B_R^{\varepsilon,\mathfrak{c}}\big)(0)\big\|_{L^{\infty}}+\varepsilon^{\frac{3}{2}}\mathcal{I}_0+\left\|\nabla_{x,p}f_R^{\varepsilon,\mathfrak{c}}(0)\right\|+\left\|\big(\nabla_x E_R^{\varepsilon,\mathfrak{c}},\nabla_x B_R^{\varepsilon,\mathfrak{c}}\big)(0)\right\| \nonumber\\
		&\quad+\left\|f_R^{\varepsilon,\mathfrak{c}}(0)\right\|+\left\|\big(E_R^{\varepsilon,\mathfrak{c}}, B_R^{\varepsilon,\mathfrak{c}}\big)(0)\right\|+1\Big].
		\end{align}
		Moreover, we have that
		\begin{align}\label{e1.17}
			\sup _{t \in [0,T]}&\Big\{\Big\|\frac{F^{\varepsilon, \mathfrak{c}}(t)-\mathbf{M}_{\mathfrak{c}}(t)}{\sqrt{J_{\mathbf{M}}}}\Big\|_{W^{1,\infty}}+\Big\|\frac{F^{\varepsilon,\mathfrak{c}}(t)-\mathbf{M}_{\mathfrak{c}}(t)}{\sqrt{\mathbf{M}_{\mathfrak{c}}(t)}}\Big\|+\left\|E^{\varepsilon, \mathfrak{c}}(t)-E_{\mathfrak{e}}(t)\right\|_{W^{1,\infty}}\nonumber\\
			+&\left\|E^{\varepsilon, \mathfrak{c}}(t)-E_{\mathfrak{e}}(t)\right\| +\left\|B^{\varepsilon, \mathfrak{c}}(t)-B_{\mathfrak{e}}(t)\right\|_{W^{1,\infty}}+\left\|B^{\varepsilon, \mathfrak{c}}(t)-B_{\mathfrak{e}}(t)\right\| \Big\}	\leq C_T \varepsilon,
		\end{align}
		the constant $C_T>0$ is independent of $\mathfrak{c}$ and $\varepsilon$.
	\end{theorem}
	
	\begin{remark}\label{r1.3}
    Theorem \ref{t1.1} yields the hydrodynamic limit from rVMB \eqref{e1.1} to rEM \eqref{e.4}, uniformly for all sufficiently large speeds of light $\mathfrak{c}\geq \mathfrak{c}_0$, without imposing any \emph{a priori} relation between $\varepsilon$ and $\mathfrak{c}$. 
    \end{remark}
    
	\begin{remark}\label{r1.4}
		When $\frac{|u_{\mathfrak{e}}|}{\mathfrak{c}}$ is suitably small, it has been shown in \cite[Lemma 1.1]{Speck-CMP-2011} that there exist positive constants $C>0, n_M>0, T_M>0$, and $\alpha \in\left(\frac{1}{2}, 1\right)$, which are independent of $\mathfrak{c}$, such that \eqref{e1.15} holds.
	\end{remark}
	Combining the Hilbert expansion in Theorem \ref{t1.1} with the Newtonian limit of rEM established in Section \ref{section 3}, we obtain the following combined hydrodynamic and Newtonian limit.
	\begin{theorem}\label{t1.5}
		Assume all the hypotheses of Theorem \ref{t1.1}, and let
    $(\rho,\mathfrak u,\phi)$ denote the corresponding smooth solution
    of the classical EP system \eqref{e0.5} constructed in Lemma \ref{l3.2}. Let $\theta:=\theta(\rho)$ be determined by \eqref{e0.6} and $\mu$ be the local Maxwellian of classical Boltzmann equation, i.e.,
		\begin{equation}\label{e1.18}
			\mu(t, x, p)=\frac{\rho}{(2 \pi \theta)^{\frac{3}{2}}} e^{-\frac{|p-\mathfrak{u}|^2}{2 \theta}} .
		\end{equation}
		Then there exist independent positive constants $\varepsilon_0 \in(0,1)$ and $\mathfrak{c}_0 \gg 1$ such that for all $0<\varepsilon \leq \varepsilon_0$ and $\mathfrak{c} \geq \mathfrak{c}_0$, the following estimate holds:
		\begin{align}\label{e1.19}
			\sup _{0 \leq t \leq T}\left\|\left(F^{\varepsilon, \mathfrak{c}}-\mu\right)(t) e^{\delta_0|p|}\right\|_{L^{\infty}} \lesssim  \varepsilon + \mathfrak{c}^{-1}, 
		\end{align}
		and 
		\begin{align}\label{e1.19*}
		    \sup _{0 \leq t \leq T}\left\|(E^{\varepsilon, \mathfrak{c}}-\nabla \phi)(t)\right\|_{L^{\infty}} \lesssim  \varepsilon + \mathfrak{c}^{-1}, \quad
			\sup _{0 \leq t \leq T}\left\|B^{\varepsilon, \mathfrak{c}}(t)\right\|_{L^{\infty}} \lesssim  \varepsilon + \mathfrak{c}^{-1},
		\end{align}
		where all the positive constants $\varepsilon_0, \mathfrak{c}_0$ and $\delta_0$ are independent of $\varepsilon$ and $\mathfrak{c}$. Moreover, the implicit constants in \eqref{e1.19}-\eqref{e1.19*} depend only on $T>0$ and the uniform bounds of the background profiles, and are independent of $\varepsilon$ and $\mathfrak{c}$.
	\end{theorem}
  
    \begin{remark}\label{r1.6}
    The compressible EP serves as an important macroscopic model in the classical Newtonian setting. However, its structure combines a hyperbolic Euler part, for which perturbations propagate at finite speeds, with an elliptic Poisson equation for the self-consistent electric field, which induces an instantaneous response and thus destroys the strict finite-speed propagation property of the system. Theorem \ref{t1.5} shows that, nevertheless, the classical EP can be rigorously derived from rVMB \eqref{e1.1} through the combined hydrodynamic and Newtonian limits, starting from a relativistic kinetic model with finite propagation speed.
    \end{remark}

\subsection{Main difficulties and strategy of the proof}

We now explain the main analytical difficulties of this paper and the strategy used to overcome them.

\medskip
\noindent
{\bf (1) Newtonian limit at the fluid level.}
The first difficulty lies in justifying the Newtonian limit for the rEM \eqref{e.4} in $\mathbb{R}^3$. 
Although \eqref{e.4} formally converges to the classical EP system \eqref{e0.5} as $\mathfrak{c}\to\infty$, a direct energy estimate for the difference is not sufficient, since the electromagnetic part remains coupled and one needs quantitative control of the first-order magnetic profile appearing in the expansion.

In \cite{Peng-CAMSB-2007}, Peng-Wang established the Newtonian limit from the classical EM to the compressible EP in $\mathbb{T}^3$ locally in time.
For rEM \eqref{e.4} in $\mathbb{R}^3$, motivated by \cite{Peng-CAMSB-2007}, we construct an asymptotic expansion in powers of $\mathfrak{c}^{-1}$. The leading-order system corresponds to the classical EP in $\mathbb{R}^3$, while the first-order correction is governed by a linearized EP system coupled with the curl-div problem \eqref{e.10-2} for $B_1$. To close the remainder estimates, in particular for the term $R_B^{\mathfrak{c}}$ in $\eqref{e.11-1}$, it is essential to establish an $L^2$-based estimate for $B_1$. We remark that, on the periodic domain $\mathbb{T}^3$, this can be easily resolved via the Poincare inequality (by similar arguments as in \cite{Peng-CAMSB-2007}); in the whole-space setting, however, this approach is no longer effective. 

Indeed, since $\operatorname{div} B_1=0$, we introduce a vector potential $\psi$ such that
\[
B_1=\nabla\times\psi,\qquad \operatorname{div}\psi=0.
\]
Then $\psi$ satisfies
\[
-\Delta\psi=\partial_t\nabla\phi-4\pi n_0u_0.
\]
However, standard elliptic regularity only yields $\nabla^2\psi\in H^s(\mathbb{R}^3)$, which does not by itself imply $\nabla\psi\in L^2(\mathbb{R}^3)$. The crucial observation is that the irrotationality condition
\[
\nabla\times u_0=0
\]
is preserved by the limiting EP flow. Hence one may write $u_0=\nabla\phi_0$, so that
\[
\partial_t\nabla\phi-4\pi n_0u_0
=\nabla(\partial_t\phi-4\pi\bar n\phi_0)-4\pi (n_0-\bar n)u_0.
\]
Therefore, the gradient part vanishes against divergence-free test functions, and the effective forcing is reduced to $-4\pi (n_0-\bar n)u_0$. Since
\[
(n_0-\bar n)u_0\in L^1(\mathbb{R}^3)\cap L^2(\mathbb{R}^3)\subset L^{6/5}(\mathbb{R}^3),
\]
the Lax-Milgram argument yields $\nabla\psi\in L^2(\mathbb{R}^3)$, and hence $B_1\in L^2(\mathbb{R}^3)$. Combining this with elliptic regularity, we further obtain the required higher-order bounds for $B_1$. Without the irrotationality condition $\nabla\times u_0=0$, the above argument no longer applies directly; it would be interesting to understand whether the whole-space Newtonian limit can still be justified in that case.

With this ingredient, together with the smooth theory for the EP \eqref{e0.5}, the standard theory of linear symmetrizable hyperbolic systems, and the curl-div analysis, we rigorously justify the expansion. For well-prepared initial data, we prove that the rEM solution converges to the EP solution at rate $\mathfrak{c}^{-1}$ on a time interval independent of $\mathfrak{c}$; see Lemma \ref{t3.6} and Theorem \ref{c3.4}.

\medskip
\noindent
{\bf (2) Uniform-in-$\mathfrak{c}$ estimates for the Hilbert expansion coefficients.}
The second difficulty is to construct the Hilbert expansion for rVMB with estimates uniform-in-$\mathfrak{c}$. Indeed, once the expansion is substituted into the equation, the coefficients $F_n^{\mathfrak{c}}$ are determined recursively, and their kinetic parts satisfy linear equations of the form
\[
\mathbf{L}_{\mathfrak{c}}F_n^{\mathfrak{c}}=G_n,
\]
see \eqref{N.33}. Therefore, a key issue is to obtain quantitative bounds on the pseudo-inverse $\mathbf{L}_{\mathfrak{c}}^{-1}$, and in particular on its derivatives, with bounds independent of $\mathfrak{c}$.

In the classical Boltzmann theory, such estimates are by now well understood. For instance, Jiang et al.~\cite{Jiang-CMP-2025} established the exponential decay estimates for the inverse of the linearized Boltzmann operator and its $(t,x)$-derivatives for both hard and soft cutoff potentials. In the relativistic setting, however, the situation is much subtler. While Wang-Xiao \cite{Wang-Xiao-JLMS-2026} proved exponential decay for $\mathbf{L}_{\mathfrak{c}}^{-1}$ itself, uniform estimates for $\partial_{t,x,p}\mathbf{L}_{\mathfrak{c}}^{-1}$ do not seem to be available in the literature. The main obstruction is that momentum differentiation of the relativistic collision operator generates singular kernels, a phenomenon that is absent, or substantially weaker, in the non-relativistic theory.

A crucial new ingredient in our analysis is a uniform-in-$\mathfrak{c}$ high-order momentum derivative estimate for the relativistic collision operator. Extending the momentum regularity framework of Guo-Strain \cite{Guo-CMP-2012} from the normalized case $\mathfrak{c}=1$ to an arbitrary speed of light $\mathfrak{c}$, we introduce a refined $\mathfrak{c}$-dependent decomposition of the momentum region and derive new derivative estimates in both the Glassey--Strauss frame and the {\it center of momentum} frame; see Lemma \ref{l2.4} and Appendix \ref{Appendix A}. This is not a simple restoration of the physical parameter $\mathfrak{c}$, since the general-speed-of-light setting brings in additional $\mathfrak{c}$-dependent scales that must be tracked carefully in all high-order momentum derivative estimates. 

With the help of Lemma \ref{l2.5}, which is established based on the above decomposition strategy, and the exponential decay estimate for $\mathbf{L}_{\mathfrak{c}}^{-1}$, we prove that for any fixed $0<\lambda_l<\lambda_{l-1}<\cdots<\lambda_1<\lambda<1$, $\mathfrak{k}>\frac32$, $1\le l\le s-2$, and $|\alpha|+|\beta|=l$,
\begin{align*}
\left|\partial_{\beta}^{\alpha}\left[\mathbf{L}_{\mathfrak{c}}^{-1} g(t, x, p)\right]\right| \leq C \sum_{i=0}^l \sum_{|\alpha^\prime|+|\beta^\prime|=i} \big\|\langle p\rangle^{\mathfrak{k}} \mathbf{M}_{\mathfrak{c}}^{-\frac{\lambda_i}{2}} \partial_{\beta^{\prime}}^{\alpha^{\prime}} g\big\|_{L^{\infty}} \cdot \mathbf{M}_{\mathfrak{c}}^{\frac{\lambda_l}{2}}(p),\quad p \in \mathbb{R}^3,
\end{align*}
where the constant $C>0$ is independent of $\mathfrak{c}$. This uniform estimate is one of the main technical ingredients of the paper. In particular, it allows us to solve the recursive equations \eqref{N.11}-\eqref{N.12} for the Hilbert coefficients in a stable way and to derive uniform bounds for $F_n^{\mathfrak{c}}$ for all $n\ge 1$; see Proposition \ref{t5.3}.

\medskip
\noindent
{\bf (3) Uniform control of the remainder and curved characteristics.}
The third major difficulty is to control the remainder terms $\big(F_R^{\varepsilon,\mathfrak{c}}, E_R^{\varepsilon,\mathfrak{c}}, B_R^{\varepsilon,\mathfrak{c}}\big)$ uniformly in both $\varepsilon$ and $\mathfrak{c}$. For this purpose, we use the $H^1$-$W^{1,\infty}$ framework developed by Guo-Xiao \cite{Guo-CMP-2021}, which may be viewed as a refinement of the classical $L^2$-$L^\infty$ method of \cite{Guo-KRM-2009,Guo-CPAM-2010}. In the present problem, this framework is especially suitable because the remainder is transported along the curved relativistic characteristics generated by the Lorentz force:
\begin{align*}
\frac{d X(\tau;t,x,p)}{d\tau}&=\hat P(\tau;t,x,p),\\
\frac{d P(\tau;t,x,p)}{d\tau}
&=-E^{\varepsilon,\mathfrak{c}}(\tau,X(\tau;t,x,p))
-\frac{P}{P^0}(\tau;t,x,p)\times
B^{\varepsilon,\mathfrak{c}}(\tau,X(\tau;t,x,p)),
\end{align*}
with $X(t;t,x,p)=x$ and $P(t;t,x,p)=p$.

To justify these characteristics and derive $L^\infty$ estimates along them, one needs a uniform Lipschitz bound for the electro-magnetic field. More precisely, for $k\ge 5$, we need
\begin{align*}
\sup_{t\in[0,T]}
\left\|
\big(E^{\varepsilon,\mathfrak{c}},B^{\varepsilon,\mathfrak{c}}\big)(t)
\right\|_{W^{1,\infty}}
\le
\sum_{i=0}^{2k-1}\varepsilon^i
\left\|
\big(E_i^{\mathfrak{c}},B_i^{\mathfrak{c}}\big)(t)
\right\|_{W^{1,\infty}}
+\varepsilon^k
\left\|
\big(E_R^{\varepsilon,\mathfrak{c}},B_R^{\varepsilon,\mathfrak{c}}\big)(t)
\right\|_{W^{1,\infty}}
\lesssim 1.
\end{align*}
Since the lower-order terms are already controlled, the main issue is therefore to estimate the remainder fields
$\left(E_R^{\varepsilon,\mathfrak{c}},B_R^{\varepsilon,\mathfrak{c}}\right)$ in $W^{1,\infty}(\mathbb{R}^3)$.
This is achieved through a bootstrap argument. We assume
\begin{align}\label{e1.20}
\begin{gathered}
			\sup _{t \in[0, T]} \varepsilon^2\big\|h_R^{\varepsilon,\mathfrak{c}}(t)\big\|_{L^{\infty}} \leq \varepsilon^{\frac{1}{4}}, \quad \sup _{t \in[0, T]} \varepsilon^3\big\|\nabla_{x,p} h_R^{\varepsilon,\mathfrak{c}}(t)\big\|_{L^{\infty}}\leq \varepsilon^{\frac{1}{8}},\\
	    \sup_{t \in[0, T]} \varepsilon^3\left\|\left(E_R^{\varepsilon,\mathfrak{c}}, B_R^{\varepsilon,\mathfrak{c}}\right)(t)\right\|_{L^{\infty}} \leq \varepsilon^{\frac{1}{4}}, \quad
		\sup_{t \in[0, T]} \varepsilon^4\left\|\big(\nabla_x  E_R^{\varepsilon,\mathfrak{c}}, \nabla_x B_R^{\varepsilon,\mathfrak{c}}\big)(t)\right\|_{L^{\infty}} \leq \varepsilon^{\frac{1}{8}},
	\end{gathered}
\end{align}
and then derive the $L^\infty$ estimate of $h_R^{\varepsilon,\mathfrak{c}}$ in Section \ref{section 5}, where the uniform-in-$\mathfrak{c}$ bounds for the relativistic collision operator from \cite{Wang-Xiao-JLMS-2026} play a key role. After that, in Section \ref{section 6}, we differentiate the remainder equation and obtain the $W^{1,\infty}$ estimate of $h_R^{\varepsilon,\mathfrak{c}}$ in a similar manner. At this stage, however, the differentiated equations \eqref{D.4} and \eqref{D.8} produce additional error terms, and therefore the $W^{1,\infty}$ argument must be supplemented by an $H^1$ estimate for the remainder triple
$\left(f_R^{\varepsilon,\mathfrak{c}},
E_R^{\varepsilon,\mathfrak{c}},
B_R^{\varepsilon,\mathfrak{c}}\right)$,
which is established in Section \ref{section 8}.

The next step is to estimate the electro-magnetic field remainders $\left(E_R^{\varepsilon,\mathfrak{c}},B_R^{\varepsilon,\mathfrak{c}}\right)$ in $W^{1,\infty}(\mathbb{R}^3)$ through the Glassey--Strauss representation formula \cite{Glassey-ARMA-1986}. 
Here, unlike the light-speed normalized setting of Guo-Xiao \cite{Guo-CMP-2021}, a direct estimate of the Glassey--Strauss kernels by absolute values produces two apparently non-uniform contributions. In the $L^\infty$ estimate of the field remainders, the term $E_{R_T,i}^{\varepsilon,\mathfrak{c}}$, if estimated directly by absolute values, seems to give
\[
\mathfrak{c}\int_0^t \|h_R^{\varepsilon,\mathfrak{c}}(\tau)\|_{L^\infty}\,d\tau,
\]
while in the $W^{1,\infty}$ estimate the singular term
\[
\iint_{|x-y|<\mathfrak{c}t}\frac{a_A(\omega,\hat p)}{|x-y|^3}
F_R^{\varepsilon,\mathfrak{c}}
\Bigl(t-\frac{|x-y|}{\mathfrak{c}},y,p\Bigr)\,d y\,d p
\]
if estimated directly by absolute values, seems to yield
\[
\ln(1+\mathfrak{c}t)\sup_{\tau\in[0,t]}
\|h_R^{\varepsilon,\mathfrak{c}}(\tau)\|_{L^\infty}.
\]
These two terms are artifacts of a crude treatment of the complete kernels rather than structural losses of the system. The key point is to separate the Coulomb and Hessian-Coulomb parts obtained at $p=0$ from the full retarded kernels. The remaining kernels carry a factor of order $\min\{|p|/\mathfrak{c},1\}$, which yields a compensating $\mathfrak{c}^{-1}$ gain after integration in momentum. The principal Coulomb parts are then compared with the instantaneous Poisson operators truncated at the same radius $\mathfrak{c} t$. The discrepancy is controlled by the continuity equation
\[
\partial_t\rho_R+\nabla_x\cdot j_R=0,
\qquad
\rho_R=\int_{\mathbb R^3}F_R\,dp,
\qquad
j_R=\int_{\mathbb R^3}\hat p F_R\,dp.
\]
This retarded-to-instantaneous comparison removes the explicit $\mathfrak{c}$-growth in both the field and gradient-field estimates.

With this refined decomposition, we close the bootstrap \eqref{e1.20} and justify the Hilbert expansion uniformly for all $\mathfrak{c}\geq\mathfrak{c}_0$, with no \emph{a priori} upper bound on $\mathfrak{c}$ in terms of $\varepsilon$; see Theorem \ref{t1.1}.
	\subsection{Organization of the paper}
	In Section \ref{section 2}, we provide two distinct representations of the relativistic collision operator in the {\it center of momentum} frame and the Glassey--Strauss frame, respectively, as well as some basic estimates. In Section \ref{section 3}, the Newtonian limit of rVMB is obtained by constructing the uniformly convergent asymptotic expansions for solutions to the rEM equations in $\mathbb{R}^3$ with respect to $\mathfrak{c}^{-1}$. Section \ref{section 4} is devoted to the uniform-in-$\mathfrak{c}$ estimates on the linear part of Hilbert expansion, along with the decay estimate of parameter derivatives of $\mathbf{L}_{\mathfrak{c}}^{-1}$. In Section \ref{section 5}, we study the characteristics and $L^{\infty}$ estimate of $h_R^{\varepsilon,\mathfrak{c}}$ under the crucial {\it a priori} assumptions \eqref{e1.20}.
	In Section \ref{section 6}, the $L^{ \infty}$ estimates for $\nabla_x h_R^{\varepsilon,\mathfrak{c}}$ and $\nabla_p h_R^{\varepsilon,\mathfrak{c}}$ are obtained.
	In Section \ref{section 7}, we estimate the $W^{1, \infty}$ norm of $\big( E_R^{\varepsilon,\mathfrak{c}}, B_R^{\varepsilon,\mathfrak{c}}\big)$ in terms of the $W^{1, \infty}$ norm of $h_R^{\varepsilon,\mathfrak{c}}$ and $H^1$ norm of $f_R^{\varepsilon,\mathfrak{c}}$, via Glassey--Strauss representation.  
	In Section \ref{section 8},  $H^1$ estimates for $\left(f_R^{\varepsilon,\mathfrak{c}}, E_R^{\varepsilon,\mathfrak{c}}, B_R^{\varepsilon,\mathfrak{c}}\right)$ are derived.
	In Section \ref{section 9}, we finally verify \eqref{e1.20} and prove our main results, Theorems \ref{t1.1} and \ref{t1.5}.
	
	\subsection{Notations}

Throughout the paper, $C$ denotes a generic positive constant which may change from line to line. Likewise, $C(a)$, $C(a,b)$, and so on denote generic positive constants depending on the indicated quantities and may also vary from line to line. We write $A \lesssim B$ if there exists a positive constant $C$ such that $A \leq C B$, and write $A \approx B$ if there exists a constant $C>1$ such that
\[
\frac{1}{C}A \le B \le CA.
\]

We use the standard notation $W^{k,2}(\mathbb R_x^3)$ (or $W^{k,2}(\mathbb R_x^3\times \mathbb R_p^3)$) and $W^{k,\infty}(\mathbb R_x^3)$ (or $W^{k,\infty}(\mathbb R_x^3\times \mathbb R_p^3)$) for Sobolev spaces, with corresponding norms $\|\cdot\|_{H^k}$ and $\|\cdot\|_{W^{k,\infty}}$, respectively. We also use $\|\cdot\|$ and $\|\cdot\|_{L^\infty}$ to denote the standard $L^2$ and $L^\infty$ norms on either $\mathbb R_x^3\times \mathbb R_p^3$ or $\mathbb R_x^3$, depending on the context. The standard inner product on $L^2(\mathbb R_x^3\times \mathbb R_p^3)$ or $L^2(\mathbb R_p^3)$ is denoted by $\langle\cdot,\cdot\rangle$.

We use $D_x$ and $D_p$ to denote arbitrary derivatives with respect to the spatial variable $x$ and the momentum variable $p$, respectively. For $\alpha=(\alpha^0,\alpha^1,\alpha^2,\alpha^3)\in\mathbb N^4$ and $\beta=(\beta^1,\beta^2,\beta^3)\in\mathbb N^3$, define
\[
\partial_\beta^\alpha
:=\partial_t^{\alpha^0}\partial_{x_1}^{\alpha^1}\partial_{x_2}^{\alpha^2}\partial_{x_3}^{\alpha^3}
\partial_{p_1}^{\beta^1}\partial_{p_2}^{\beta^2}\partial_{p_3}^{\beta^3}.
\]
For multi-indices $\alpha_1,\alpha_2\in\mathbb N^4$, we write $\alpha_1\le \alpha_2$ if $\alpha_1^i\le \alpha_2^i$ for each component, and $\alpha_1<\alpha_2$ if $\alpha_1\le \alpha_2$ and $|\alpha_1|<|\alpha_2|$, where
$|\alpha|:=\alpha^0+\alpha^1+\alpha^2+\alpha^3$.
The same convention applies to multi-indices in $\mathbb N^3$.

For later use, for any integer $k\ge 0$, we define the mixed space-time Sobolev norm
\[
\|f(t)\|_{\mathcal H^k}^2
:=\sum_{l+|\alpha|\le k}\|\partial_t^l\partial_x^\alpha f(t)\|^2,
\]
where $\alpha\in\mathbb N^3$ is a spatial multi-index and $\|\cdot\|$ denotes the standard $L^2(\mathbb R^3)$ norm.
	\
	\section{Preliminaries}\label{section 2}
	
	\subsection{Bessel function}
	The modified Bessel function of the second kind is given by
    \begin{align}\label{2.1}
    K_j(z)=\left(\frac{z}{2}\right)^j \frac{\Gamma\left(\frac{1}{2}\right)}{\Gamma\left(j+\frac{1}{2}\right)} \int_1^{\infty} e^{-z t}\left(t^2-1\right)^{j-\frac{1}{2}} d t, \quad j \geq 0, z>0.
    \end{align}

    Some of the properties of the modified Bessel function $K_j$ are needed.
    \begin{lemma}\emph{(\cite{Olver-1997,Watson-1980})}\label{l2.0}
     It holds that
    $$
    K_{j+1}(z)=\frac{2 j}{z} K_j(z)+K_{j-1}(z), \quad j \geq 1,
    $$
    and
    $$
    \frac{d}{d z}\left(\frac{K_j(z)}{z^j}\right)=-\left(\frac{K_{j+1}(z)}{z^j}\right), \quad j \geq 0 .
    $$
    The asymptotic expansion for $K_j(z)$ takes the form
    $$
    K_j(z)=\sqrt{\frac{\pi}{2 z}} \frac{1}{e^z}\left[\sum_{m=0}^{n-1} A_{j, m} z^{-m}+\gamma_{j, n}(z) z^{-n}\right], \quad j \geq 0, n \geq 1,
    $$
    where the following additional identities and inequalities also hold:
    \begin{align*}
    A_{j, 0} & =1, \\
    A_{j, m} & =\frac{1}{m!8^m}\left(4 j^2-1\right)\left(4 j^2-3^2\right) \cdots\left(4 j^2-(2 m-1)^2\right), \quad j \geq 0, m \geq 1, \\
    \left|\gamma_{j, n}(z)\right| & \leq 2\left|A_{j, n}\right| \exp \left(\left[j^2-\frac{1}{4}\right] z^{-1}\right), \quad j \geq 0, n \geq 1, \\
    K_j(z) & <K_{j+1}(z), \quad j \geq 0.
    \end{align*}
    Furthermore, for $j \leq n+\frac{1}{2}$, one has a more exact estimate
    \begin{align*}
    \left|\gamma_{j, n}(z)\right| \leq\left|A_{j, n}\right|.
   \end{align*}
    \end{lemma}
    Using the asymptotic expansions for $K_2(z)$ and $K_3(z)$ in Lemma \ref{l2.0}, one can obtain
    \begin{gather}
    K_2(z)=\sqrt{\frac{\pi}{2 z}} \frac{1}{e^z}\left[1+\frac{15}{8 z}+\frac{105}{2!(8 z)^2}-\frac{945}{3!(8 z)^3}+\frac{31185}{4!(8 z)^4}+O\left(z^{-5}\right)\right],\label{2.2}\\
     K_3(z)=\sqrt{\frac{\pi}{2 z}} \frac{1}{e^z}\left[1+\frac{35}{8 z}+\frac{945}{2!(8 z)^2}+\frac{10395}{3!(8 z)^3}-\frac{135135}{4!(8 z)^4}+O\left(z^{-5}\right)\right], \label{2.3}\\
        \frac{K_3\left(z\right)}{K_2\left(z\right)}-1=\frac{5}{2 z}+\frac{15}{8 z^2}-\frac{15}{8 z^3}+O\left(z^{-4}\right),\label{2.4}\\
    \left(\frac{K_3(z)}{K_2(z)}\right)^2-\frac{5}{z} \frac{K_3(z)}{K_2(z)}+\frac{1}{z^2}-1=-\frac{3}{2 z^2}-\frac{15}{4 z^3}+\frac{45}{8 z^4}+O\left(z^{-5}\right).\label{2.5}
    \end{gather}
	\subsection{Estimates for parameter derivatives of the relativistic local Maxwellian}
	For any $l=0,1,\cdots,s-2$($s$ is a fixed integer defined in Lemma \ref{l3.2}), denote
	\begin{align}\label{e2.34}
		\mathfrak{e}_l:=\frac{1}{\inf_{t, x}n_{\mathfrak{e}}}+\frac{1}{\inf_{t, x}T_\mathfrak{e}}+\sum_{|\alpha|\leq l}\sup_{t, x}\left|\partial^{\alpha}_{t,x}\left(n_{\mathfrak{e}},u_{\mathfrak{e}},T_\mathfrak{e}\right) \right|,
	\end{align}
	where the derivative operator $\partial_{t, x}^\alpha$ is defined by
	$\partial_{t, x}^\alpha=\partial_t^{\alpha^0} \partial_{x_1}^{\alpha^1} \partial_{x_2}^{\alpha^2} \partial_{x_3}^{\alpha^3}$ for $\alpha=\left(\alpha^0, \alpha^1, \alpha^2, \alpha^3\right)\in\mathbb{N}^4$.
	From Theorem \ref{c3.4}, we obtain that $n_{\mathfrak{e}}$, $T_\mathfrak{e}$ have both positive lower and upper bounds. Moreover, all derivatives of $\left(n_{\mathfrak{e}},u_{\mathfrak{e}},T_\mathfrak{e}\right)$ up to order $s-2$ are also bounded. Consequently, $\mathfrak{e}_l<\infty$. And it's obvious that $\mathfrak{e}_l\leq\mathfrak{e}_{l+1}$ for $l=0,1,\cdots,s-2$. Recall the relativistic local Maxwellian \eqref{e1.11}, we have the following estimate:
	
	\begin{lemma}\label{l2.12}
	    For any $0\leq l\leq N$($N$ is a fixed integer), let $\alpha=\left(\alpha_0, \alpha_1, \alpha_2, \alpha_3\right)\in \mathbb{N}^4$ and $\beta=\left(\beta_1, \beta_2, \beta_3\right) \in \mathbb{N}^3$ with $|\alpha|+|\beta|=l$. Then there holds
	    \begin{align}\label{e2.35}
	        \left|\partial_{\beta}^{\alpha} \mathbf{M}_{\mathfrak{c}}(t, x, p)\right|\leq C(\mathfrak{e}_l)\left\langle p\right\rangle^{2|\alpha|+|\beta|} \mathbf{M}_{\mathfrak{c}}(t, x, p),
	    \end{align}
	    where $C(\mathfrak{e}_l)$ depends only on $\mathfrak{e}_l$ and is independent of $\mathfrak{c}$.
	\end{lemma}
    The proof follows from a direct calculation and is omitted for brevity.
	\subsection{Hilbert-Schmidt formulation}
	Using Lorentz transformations in \cite{Groot-RKT-1980, Strain-CMP-2010}, in the \emph{center of momentum system}, $Q_{\mathfrak{c}}(f_1, f_2)$ can be written as
	\begin{align}\label{04}
		Q_{\mathfrak{c}}(f_1, f_2)&=\iint_{\mathbb{R}^3\times \mathbb{S}^2} v_\phi \left[f_1\left(p^{\prime}\right) f_2\left(q^{\prime}\right)-f_1(p) f_2(q)\right] d \omega d q	\nonumber\\
        &:=Q_{\mathfrak{c}}^{+}(f_1, f_2)-Q_{\mathfrak{c}}^{-}(f_1, f_2),
	\end{align}
	where $v_\phi=v_\phi(p, q)$ is the M{\o}ller velocity
	$$
	v_\phi(p, q):=\frac{\mathfrak{c}}{2} \sqrt{\Big|\frac{p}{p^0}-\frac{q}{q^0}\Big|^2-\Big|\frac{p}{p^0} \times \frac{q}{q^0}\Big|^2}=\frac{\mathfrak{c}}{4} \frac{g \sqrt{\mathfrak{s}}}{p^0 q^0} .
	$$
	The pre-post collisional momentum in \eqref{04} satisfies
	\begin{equation}\label{05}
		\left\{
		\begin{aligned}
			&p^{\prime}=\frac12(p+q)+\frac12 g\Big(\omega+(\gamma_0-1)(p+q)\frac{(p+q)\cdot \omega}{|p+q|^2}\Big),\\
			&q^{\prime}=\frac12(p+q)-\frac12 g\Big(\omega+(\gamma_0-1)(p+q)\frac{(p+q)\cdot \omega}{|p+q|^2}\Big),
		\end{aligned}
		\right.
	\end{equation}
	where $\gamma_0:=\left(p^0+q^0\right) / \sqrt{\mathfrak{s}}$. The pre-post collisional energy is given by
	\begin{equation}\label{06}
		\left\{
		\begin{aligned}
			&p^{\prime0}=\frac12(p^0+q^0)+\frac{1}{2}\frac{g}{\sqrt{\mathfrak{s}}}(p+q)\cdot \omega,\\
			&q^{\prime0}=\frac12(p^0+q^0)-\frac{1}{2}\frac{g}{\sqrt{\mathfrak{s}}}(p+q)\cdot \omega.
		\end{aligned}
		\right.
	\end{equation}
	
	We define
	$$
	\begin{aligned}
		\mathbf{L}_{\mathfrak{c}} f&:=-\frac{1}{\sqrt{\mathbf{M}_{\mathfrak{c}}}}\left[Q_{\mathfrak{c}}\left(\sqrt{\mathbf{M}_{\mathfrak{c}}} f, \mathbf{M}_{\mathfrak{c}}\right)+Q_{\mathfrak{c}}\left(\mathbf{M}_{\mathfrak{c}}, \sqrt{\mathbf{M}_{\mathfrak{c}}} f\right)\right]=\nu_{\mathfrak{c}} f-\mathbf{K}_{\mathfrak{c}} f,\\
		\Gamma_{\mathfrak{c}}\left(f_1, f_2\right)&:=\frac{1}{\sqrt{\mathbf{M}_{\mathfrak{c}}}}Q_{\mathfrak{c}}^+\left(\sqrt{\mathbf{M}_{\mathfrak{c}}} f_1, \sqrt{\mathbf{M}_{\mathfrak{c}}} f_2\right)-\frac{1}{\sqrt{\mathbf{M}_{\mathfrak{c}}}}Q_{\mathfrak{c}}^-\left(\sqrt{\mathbf{M}_{\mathfrak{c}}} f_1, \sqrt{\mathbf{M}_{\mathfrak{c}}} f_2\right)\\
		&:=\Gamma_{\mathfrak{c}}^+\left(f_1, f_2\right)-\Gamma_{\mathfrak{c}}^-\left(f_1, f_2\right),
	\end{aligned}
	$$
	where the collision frequency $\nu_{\mathfrak{c}}=\nu_{\mathfrak{c}}(t, x, p)$ is defined as
	\begin{align}\label{e2.1}
		\nu_{\mathfrak{c}}(t,x,p)&:=\int_{\mathbb{R}^3} \int_{\mathbb{S}^2} v_\phi(p, q) \mathbf{M}_{\mathfrak{c}}(q) d \omega d q\nonumber\\
        &=\frac{\mathfrak{c}}{2} \frac{1}{p^0} \int_{\mathbb{R}^3} \frac{d q}{q^0} \int_{\mathbb{R}^3} \frac{d q^{\prime}}{q^{\prime 0}} \int_{\mathbb{R}^3} \frac{d p^{\prime}}{p^{\prime 0}} W\left(p, q \mid p^{\prime}, q^{\prime}\right) \mathbf{M}_{\mathfrak{c}}(q),
	\end{align}
	and $\mathbf{K}_{\mathfrak{c}} f$ takes the following form:
	$$
	\begin{aligned}
		\mathbf{K}_{\mathfrak{c}} f:= & \int_{\mathbb{R}^3} \int_{\mathbb{S}^2} v_\phi(p, q) \sqrt{\mathbf{M}_{\mathfrak{c}}(q)}\left[\sqrt{\mathbf{M}_{\mathfrak{c}}\left(q^{\prime}\right)} f\left(p^{\prime}\right)+\sqrt{\mathbf{M}_{\mathfrak{c}}\left(p^{\prime}\right)} f\left(q^{\prime}\right)\right] d \omega d q \\
		& \quad-\int_{\mathbb{R}^3} \int_{\mathbb{S}^2} v_\phi(p, q) \sqrt{\mathbf{M}_{\mathfrak{c}}(p) \mathbf{M}_{\mathfrak{c}}(q)} f(q) d \omega d q \\
		= & \mathbf{K}_{\mathfrak{c} 2} f-\mathbf{K}_{\mathfrak{c} 1} f .
	\end{aligned}
	$$
	Then it is clear that the kernel of $\mathbf{K}_{\mathfrak{c}1}$ takes the form 
	\begin{align}\label{e2.2}
		k_{\mathfrak{c}1}(p,q)&=\int_{\mathbb{S}^2} v_{\phi}(p,q) \sqrt{\mathbf{M}_{\mathfrak{c}}(p)\mathbf{M}_{\mathfrak{c}}(q)} d\omega
		=\frac{\pi \mathfrak{c}g\sqrt{\mathfrak{s}}}{p^0q^0}\sqrt{\mathbf{M}_{\mathfrak{c}}(p)\mathbf{M}_{\mathfrak{c}}(q)}.
	\end{align}
	By similar arguments as in \cite{Strain-CMP-2010}, we can deduce that the kernel of $\mathbf{K}_{\mathfrak{c}2}$ is
	\begin{align}\label{e2.3}
		k_{\mathfrak{c}2}(p,q)=\frac{\mathfrak{c}}{p^0q^0}\int_{\mathbb{R}^3}\frac{d q^{\prime}}{q^{\prime 0}}\int_{\mathbb{R}^3}\frac{d p^{\prime}}{p^{\prime 0}}\bar{\mathfrak{s}}\delta^{(4)}(p^\mu+p^{\prime \mu}-q^{\mu}-q^{\prime \mu})\sqrt{\mathbf{M}_{\mathfrak{c}}(p^{\prime})\mathbf{M}_{\mathfrak{c}}(q^{\prime})},
	\end{align}
	where
	\begin{align*}
	 \bar{g}^2:=g^2-\frac{1}{2}(p^{\mu}+q^{\mu})(p^{\prime}_{\mu}+q^{\prime}_{\mu}-p_{\mu}-q_{\mu}), \quad \bar{\mathfrak{s}}:=\bar{g}^2+4\mathfrak{c}^2.
	\end{align*}
	
	Denote $k_{\mathfrak{c}}(p,q):=k_{\mathfrak{c}2}(p,q)-k_{\mathfrak{c}1}(p,q)$ and
	\begin{align*}
		k_{\mathfrak{c}w}(p,q):=k_{\mathfrak{c}}(p,q)\frac{w_\ell(p)}{w_\ell(q)}.
	\end{align*}
	
	It is well-known that $\mathbf{L}_{\mathfrak{c}}$ is a self-adjoint non-negative definite operator in $L_p^2$ space with the kernel
	$$
	\mathcal{N}_{\mathfrak{c}}=\operatorname{span}\left\{\sqrt{\mathbf{M}_{\mathfrak{c}}}, p_i \sqrt{\mathbf{M}_{\mathfrak{c}}}(i=1,2,3), \frac{p^0}{\mathfrak{c}} \sqrt{\mathbf{M}_{\mathfrak{c}}}\right\} .
	$$
	Let $\mathbf{P}_{\mathfrak{c}}$ be the orthogonal projection from $L_p^2$ onto $\mathcal{N}_{\mathfrak{c}}$. For given $f$, we denote the macroscopic part $\mathbf{P}_{\mathfrak{c}} f$ as
    \begin{align*}
    \mathbf{P}_{\mathfrak{c}} f=\left\{a_f+b_f \cdot p+c_f \frac{p^0}{\mathfrak{c}}\right\} \sqrt{\mathbf{M}_{\mathfrak{c}}},
    \end{align*}
    and further denote $\left\{\mathbf{I}-\mathbf{P}_{\mathfrak{c}}\right\} f$ to be the microscopic part of $f$.
	
	Indeed, the orthonormal basis of $\mathcal{N}_{\mathfrak{c}}$ for the relativistic Boltzmann equation has the form
	\begin{align}\label{e2.36-0}
            \chi_0^{\mathfrak{c}}=\mathfrak{a}_0 \sqrt{\mathbf{M}_{\mathfrak{c}}}, \quad \chi_j^{\mathfrak{c}}=\frac{p_j-\mathfrak{a}_j}{\mathfrak{b}_j} \sqrt{\mathbf{M}_{\mathfrak{c}}}(j=1,2,3), \quad \chi_4^{\mathfrak{c}}=\frac{p^0 / \mathfrak{c}+\sum_{i=1}^3 \lambda_i p_i+\mathfrak{e}}{\zeta} \sqrt{\mathbf{M}_{\mathfrak{c}}},
        \end{align}
    where $\mathfrak{a}_0$, $\mathfrak{a}_j(j=1,2,3)$, $\mathfrak{b}_j(j=1,2,3)$, $\lambda_i(i=1,2,3)$, $\mathfrak{e}$ and $\zeta$ are functions of $\left(n_{\mathfrak{e}}, u_{\mathfrak{e}}, T_\mathfrak{e}\right)$. One can refer to the appendix of \cite{Wang-Xiao-JLMS-2026} for detailed expressions. Furthermore, a similar argument as in Lemma \ref{l2.12} yields 
    \begin{align}\label{e2.36}
	        \left|\partial_{\beta}^{\alpha} \chi_j^{\mathfrak{c}}(t, x, p)\right|\leq C(\mathfrak{e}_l)\left\langle p\right\rangle^{2|\alpha|+|\beta|+1} \mathbf{M}^{\frac{1}{2}}_{\mathfrak{c}}(t, x, p), \quad j=0,1,\cdots,4,
	    \end{align}
	where $C(\mathfrak{e}_l)$ depends only on $\mathfrak{e}_l$ and is independent of $\mathfrak{c}$.
	\begin{definition}[{Pseudo-inverse operator of $\mathbf{L}_{\mathfrak{c}}$}]
		The inverse operator
		\begin{align*}
			(\mathbf{L}_{\mathfrak{c}}|_{\mathcal{N}_{\mathfrak{c}}^\perp})^{-1}:\mathcal{N}_{\mathfrak{c}}^\perp\rightarrow \mathcal{N}_{\mathfrak{c}}^\perp,
		\end{align*}	
		is called the pseudo-inverse operator of the linearized relativistic Boltzmann collision operator $\mathbf{L}_{\mathfrak{c}}$, which
		is briefly denoted by $\mathbf{L}_{\mathfrak{c}}^{-1}$.
	\end{definition}
	\subsection{Some useful estimates in the center of momentum system}
	As shown in \cite{Wang-Xiao-JLMS-2026}, using expression \eqref{04}, one can obtain a series of uniform-in-$\mathfrak{c}$ estimates for the linear and nonlinear collision operators, which are as follows: 
	\begin{lemma}\emph{(\cite{Wang-Xiao-JLMS-2026})}\label{l2.1}
		Recall $w_\ell(p)$ in \eqref{e1.14}. Then it holds that 
		\item [(1)]
		$\displaystyle 
		 \frac{\sqrt{|p \times q|^2+\mathfrak{c}^2|p-q|^2}}{\sqrt{p^0 q^0}} \leq g \leq|p-q| \quad \text { and } \quad  g^2<\mathfrak{s} \leq 4 p^0 q^0 $.\\ 
		 \item [(2)]
		 $\displaystyle
		  v_\phi=\frac{\mathfrak{c}}{4} \frac{g \sqrt{\mathfrak{s}}}{p^0 q^0} \leq \min \left\{\mathfrak{c}, \frac{|p-q|}{2}\right\}$.\\
		\item [(3)]
		$\displaystyle
		 \mathbf{M}_{\mathfrak{c}}(t, x, p) \lesssim e^{-\frac{c_0}{2}|p|}, \quad \text{$c_0$ is a positive constant defined in \eqref{e.30}}.$\\
		 \item [(4)]
		 $\displaystyle
		 \int_{\mathbb{R}^3} \int_{\mathbb{S}^2} v_\phi \mathbf{M}_{\mathfrak{c}}^\alpha(q) d \omega d q \cong \nu_{\mathfrak{c}}(p), \quad \text { for } \alpha>0.$\\
		\item [(5)]
		$\displaystyle 
		\nu_{\mathfrak{c}}(p)\simeq\begin{cases}
			1+|p|,\ |p|\leq \mathfrak{c},\\
			\mathfrak{c},\ |p|\geq\mathfrak{c}.
		\end{cases}$ \\
		\item [(6)]
		$\displaystyle \int_{\R^3}k^2_{\mathfrak{c}w}(p,q)dq\lesssim
		\left\{
		\begin{aligned}
			&\frac{1}{1+|p|}, \ |p|\le \mathfrak{c},\\
			&\frac{1}{\mathfrak{c}}, \  |p|\ge \mathfrak{c}.
		\end{aligned}
		\right.$\\
		\item [(7)]
		$\displaystyle 
		\int_{\mathbb{R}^3}|k_{\mathfrak{c}w}(p,q)|e^{\frac{c_0}{16}|p-q|}dq\lesssim
		\left\{
		\begin{aligned}
			&\frac{1}{1+|p|}, \ |p|\le \mathfrak{c},\\
			&\frac{1}{\mathfrak{c}}, \  |p|\ge \mathfrak{c},
		\end{aligned}
		\right. \quad \text{$c_0$ is a positive constant defined in \eqref{e.30}}.$\\
		\item [(8)]
		$\displaystyle 
		\langle\mathbf{L}_{\mathfrak{c}} g,g\rangle\geq \zeta_0\|(\mathbf{I-P_{\mathfrak{c}}})g\|_{\nu_{\mathfrak{c}}}^2,\quad \zeta_0>0\, \text{and $\zeta_0$ is independent of $\mathfrak{c}$} $.\\
		\item [(9)]
		$\displaystyle 
		w_\ell(p)\lesssim w_\ell(p^{\prime})w_\ell(q^{\prime}), \quad w_\ell(q)\lesssim w_\ell(p^{\prime})w_\ell(q^{\prime})$, \quad for $\ell\geq1$.
	\end{lemma}
	
	Furthermore, we provide some commonly used estimates.
	\begin{lemma}\label{l2.05}
		Recall that $\left(p,q\right)$, $\left(p^{\prime},q^{\prime}\right)$ are the pre-post collisional momentum pairs which satisfy the principles of conservation of energy and momentum. Then it holds that
		\begin{align}\label{bc.1}
			\nu_{\mathfrak{c}}(p)\lesssim\nu_{\mathfrak{c}}(p^{\prime})\nu_{\mathfrak{c}}(q^{\prime}),\quad\nu_{\mathfrak{c}}(q)\lesssim\nu_{\mathfrak{c}}(p^{\prime})\nu_{\mathfrak{c}}(q^{\prime}).
		\end{align}
	\end{lemma}
	\begin{proof}
		We divide the proof into two cases.
		
		{\it Case 1:} $|p^{\prime}|\geq \mathfrak{c}$ or $|q^{\prime}|\geq \mathfrak{c}$. It holds that $\nu_{\mathfrak{c}}(p^{\prime})\nu_{\mathfrak{c}}(q^{\prime})\geq \mathfrak{c}$.
		
		For $|p|\geq \mathfrak{c}$, it is easy to see that $\nu_{\mathfrak{c}}(p)=\mathfrak{c}\leq \nu_{\mathfrak{c}}(p^{\prime})\nu_{\mathfrak{c}}(q^{\prime})$;
		
		For $|p|\leq \mathfrak{c}$, there holds
		\begin{align*}
			\nu_{\mathfrak{c}}(p)\lesssim1+|p|\leq 2\mathfrak{c} \leq2\nu_{\mathfrak{c}}(p^{\prime})\nu_{\mathfrak{c}}(q^{\prime}).
		\end{align*}
		
		{\it Case 2:} $|p^{\prime}|\leq \mathfrak{c}$ and $|q^{\prime}|\leq \mathfrak{c}$. It holds that $\nu_{\mathfrak{c}}(p^{\prime})\nu_{\mathfrak{c}}(q^{\prime})= \left(1+\left|p^{\prime}\right|\right)\left(1+\left|q^{\prime}\right|\right)$.
		
		For $|p|\leq \mathfrak{c}$, it's clear that by the ninth inequality in Lemma \ref{l2.1} for $\ell=1$
		\begin{align*}
			\nu_{\mathfrak{c}}(p)\lesssim1+|p|\lesssim\left(1+\left|p^{\prime}\right|\right)\left(1+\left|q^{\prime}\right|\right)= \nu_{\mathfrak{c}}(p^{\prime})\nu_{\mathfrak{c}}(q^{\prime});
		\end{align*}

	 	For $|p|\geq \mathfrak{c}$, we have
	 	\begin{align*}
	 		\nu_{\mathfrak{c}}(p)\lesssim\mathfrak{c}\leq 1+|p| \lesssim\left(1+\left|p^{\prime}\right|\right)\left(1+\left|q^{\prime}\right|\right)= \nu_{\mathfrak{c}}(p^{\prime})\nu_{\mathfrak{c}}(q^{\prime}).
	 	\end{align*}
	 	 Thus we complete the proof of \eqref{bc.1} for $\nu_{\mathfrak{c}}(p)$. The estimate \eqref{bc.1} for $\nu_{\mathfrak{c}}(q)$ can be obtained by the same way. Then we complete the proof of Lemma \ref{l2.05}.
	\end{proof}
	\begin{lemma}\emph{(\cite{Wang-Xiao-JLMS-2026})}\label{l2.2}
		It holds that
		\begin{align}\label{e2.4}
			\left|\frac{w_{\ell}}{\sqrt{\mathbf{M}_{\mathfrak{c}}}} Q_{\mathfrak{c}}\left(f_1 \sqrt{\mathbf{M}_{\mathfrak{c}}}, f_2 \sqrt{\mathbf{M}_{\mathfrak{c}}}\right)\right| \lesssim \nu_{\mathfrak{c}}(p)\left\|f_1\right\|_{\infty, \ell}\left\|f_2\right\|_{\infty, \ell},
		\end{align}
		and 
		\begin{align}\label{e2.4-1}
			\left|\frac{w_{\ell}}{\sqrt{J_{\mathbf{M}}}} Q_{\mathfrak{c}}\left(f_1 \sqrt{J_{\mathbf{M}}}, f_2 \sqrt{J_{\mathbf{M}}}\right)\right| \lesssim \nu_{\mathfrak{c}}(p)\left\|f_1\right\|_{\infty, \ell}\left\|f_2\right\|_{\infty, \ell},
		\end{align}
		where the constants are independent of $\mathfrak{c}$.
	\end{lemma}
	\begin{lemma}\label{l2.3}\emph{(\cite{Speck-CMP-2011})}
		For any $\ell \geq 9$, it holds that
		\begin{equation}\label{e2.5}
			\left|\left\langle\Gamma_{\mathfrak{c}}\left(f_1, f_2\right), f_3\right\rangle\right| \lesssim\left\|f_3\right\|_{\infty, \ell}\left\|f_2\right\|\left\|f_1\right\| .
		\end{equation}
		Furthermore, if $\chi(p)$ satisfies $|\chi(p)| \lesssim e^{-\delta_1|p|}$ for some positive constant $\delta_1>0$, then we have
		\begin{equation}\label{e2.6}
			\left|\left\langle\Gamma_{\mathfrak{c}}\left(f_1, \chi\right), f_3\right\rangle\right|+\left|\left\langle\Gamma_{\mathfrak{c}}\left(\chi, f_1\right), f_3\right\rangle\right| \lesssim\left\|f_3\right\|_{\nu_{\mathfrak{c}}}\left\|f_1\right\|_{\nu_{\mathfrak{c}}},
		\end{equation}
		where the constants are independent of $\mathfrak{c}$.
	\end{lemma}
	\begin{remark}
    We point out that Lemma~\ref{l2.3} remains valid when the bracket 
    $\langle\cdot,\cdot\rangle$ is understood as either the 
    $L^2(\mathbb{R}^3_p)$ inner product or the 
    $L^2(\mathbb{R}^3_x\times\mathbb{R}^3_p)$ inner product. 
    Although \cite[Lemma 2.3]{Wang-Xiao-JLMS-2026} only states and proves the 
    estimate in the $L^2(\mathbb{R}^3_p)$ setting, the same conclusion in the 
    $L^2(\mathbb{R}^3_x\times\mathbb{R}^3_p)$ framework follows by applying exactly 
    the same argument pointwise in the spatial variable $x$, and then integrating 
    over $x$.
    \end{remark}
	
	\subsection{Glassey--Strauss expression}
	Glassey and Strauss illustrated in \cite{Glassey-PRIMS-1993} that a reduction of the collision integrals can be performed, without using Lorentz transformations, to obtain
	\begin{align}\label{01}
		Q_{\mathfrak{c}}(f_1, f_2)=\iint_{\mathbb{R}^3\times \mathbb{S}^2}\f{\mathfrak{s}\mathfrak{B}(p,q,\omega)}{p^0q^0}[f_1(p^{\prime\prime})f_2(q^{\prime\prime})-f_1(p)f_2(q)]d\omega d q
	\end{align}
	with
	\begin{align}\label{01-0}
		\mathfrak{B}(p,q,\omega):=\mathfrak{c}\frac{(p^0+q^0)^2p^0q^0|\omega\cdot(\f{p}{p^0}-\f{q}{q^0})|}{[(p^0+q^0)^2-(\omega\cdot(p+q))^2]^2}.
	\end{align} 
	
	In this expression, the post-collision momentum in \eqref{01} is given by
	\begin{align}\label{02}
		p^{\prime\prime}=p+a(p,q,\omega)\omega,\quad q^{\prime\prime}=q-a(p,q,\omega)\omega,
	\end{align}
	where
	\begin{align*}
		a(p,q,\omega):=\f{2p^0q^0(p^0+q^0)\{\omega\cdot(\f{q}{q^0}-\f{p}{p^0})\}}{(p^0+q^0)^2-\{\omega\cdot(p+q)\}^2}.
	\end{align*} 
		
	\begin{lemma}\cite[Corollary 1.5]{Strain-KRM-2011}\label{l-new}
	Recall $v_\phi(p,q)$ in \eqref{04}. For any suitable integrable function $G: \mathbb{R}^3\times \mathbb{R}^3\times \mathbb{R}^3\times \mathbb{R}^3\rightarrow \mathbb{R}$, it holds that
	\begin{align}\label{07}
		\int_{\mathbb{S}^{2}} v_\phi G\left(p, q, p^{\prime}, q^{\prime}\right)d \omega=\int_{\mathbb{S}^2}  \f{\mathfrak{s}\mathfrak{B}(p,q,\omega)}{p^0q^0} G\left(p, q, p^{\prime\prime}, q^{\prime\prime}\right)d \omega,
	\end{align}
	where the post-collision momentum on the left are defined as in \eqref{05} and the similar expressions on the right are the ones given  in \eqref{02}.
	\end{lemma}

	In the \emph{center of momentum system}, it is difficult to compute the Jacobian determinant
\[
\frac{\partial(p',q')}{\partial(p,q)},
\]
since the post-collision momenta $(p',q')$ are given by the complicated expressions in \eqref{05}. For more details on this issue, we refer to \cite{Chapman-CMP-2021}. On the other hand, for the Glassey--Strauss representation \eqref{01}, it was proved in \cite{Glassey-TTSP-1991} that
\begin{align}\label{03}
\frac{\partial (p^{\prime\prime},q^{\prime\prime})}{\partial (p,q)}
=-\frac{p^{\prime\prime0}q^{\prime\prime0}}{p^0q^0},
\end{align}
which is useful for changes of variables in integrals involving both $(p,q)$ and $(p'',q'')$. Therefore, although the Jacobian for the \emph{center of momentum system} parametrization is not explicitly available, combining Lemma \ref{l-new} with the Glassey--Strauss change of variables yields the following pre-post collisional symmetry.
\begin{corollary}\label{c-new}
For any suitable integrable function
$G:\mathbb{R}^3\times\mathbb{R}^3\times\mathbb{R}^3\times\mathbb{R}^3\to\mathbb{R}$,
it holds that
\begin{equation}\label{08}
\int_{\mathbb{R}^3}dp\int_{\mathbb{R}^3}dq\int_{\mathbb{S}^2}d\omega\,
v_\phi\, G(p,q,p',q')
=
\int_{\mathbb{R}^3}dp\int_{\mathbb{R}^3}dq\int_{\mathbb{S}^2}d\omega\,
v_\phi\, G(p',q',p,q).
\end{equation}
\end{corollary}
	\subsection{High-order momentum derivative estimates for the collision operator}
	Due to the presence of the Lorentz force term $E+\frac{p}{p^0}\times B$,
    one must study momentum regularity for rVMB. A basic difficulty is that neither of the two standard representations of the relativistic Boltzmann collision operator is suitable for taking high momentum derivatives over the entire integration region.

In the \emph{center of momentum system} \eqref{04}-\eqref{05}, high $p$-derivatives of the post-collisional momenta $p'$ and $q'$ become singular near two distinct degeneracy sets. On the one hand, derivatives of the relative momentum $g$ degenerate near $p=q$. On the other hand, derivatives of the angular term
\[
(p+q)\frac{(p+q)\cdot \omega}{|p+q|^2},
\]
degenerate near $p+q=0$, that is, near $p=-q$. Thus, the \emph{center of momentum system} variables are not suitable for global high-order momentum estimates.

By contrast, in the Glassey--Strauss representation \eqref{01}-\eqref{02}, momentum derivatives of $p''$ and $q''$ do not exhibit the above singular behavior, but instead generate unfavorable polynomial growth in $p$. Such growth is too strong to be controlled directly in a global energy estimate.

Following \cite{Guo-CMP-2012}, we overcome these two difficulties by splitting the integration region into
\begin{align}\label{A.1}
A:=\{|p|\le \mathfrak{c}\}\cup\{|p|\ge \mathfrak{c},\ |p|\le 2q^0\},
\qquad
A^c:=\{|p|> 2q^0\}.
\end{align}
On the set $A$, we employ the Glassey--Strauss representation \eqref{01}. Since $|p|\lesssim q^0$ on $A$, the polynomial growth in $p$ can be absorbed by the factor $\mathbf M_{\mathfrak{c}}^{1/2}(q)$. On the complementary region $A^c$, we use the \emph{center of momentum system} \eqref{04}. Indeed, if $|p|\ge 2q^0$, then
\[
|p-q|\ge |p|-|q|\ge \frac{|p|}{2},
\qquad
|p+q|\ge |p|-|q|\ge \frac{|p|}{2},
\]
so this region stays uniformly away from both degeneracy sets $p=q$ and $p=-q$. Hence the singular behavior of high momentum derivatives in the \emph{center of momentum system} variables is avoided on $A^c$.
	
	Based on the region-dependent frame decomposition above, we establish the following estimates:
	\begin{lemma}\label{l2.4}
		For any $\ell\geq 9$ and $|\alpha|+|\beta|>0$, we have the following estimate for the collision operator
		\begin{equation}\label{e2.7}
			\left|\left\langle\partial_{\beta}^{\alpha} \Gamma_{\mathfrak{c}}\left(f_1, f_2\right), \partial_{\beta}^{\alpha} f_3\right\rangle\right| \lesssim\left\|\partial_{\beta}^{\alpha} f_3\right\|_{\infty, \ell}\sum_{\substack{\alpha_1+\alpha_2 \leq \alpha\\\beta_1+\beta_2 \leq \beta}}\big\|\partial_{\beta_1}^{\alpha_1} f_1\big\|\big\|\partial_{\beta_2}^{\alpha_2} f_2\big\| ,
		\end{equation}
		where $\partial_{\beta}^{\alpha} \stackrel{\text { def }}{=} \partial_t^{\alpha^0} \partial_{x_1}^{\alpha^1} \partial_{x_2}^{\alpha^2} \partial_{x_3}^{\alpha^3} \partial_{p_1}^{\beta^1} \partial_{p_2}^{\beta^2} \partial_{p_3}^{\beta^3}$, $\alpha=\left(\alpha^0, \alpha^1, \alpha^2, \alpha^3\right)$, $\beta=\left(\beta^1, \beta^2, \beta^3\right)$. If  $\chi$ is any rapidly decaying function$(|\chi(p)| \leq C\left(p^0\right)^{-m}$, with $m>3)$, then we have
		\begin{equation}\label{e2.8}
			\left|\left\langle\partial_{\beta}^{\alpha} \Gamma_{\mathfrak{c}}\left(f_1, \chi\right), \partial_{\beta}^{\alpha} f_3\right\rangle \right|+\left|\left\langle\partial_{\beta}^{\alpha} \Gamma_{\mathfrak{c}}\left(\chi, f_1\right),\partial_{\beta}^{\alpha} f_3\right\rangle \right| \lesssim\sum_{\substack{\alpha_1\leq \alpha\\\beta_1\leq \beta}}\big\|\partial_{\beta_1}^{\alpha_1} f_1\big\|_{\nu_{\mathfrak{c}}}\big\|\partial_{\beta}^{\alpha} f_3\big\|_{\nu_{\mathfrak{c}}},
		\end{equation}
		where the constants are independent of $\mathfrak{c}$.
	\end{lemma}
\begin{proof}
We only prove the case $\alpha=0$. In fact, the space-time derivatives $\partial^\alpha$
do not hit the collision kernel in $\Gamma_{\mathfrak{c}}$, but only the local Maxwellian and the distribution functions
appearing in the integrand. Therefore, they do not generate any new singularity, and the
general case follows immediately from Leibniz's rule together with the same argument as in
\cite{Speck-CMP-2011}.

For pure momentum derivatives
$\partial_\beta$ on $\Gamma_{\mathfrak{c}}$,
we use the decomposition \eqref{A.1} and treat the two regions separately.
On $A$, we apply the Glassey--Strauss representation \eqref{01}. The differentiated kernel estimates in this region are provided by Lemma \ref{lA.1}. On $A^c$, we use the {\it center of momentum} representation \eqref{04}. Since this region stays uniformly away from the degeneracy sets $p=q$ and $p=-q$, the corresponding estimates follow from Lemma \ref{lA.3}.

Starting from these kernel bounds, and following the argument of \cite[Lemma 3.12]{Speck-CMP-2011}, we derive Lemmas \ref{lA.2} and \ref{lA.5} on $A$ and $A^c$, respectively. Combining the contributions from the two regions then yields \eqref{e2.7}--\eqref{e2.8}. We omit the details and refer to Appendix \ref{Appendix A}.
\end{proof}
	\begin{remark}
The analogue of Lemma \ref{lA.1} and Lemma \ref{lA.3} in the normalized case $\mathfrak{c}=1$ go back to \cite{Guo-CMP-2012}. In the normalized case $\mathfrak{c}=1$, one may use the lower bound $g\gtrsim 1$ after staying away from the degeneracy sets. In the general-speed-of-light setting, however, this becomes the scale-dependent bound $g\gtrsim \mathfrak{c}$, and all high-order momentum derivative estimates must keep track of the resulting $\mathfrak{c}$-dependence. 

More precisely, on the region $A$ we rework the Glassey--Strauss representation so that the polynomial growth generated by momentum derivatives is measured in a way compatible with the $\mathfrak{c}$-scale and can still be absorbed by the relativistic Maxwellian. On the region $A^c$, we revisit the {\it center of momentum} representation and establish high-order derivative bounds for $g$, $\sqrt{\mathfrak s}$, $v_\phi$, $p'$, and $q'$ uniformly in $\mathfrak{c}$. These refined estimates are proved in Appendix \ref{Appendix A} and form the key input in the proof of Lemma \ref{l2.4}.
\end{remark}

    For later use, we also need the following pointwise estimate, see the proof in Appendix \ref{appendix A.3.}.
	\begin{lemma}\label{l2.5}
For any fixed $0<\lambda<1$, let $\alpha\in\mathbb N^4$ and $\beta\in\mathbb N^3$.

Suppose that for $i=1,2$,
\[
\sum_{\bar\alpha\le \alpha, \bar\beta\le \beta}
\left|\partial_{\bar\beta}^{\bar\alpha}f_i(p)\right|
\lesssim \mathbf M_{\mathfrak{c}}^{\frac{\lambda}{2}}(p).
\]
Then
\begin{align}\label{A.16}
\left|\partial_\beta^\alpha \Gamma_{\mathfrak{c}}(f_1,f_2)(p)\right|
\lesssim
\nu_{\mathfrak{c}}(p)\mathbf M_{\mathfrak{c}}^{\frac{\lambda}{2}}(p).
\end{align}

Furthermore, suppose that
\[
\sum_{\substack{\bar\alpha\le \alpha,\ \bar\beta\le \beta\\
|\bar\alpha|+|\bar\beta|<|\alpha|+|\beta|}}
\left|\partial_{\bar\beta}^{\bar\alpha}f(p)\right|
\lesssim \mathbf M_{\mathfrak{c}}^{\frac{\lambda}{2}}(p).
\]
Then
\begin{align}\label{19}
\left|[\partial_\beta^\alpha,\mathbf L_{\mathfrak{c}}]f(p)\right|
\lesssim
\nu_{\mathfrak{c}}(p)\mathbf M_{\mathfrak{c}}^{\frac{\lambda}{2}}(p).
\end{align}
Moreover,
\begin{align}\label{20}
\left|\partial_{\beta}\nu_{\mathfrak{c}}(p)\right|
\lesssim \nu_{\mathfrak{c}}(p),
\end{align}
where the implicit constant is independent of $\mathfrak{c}$.
\end{lemma}
	\subsection{Estimates on the operator $\mathcal{K}_{\mathfrak{c}}$ and its derivatives}
	Recall $J_{\mathbf{M}}(p)$ in \eqref{e1.15}, we define
	$$
	\mathcal{L}_{\mathfrak{c}} f:=-\frac{1}{\sqrt{J_{\mathbf{M}}}}\left\{Q_{\mathfrak{c}}\left(\mathbf{M}_{\mathfrak{c}}, \sqrt{J_{\mathbf{M}}} f\right)+Q_{\mathfrak{c}}\left(\sqrt{J_{\mathbf{M}}}f,\mathbf{M}_{\mathfrak{c}}\right)\right\}=\nu_{\mathfrak{c}} f-\mathcal{K}_{\mathfrak{c}} f,
	$$
	where $\mathcal{K}_{\mathfrak{c}}=\mathcal{K}_{\mathfrak{c}_2}-\mathcal{K}_{\mathfrak{c}_1}$. More specifically, $\nu_{\mathfrak{c}}$ is defined in \eqref{e2.1} and operators $\mathcal{K}_{\mathfrak{c}_1} f$ and $\mathcal{K}_{\mathfrak{c}_2} f$ are defined as
	\begin{align*}
		\mathcal{K}_{\mathfrak{c}_1} f & :=\frac{1}{\sqrt{J_{\mathbf{M}}}} Q_{\mathfrak{c}}^{-}\left(\mathbf{M}_{\mathfrak{c}}, \sqrt{J_{\mathbf{M}}} f\right)=\int_{\mathbb{R}^3} \int_{\mathbb{S}^2} v_\phi\left\{\sqrt{J_{\mathbf{M}}(q)}\frac{\mathbf{M}_{\mathfrak{c}}(p)}{\sqrt{J_{\mathbf{M}}(p)}} f(q)\right\} d \omega d q, \\
		\mathcal{K}_{\mathfrak{c}_2} f & :=\frac{1}{\sqrt{J_{\mathbf{M}}}}\left\{Q_{\mathfrak{c}}^{+}\left(\mathbf{M}_{\mathfrak{c}},\sqrt{J_{\mathbf{M}}}f\right)+Q_{\mathfrak{c}}^{+}\left(\sqrt{J_{\mathbf{M}}}f, \mathbf{M}_{\mathfrak{c}}\right)\right\} \\
		& =\int_{\mathbb{R}^3} \int_{\mathbb{S}^2} v_\phi\left\{\mathbf{M}_{\mathfrak{c}}\left(p^{\prime}\right) \frac{\sqrt{J_{\mathbf{M}}(q^{\prime})}}{\sqrt{J_{\mathbf{M}}(p)}} f\left(q^{\prime}\right)+\mathbf{M}_{\mathfrak{c}}\left(q^{\prime}\right) \frac{\sqrt{J_{\mathbf{M}}(p^{\prime})}}{\sqrt{J_{\mathbf{M}}(p)}} f\left(p^{\prime}\right)\right\} d \omega d q.
	\end{align*}
	Noting \eqref{e1.15}, by similar arguments as in \cite{Strain-CMP-2010}, one can show that
    \begin{align}\label{e2.42}
    \left|\mathcal{K}_{\mathfrak{c}_i}(f)\right| \lesssim \int_{\mathbb{R}^3} k_i(p, q)|f(q)| d q, \quad i=1,2
    \end{align}
    where
    \begin{align}\label{e2.43}
    k_1(p, q)=|p-q| e^{-\frac{\delta}{T_M}|p|} e^{-\frac{\delta}{T_M}|q|} \quad \text{and} \quad k_2(p, q)=\frac{1}{|p-q|} e^{-\frac{\delta}{2 T_M}|p-q|},
    \end{align}
    with $\delta:=\alpha-\frac{1}{2}>0$. We denote $k(p, q):=k_1(p, q)+k_2(p, q)$, then it holds that
    $$
    \left|\mathcal{K}_{\mathfrak{c}}(f)\right| \lesssim \int_{\mathbb{R}^3} k(p, q)|f(q)| d q, \quad   i=1,2 .
    $$
    Denote
    $$
    k_w(p, q):=k(p, q) \frac{w_{\ell}(p)}{w_{\ell}(q)} .
    $$
    By similar arguments as in \cite[Lemmas 4.4-4.5]{Wang-Xiao-JLMS-2026}, one has
    \begin{align}\label{e2.39}
    \int_{\mathbb{R}^3} k_w(p, q) e^{\frac{\delta}{4T_M}|p-q|} d q+\int_{\mathbb{R}^3} k_w^2(p, q) d q \lesssim \max \left\{\frac{1}{\mathfrak{c}}, \frac{1}{1+|p|}\right\} .
    \end{align}
    
     Next, we provide several estimates for the derivatives of $\mathcal{K}_{\mathfrak{c}}f$ with respect to $x$,$p$. 
	\begin{lemma}\label{l2.13}
	    For $j=1,2,3$, there hold
	    \begin{align}\label{e2.44}
	    \left|\partial_{x_j}\mathcal{K}_{\mathfrak{c}_i}(f)\right| \lesssim \int_{\mathbb{R}^3} \hat{k}_i(p, q)|f(q)|d q+\int_{\mathbb{R}^3} k_i(p, q)|\partial_{x_j}f(q)|d q, \quad i=1,2,
	    \end{align}
        and
        \begin{align}
        &\left|\partial_{p_j}\mathcal{K}_{\mathfrak{c}_1}(f)\right| \lesssim \int_{\mathbb{R}^3} \hat{k}_1(p, q)|f(q)|d q,\label{e2.40-1}\\
        &\left|\partial_{p_j}\mathcal{K}_{\mathfrak{c}_2}(f)\right| \lesssim \int_{\mathbb{R}^3} \hat{k}_2(p, q)|f(q)|d q+\int_{\mathbb{R}^3} k_2(p, q)|\partial_{p_j}f(q)|d q,\label{e2.40-2}
        \end{align}
	    where $k_i(p, q)(i=1,2)$ are defined in \eqref{e2.43} and
	    \begin{gather*}
	        \hat{k}_1(p, q):=|p-q| e^{-\frac{\delta}{2T_M}|p|} e^{-\frac{\delta}{2T_M}|q|}, \quad \hat{k}_2(p, q):=\frac{1}{|p-q|} e^{-\frac{\delta}{4T_M}|p-q|}
	    \end{gather*}
	    with $\delta=\alpha-\frac{1}{2}>0$.
	\end{lemma}
	\begin{proof}
	The proof is deferred to appendix \ref{appendix A.4} for brevity.
	\end{proof}
	Following the same procedure as before, we denote $\hat{k}(p, q):=\hat{k}_1(p, q)+\hat{k}_2(p, q)$, then there holds
	\begin{align}\label{e2.41}
	\left|\nabla_{x,p}\mathcal{K}_{\mathfrak{c}}(f)\right| \lesssim \int_{\mathbb{R}^3} \hat{k}(p, q)|f(q)| d q+\int_{\mathbb{R}^3} k(p, q)|\nabla_{x,p}f(q)|d q.
	\end{align}
	Denote 
	$$
	\hat{k}_w(p, q):=\hat{k}(p, q)\frac{w_\ell(p)}{w_\ell(q)}.
	$$
	By similar arguments as in \cite[Lemmas 4.4-4.5]{Wang-Xiao-JLMS-2026}, one has
	\begin{equation}\label{e2.41-1}
		\int_{\mathbb{R}^3} \hat{k}_w(p, q) e^{\frac{\delta}{8T_M}|p-q|} d q+\int_{\mathbb{R}^3} \hat{k}_w^2(p, q) d q \lesssim \max \left\{\frac{1}{\mathfrak{c}}, \frac{1}{1+|p|}\right\}.
	\end{equation}
	
\section{The Newtonian Limit for the Cauchy Problem of rEM}\label{section 3}

In this section, we study the Cauchy problem for the relativistic Euler-Maxwell system in the whole space $\mathbb{R}^3$ and justify its Newtonian limit as $\mathfrak{c}\to\infty$. 

 We first derive a formal asymptotic expansion in powers of $\mathfrak{c}^{-1}$, where the leading-order system is the classical isentropic compressible Euler-Poisson system and the first-order correction is determined by a linearized Euler-Poisson system together with a curl-div system for the magnetic field. We then rigorously justify this expansion and establish uniform estimates for the remainder terms.

\subsection{Formal asymptotic expansions}
	As a starting point, we assume that the initial data defined in \eqref{e.4-1} admit the following asymptotic expansion with respect to the speed of light $\mathfrak{c}$:
	$$
	(n_{\mathfrak{e}}, u_{\mathfrak{e}}, T_\mathfrak{e}, E_{\mathfrak{e}}, B_{\mathfrak{e}})(0,x)=\sum_{j=0}^{m} \mathfrak{c}^{-j}\left(n_{j}, u_{j}, T_{j}, E_{j}, B_{j}\right)(0,x)+O\left(\mathfrak{c}^{-m}\right).
    $$
	Then we consider the following expansion for solutions of \eqref{e.4}:
	\begin{equation}\label{e.6}
		\left(n_{\mathfrak{e}}, u_{\mathfrak{e}}, T_\mathfrak{e}, E_{\mathfrak{e}}, B_{\mathfrak{e}}\right)=\sum_{j \geq 0} \mathfrak{c}^{-j}\left(n_j, u_j, T_j, E_j, B_j\right).
	\end{equation}
	
	Under the isentropic assumption \eqref{e.1}, \(S_{\mathfrak e}\) is a constant, so \(T_{\mathfrak e}\) is uniquely determined by \(n_{\mathfrak e}\).
	Inserting \eqref{e.6} into \eqref{e.7}, we obtain
	\begin{align*}
	 n_0+\frac{1}{\mathfrak{c}}n_1+O\left(\mathfrak{c}^{-2}\right)=e^{\frac{5}{2}-S_\mathfrak{e}}(2 \pi)^{\frac{3}{2}}\Big[T_0^{\frac{3}{2}}+\frac{1}{\mathfrak{c}} \frac{3}{2} T_0^{\frac{1}{2}} T_1+O\left(\mathfrak{c}^{-2}\right)\Big].
	\end{align*}
	Comparing the order of $\mathfrak{c}$, one obtains
	\begin{align}\label{e.8}
        T_0=\frac{1}{2 \pi} e^{\frac{2}{3} S_\mathfrak{e}-\frac{5}{3}} n_0^{\frac{2}{3}}, \quad T_1=\frac{1}{3 \pi} e^{\frac{2}{3} S_\mathfrak{e}-\frac{5}{3}} n_0^{-\frac{1}{3}} n_1=T_0^{\prime}\left(n_0\right) n_1 .
	\end{align}
	
	For later use, we now address the two complex terms in \eqref{e.4}. Using \eqref{2.4}, we have from \eqref{e0.3} and \eqref{e.1} that
	\begin{align}\label{e.8-1}
	    \frac{h}{\mathfrak{c}^2}=\frac{K_3(\gamma)}{K_2(\gamma)}&=1+\frac{5}{2 \gamma}+\frac{15}{8 \gamma^2}-\frac{15}{8 \gamma^3}+O\left(\gamma^{-4}\right)\nonumber\\
	    &=1+\frac{5}{2}\frac{T_\mathfrak{e}}{\mathfrak{c}^2}+\frac{15}{8}\frac{T_\mathfrak{e}^2}{\mathfrak{c}^4}-\frac{15}{8}\frac{T_\mathfrak{e}^3}{\mathfrak{c}^6}+O\Big(\frac{T_\mathfrak{e}^4}{\mathfrak{c}^8}\Big).
	\end{align}
	By Taylor's theorem, it's clear that
	\begin{align}\label{e.8-2}
	\frac{u_{\mathfrak{e}}^0}{\mathfrak{c}}
	=1+\frac{1}{2}\left(\frac{u_{\mathfrak{e}}}{\mathfrak{c}}\right)^2-\frac{1}{8}\left(\frac{u_{\mathfrak{e}}}{\mathfrak{c}}\right)^4+\cdots+(-1)^{n-1} \frac{(2 n-3)!!}{(2 n)!!}\left(\frac{u_{\mathfrak{e}}}{\mathfrak{c}}\right)^{2 n}+ O\Big(\left(\frac{u_{\mathfrak{e}}}{\mathfrak{c}}\right)^{2 n+2}\Big).
	\end{align}
    
	Substituting \eqref{e.6}-\eqref{e.8-2} into \eqref{e.4}, and comparing the coefficients to each order of $\mathfrak{c}$, we obtain
	
	(1) The $\mathfrak{c}^0$-th order equations:
	\begin{equation}\label{e.9}
		\left\{\begin{aligned}
			&\partial_t n_0+\operatorname{div}\left(n_0 u_0\right)=0, \\
			&\partial_t u_0+\left(u_0 \cdot \nabla\right) u_0+\frac{5}{2}\nabla T_0(n_0)=-E_0, \\
			&\operatorname{div} E_0=4\pi\left(\bar{n}-n_0\right), \quad \nabla \times E_0=0, \\
			&\operatorname{div} B_0=0, \quad \nabla \times B_0=0, \\
			&\left.\left(n_0, u_0\right)\right|_{t=0}=\left(n_0(0), u_0(0)\right).
		\end{aligned}\right.
	\end{equation}
	
    To uniquely determine the electromagnetic fields, we impose the following far-field boundary conditions:
\begin{align}\label{e.9-1}
    E_0(t,x) \to 0, \quad B_0(t,x) \to 0 \quad \text{as } |x| \to \infty.
\end{align}
Since $\operatorname{div} B_0 = 0$ and $\nabla \times B_0 = 0$, the boundary condition \eqref{e.9-1} uniquely implies $B_0 \equiv 0$ in $\mathbb{R}^3$. Furthermore, the condition $\nabla \times E_0 = 0$ guarantees the existence of a scalar potential $\phi$ satisfying $E_0 = \nabla \phi$. Consequently, in view of \eqref{e.8}, the system \eqref{e.9} is exactly the classical isentropic compressible EP \eqref{e0.5}.

For smooth, irrotational initial data that are small perturbations of the constant equilibrium state $(\bar{n},0)$, Guo \cite{Guo-CMP-1998} established the global existence of smooth solutions.

\begin{lemma}[\cite{Guo-CMP-1998}]\label{l3.2}
Let $\varrho(x) \in C_c^\infty(\mathbb{R}^3)$ and $\upsilon(x) \in C_c^\infty(\mathbb{R}^3;\mathbb{R}^3)$ satisfy
\begin{equation}\label{e.9-2}
    \nabla \times \upsilon(x) = 0,
    \qquad
    \int_{\mathbb{R}^3} \varrho(x)\,dx = 0.
\end{equation}
Assume that the initial data are given by the compactly supported small perturbation
\[
    n_0(0,x)=\bar n+\epsilon \varrho(x),
    \qquad
    u_0(0,x)=\epsilon \upsilon(x).
\]
Then there exists $\epsilon_*>0$ such that, if $0<\epsilon<\epsilon_*$, the EP
\eqref{e.9} admits a unique global smooth solution $(n_0,u_0,\nabla \phi)$ on $[0,\infty)$.

Furthermore, for any $1<p<\frac32$ and every sufficiently large integer $s>0$, one has
\begin{equation}\label{e.9-3}
    \sup_{t\ge0}\|(n_0-\bar n,u_0,\nabla \phi)(t)\|_{H^s}
    +\sup_{t\ge0}\big[(1+t)^p\|(n_0-\bar n,u_0,\nabla \phi)(t)\|_{W^{s,\infty}}\big]
    \le C\epsilon,
\end{equation}
where $C>0$ is independent of $\epsilon$.
\end{lemma}
    
(2) The $\mathfrak{c}^{-1}$-th order equations:
	\begin{align}\label{e.10}
		\left\{\begin{aligned}
			&\partial_t n_1+\operatorname{div}\left(n_0 u_1+n_1 u_0\right)=0, \\
			&\partial_t u_1+\left(u_0 \cdot \nabla\right) u_1+\left(u_1 \cdot \nabla\right) u_0+\frac{5}{2}\nabla\left(T_0^{\prime}\left(n_0\right) n_1\right)=-E_1-u_0 \times B_0, \\
			&\operatorname{div} E_1=-4\pi n_1, \quad \nabla \times E_1=-\partial_t B_0, \\
			&\operatorname{div} B_1=0, \quad \nabla \times B_1=\partial_t E_0-4\pi n_0 u_0, \\
			&\left.\left(n_1, u_1\right)\right|_{t=0}=\left(n_1(0), u_1(0)\right).
		\end{aligned}\right.
	\end{align}
	
	Similarly, to uniquely determine the electromagnetic fields, we impose the following far-field boundary conditions:
    \begin{align}\label{e.10-0}
    E_1(t,x) \to 0, \quad B_1(t,x) \to 0 \quad \text{as } |x| \to \infty.
    \end{align}
	Since $B_0 \equiv 0$, the subsystem for $(n_1,u_1,E_1)$ decouples into the following linearized EP:
	\begin{align}\label{e.10-1}
		\left\{\begin{aligned}
			&\partial_t n_1+\operatorname{div}\left(n_0 u_1+n_1 u_0\right)=0, \\
			&\partial_t u_1+\left(u_0 \cdot \nabla\right) u_1+\left(u_1 \cdot \nabla\right) u_0+\frac{5}{2}\nabla\left(T_0^{\prime}\left(n_0\right) n_1\right)=-E_1, \\
			&\operatorname{div} E_1=-4\pi n_1, \quad \nabla \times E_1=0, \\
			&\left.\left(n_1, u_1\right)\right|_{t=0}=\left(n_1(0), u_1(0)\right),
		\end{aligned}\right.
	\end{align}
	while $B_1$ satisfies the linear curl-div system:
	\begin{align}\label{e.10-2}
	    \operatorname{div} B_1=0, \quad \nabla \times B_1=\partial_t E_0-4\pi n_0 u_0.
	\end{align}
	Relying on the existence theory for linear symmetrizable hyperbolic systems (see \cite{Friedrichs-CPAM-1954,Dafermos-2016}) and properties of curl-div systems in $\mathbb{R}^3$, we obtain the following lemma:
	
	\begin{lemma}\label{l3.3}
    Let $s$ be the fixed integer from \eqref{e.9-3}. Suppose that $(n_1, u_1)(0,x) \in H^{s-1}(\mathbb{R}^3)$ and $E_1(0,x)$ is determined by the compatibility conditions $\eqref{e.10-1}_{3}$. Then the linear system \eqref{e.10} admits a unique global classical solution $(n_1, u_1, E_1, B_1)$ on $[0,\infty) \times\mathbb{R}^3$, satisfying:
    $$
	(n_1, u_1, E_1) \in C \big([0, \infty); H^{s-1}(\mathbb{R}^3)\big)\cap C^1 \big([0, \infty); H^{s-2}(\mathbb{R}^3)\big), \quad  B_1 \in C \big([0, \infty); H^{s+1}(\mathbb{R}^3)\big).
    $$
    Moreover, the solution satisfies the following uniform energy estimate:
    \begin{align*}
	\sup_{t \geq 0}  \left\|(n_1, u_1, E_1, B_1)(t)\right\|_{\mathcal{H}^{s-1}} \leq C,
    \end{align*}
    where $C>0$ depends only on the initial data and the corresponding uniform bounds of the background solution $(n_0,u_0,E_0)$.
    \end{lemma}
	
	\begin{proof}
    The well-posedness of \eqref{e.10-1} follows directly from standard linear theory (see, for instance, \cite{Evans-1998,Dafermos-2016}). Here, we focus on the \textit{a priori} estimates. Applying the standard energy method for linear symmetrizable hyperbolic systems (cf. \cite{Friedrichs-CPAM-1954,Majda-AMS-1984,Guo-CMP-2010}), we obtain the following differential inequality:
    \begin{align*}
        \frac{d}{d t}\big\|(n_1,u_1,E_1)\big\|_{H^{s-1}}^2 \leq \frac{C}{(1+t)^p}\big\|(n_1,u_1,E_1)\big\|_{H^{s-1}}^2+\frac{C}{(1+t)^p}\big\|(n_1,u_1,E_1)\big\|_{H^{s-1}}.
    \end{align*}
    Applying Gronwall's inequality yields:
    \begin{align}\label{e.10-3}
        \big\|(n_1,u_1,E_1)(t)\big\|_{H^{s-1}}\leq C \left[\big\|(n_1,u_1,E_1)(0)\big\|_{H^{s-1}} + 1 \right],
    \end{align}
    where $C > 0$ depends only on $\|(n_0, u_0)\|_{H^s}$ and is uniform in $t$. The higher-order time estimates follow by repeated time differentiation of the system together with the same linear hyperbolic estimates.

Now we solve the linear curl-div system \eqref{e.10-2}. 
Fix $t\ge 0$ and write
\[
f:=\partial_t \nabla\phi-4\pi n_0u_0.
\]

\smallskip
\noindent\textit{Step 1. Regularity of $f$.}
By \eqref{e.9-2} and the preservation of irrotationality, there exists a scalar potential $\phi_0$ such that
\[
u_0=\nabla\phi_0,
\]
which yields that
\begin{align}\label{e.10-5}
f=\nabla\Phi-4\pi (n_0-\bar n)u_0,
\qquad
\Phi:=\partial_t\phi-4\pi\bar n\,\phi_0.
\end{align}

Since $n_0-\bar n\in H^s(\mathbb R^3),u_0\in H^s(\mathbb R^3;\mathbb R^3)$ and $s>\frac32$, the Sobolev algebra property yields
\begin{align}\label{e.10-15}
n_0u_0\in H^s(\mathbb R^3;\mathbb R^3),
\end{align}
and also
\[
n_0-\bar n,\ u_0\in L^2(\mathbb R^3)\cap L^\infty(\mathbb R^3).
\]
Hence
\begin{align}\label{e.10-6}
(n_0-\bar n)u_0\in L^1(\mathbb R^3)\cap L^2(\mathbb R^3)\subset L^{6/5}(\mathbb R^3).
\end{align}

On the other hand, since $-\Delta \phi=4\pi(n_0-\bar n)$,
\[
\partial_t \nabla\phi
=
-4\pi \nabla(-\Delta)^{-1}\operatorname{div}(n_0u_0).
\]
As $\nabla(-\Delta)^{-1}\operatorname{div}$ is a zero-order Fourier multiplier, it is bounded on
$H^s(\mathbb R^3)$. Hence
\[
\partial_t \nabla\phi\in H^s(\mathbb R^3;\mathbb R^3),
\qquad
f\in H^s(\mathbb R^3;\mathbb R^3).
\]
Moreover, by \eqref{e.9}, we have
\begin{align}\label{e.10-4}
\operatorname{div}f=0
\qquad  \text{in } \mathbb R^3.
\end{align}
\smallskip
\noindent\textit{Step 2. Construction of $\psi$.}
Define
\[
C_{c,\sigma}^\infty(\mathbb R^3)
:=
\bigl\{
w\in C_c^\infty(\mathbb R^3;\mathbb R^3):\operatorname{div}w=0
\bigr\},
\]
and
\[
\dot H^1_\sigma(\mathbb R^3)
:=
\overline{C_{c,\sigma}^\infty(\mathbb R^3)}^{\|\nabla(\cdot)\|_{L^2(\mathbb R^3)}}
\]
to be the completion of $C_{c,\sigma}^\infty(\mathbb R^3)$ with respect to the norm $\|\nabla(\cdot)\|_{L^2(\mathbb R^3)}$.

Define the linear functional
\[
\Lambda(w)
:=
-4\pi\int_{\mathbb R^3}(n_0-\bar n)u_0\cdot w\,dx,
\qquad
w\in \dot H^1_\sigma(\mathbb R^3).
\]
By \eqref{e.10-6} and H\"older's inequality,
\[
|\Lambda(w)|
\le
4\pi \|(n_0-\bar n)u_0\|_{L^{6/5}}\|w\|_{L^6}
\le
C\|(n_0-\bar n)u_0\|_{L^{6/5}}\|\nabla w\|_{L^2},
\]
where we have used the Sobolev inequality
\[
\|w\|_{L^6(\mathbb R^3)}
\le C\|\nabla w\|_{L^2(\mathbb R^3)},
\qquad
\forall\, w\in \dot H^1_\sigma(\mathbb R^3).
\]
Hence $\Lambda$ is continuous on $\dot H^1_\sigma(\mathbb R^3)$. 
Next define the bilinear form
\[
a(\psi,w)
:=
\int_{\mathbb R^3}\nabla\psi:\nabla w\,dx,
\qquad
\psi,w\in \dot H^1_\sigma(\mathbb R^3).
\]
Clearly $a$ is continuous and coercive. The Lax-Milgram theorem yields a unique
\[
\psi\in \dot H^1_\sigma(\mathbb R^3)
\]
such that
\begin{align}\label{e.10-7}
\int_{\mathbb R^3}\nabla\psi:\nabla w\,d x
=
-4\pi\int_{\mathbb R^3}(n_0-\bar n)u_0\cdot w\,d x,
\qquad
\forall\, w\in \dot H^1_\sigma(\mathbb R^3).
\end{align}
In particular,
\begin{align}\label{e.10-8}
\|\nabla\psi\|_{L^2(\mathbb R^3)}
\le
C\|(n_0-\bar n)u_0\|_{L^{6/5}(\mathbb R^3)}.
\end{align}
By definition of $\dot H^1_\sigma(\mathbb R^3)$,
\begin{align}\label{e.10-9}
\operatorname{div}\psi=0
\qquad \text{in }\mathcal D'(\mathbb R^3).
\end{align}

\smallskip
\noindent\textit{Step 3. Identification of the equation $-\Delta\psi=f$.}
For $w\in C_{c,\sigma}^\infty(\mathbb R^3)$, using \eqref{e.10-5},
\[
\langle \nabla\Phi,w\rangle
=
-\langle \Phi,\operatorname{div}w\rangle
=
0,
\]
and therefore
\[
\langle f,w\rangle
=
-4\pi\int_{\mathbb R^3}(n_0-\bar n)u_0\cdot w\,dx.
\]
Substituting this into \eqref{e.10-7}, we obtain
\begin{align}\label{e.10-10}
\int_{\mathbb R^3}\nabla\psi:\nabla w\,dx
=
\langle f,w\rangle,
\qquad
\forall\, w\in C_{c,\sigma}^\infty(\mathbb R^3).
\end{align}

Define
\[
g:=-\Delta\psi-f.
\]
Since $\nabla\psi\in L^2(\mathbb R^3)$, we have $-\Delta\psi\in H^{-1}(\mathbb R^3;\mathbb R^3)$; since $f\in H^s(\mathbb R^3)\subset L^2(\mathbb R^3)\subset H^{-1}(\mathbb R^3)$, it follows that
\[
g\in H^{-1}(\mathbb R^3;\mathbb R^3)\subset \mathcal S'(\mathbb R^3;\mathbb R^3).
\]
By \eqref{e.10-10},
\[
\langle g,w\rangle=0,
\qquad
\forall\, w\in C_{c,\sigma}^\infty(\mathbb R^3).
\]
Hence, by the distributional de Rham theorem, there exists $p\in \mathcal D'(\mathbb R^3)$ such that
\[
g=\nabla p.
\]
On the other hand, by \eqref{e.10-4} and \eqref{e.10-9},
\[
\operatorname{div}g
=
-\Delta(\operatorname{div}\psi)-\operatorname{div}f
=
0
\qquad \text{in }\mathcal D'(\mathbb R^3),
\]
while $g=\nabla p$ implies
\[
\nabla\times g=0
\qquad \text{in }\mathcal D'(\mathbb R^3).
\]

Taking Fourier transforms yields
\[
\xi\cdot \widehat g=0,
\qquad
\xi\times \widehat g=0
\qquad \text{in }\mathcal S'(\mathbb R^3).
\]
Let $\chi\in C_c^\infty(\mathbb R^3\setminus\{0\})$. Since $|\xi|^{-2}$ is smooth on $\operatorname{supp}\chi$, the vector identity
\[
|\xi|^2 v=\xi(\xi\cdot v)-\xi\times(\xi\times v)
\]
gives
\[
|\xi|^2 \chi \widehat g
=
\xi(\xi\cdot \chi \widehat g)-\xi\times(\xi\times \chi \widehat g)
=
0.
\]
Thus $\chi \widehat g=0$ for every $\chi\in C_c^\infty(\mathbb R^3\setminus\{0\})$, and hence
\[
\operatorname{supp}\widehat g\subset\{0\}.
\]
Therefore $g$ is a polynomial distribution. Since $g\in H^{-1}(\mathbb R^3)$, 
\[
g=0.
\]
Consequently,
\begin{align}\label{e.10-11}
-\Delta\psi=f
\qquad \text{in }\mathcal D'(\mathbb R^3;\mathbb R^3).
\end{align}

\smallskip
\noindent\textit{Step 4. Construction of $B_1$.}
Noting \eqref{e.10-11}, the standard global elliptic estimate on $\mathbb R^3$ yields
\begin{align}\label{e.10-13}
\|\nabla^2\psi\|_{H^s(\mathbb R^3)}
\le
C\|f\|_{H^s(\mathbb R^3)}.
\end{align}
Combining \eqref{e.10-8} and \eqref{e.10-13}, we obtain
\begin{align}\label{e.10-14}
\|\nabla\psi\|_{H^{s+1}(\mathbb R^3)}
\le
C\Bigl(
\|(n_0-\bar n)u_0\|_{L^{6/5}(\mathbb R^3)}
+
\|f\|_{H^s(\mathbb R^3)}
\Bigr).
\end{align}

Define
\begin{align*}
B_1:=\nabla\times\psi.
\end{align*}
Since $B_1$ is a first-order linear combination of $\nabla\psi$, \eqref{e.10-14} implies
\begin{align}\label{e.10-16}
B_1\in H^{s+1}(\mathbb R^3;\mathbb R^3),
\qquad
\|B_1\|_{H^{s+1}(\mathbb R^3)}
\le
C\Bigl(
\|(n_0-\bar n)u_0\|_{L^{6/5}(\mathbb R^3)}
+
\|f\|_{H^s(\mathbb R^3)}
\Bigr).
\end{align}
Moreover, by \eqref{e.10-9} and \eqref{e.10-11},
\[
\operatorname{div}B_1
=
\operatorname{div}(\nabla\times\psi)
=
0,
\]
and
\[
\nabla\times B_1
=
\nabla\times(\nabla\times\psi)
=
\nabla(\operatorname{div}\psi)-\Delta\psi
=f.
\]
Hence $B_1$ solves \eqref{e.10-2}. The uniqueness is standard: if $\widetilde B_1$ is another solution to \eqref{e.10-2} in $H^{s+1}(\mathbb R^3;\mathbb R^3)$, then $D:=B_1-\widetilde B_1\in L^2(\mathbb R^3)$ is harmonic in $\mathbb R^3$, and thus $D\equiv0$.

The time regularity of $B_1$ can be obtained in the same way: for each integer $0\le l\le s-1$, one differentiates \eqref{e.10-2} $l$ times with respect to $t$ and applies the above construction to $\partial_t^l B_1$, which yields the stated regularity. This completes the proof of Lemma \ref{l3.3}.
\end{proof}
\subsection{Justification of the expansions}
In this subsection, we rigorously justify the asymptotic expansions of solutions $\left(n_{\mathfrak{e}}, u_{\mathfrak{e}}, E_{\mathfrak{e}}, B_{\mathfrak{e}}\right)$ of \eqref{e.4} under the assumption of well-prepared initial data. As a consequence, we obtain the convergence of $\left(n_{\mathfrak{e}}, u_{\mathfrak{e}}, E_{\mathfrak{e}}, B_{\mathfrak{e}}\right)$ to the solution $\left(n_0, u_0, E_0, B_0\right)$ of the compressible EP equations \eqref{e.9} as $\mathfrak{c}\rightarrow \infty$.

For later use, we denote
\[
\left(n_{a}^{\mathfrak{c}}, u_{a}^{\mathfrak{c}}, E_{a}^{\mathfrak{c}}, B_{a}^{\mathfrak{c}}\right):=\sum_{j=0}^1 \mathfrak{c}^{-j}\left(n_j, u_j, E_j, B_j\right)
\]
with $(n_j, u_j, E_j, B_j)$ ($j=0,1$) obtained in Lemmas \ref{l3.2}--\ref{l3.3}. One can verify that $\left(n_{a}^{\mathfrak{c}}, u_{a}^{\mathfrak{c}}, E_{a}^{\mathfrak{c}}, B_{a}^{\mathfrak{c}}\right)$ satisfy:
\begin{equation}\label{e.11}
	\left\{\begin{aligned}
		&\partial_t n_{a}^{\mathfrak{c}}+\operatorname{div}\left(n_{a}^{\mathfrak{c}} u_{a}^{\mathfrak{c}}\right)=R_n^{\mathfrak{c}}, \\
		&\partial_t u_{a}^{\mathfrak{c}}+\left(u_{a}^{\mathfrak{c}}\cdot \nabla\right) u_{a}^{\mathfrak{c}}+\frac{5}{2}\nabla T_0(n_0)+\frac{1}{\mathfrak{c}}\frac{5}{2}\nabla\left(T_0^{\prime}\left(n_0\right) n_1\right)=-E_{a}^{\mathfrak{c}}-\frac{1}{\mathfrak{c}} u_{a}^{\mathfrak{c}} \times B_{a}^{\mathfrak{c}}+R_u^{\mathfrak{c}}, \\
		&\partial_t E_{a}^{\mathfrak{c}}-\mathfrak{c} \nabla \times B_{a}^{\mathfrak{c}}=4\pi n_{a}^{\mathfrak{c}} u_{a}^{\mathfrak{c}}+R_E^{\mathfrak{c}}, \quad \operatorname{div} E_{a}^{\mathfrak{c}}=4\pi\left(\bar{n}-n_{a}^{\mathfrak{c}}\right), \\
		&\partial_t B_{a}^{\mathfrak{c}}+\mathfrak{c} \nabla \times E_{a}^{\mathfrak{c}}=R_B^{\mathfrak{c}}, \quad \operatorname{div} B_{a}^{\mathfrak{c}}=0, \\
		&\left.\left(n_{a}^{\mathfrak{c}}, u_{a}^{\mathfrak{c}}, E_{a}^{\mathfrak{c}}, B_{a}^{\mathfrak{c}}\right)\right|_{t=0}=\sum_{j=0}^1 \mathfrak{c}^{-j}\left(n_j, u_j, E_j, B_j\right)(0),
	\end{aligned}\right.
\end{equation}
where 
\begin{align}\label{e.11-1}
	&R_n^{\mathfrak{c}}=\frac{1}{\mathfrak{c}^2}\operatorname{div}\left(n_1u_1\right), \quad R_u^{\mathfrak{c}}=\frac{1}{\mathfrak{c}^2}\left(u_1\cdot \nabla\right)u_1+\frac{1}{\mathfrak{c}^2}\left(u_1 \times B_0+u_0 \times B_1\right)+\frac{1}{\mathfrak{c}^3}u_1\times B_1,\nonumber\\
	&R_E^{\mathfrak{c}}=\frac{1}{\mathfrak{c}}\partial_t E_1-4\pi \frac{1}{\mathfrak{c}}\left(n_1u_0+n_0u_1\right)-4\pi\frac{1}{\mathfrak{c}^2}n_1u_1, \quad R_B^{\mathfrak{c}}=\frac{1}{\mathfrak{c}}\partial_t B_1.
\end{align}
It is clear that
\begin{align}\label{e.12}
	\operatorname{div} R_E^{\mathfrak{c}}=-4\pi R_n^{\mathfrak{c}}, \quad \operatorname{div} R_B^{\mathfrak{c}}=0.
\end{align}
According to Lemmas \ref{l3.2}--\ref{l3.3}, we have
\begin{align}\label{e.13}
	\sup_{t\geq 0}\left\|\left(R_n^{\mathfrak{c}}, R_u^{\mathfrak{c}}, R_E^{\mathfrak{c}}, R_B^{\mathfrak{c}}\right)\right\|_{\mathcal{H}^{s-1}} \leq C_0 \mathfrak{c}^{-1},
\end{align}
where $C_0>0$ is a constant independent of $\mathfrak{c}$.

Now we proceed to derive the remainder equations. Denote
\begin{equation}\label{e.14}
	\left(\mathcal{N}^{\mathfrak{c}}, \mathcal{U}^{\mathfrak{c}}, \mathcal{E}^{\mathfrak{c}}, \mathcal{B}^{\mathfrak{c}}\right):=\left(n_{\mathfrak{e}}-n_{a}^{\mathfrak{c}}, u_{\mathfrak{e}}-u_{a}^{\mathfrak{c}}, E_{\mathfrak{e}}-E_{a}^{\mathfrak{c}}, B_{\mathfrak{e}}-B_{a}^{\mathfrak{c}}\right).
\end{equation}
Subtracting \eqref{e.11} from \eqref{e.4}, one has:
\begin{align}\label{e.15}
	\begin{cases}
		h^{\prime}\partial_t \mathcal{N}^{\mathfrak{c}}+\frac{h^{\prime}(n_{a}^{\mathfrak{c}}+\mathcal{N}^{\mathfrak{c}})}{\left(u_{\mathfrak{e}}^0\right)^2}\left(u_{a}^{\mathfrak{c}}+\mathcal{U}^{\mathfrak{c}}\right)\cdot \partial_t \mathcal{U}^{\mathfrak{c}}+\frac{\mathfrak{c}h^{\prime}}{u_{\mathfrak{e}}^0}\left(u_{a}^{\mathfrak{c}}+\mathcal{U}^{\mathfrak{c}}\right)\cdot\nabla \mathcal{N}^{\mathfrak{c}}+\frac{\mathfrak{c}h^{\prime}}{u_{\mathfrak{e}}^0}\left(n_{a}^{\mathfrak{c}}+\mathcal{N}^{\mathfrak{c}}\right)\operatorname{div} \mathcal{U}^{\mathfrak{c}}\\
		\qquad\qquad\qquad\qquad \qquad \qquad \qquad \qquad \qquad +B_{11}\mathcal{N}^{\mathfrak{c}}+B_{12}\cdot\mathcal{U}^{\mathfrak{c}}=H_1^{\mathfrak{c}},\\
		\frac{h^{\prime}(n_{a}^{\mathfrak{c}}+\mathcal{N}^{\mathfrak{c}})}{\left(u_{\mathfrak{e}}^0\right)^2}(u_{a}^{\mathfrak{c}}+\mathcal{U}^{\mathfrak{c}})\partial_t \mathcal{N}^{\mathfrak{c}}+\frac{h}{\mathfrak{c}^2}\left(n_{a}^{\mathfrak{c}}+\mathcal{N}^{\mathfrak{c}}\right)\Big(I-\frac{u_{\mathfrak{e}} \otimes u_{\mathfrak{e}}}{\left(u_{\mathfrak{e}}^0\right)^2}\Big)\partial_t\mathcal{U}^{\mathfrak{c}}+\frac{\mathfrak{c}h^{\prime}}{u_{\mathfrak{e}}^0}\left(n_{a}^{\mathfrak{c}}+\mathcal{N}^{\mathfrak{c}}\right) \nabla_x \mathcal{N}^{\mathfrak{c}} \\
		\qquad \qquad 
		+\frac{h}{\mathfrak{c} u_{\mathfrak{e}}^0}\left(n_{a}^{\mathfrak{c}}+\mathcal{N}^{\mathfrak{c}}\right)\Big(I-\frac{u_{\mathfrak{e}} \otimes u_{\mathfrak{e}}}{\left(u_{\mathfrak{e}}^0\right)^2}\Big)\left(u_{a}^{\mathfrak{c}}+\mathcal{U}^{\mathfrak{c}}\right)\cdot \nabla \mathcal{U}^{\mathfrak{c}}+\mathcal{N}^{\mathfrak{c}}B_{21}+B_{22}\mathcal{U}^{\mathfrak{c}}\\
		\qquad\qquad \qquad  =-\left(I-\frac{u_{\mathfrak{e}} \otimes u_{\mathfrak{e}}}{(u_{\mathfrak{e}}^0)^2}\right)\left(\mathcal{N}^{\mathfrak{c}}+n_{a}^{\mathfrak{c}}\right) \mathcal{E}^{\mathfrak{c}} -\frac{\mathcal{N}^{\mathfrak{c}}+n_{a}^{\mathfrak{c}}}{u_{\mathfrak{e}}^0} \left(\mathcal{U}^{\mathfrak{c}}+u_{a}^{\mathfrak{c}}\right)  \times \mathcal{B}^{\mathfrak{c}}+H_2^{\mathfrak{c}},\\
		\partial_t \mathcal{E}^{\mathfrak{c}}-\mathfrak{c} \nabla \times \mathcal{B}^{\mathfrak{c}}=4 \pi \left(\mathcal{N}^{\mathfrak{c}}+n_{a}^{\mathfrak{c}}\right) \mathcal{U}^{\mathfrak{c}}+4 \pi\mathcal{N}^{\mathfrak{c}} u_{a}^{\mathfrak{c}}-R_E^{\mathfrak{c}}, \quad \operatorname{div} \mathcal{B}^{\mathfrak{c}}=0,\\
		\partial_t \mathcal{B}^{\mathfrak{c}}+\mathfrak{c} \nabla \times \mathcal{E}^{\mathfrak{c}}=-R_B^{\mathfrak{c}}, \quad \operatorname{div} \mathcal{E}^{\mathfrak{c}}=4\pi\Big(n_{a}^{\mathfrak{c}}-\frac{u_{\mathfrak{e}}^0}{\mathfrak{c}}\left(\mathcal{N}^{\mathfrak{c}}+n_{a}^{\mathfrak{c}}\right)\Big), 
	\end{cases}
\end{align}
where
	\begin{align*}
	    B_{11}=&\frac{\mathfrak{c}h^{\prime}}{u_{\mathfrak{e}}^0}\operatorname{div}u_{a}^{\mathfrak{c}}+\frac{h^{\prime}}{(u_{\mathfrak{e}}^0)^2}(u_{a}^{\mathfrak{c}}+ \mathcal{U}^{\mathfrak{c}})\cdot\partial_t u_{a}^{\mathfrak{c}}, \quad  B_{12}=\frac{\mathfrak{c}h^{\prime}}{u_{\mathfrak{e}}^0}\nabla n_{a}^{\mathfrak{c}}+\frac{h^{\prime}}{(u_{\mathfrak{e}}^0)^2}n_{a}^{\mathfrak{c}}\partial_t u_{a}^{\mathfrak{c}}, \nonumber\\
	    B_{21}=&  \frac{h}{\mathfrak{c}^2}\Big(I-\frac{u_{\mathfrak{e}} \otimes u_{\mathfrak{e}}}{\left(u_{\mathfrak{e}}^0\right)^2}\Big) \partial_t u_a^{\mathfrak{c}}-\frac{\mathfrak{c}}{u_{\mathfrak{e}}^0}\partial_t u_a^{\mathfrak{c}}+\frac{\mathfrak{c} h^{\prime}}{u_{\mathfrak{e}}^0} \nabla_x n_a^{\mathfrak{c}}-\frac{\mathfrak{c}}{u_{\mathfrak{e}}^0}\Big(\frac{5}{2} \nabla T_0\left(n_0\right)+\frac{1}{\mathfrak{c}} \frac{5}{2} \nabla\left(T_0^{\prime}\left(n_0\right) n_1\right)\Big) \nonumber\\ 
	    & +\frac{h^{\prime}}{\left(u_{\mathfrak{e}}^0\right)^2}\left(u_a^{\mathfrak{c}}+\mathcal{U}^{\mathfrak{c}}\right) \partial_t n_a^{\mathfrak{c}}+\Big(I-\frac{u_{\mathfrak{e}} \otimes u_{\mathfrak{e}}}{\left(u_{\mathfrak{e}}^0\right)^2}\Big) E_a^{\mathfrak{c}}-\frac{\mathfrak{c}}{u_{\mathfrak{e}}^0} E_a^{\mathfrak{c}}+\frac{h}{\mathfrak{c} u_{\mathfrak{e}}^0}\Big(I-\frac{u_{\mathfrak{e}} \otimes u_{\mathfrak{e}}}{\left(u_{\mathfrak{e}}^0\right)^2}\Big)\left(\mathcal{U}^{\mathfrak{c}} \cdot \nabla \right) u_a^{\mathfrak{c}} \nonumber\\ 
	    & +\frac{h}{\mathfrak{c} u_{\mathfrak{e}}^0}\Big(I-\frac{u_{\mathfrak{e}} \otimes u_{\mathfrak{e}}}{\left(u_{\mathfrak{e}}^0\right)^2}\Big)\left(u_a^{\mathfrak{c}} \cdot \nabla \right)u_a^{\mathfrak{c}}-\frac{\mathfrak{c}}{ u_{\mathfrak{e}}^0}\left(u_a^{\mathfrak{c}} \cdot \nabla \right)u_a^{\mathfrak{c}} +\frac{\mathfrak{c}}{u_{\mathfrak{e}}^0} R_u^{\mathfrak{c}},\nonumber \\
	    B_{22}=&\frac{h^{\prime}}{\left(u_{\mathfrak{e}}^0\right)^2} n_a^{\mathfrak{c}}\partial_t n_a^{\mathfrak{c}} I+\frac{h}{\mathfrak{c} u_{\mathfrak{e}}^0} n_a^{\mathfrak{c}}\Big(I-\frac{u_{\mathfrak{e}}\otimes u_{\mathfrak{e}}}{\left(u_{\mathfrak{e}}^0\right)^2}\Big) \nabla u_a^{\mathfrak{c}}+\frac{n_a^{\mathfrak{c}}+\mathcal{N}^{\mathfrak{c}}}{u_{\mathfrak{e}}^0} [\times B_{a}^{\mathfrak{c}}], \nonumber\\
		H_1^{\mathfrak{c}}=&h^{\prime}\Big(\frac{\mathfrak{c}}{u_{\mathfrak{e}}^0}-1\Big) \partial_t n_{a}^{\mathfrak{c}}-\frac{h^{\prime}n_{a}^{\mathfrak{c}}}{(u_{\mathfrak{e}}^0)^2}u_{a}^{\mathfrak{c}}\cdot\partial_t u_{a}^{\mathfrak{c}}-\frac{\mathfrak{c}h^{\prime}}{u_{\mathfrak{e}}^0}R_n^{\mathfrak{c}},\\
		H_2^{\mathfrak{c}}=& \Big(\frac{\mathfrak{c}}{u_{\mathfrak{e}}^0}-\frac{h}{\mathfrak{c}^2}\Big) n_a^{\mathfrak{c}} \partial_t u_a^{\mathfrak{c}}+\frac{h}{\mathfrak{c}^2} \frac{u_{\mathfrak{e}}\otimes u_{\mathfrak{e}}}{\left(u_{\mathfrak{e}}^0\right)^2} n_a^{\mathfrak{c}} \partial_t u_a^{\mathfrak{c}}-\frac{h^{\prime} n_a^{\mathfrak{c}} u_a^{\mathfrak{c}}}{\left(u_{\mathfrak{e}}^0\right)^2} \partial_t n_a^{\mathfrak{c}}+\Big(\mathfrak{c}-\frac{h}{\mathfrak{c}}\Big) \frac{n_a^{\mathfrak{c}}}{u_{\mathfrak{e}}^0}\left(u_a^{\mathfrak{c}} \cdot \nabla\right) u_a^{\mathfrak{c}} \\ 
		&-\frac{\mathfrak{c} h^{\prime}}{u_{\mathfrak{e}}^0} n_a^{\mathfrak{c}} \nabla_x n_a^{\mathfrak{c}} +\frac{\mathfrak{c}}{u_{\mathfrak{e}}^0} n_a^{\mathfrak{c}}\Big(\frac{5}{2} \nabla T_0\left(n_0\right)+\frac{1}{\mathfrak{c}} \frac{5}{2} \nabla\left(T_0^{\prime}\left(n_0\right) n_1\right)\Big)  \\ 
		& +\frac{h}{\mathfrak{c} u_{\mathfrak{e}}^0} \frac{u_{\mathfrak{e}} \otimes u_{\mathfrak{e}}}{\left(u_{\mathfrak{e}}^0\right)^2} n_a^{\mathfrak{c}}\left(u_a^{\mathfrak{c}} \cdot \nabla\right) u_a^{\mathfrak{c}}+\Big(\frac{\mathfrak{c}}{u_{\mathfrak{e}}^0}-1\Big) n_a^{\mathfrak{c}} E_a^{\mathfrak{c}}+\frac{u_{\mathfrak{e}}\otimes u_{\mathfrak{e}}}{\left(u_{\mathfrak{e}}^0\right)^2} n_a^{\mathfrak{c}} E_a^{\mathfrak{c}}-\frac{\mathfrak{c}}{u_{\mathfrak{e}}^0} n_a^{\mathfrak{c}} R_u^{\mathfrak{c}},
	\end{align*}
    and 
    \begin{align*}
        [\times B_{a}^{\mathfrak{c}}]:=\left(\begin{array}{ccc}0 & B_{a, 3}^{\mathfrak{c}} & -B_{a, 2}^{\mathfrak{c}} \\ -B_{a, 3}^{\mathfrak{c}} & 0 & B_{a, 1}^{\mathfrak{c}} \\ B_{a, 2}^{\mathfrak{c}} & -B_{a, 1}^{\mathfrak{c}} & 0\end{array}\right).
    \end{align*}
    
Set
\begin{align*}
	& \mathcal{W}^{\mathfrak{c}}=\binom{\mathcal{N}^{\mathfrak{c}}}{\mathcal{U}^{\mathfrak{c}}}, \quad  B=\left(\begin{array}{cc}B_{11} & B_{12}^t \\ B_{21} & B_{22} \end{array}\right),\quad H^{\mathfrak{c}}=\binom{H_1^{\mathfrak{c}}}{H_2^{\mathfrak{c}}},\\
	& A_0=\left(\begin{array}{cc}
		h^{\prime} & \frac{h^{\prime}(n_{a}^{\mathfrak{c}}+\mathcal{N}^{\mathfrak{c}})}{\left(u_{\mathfrak{e}}^0\right)^2}\left(u_{a}^{\mathfrak{c}}+\mathcal{U}^{\mathfrak{c}}\right)^t \\
		\frac{h^{\prime}(n_{a}^{\mathfrak{c}}+\mathcal{N}^{\mathfrak{c}})}{\left(u_{\mathfrak{e}}^0\right)^2}\left(u_{a}^{\mathfrak{c}}+\mathcal{U}^{\mathfrak{c}}\right) \quad& \frac{ h}{\mathfrak{c}^2}(n_{a}^{\mathfrak{c}}+\mathcal{N}^{\mathfrak{c}})\Big(I-\frac{u_{\mathfrak{e}} \otimes u_{\mathfrak{e}}}{\left(u_{\mathfrak{e}}^0\right)^2}\Big)
	\end{array}\right),\\
	&A_i=\left(\begin{array}{cc}
		\frac{\mathfrak{c}h^{\prime}}{u_{\mathfrak{e}}^0}\left(u_{a}^{\mathfrak{c}}+\mathcal{U}^{\mathfrak{c}}\right)_i & \frac{\mathfrak{c}h^{\prime}}{u_{\mathfrak{e}}^0} \left(n_{a}^{\mathfrak{c}}+\mathcal{N}^{\mathfrak{c}}\right)\mathbf{e}_i^t \\
		\frac{\mathfrak{c}h^{\prime}}{u_{\mathfrak{e}}^0}\left(n_{a}^{\mathfrak{c}}+\mathcal{N}^{\mathfrak{c}}\right) \mathbf{e}_i \quad& \frac{h}{\mathfrak{c} u_{\mathfrak{e}}^0}\left(n_{a}^{\mathfrak{c}}+\mathcal{N}^{\mathfrak{c}}\right)\Big(I-\frac{u_{\mathfrak{e}} \otimes u_{\mathfrak{e}}}{\left(u_{\mathfrak{e}}^0\right)^2}\Big)\left(u_{a}^{\mathfrak{c}}+\mathcal{U}^{\mathfrak{c}}\right)_i
	\end{array}\right),\\
	& R^{\mathfrak{c}}=\left(\begin{array}{c}
		0 \\
		-\left(I-\frac{u_{\mathfrak{e}} \otimes u_{\mathfrak{e}}}{(u_{\mathfrak{e}}^0)^2}\right)\left(\mathcal{N}^{\mathfrak{c}}+n_{a}^{\mathfrak{c}}\right) \mathcal{E}^{\mathfrak{c}} -\frac{n_{a}^{\mathfrak{c}}+\mathcal{N}^{\mathfrak{c}}}{u_{\mathfrak{e}}^0} \left(\mathcal{U}^{\mathfrak{c}}+u_{a}^{\mathfrak{c}}\right)  \times \mathcal{B}^{\mathfrak{c}}
	\end{array}\right),
\end{align*}
where $\left(\mathbf{e}_1, \mathbf{e}_2, \mathbf{e}_3\right)$ is the canonical basis of $\mathbb{R}^3$ and $y_i$ denotes the $i$-th component of $y \in \mathbb{R}^3$.

It is clear that $A_0$ and $A_i$ are symmetric matrices. Furthermore, $A_0$ is positive definite. The proof of this statement will be given in Appendix \ref{Appendix0}. Thus, we can rewrite \eqref{e.15} as the following symmetric hyperbolic system:
\begin{align}\label{e.16}
	\left\{\begin{aligned}
		&A_0\partial_t \mathcal{W}^{\mathfrak{c}}+\sum_{i=1}^3 A_i \partial_{x_i} \mathcal{W}^{\mathfrak{c}}+B\mathcal{W}^{\mathfrak{c}}= R^{\mathfrak{c}}+H^{\mathfrak{c}}, \\
		&\partial_t \mathcal{E}^{\mathfrak{c}}-\mathfrak{c} \nabla \times \mathcal{B}^{\mathfrak{c}}=4\pi\left(\mathcal{N}^{\mathfrak{c}}+n_{a}^{\mathfrak{c}}\right) \mathcal{U}^{\mathfrak{c}}+4\pi\mathcal{N}^{\mathfrak{c}} u_{a}^{\mathfrak{c}}-R_E^{\mathfrak{c}}, \\
		&\partial_t \mathcal{B}^{\mathfrak{c}}+\mathfrak{c} \nabla \times \mathcal{E}^{\mathfrak{c}}=-R_B^{\mathfrak{c}}.
	\end{aligned}\right.
\end{align}
Also, note that the redundant equations $\operatorname{div} \mathcal{E}^{\mathfrak{c}}=4\pi\Big(n_{a}^{\mathfrak{c}}-\frac{u_{\mathfrak{e}}^0}{\mathfrak{c}}\left(\mathcal{N}^{\mathfrak{c}}+n_{a}^{\mathfrak{c}}\right)\Big)$ and $\operatorname{div} \mathcal{B}^{\mathfrak{c}}=0$ in the system \eqref{e.15} hold globally if they are satisfied by the initial conditions. This is ensured by \eqref{e.19-1} and \eqref{e.9}--\eqref{e.10} via the transform \eqref{e.14}.

The existence and uniqueness of local smooth solutions to \eqref{e.4} or \eqref{e.16} are addressed by the following lemma:

\begin{lemma}\label{t3.6}
Let $(n_0, u_0, E_0)$ and $\left(n_1, u_1, E_1, B_1\right)$ be the global smooth solutions given by Lemma \ref{l3.2} and Lemma \ref{l3.3}, respectively. For any sufficiently large integer $s>0$, we suppose that the initial data satisfy
\begin{equation}\label{e.23}
	\left\|\left(\mathcal{W}^{\mathfrak{c}}, \mathcal{E}^{\mathfrak{c}}, \mathcal{B}^{\mathfrak{c}}\right)(0)\right\|_{\mathcal{H}^{s-2}} \leq C_1 \mathfrak{c}^{-1}
\end{equation}
for some constant $C_1 > 0$ independent of $\mathfrak{c}$. Then for any given finite time $T > 0$, there exists a sufficiently large constant $\mathfrak{c}_0 = \mathfrak{c}_0(T, C_0, C_1) > 0$ such that for all $\mathfrak{c} \geq \mathfrak{c}_0$, the remainder equations \eqref{e.16} admit a unique classical solution 
\[
	\left(\mathcal{W}^{\mathfrak{c}}, \mathcal{E}^{\mathfrak{c}}, \mathcal{B}^{\mathfrak{c}}\right) \in C \big([0, T]; H^{s-2}(\mathbb{R}^3)\big)\cap C^1 \big([0, T]; H^{s-3}(\mathbb{R}^3)\big),
\]
and the following uniform energy estimate holds:
\begin{equation}\label{e.27}
	\sup_{t \in [0,T]} \left\|\left(\mathcal{W}^{\mathfrak{c}}, \mathcal{E}^{\mathfrak{c}}, \mathcal{B}^{\mathfrak{c}}\right)(t)\right\|_{\mathcal{H}^{s-2}} \leq C(C_0, C_1, T) \mathfrak{c}^{-1},
\end{equation}
where $C(C_0, C_1, T) > 0$ is a constant independent of $\mathfrak{c}$.
\end{lemma}

\begin{proof}
 According to the classical theory of symmetric hyperbolic systems (see, e.g., \cite{Majda-AMS-1984, Friedrichs-CPAM-1954}), the nonlinear IVP \eqref{e.16} admits a local smooth solution. To extend the solution up to any given finite time $T>0$, we employ a bootstrap argument. First, we make the following \textit{a priori} assumption:
\begin{align}\label{e.24}
	\left\|\left(\mathcal{W}^{\mathfrak{c}}, \mathcal{E}^{\mathfrak{c}}, \mathcal{B}^{\mathfrak{c}}\right)\right\|_{\mathcal{H}^{s-2}}\leq 1.
\end{align}
Under this assumption, we prove that estimate \eqref{e.27} holds. When $\mathfrak{c}\gg 1$, \eqref{e.27} implies a bound that is stricter than \eqref{e.24}. Then, by a standard continuity argument, the solution can be extended up to $T$.
    
For any integer $k\geq 1$, it follows from the Fundamental Theorem of Calculus that
\begin{align}\label{e.21}
	&\left\|\left(\mathcal{W}^{\mathfrak{c}}, \mathcal{E}^{\mathfrak{c}}, \mathcal{B}^{\mathfrak{c}}\right)(t)\right\|_{\mathcal{H}^{k-1}}^2 \nonumber\\
	\leq&  \left\|\left(\mathcal{W}^{\mathfrak{c}}, \mathcal{E}^{\mathfrak{c}}, \mathcal{B}^{\mathfrak{c}}\right)(0)\right\|_{\mathcal{H}^{k-1}}^2+2 \int_0^t\left\|\partial_\tau \left(\mathcal{W}^{\mathfrak{c}}, \mathcal{E}^{\mathfrak{c}}, \mathcal{B}^{\mathfrak{c}}\right)(\tau)\right\|_{\mathcal{H}^{k-1}}\|\left(\mathcal{W}^{\mathfrak{c}}, \mathcal{E}^{\mathfrak{c}}, \mathcal{B}^{\mathfrak{c}}\right)(\tau)\|_{\mathcal{H}^{k-1}} d \tau \nonumber\\
	\leq&  \left\|\left(\mathcal{W}^{\mathfrak{c}}, \mathcal{E}^{\mathfrak{c}}, \mathcal{B}^{\mathfrak{c}}\right)(0)\right\|_{\mathcal{H}^{k-1}}^2+2\int_0^t\|\left(\mathcal{W}^{\mathfrak{c}}, \mathcal{E}^{\mathfrak{c}}, \mathcal{B}^{\mathfrak{c}}\right)(\tau)\|_{\mathcal{H}^k}^2 d \tau.
\end{align}
Hence, we only need to perform the energy estimates for the highest-order derivatives.
	
Let $l+|\alpha|=s-2$. Applying the operator $\partial_t^l\partial_x^{\alpha}$ directly to \eqref{e.16}, we obtain:
\begin{equation}\label{B.2}
	\left\{\begin{aligned}
		&A_0 \partial_t (\partial_t^l\partial_x^{\alpha} \mathcal{W}^{\mathfrak{c}}) + \sum_{i=1}^3 A_i \partial_{x_i} (\partial_t^l\partial_x^{\alpha} \mathcal{W}^{\mathfrak{c}}) + B (\partial_t^l\partial_x^{\alpha} \mathcal{W}^{\mathfrak{c}}) = \partial_t^l\partial_x^{\alpha}(H^{\mathfrak{c}}+R^{\mathfrak{c}}) + \mathcal{C}_{l,\alpha}^{\mathfrak{c}}, \\
		&\partial_t (\partial_t^l\partial_x^{\alpha} \mathcal{E}^{\mathfrak{c}})-\mathfrak{c} \nabla \times (\partial_t^l\partial_x^{\alpha} \mathcal{B}^{\mathfrak{c}})=\partial_t^l\partial_x^{\alpha}[\left(\mathcal{N}^{\mathfrak{c}}+n_{a}^{\mathfrak{c}}\right) \mathcal{U}^{\mathfrak{c}}+\mathcal{N}^{\mathfrak{c}} u_{a}^{\mathfrak{c}}-R_E^{\mathfrak{c}}], \\
		&\partial_t (\partial_t^l\partial_x^{\alpha} \mathcal{B}^{\mathfrak{c}})+\mathfrak{c} \nabla \times (\partial_t^l\partial_x^{\alpha} \mathcal{E}^{\mathfrak{c}})=-\partial_t^l\partial_x^{\alpha}R_B^{\mathfrak{c}},
	\end{aligned}\right.
\end{equation}
where the commutator $\mathcal{C}_{l,\alpha}^{\mathfrak{c}}$ is given by
\[
\mathcal{C}_{l,\alpha}^{\mathfrak{c}} =-[\partial_t^l\partial_x^{\alpha},A_0]\partial_t \mathcal{W}^{\mathfrak{c}} - \sum_{i=1}^3 [\partial_t^l\partial_x^{\alpha},A_i]\partial_{x_i} \mathcal{W}^{\mathfrak{c}} - [\partial_t^l\partial_x^{\alpha}, B] \mathcal{W}^{\mathfrak{c}}.
\]
By Lemmas \ref{l3.2}--\ref{l3.3}, the \textit{a priori} assumption \eqref{e.24}, and standard Moser-type commutator estimates, we can bound the commutator as:
\[
\left\|\mathcal{C}_{l,\alpha}^{\mathfrak{c}}\right\| \leq C(c_0) \left\|\mathcal{W}^{\mathfrak{c}}\right\|_{\mathcal{H}^{s-2}}.
\]
Multiplying $\eqref{B.2}_1$ by $\partial_t^l\partial_x^{\alpha} \mathcal{W}^{\mathfrak{c}}$, integrating over $\mathbb{R}^3$, and utilizing the symmetry of $A_0$ and $A_i$ alongside integration by parts, we deduce:
\begin{align}\label{B.9}
	\frac{d}{d t}\left\|\partial_t^l\partial_x^{\alpha} \mathcal{W}^{\mathfrak{c}}\right\|^2\leq C(C_0)\left\|(\mathcal{W}^{\mathfrak{c}},\mathcal{E}^{\mathfrak{c}},\mathcal{B}^{\mathfrak{c}})\right\|_{\mathcal{H}^{s-2}}^2+\big\|\partial_t^l\partial_x^{\alpha} H^{\mathfrak{c}}\big\|^2,
\end{align}
where we have used $\left\langle A_0\partial_t^l\partial_x^{\alpha} \mathcal{W}^{\mathfrak{c}}, \partial_t^l\partial_x^{\alpha} \mathcal{W}^{\mathfrak{c}}\right\rangle\cong \|\partial_t^l\partial_x^{\alpha} \mathcal{W}^{\mathfrak{c}}\|^2$ since $A_0$ is uniformly positive definite.

Multiplying $\eqref{B.2}_{2,3}$ by $\partial_t^l\partial_x^{\alpha} \mathcal{E}^{\mathfrak{c}}$ and $\partial_t^l\partial_x^{\alpha}\mathcal{B}^{\mathfrak{c}}$ respectively, and integrating the resultant equations over $\mathbb{R}^3$, we obtain:
\begin{equation}\label{B.10}
\frac{d}{d t}\left\|\partial_t^l\partial_x^{\alpha}(\mathcal{E}^{\mathfrak{c}},\mathcal{B}^{\mathfrak{c}})\right\|^2\leq C(C_0)\left\|(\mathcal{W}^{\mathfrak{c}},\mathcal{E}^{\mathfrak{c}},\mathcal{B}^{\mathfrak{c}})\right\|_{\mathcal{H}^{s-2}}^2+\big\|\partial_t^l\partial_x^{\alpha} R_E^{\mathfrak{c}}\big\|^2+\big\|\partial_t^l\partial_x^{\alpha} R_B^{\mathfrak{c}}\big\|^2,
\end{equation}
where we have used the exact cancellation of the curl operators:
\[
\int_{\mathbb{R}^3}\left(-\mathfrak{c} \nabla \times \partial_t^l\partial_x^\alpha \mathcal{B}^{\mathfrak{c}} \cdot \partial_t^l\partial_x^\alpha \mathcal{E}^{\mathfrak{c}}+\mathfrak{c} \nabla \times \partial_t^l\partial_x^\alpha \mathcal{E}^{\mathfrak{c}} \cdot \partial_t^l\partial_x^\alpha \mathcal{B}^{\mathfrak{c}}\right) d x=\mathfrak{c} \int_{\mathbb{R}^3} \operatorname{div}\left(\partial_t^l\partial_x^\alpha \mathcal{E}^{\mathfrak{c}} \times \partial_t^l\partial_x^\alpha \mathcal{B}^{\mathfrak{c}}\right) d x=0.
\]
Combining \eqref{B.9} and \eqref{B.10}, and applying Gronwall's inequality, we get:
\begin{align}\label{B.11}
	\sup _{t \in [0,T]}&\left\|(\mathcal{W}^{\mathfrak{c}},\mathcal{E}^{\mathfrak{c}},\mathcal{B}^{\mathfrak{c}})(t)\right\|_{\mathcal{H}^{s-2}}^2
	\leq e^{C(C_0)T}\left\|(\mathcal{W}^{\mathfrak{c}},\mathcal{E}^{\mathfrak{c}},\mathcal{B}^{\mathfrak{c}})(0)\right\|_{\mathcal{H}^{s-2}}^2\nonumber\\
	&+\int_0^{T}e^{C(C_0)(T-\tau)}\big\|(H^{\mathfrak{c}},R_E^{\mathfrak{c}},R_B^{\mathfrak{c}})(\tau)\big\|_{\mathcal{H}^{s-2}}^2 d \tau .
\end{align}
By the definition of $(H^{\mathfrak{c}}, R_E^{\mathfrak{c}}, R_B^{\mathfrak{c}})$, \eqref{e.13}, and \eqref{e.24}, we deduce:
\begin{align}\label{B.12}
	 \int_0^{T}\big\|(H^{\mathfrak{c}},R_E^{\mathfrak{c}},R_B^{\mathfrak{c}})(\tau)\big\|_{\mathcal{H}^{s-2}}^2 d \tau \leq C\left(C_0\right) T \mathfrak{c}^{-2}.
\end{align}
Substituting \eqref{B.12} into \eqref{B.11} and using \eqref{e.23}, we get:
\begin{align*}
	\sup _{t \in [0,T]}&\left\|(\mathcal{W}^{\mathfrak{c}},\mathcal{E}^{\mathfrak{c}},\mathcal{B}^{\mathfrak{c}})(t)\right\|_{\mathcal{H}^{s-2}}
	\leq C(C_0,C_1,T)\mathfrak{c}^{-1}.
\end{align*}
To strictly close the \textit{a priori} assumption \eqref{e.24}, we choose $\mathfrak{c}$ large enough such that
$$
 C(C_0, C_1,T)\mathfrak{c}^{-1} \leq \frac{1}{2}.
$$
Denote $\mathfrak{c}_0(T, C_0, C_1) := 2 C(C_0, C_1, T)$, therefore, for all $\mathfrak{c} \geq \mathfrak{c}_0$, the \textit{a priori} assumption \eqref{e.24} is fully justified, and the local solution is successfully extended to the whole time interval $[0, T]$. Thus we complete the proof of Lemma \ref{t3.6}.
\end{proof}

Combining the background profiles from Lemmas \ref{l3.2}-\ref{l3.3} and the uniform remainder estimates from Lemma \ref{t3.6}, we obtain the rigorous Newtonian limit of rEM \eqref{e.4}.
\begin{theorem}\label{c3.4}
Suppose the assumptions in Lemmas \ref{l3.2}-\ref{l3.3} and Lemma \ref{t3.6} hold. Then, for any fixed finite time $T>0$, there exists a sufficiently large constant $\mathfrak{c}_0 > 0$ such that for all $\mathfrak{c} \geq \mathfrak{c}_0$, rEM \eqref{e.4} admits a unique classical solution $(n_{\mathfrak{e}}, u_{\mathfrak{e}}, E_{\mathfrak{e}}, B_{\mathfrak{e}})$ on $[0,T] \times \mathbb{R}^3$ satisfying: 
\begin{align}\label{e.29-1}
	\sup_{t \in [0,T]} \left\|(n_{\mathfrak{e}}-\bar n, u_{\mathfrak{e}}, E_{\mathfrak{e}}, B_{\mathfrak{e}})(t)\right\|_{\mathcal{H}^{s-2}} \leq C\epsilon+C\mathfrak{c}^{-1}.
\end{align}

Moreover, as $\mathfrak{c} \to \infty$, the rEM solution converges uniformly to the solution $(n_0, u_0, E_0)$ of the classical EP \eqref{e.9} on $[0,T] \times \mathbb{R}^3$, with the convergence rate:
\begin{equation}\label{e.29}
	\sup_{t \in [0,T]} \Big[ \|n_{\mathfrak{e}} - n_0\|_{L^\infty} + \|u_{\mathfrak{e}} - u_0\|_{L^\infty} + \|E_{\mathfrak{e}} - E_0\|_{L^\infty} + \|B_{\mathfrak{e}}\|_{L^\infty} \Big] \leq C\mathfrak{c}^{-1}.
\end{equation}

In addition, denote $\bar T:=\frac{1}{2 \pi} e^{\frac{2}{3} S_\mathfrak{e}-\frac{5}{3}} \bar{n}^{\frac{2}{3}}$, then it holds that
\begin{equation}\label{e.26}
\sup_{t\in[0,T]}\|T_{\mathfrak{e}}-\bar T\|_{\mathcal{H}^{s-2}}\le C\epsilon+C\mathfrak{c}^{-1},
\qquad
\sup_{t\in[0,T]}\|T_{\mathfrak{e}}-T_0\|_{L^\infty}\le C\mathfrak{c}^{-1}.
\end{equation}
Here the constants $C>0$ above are all uniform in $\mathfrak{c}$.
\end{theorem}

\begin{proof}
By the definition \eqref{e.14} of the remainder, the exact solution takes the form
\[
	(n_{\mathfrak{e}}, u_{\mathfrak{e}}, E_{\mathfrak{e}}, B_{\mathfrak{e}}) = (n_0, u_0, E_0, 0) + \frac{1}{\mathfrak{c}}(n_1, u_1, E_1, B_1) + (\mathcal{N}^{\mathfrak{c}}, \mathcal{U}^{\mathfrak{c}}, \mathcal{E}^{\mathfrak{c}}, \mathcal{B}^{\mathfrak{c}}).
\]
Then \eqref{e.29-1} follows directly from Lemmas \ref{l3.2}--\ref{l3.3} and the remainder estimate \eqref{e.27} in Lemma \ref{t3.6}, while \eqref{e.29} follows from Lemma \ref{l3.3}, \eqref{e.27}, and the Sobolev's embedding $H^{s-2}(\mathbb R^3)\hookrightarrow L^\infty(\mathbb R^3)$. 

The bounds \eqref{e.26} follow from the fact that $T_{\mathfrak{e}}$ depends smoothly on $n_{\mathfrak{e}}$ through \eqref{e0.4}. Indeed, the uniform $H^{s-2}$ bound of $T_{\mathfrak{e}}$ is a consequence of the corresponding bound for $n_{\mathfrak{e}}$ and standard composition estimates, while the convergence estimate
follows from the convergence of $n_{\mathfrak{e}}$ to $n_0$ together with \eqref{e.7} and \eqref{e.8}.
\end{proof}

As a direct consequence of the estimates \eqref{e.29-1}-\eqref{e.26}, the density and temperature stay uniformly away from vacuum, and the velocity remains uniformly bounded. Specifically, there exists a uniform constant $c_0 \in (0, 1)$, independent of $\mathfrak{c}$, such that for all $(t,x) \in [0,T] \times \mathbb{R}^3$, there holds:
\begin{equation}\label{e.30}
	|u_{\mathfrak{e}}(t, x)| \leq \frac{1}{c_0}, \quad c_0 \leq n_{\mathfrak{e}}(t, x) \leq \frac{1}{c_0}, \quad c_0 \leq T_{\mathfrak{e}}(t, x) \leq \frac{1}{c_0}.
\end{equation}

	\section{Uniform-in-$\mathfrak{c}$ Estimates on the Linear Part of Hilbert Expansion} \label{section 4}
	
	Denote $\overline{\mathbf{M}}_{\mathfrak{c}}$ as the local Maxwellian in the rest frame:
	$$
	\overline{\mathbf{M}}_{\mathfrak{c}}(t, x, p):=\frac{n_{\mathfrak{e}} \gamma}{4 \pi \mathfrak{c}^3 K_2(\gamma)} \exp \left\{\frac{-\mathfrak{c} p^0}{T_\mathfrak{e}}\right\} .
	$$
	The corresponding Lorentz transformation is
	\begin{align*}
	\bar{\Lambda}=\left(\bar{\Lambda}_\nu^\mu\right)=\left(\begin{array}{cccc}
    \tilde{r} & \frac{\tilde{r} v_1}{\mathfrak{c}} & \frac{\tilde{r} v_2}{\mathfrak{c}} & \frac{\tilde{r} v_3}{\mathfrak{c}} \\
    \frac{\tilde{r} v_1}{\mathfrak{c}} & 1+(\tilde{r}-1) \frac{v_1^2}{|v|^2} & (\tilde{r}-1) \frac{v_1 v_2}{|v|^2} & (\tilde{r}-1) \frac{v_1 v_3}{|v|^2} \\
    \frac{\tilde{r} v_2}{\mathfrak{c}} & (\tilde{r}-1) \frac{v_1 v_2}{|v|^2} & 1+(\tilde{r}-1) \frac{v_2^2}{|v|^2} & (\tilde{r}-1) \frac{v_2 v_3}{|v|^2} \\
    \frac{\tilde{r} v_3}{\mathfrak{c}} & (\tilde{r}-1) \frac{v_1 v_3}{|v|^2} & (\tilde{r}-1) \frac{v_2 v_3}{|v|^2} & 1+(\tilde{r}-1) \frac{v_3^2}{|v|^2}
    \end{array}\right),
	\end{align*}
	where $\tilde{r}=\frac{u_{\mathfrak{e}}^0}{\mathfrak{c}}, v_i=\frac{\mathfrak{c} u_{\mathfrak{e},i}}{u_{\mathfrak{e}}^0}$.
	A direct calculation shows that
    $$
    \left(u_{\mathfrak{e}}^0, u_{\mathfrak{e}}^1, u_{\mathfrak{e}}^2, u_{\mathfrak{e}}^3\right)^t=\bar{\Lambda}(\mathfrak{c}, 0,0,0)^t .
    $$
    
	Define the third momentum
	$$
	T^{\alpha \beta \gamma}\left[\mathbf{M}_{\mathfrak{c}}\right]:=\int_{\mathbb{R}^3} \frac{p^\alpha p^\beta p^\gamma}{p^0} \mathbf{M}_{\mathfrak{c}} d p, \quad \bar{T}^{\alpha \beta \gamma}:=\int_{\mathbb{R}^3} \frac{p^\alpha p^\beta p^\gamma}{p^0} \overline{\mathbf{M}}_{\mathfrak{c}} d p .
	$$
	We first give the expression of $\bar{T}^{\alpha \beta \gamma}$ which can be proved directly and we omit the details here for brevity.
	\begin{lemma}\label{l5.1}\emph{(\cite{Wang-Xiao-JLMS-2026})}
		Let $i, j, k \in\{1,2,3\}$. It holds that
		$$
		\begin{aligned}
			\bar{T}^{000} & =\frac{n_{\mathfrak{e}} \mathfrak{c}^2\left[3 K_3(\gamma)+\gamma K_2(\gamma)\right]}{\gamma K_2(\gamma)}, \\
			\bar{T}^{0 i i} & =\bar{T}^{i i 0}=\bar{T}^{i 0 i}=\frac{n_{\mathfrak{e}} \mathfrak{c}^2 K_3(\gamma)}{\gamma K_2(\gamma)}, \\
			\bar{T}^{\alpha \beta \gamma} & =0, \quad \text { if }(\alpha, \beta, \gamma) \neq(0,0,0),(0, i, i),(i, i, 0),(i, 0, i) .
		\end{aligned}
		$$
	\end{lemma}
	Observing that
	$$
	T^{\alpha \beta \gamma}\left[\mathbf{M}_{\mathfrak{c}}\right]=\bar{\Lambda}_{\alpha^{\prime}}^\alpha \bar{\Lambda}_{\beta^{\prime}}^\beta \bar{\Lambda}_{\gamma^{\prime}}^\gamma \bar{T}^{\alpha^{\prime} \beta^{\prime} \gamma^{\prime}},
	$$
	we can obtain the expression of $T^{\alpha \beta \gamma}\left[\mathbf{M}_{\mathfrak{c}}\right]$ from Lemma \ref{l5.1}.
	\begin{lemma}\label{l5.2}\emph{(\cite{Wang-Xiao-JLMS-2026})}
		For $i, j, k \in\{1,2,3\}$, there hold
		$$
		\begin{aligned}
			T^{000}\left[\mathbf{M}_{\mathfrak{c}}\right] & =\frac{n_{\mathfrak{e}}}{\mathfrak{c} \gamma K_2(\gamma)}\left[\left(3 K_3(\gamma)+\gamma K_2(\gamma)\right)\left(u_{\mathfrak{e}}^0\right)^3+3 K_3(\gamma) u_{\mathfrak{e}}^0|u_{\mathfrak{e}}|^2\right], \\
			T^{00 i}\left[\mathbf{M}_{\mathfrak{c}}\right] & =\frac{n_{\mathfrak{e}}}{\mathfrak{c} \gamma K_2(\gamma)}\left[\left(5 K_3(\gamma)+\gamma K_2(\gamma)\right)\left(u_{\mathfrak{e}}^0\right)^2 u_{\mathfrak{e},i}+K_3(\gamma)|u_{\mathfrak{e}}|^2 u_{\mathfrak{e},i}\right], \\
			T^{0 i j}\left[\mathbf{M}_{\mathfrak{c}}\right] & =\frac{n_{\mathfrak{e}}}{\mathfrak{c} \gamma K_2(\gamma)}\left[\left(6 K_3(\gamma)+\gamma K_2(\gamma)\right) u_{\mathfrak{e}}^0 u_{\mathfrak{e},i} u_{\mathfrak{e},j}+\mathfrak{c}^2 K_3(\gamma) u_{\mathfrak{e}}^0 \delta_{i j}\right], \\
			T^{i j k}\left[\mathbf{M}_{\mathfrak{c}}\right] & =\frac{n_{\mathfrak{e}}}{\mathfrak{c} \gamma K_2(\gamma)}\left[\left(6 K_3(\gamma)+\gamma K_2(\gamma)\right) u_{\mathfrak{e},i} u_{\mathfrak{e},j} u_{\mathfrak{e},k}+\mathfrak{c}^2 K_3(\gamma)\left(u_{\mathfrak{e},i} \delta_{j k}+u_{\mathfrak{e},j} \delta_{i k}+u_{\mathfrak{e},k} \delta_{i j}\right)\right] .
		\end{aligned}
		$$
	\end{lemma}
	\subsection{Reformulation of $F_{n+1}^{\mathfrak{c}}$}
	For $n=0,1, \cdots, 2 k-2$, we decompose $F_{n+1}^{\mathfrak{c}}$ as
	$$
	\frac{F_{n+1}^{\mathfrak{c}}}{\sqrt{\mathbf{M}_{\mathfrak{c}}}}=\mathbf{P}_{\mathfrak{c}}\Big(\frac{F_{n+1}^{\mathfrak{c}}}{\sqrt{\mathbf{M}_{\mathfrak{c}}}}\Big)+\left\{\mathbf{I}-\mathbf{P}_{\mathfrak{c}}\right\}\Big(\frac{F_{n+1}^{\mathfrak{c}}}{\sqrt{\mathbf{M}_{\mathfrak{c}}}}\Big),
	$$
	where
	\begin{equation}\label{N.1}
		\mathbf{P}_{\mathfrak{c}}\Big(\frac{F_{n+1}^{\mathfrak{c}}}{\sqrt{\mathbf{M}_{\mathfrak{c}}}}\Big)=\Big[a_{n+1}+b_{n+1} \cdot p+c_{n+1} \frac{p^0}{\mathfrak{c}}\Big] \sqrt{\mathbf{M}_{\mathfrak{c}}} .
	\end{equation}	
	
	 Notice from \eqref{e1.9} that
	\begin{align}\label{N.2}
		\left\{\begin{aligned}
			&\partial_t F_{n+1}^{\mathfrak{c}}+\hat{p} \cdot \nabla_x F_{n+1}^{\mathfrak{c}}=\sum_{\substack{i+j=n+2 \\
					i, j \geq 0}} Q_{\mathfrak{c}}\left(F_i^{\mathfrak{c}}, F_j^{\mathfrak{c}}\right)+ \sum_{\substack{i+j=n+1 \\
					i,j \geq 0}}\Big(E_i^{\mathfrak{c}}+\frac{p}{p^0} \times B_i^{\mathfrak{c}}\Big) \cdot \nabla_p F_j^{\mathfrak{c}}, \\
			&\partial_t E_{n+1}^{\mathfrak{c}}-\mathfrak{c} \nabla_x \times B_{n+1}^{\mathfrak{c}}=4 \pi\int_{\mathbb{R}^3} \hat{p} F_{n+1}^{\mathfrak{c}} d p, \\
			&\partial_t B_{n+1}^{\mathfrak{c}}+\mathfrak{c} \nabla_x \times E_{n+1}^{\mathfrak{c}}=0, \\
			&\operatorname{div} E_{n+1}^{\mathfrak{c}}=-4\pi \int_{\mathbb{R}^3} F_{n+1}^{\mathfrak{c}} d p, \\
			&\operatorname{div} B_{n+1}^{\mathfrak{c}}=0.
		\end{aligned}\right.
	\end{align}
	With the help of \eqref{N.8}-\eqref{N.10} governing the macroscopic quantities $(a_{n+1},b_{n+1},c_{n+1})$, we have
	\begin{equation}\label{N.11}
		\mathbf{A}_0 \partial_t U_{n+1}+\sum_{i=1}^3 \mathbf{A}_i \partial_i U_{n+1}+\mathbf{B_1} U_{n+1}+\mathbf{B_2} \bar{U}_{n+1}=\mathbf{S}_{n+1},
	\end{equation}
	where
	$$
	U_{n+1}=\left(\begin{array}{c}
		a_{n+1} \\
		b_{n+1} \\
		c_{n+1}
	\end{array}\right), \quad \bar{U}_{n+1}=\left(\begin{array}{c}
		E_{n+1}^{\mathfrak{c}} \\
		B_{n+1}^{\mathfrak{c}}
	\end{array}\right) , \quad
	\mathbf{S}_{n+1}=\left(\begin{array}{c}
		\mathbf{S}_{n+1}^{1}\\
		\mathbf{S}_{n+1}^2 \\
		\mathbf{S}_{n+1}^3
	\end{array}\right),
	$$
    \begin{align*}
        \mathbf{S}_{n+1}^{1}=-\operatorname{div} \Big\{\int_{\mathbb{R}^{3}} \hat{p} \sqrt{\mathbf{M}_{\mathfrak{c}}}\left\{\mathbf{I}-\mathbf{P}_{\mathfrak{c}}\right\}\Big(\frac{F_{n+1}^{\mathfrak{c}}}{\sqrt{\mathbf{M}_{\mathfrak{c}}}}\Big) d p \Big\},
    \end{align*}
	{\small
	\begin{align*}
		\mathbf{S}_{n+1,j}^2=&-\sum_{\substack{k+l=n+1 \\
				k, l \geq 1}} E_{k, j}^{\mathfrak{c}}\Bigg\{\frac{n_{\mathfrak{e}} u_{\mathfrak{e}}^0}{\mathfrak{c}} a_l+\frac{e_{\mathfrak{e}}+P_{\mathfrak{e}}}{\mathfrak{c}^3} u_{\mathfrak{e}}^0\left(u_{\mathfrak{e}} \cdot b_l\right)+\frac{e_{\mathfrak{e}}\left(u_{\mathfrak{e}}^0\right)^2+P_{\mathfrak{e}}|u_{\mathfrak{e}}|^2}{\mathfrak{c}^4} c_l\Bigg\} \\
		&-\sum_{\substack{k+l=n+1 \\
				k, l \geq 1}}\Big\{\Big(\frac{n_{\mathfrak{e}} u_{\mathfrak{e}}}{\mathfrak{c}} a_l+\frac{e_{\mathfrak{e}}+P_{\mathfrak{e}}}{\mathfrak{c}^3} u_{\mathfrak{e}}\left(u_{\mathfrak{e}} \cdot b_l\right)+\frac{P_{\mathfrak{e}}}{\mathfrak{c}} b_l+\frac{e_{\mathfrak{e}}+P_{\mathfrak{e}}}{\mathfrak{c}^4} u_{\mathfrak{e}}^0 u_{\mathfrak{e}} c_l\Big) \times B_k^{\mathfrak{c}}\Big\}_j \\
		&-\operatorname{div} \int_{\mathbb{R}^3} p_j \hat{p} \sqrt{\mathbf{M}_{\mathfrak{c}}}\left\{\mathbf{I}-\mathbf{P}_{\mathfrak{c}}\right\}\Big(\frac{F_{n+1}^{\mathfrak{c}}}{\sqrt{\mathbf{M}_{\mathfrak{c}}}}\Big) d p+\int_{\mathbb{R}^3}\Big(\frac{p}{p^0} \times B_0^{\mathfrak{c}}\Big)_j \sqrt{\mathbf{M}_{\mathfrak{c}}}\left\{\mathbf{I}-\mathbf{P}_{\mathfrak{c}}\right\}\Big(\frac{F_{n+1}^{\mathfrak{c}}}{\sqrt{\mathbf{M}_{\mathfrak{c}}}}\Big) d p \\
		&-\sum_{\substack{k+l=n+1 \\
				k, l \geq 1}} \int_{\mathbb{R}^3}\Big(\frac{p}{p^0} \times B_k^{\mathfrak{c}}\Big)_j \sqrt{\mathbf{M}_{\mathfrak{c}}}\left\{\mathbf{I}-\mathbf{P}_{\mathfrak{c}}\right\}\Big(\frac{F_l^{\mathfrak{c}}}{\sqrt{\mathbf{M}_{\mathfrak{c}}}}\Big) d p , \quad j=1,2,3,
	\end{align*}}
	and 
	{\small
	\begin{align*}
		\mathbf{S}_{n+1}^3=&-\int_{\mathbb{R}^3} \frac{p}{\mathfrak{c} p^0} \cdot E_0^{\mathfrak{c}} \sqrt{\mathbf{M}_{\mathfrak{c}}}\left\{\mathbf{I}-\mathbf{P}_{\mathfrak{c}}\right\}\Big(\frac{F_{n+1}^{\mathfrak{c}}}{\sqrt{\mathbf{M}_{\mathfrak{c}}}}\Big) d p \\
		& -\sum_{\substack{k+l=n+1 \\
				k, l \geq 1}}\Big\{\frac{n_{\mathfrak{e}} a_l}{\mathfrak{c}^2} u_{\mathfrak{e}} \cdot E_k^{\mathfrak{c}}+\frac{e_{\mathfrak{e}}+P_{\mathfrak{e}}}{\mathfrak{c}^4}\left(u_{\mathfrak{e}} \cdot b_l\right)\left(E_k^{\mathfrak{c}} \cdot u_{\mathfrak{e}}\right)+\frac{P_{\mathfrak{e}}}{\mathfrak{c}^2} E_k^{\mathfrak{c}} \cdot b_l+\frac{e_{\mathfrak{e}}+P_{\mathfrak{e}}}{\mathfrak{c}^5} u_{\mathfrak{e}}^0 c_l\left(u_{\mathfrak{e}} \cdot E_k^{\mathfrak{c}}\right) \\
		& -\int_{\mathbb{R}^3} \frac{p \cdot E_k^{\mathfrak{c}}}{\mathfrak{c} p^0} \sqrt{\mathbf{M}_{\mathfrak{c}}}\left\{\mathbf{I}-\mathbf{P}_{\mathfrak{c}}\right\}\Big(\frac{F_l^{\mathfrak{c}}}{\sqrt{\mathbf{M}_{\mathfrak{c}}}}\Big) d p\Big\}.
	\end{align*}}
	
	Recall that $h=\frac{e_{\mathfrak{e}}+P_{\mathfrak{e}}}{n_{\mathfrak{e}}}=\mathfrak{c}^2 \frac{K_3(\gamma)}{K_2(\gamma)}$ in \eqref{e.1}. We denote
	$$
	h_1(t, x):=\frac{n_{\mathfrak{e}}}{\gamma K_2(\gamma)}\left(6 K_3(\gamma)+\gamma K_2(\gamma)\right), \quad h_2(t, x):=\frac{n_{\mathfrak{e}} K_3(\gamma)}{\gamma K_2(\gamma)} .
	$$	
	Then the matrices $\mathbf{A}_i(i=0,1,2,3)$ in \eqref{N.11} are
	$$
	\mathbf{A}_0=\left(\begin{array}{ccc}
		\frac{n_{\mathfrak{e}} u_{\mathfrak{e}}^0}{\mathfrak{c}} & \frac{n_{\mathfrak{e}} u_{\mathfrak{e}}^0 h u_{\mathfrak{e}}^t}{\mathfrak{c}^3} & \frac{e_{\mathfrak{e}}\left(u_{\mathfrak{e}}^0\right)^2+P_{\mathfrak{e}}|u_{\mathfrak{e}}|^2}{\mathfrak{c}^4} \\
		\frac{n_{\mathfrak{e}} u_{\mathfrak{e}}^0 h u_{\mathfrak{e}}^t}{\mathfrak{c}^3} & \left(\frac{h_1}{\mathfrak{c}} u_{\mathfrak{e}} \otimes u_{\mathfrak{e}}+\mathfrak{c} h_2 \mathbf{I}\right) u_{\mathfrak{e}}^0 & \left(\frac{h_1}{\mathfrak{c}^2}\left(u_{\mathfrak{e}}^0\right)^2-h_2\right) u_{\mathfrak{e}} \\
		\frac{e_{\mathfrak{e}}\left(u_{\mathfrak{e}}^0\right)^2+P_{\mathfrak{e}}|u_{\mathfrak{e}}|^2}{\mathfrak{c}^4} & \left(\frac{h_1}{\mathfrak{c}^2}\left(u_{\mathfrak{e}}^0\right)^2-h_2\right) u_{\mathfrak{e}}^t & \left(\frac{h_1}{\mathfrak{c}^3}\left(u_{\mathfrak{e}}^0\right)^2-\frac{3 h_2}{\mathfrak{c}}\right) u_{\mathfrak{e}}^0
	\end{array}\right)
	$$
	and
	$$
	\mathbf{A}_i=\left(\begin{array}{ccc}
		n_{\mathfrak{e}} u_{\mathfrak{e},i} & \frac{1}{\mathfrak{c}^2} n_{\mathfrak{e}} h u_{\mathfrak{e},i} u_{\mathfrak{e}}^t+P_{\mathfrak{e}} \mathbf{e}_i^t & \frac{1}{\mathfrak{c}^3} n_{\mathfrak{e}} h u_{\mathfrak{e}}^0 u_{\mathfrak{e},i} \\
		\frac{1}{\mathfrak{c}^2} n_{\mathfrak{e}} h u_{\mathfrak{e},i} u_{\mathfrak{e}}+P_{\mathfrak{e}} \mathbf{e}_i & h_1 u_{\mathfrak{e},i} u_{\mathfrak{e}} \otimes u_{\mathfrak{e}}+\mathfrak{c}^2 h_2\left(u_{\mathfrak{e},i} \mathbf{I}+\tilde{\mathbf{A}}_i\right) & \left(\frac{h_1}{\mathfrak{c}} u_{\mathfrak{e},i} u_{\mathfrak{e}}+\mathfrak{c} h_2 \mathbf{e}_i\right) u_{\mathfrak{e}}^0 \\
		\frac{1}{\mathfrak{c}^3} n_{\mathfrak{e}} h u_{\mathfrak{e}}^0 u_{\mathfrak{e},i} & \left(\frac{h_1}{\mathfrak{c}} u_{\mathfrak{e},i} u_{\mathfrak{e}}^t+\mathfrak{c} h_2 \mathbf{e}_i^t\right) u_{\mathfrak{e}}^0 & \left(\frac{h_1}{\mathfrak{c}^2}\left(u_{\mathfrak{e}}^0\right)^2-h_2\right) u_{\mathfrak{e},i}
	\end{array}\right),
	$$
	where
	$$
	\left(\tilde{\mathbf{A}}_i\right)_{j k}=\delta_{i j} u_{\mathfrak{e},k}+\delta_{i k} u_{\mathfrak{e},j}, \quad 1 \leq j, k \leq 3 .
	$$
	The matrix $\mathbf{B_1}=\left(b_{i j}\right)$ has the form
	{\small
	\begin{align*}
		& b_{11}=0, \quad\left(b_{12}, b_{13}, b_{14}\right)=\frac{n_{\mathfrak{e}} u_{\mathfrak{e}}^0}{\mathfrak{c}} \partial_t\Big(\frac{h u_{\mathfrak{e}}^t}{\mathfrak{c}^2}\Big)+n_{\mathfrak{e}} u_{\mathfrak{e}}^t\Big[\nabla_x\Big(\frac{h u_{\mathfrak{e}}}{\mathfrak{c}^2}\Big)\Big]^t+\left(\nabla_x P_{\mathfrak{e}}\right)^t, \\
		& b_{15}=\frac{n_{\mathfrak{e}} u_{\mathfrak{e}}^0}{\mathfrak{c}^2} \partial_t\Big(\frac{h u_{\mathfrak{e}}^0}{\mathfrak{c}^2}\Big)+\frac{n_{\mathfrak{e}} u_{\mathfrak{e}}}{\mathfrak{c}} \cdot \nabla_x\Big(\frac{h u_{\mathfrak{e}}^0}{\mathfrak{c}^2}\Big)-\partial_t\Big(\frac{P_{\mathfrak{e}}}{\mathfrak{c}^2}\Big), \\
		& \left(b_{21}, b_{31}, b_{41}\right)=\frac{n_{\mathfrak{e}} u_{\mathfrak{e}}^0}{\mathfrak{c}} \partial_t\Big(\frac{h u_{\mathfrak{e}}}{\mathfrak{c}^2}\Big)+\nabla_x P_{\mathfrak{e}}+\nabla_x\Big(\frac{h u_{\mathfrak{e}}}{\mathfrak{c}^2}\Big) n_{\mathfrak{e}} u_{\mathfrak{e}}+\frac{n_{\mathfrak{e}} u_{\mathfrak{e}}^0}{\mathfrak{c}} E_0^{\mathfrak{c}}+\frac{n_{\mathfrak{e}} u_{\mathfrak{e}}}{\mathfrak{c}}\times B_0^{\mathfrak{c}}, \\
		& \left(b_{j+1,2}, b_{j+1,3}, b_{j+1, 4}\right)=\frac{n_{\mathfrak{e}} u_{\mathfrak{e}}^0}{\mathfrak{c}} \partial_t\Big[\frac{h_1}{n_{\mathfrak{e}}} u_{\mathfrak{e},j} u_{\mathfrak{e}}^t+\frac{\mathfrak{c}^2 h_2}{n_{\mathfrak{e}}} \mathbf{e}_j^t\Big]+n_{\mathfrak{e}}\left(u_{\mathfrak{e}} \cdot \nabla_x\right)\Big(\frac{h_1}{n_{\mathfrak{e}}} u_{\mathfrak{e},j} u_{\mathfrak{e}}^t\Big)+n_{\mathfrak{e}} u_{\mathfrak{e}}^t \nabla_x\Big(\frac{\mathfrak{c}^2 h_2}{n_{\mathfrak{e}}}\Big) \mathbf{e}_j^t \\
		& \qquad\qquad+\left[\nabla_x\left(\mathfrak{c}^2 h_2 u_{\mathfrak{e},j}\right)\right]^t+\partial_{x_j}\left(\mathfrak{c}^2 h_2 u_{\mathfrak{e}}^t\right)+\frac{n_{\mathfrak{e}} h u_{\mathfrak{e}}^0}{\mathfrak{c}^3}E_{0,j}^{\mathfrak{c}} u_{\mathfrak{e}}^t+\frac{n_{\mathfrak{e}} h (u_{\mathfrak{e}}\times B_0^{\mathfrak{c}})_j u_{\mathfrak{e}}^t}{\mathfrak{c}^3}+\frac{P_{\mathfrak{e}} \beta_j}{\mathfrak{c}}, \quad j=1,2,3,\\
		& b_{j+1,5}=-\partial_t\left(h_2 u_{\mathfrak{e},j}\right)+\frac{n_{\mathfrak{e}} u_{\mathfrak{e}}^0}{\mathfrak{c}^2} \partial_t\Big(\frac{h_1}{n_{\mathfrak{e}}} u_{\mathfrak{e},j} u_{\mathfrak{e}}^0\Big)+\frac{n_{\mathfrak{e}}}{\mathfrak{c}} u_{\mathfrak{e}}^t \nabla_x\Big(\frac{h_1}{n_{\mathfrak{e}}} u_{\mathfrak{e},j} u_{\mathfrak{e}}^0\Big)+\partial_{x_j}\left(\mathfrak{c} h_2 u_{\mathfrak{e}}^0\right)\\
		&\qquad\qquad+\Big( \frac{n_{\mathfrak{e}} h (u_{\mathfrak{e}}^0)^2}{\mathfrak{c}^4}-\frac{P_{\mathfrak{e}}}{\mathfrak{c}^2}\Big)E_{0,j}^{\mathfrak{c}}+\Big[\frac{n_{\mathfrak{e}} h u_{\mathfrak{e}}^0(u_{\mathfrak{e}}\times B_0^{\mathfrak{c}})}{\mathfrak{c}^4}\Big]_j, \quad j=1,2,3,\\
		& b_{51}=\frac{n_{\mathfrak{e}} u_{\mathfrak{e}}^0}{\mathfrak{c}^2} \partial_t\Big(\frac{h u_{\mathfrak{e}}^0}{\mathfrak{c}^2}\Big)+\frac{n_{\mathfrak{e}} u_{\mathfrak{e}}}{\mathfrak{c}} \cdot \nabla_x\Big(\frac{h u_{\mathfrak{e}}^0}{\mathfrak{c}^2}\Big)-\partial_t\Big(\frac{P_{\mathfrak{e}}}{\mathfrak{c}^2}\Big)+\frac{n_{\mathfrak{e}}}{\mathfrak{c}^2}(u_{\mathfrak{e}}\cdot E_0^{\mathfrak{c}}), \\
		& \left(b_{52}, b_{53}, b_{54}\right)=\frac{n_{\mathfrak{e}} u_{\mathfrak{e}}^0}{\mathfrak{c}^2} \partial_t\Big(\frac{h_1}{n_{\mathfrak{e}}} u_{\mathfrak{e}}^0 u_{\mathfrak{e}}^t\Big)-\partial_t\left(h_2 u_{\mathfrak{e}}^t\right)+\frac{n_{\mathfrak{e}}}{\mathfrak{c}} u_{\mathfrak{e}}^t\Big[\nabla_x\left(\frac{h_1}{n_{\mathfrak{e}}} u_{\mathfrak{e}}^0 u_{\mathfrak{e}}\right)\Big]^t+\left(\nabla\left(\mathfrak{c} h_2 u_{\mathfrak{e}}^0\right)\right)^t \\
		&\qquad\qquad+\frac{n_{\mathfrak{e}} h}{\mathfrak{c}^4}(E_0^{\mathfrak{c}} \cdot u_{\mathfrak{e}})u_{\mathfrak{e}}^t+\frac{P_{\mathfrak{e}} (E_0^{\mathfrak{c}})^t}{\mathfrak{c}^2},\\
		& b_{55}=\frac{n_{\mathfrak{e}} u_{\mathfrak{e}}^0}{\mathfrak{c}^3} \partial_t\Big(\frac{h_1}{n_{\mathfrak{e}}}\left(u_{\mathfrak{e}}^0\right)^2\Big)-\partial_t\Big(\frac{3h_2u_{\mathfrak{e}}^0}{\mathfrak{c}}\Big)+\frac{n_{\mathfrak{e}}}{\mathfrak{c}^2} u_{\mathfrak{e}} \cdot \nabla_x\Big(\frac{h_1}{n_{\mathfrak{e}}}\left(u_{\mathfrak{e}}^0\right)^2\Big)-\operatorname{div}\left(h_2 u_{\mathfrak{e}}\right)+\frac{n_{\mathfrak{e}} h u_{\mathfrak{e}}^0}{\mathfrak{c}^5}(u_{\mathfrak{e}}\cdot E_0^{\mathfrak{c}}) ,
	\end{align*}}
	where $\beta_j:=\operatorname{row}_{j}([\times B_0^{\mathfrak{c}}])$, $j=1,2,3$ and 
	$[\times B_0^{\mathfrak{c}}]:=\left(\begin{array}{ccc}
		0 & B_{0,3}^{\mathfrak{c}} & -B_{0,2}^{\mathfrak{c}}\\
		-B_{0,3}^{\mathfrak{c}} & 0 & B_{0,1}^{\mathfrak{c}} \\
		B_{0,2}^{\mathfrak{c}} & -B_{0,1}^{\mathfrak{c}} & 0
	\end{array}\right)$.
	
	The matrix $\mathbf{B}_2$ in \eqref{N.11} is
	\begin{align*}
	\mathbf{B}_2=\left(\begin{array}{cc}
		0 & 0 \\
		\frac{n_{\mathfrak{e}} u_{\mathfrak{e}}^0}{\mathfrak{c}}\mathbf{I} & G \\
		\frac{n_{\mathfrak{e}} u_{\mathfrak{e}}^t}{\mathfrak{c}^2} & \mathbf{0}^t 
	\end{array}\right), 
	\qquad
	G=\left(\begin{array}{ccc}
		0 & -\frac{n_{\mathfrak{e}} u_{\mathfrak{e},3}}{\mathfrak{c}} & \frac{n_{\mathfrak{e}} u_{\mathfrak{e},2}}{\mathfrak{c}} \\
		\frac{n_{\mathfrak{e}} u_{\mathfrak{e},3}}{\mathfrak{c}} & 0 & -\frac{n_{\mathfrak{e}} u_{\mathfrak{e},1}}{\mathfrak{c}} \\
		-\frac{n_{\mathfrak{e}} u_{\mathfrak{e},2}}{\mathfrak{c}} & \frac{n_{\mathfrak{e}} u_{\mathfrak{e},1}}{\mathfrak{c}} & 0
	\end{array}\right).
	\end{align*}
	By the same argument as in \cite{Wang-Xiao-JLMS-2026}, $\mathbf{A}_0$ is a positive definite matrix.
	
	Now we consider the Maxwell system \eqref{N.7}, and we write the system as a linear symmetric hyperbolic system of $\left(E_{n+1}^{\mathfrak{c}}, B_{n+1}^{\mathfrak{c}}\right)$ :
	\begin{equation}\label{N.12}
		\partial_t \bar{U}_{n+1}+\sum_{i=1}^3 \bar{\mathbf{A}}_i \partial_i \bar{U}_{n+1}+\mathbf{B} U_{n+1}=\bar{S}_{n+1} .
	\end{equation}
	where 
	$$
	\mathbf{B}=-\left(\begin{array}{ccc}
		n_{\mathfrak{e}} u_{\mathfrak{e}} \quad& \frac{e_{\mathfrak{e}}+P_{\mathfrak{e}}}{\mathfrak{c}^2} u_{\mathfrak{e}}\otimes u_{\mathfrak{e}}+P_{\mathfrak{e}}I \quad&\frac{e_{\mathfrak{e}}+P_{\mathfrak{e}}}{\mathfrak{c}^3} u_{\mathfrak{e}}^0 u_{\mathfrak{e}} \\
		\mathbf{0} & \mathbf{0}_{3 \times 3} & \mathbf{0}
	\end{array}\right), \quad
	\bar{\mathbf{S}}_{n+1}=\left(\begin{array}{c}
		4\pi \int_{\mathbb{R}^3} \hat{p} \sqrt{\mathbf{M}_{\mathfrak{c}}}\left\{\mathbf{I}-\mathbf{P}_{\mathfrak{c}}\right\}\Big(\frac{F_{n+1}^{\mathfrak{c}}}{\sqrt{\mathbf{M}_{\mathfrak{c}}}}\Big) d p \\
		\mathbf{0} 
	\end{array}\right),	
	$$
	Let $\mathbf{O}$ be the $3 \times 3$ zero matrix, and define matrices:
	\begin{align*}
		&\mathbf{\bar{A}_{11}} =\left(\begin{array}{ccc}
			0 & 0 & 0 \\
			0 & 0 & \mathfrak{c} \\
			0 & -\mathfrak{c} & 0
		\end{array}\right),\quad \mathbf{\bar{A}_{12}}=\left(\begin{array}{ccc}
			0 & 0 & 0 \\
			0 & 0 & -\mathfrak{c} \\
			0 & \mathfrak{c} & 0
		\end{array}\right), \quad \mathbf{\bar{A}_{21}}=\left(\begin{array}{ccc}
			0 & 0 & -\mathfrak{c} \\
			0 & 0 & 0 \\
			\mathfrak{c} & 0 & 0
		\end{array}\right), \\
		& \mathbf{\bar{A}_{22}}=\left(\begin{array}{ccc}
			0 & 0 & \mathfrak{c} \\
			0 & 0 & 0 \\
			-\mathfrak{c} & 0 & 0
		\end{array}\right), \quad \mathbf{\bar{A}_{31}}=\left(\begin{array}{ccc}
			0 & \mathfrak{c} & 0 \\
			-\mathfrak{c} & 0 & 0 \\
			0 & 0 & 0
		\end{array}\right), \quad \mathbf{\bar{A}_{32}}=\left(\begin{array}{ccc}
			0 & -\mathfrak{c} & 0 \\
			\mathfrak{c} & 0 & 0 \\
			0 & 0 & 0
		\end{array}\right) .
	\end{align*}
	Then $\mathbf{\bar{A}_i}$ can be expressed as
	$$
	\mathbf{\bar{A}_i}=\left(\begin{array}{cc}
		\mathbf{O} & \mathbf{\bar{A}_{i1}} \\
		\mathbf{\bar{A}_{i2}} & \mathbf{O}
	\end{array}\right).
	$$
	
	\subsection{Uniform estimate on $\mathbf{L}_{\mathfrak{c}}^{-1}$.}
	
	To apply the Hilbert expansion procedure, the exponential decay property of $\mathbf{L}_{\mathfrak{c}}^{-1}$ is required.  In the non-relativistic context, the foundational results for exponential decay of $\mathbf{L}_{\mathfrak{c}}^{-1}$ was established in \cite{Jiang-CMP-2025}. For the relativistic case, Wang-Xiao \cite{Wang-Xiao-JLMS-2026} first established the following estimate on $\mathbf{L}_{\mathfrak{c}}^{-1}$. 
	\begin{lemma}\emph{(Exponential decay of $\mathbf{L}_{\mathfrak{c}}^{-1}$,\cite{Wang-Xiao-JLMS-2026})}\label{l2.10}
	For any fixed $0 \leq \lambda<1, \mathfrak{k}>\frac{3}{2}$, we assume 
		\begin{align*}
		    g \in \mathcal{N}_{\mathfrak{c}}^{\perp}
		    \quad \text{and} \quad
		    \big\|(1+|p|)^\mathfrak{k} \mathbf{M}_{\mathfrak{c}}^{-\frac{\lambda}{2}} g\big\|_{L^{\infty}}<\infty,
		\end{align*}
		then it holds
		\begin{equation}\label{e2.32}
			\left|\mathbf{L}_{\mathfrak{c}}^{-1} g(p)\right| \lesssim\big\|(1+|p|)^\mathfrak{k} \mathbf{M}_{\mathfrak{c}}^{-\frac{\lambda}{2}} g\big\|_{L^{\infty}} \cdot \mathbf{M}_{\mathfrak{c}}^{\frac{\lambda}{2}}(p), \quad p \in \mathbb{R}^3,
		\end{equation}
		where the constant is independent of $\mathfrak{c}$.
	\end{lemma}
	We extend this result by proving the exponential decay of parameter derivatives of  $\mathbf{L}_{\mathfrak{c}}^{-1}$ with respect to $t$, $x$, and $p$. Our analysis is based on Lemma \ref{l2.10} and the previously established commutator estimates \eqref{19} in Lemma \ref{l2.5} concerning $\mathbf{L}_{\mathfrak{c}}$ and its derivatives.
    \begin{lemma}\emph{(Exponential decay of parameter derivatives of $\mathbf{L}_{\mathfrak{c}}^{-1}$).}\label{l2.11}
		For any fixed $0<\lambda_l<\lambda_{l-1}<\cdots<\lambda_1<\lambda<1$, $\mathfrak{k}>\frac{3}{2}$, $l=1,2,\cdots,s-2$\textup{(}$s$ is a fixed integer defined in Lemma \ref{l3.2}\textup{)}, let $\alpha=\left(\alpha_0, \alpha_1, \alpha_2, \alpha_3\right)\in \mathbb{N}^4$ and $\beta=\left(\beta_1, \beta_2, \beta_3\right) \in \mathbb{N}^3$ with $|\alpha|+|\beta|=l$. Suppose $g \in \mathcal{N}_{\mathfrak{c}}^\perp$ and $\sum_{i=0}^l \sum_{\left|\alpha^{\prime}\right|+\left|\beta^{\prime}\right|=i}\big\|\langle p\rangle^{\mathfrak{k}} \mathbf{M}_{\mathfrak{c}}^{-\frac{\lambda_i}2}\partial_{\beta^{\prime}}^{\alpha^{\prime}} g\big\|_{L^{\infty}}<\infty$, then it holds
		\begin{equation}\label{e2.33}
			\left|\partial_{\beta}^{\alpha}\left[\mathbf{L}_{\mathfrak{c}}^{-1} g(t, x, p)\right]\right| \lesssim \sum_{i=0}^l \sum_{|\alpha^\prime|+|\beta^\prime|=i} \big\|\langle p\rangle^{\mathfrak{k}} \mathbf{M}_{\mathfrak{c}}^{-\frac{\lambda_i}{2}} \partial_{\beta^{\prime}}^{\alpha^{\prime}} g\big\|_{L^{\infty}} \cdot \mathbf{M}_{\mathfrak{c}}^{\frac{\lambda_l}{2}}(p),\quad p \in \mathbb{R}^3,
		\end{equation}
		where the constant is independent of $\mathfrak{c}$.
	\end{lemma}
	\begin{proof}
	Denote $\mathbf{L}_{\mathfrak{c}}^{-1} g=f(t, x, p)$, $f \in \mathcal{N}_{\mathfrak{c}}^\perp$. It follows from Lemma \ref{l2.10} that
	\begin{align}\label{F.1}
	|f(t, x, p)| \lesssim\big\|\langle p\rangle^{\mathfrak{k}} \mathbf{M}_{\mathfrak{c}}^{-\frac{\lambda}{2}} g\big\|_{L^{\infty}} \cdot \mathbf{M}_{\mathfrak{c}}^{\frac{\lambda}{2}}(p).
	\end{align}
	which yields Lemma \ref{l2.11} for $l=0$.
	We proceed to prove the bounds \eqref{e2.33} by induction on $l \geq 1$ .
	
	{\it Step 1.} $l=1$. 
	
	 For the microscopic part, note that $\partial_{\beta}^{\alpha}(\mathbf{L}_{\mathfrak{c}}f)=\partial_{\beta}^{\alpha} g$, hence
	\begin{align*}
	\mathbf{L}_{\mathfrak{c}}\left\{\left(\mathbf{I}-\mathbf{P}_{\mathfrak{c}}\right) \partial_{\beta}^{\alpha} f\right\}=\mathbf{L}_{\mathfrak{c}}\left\{\partial_{\beta}^{\alpha} f\right\}=\partial_{\beta}^{\alpha} g-\left[\partial_{\beta}^{\alpha}, \mathbf{L}_{\mathfrak{c}}\right] f.
	\end{align*}
	Noting \eqref{F.1}, Lemma \ref{l2.5},
	\begin{align*}
		\left|\left[\partial_{\beta}^{\alpha}, \mathbf{L}_{\mathfrak{c}}\right]f\right|\lesssim\big\|\langle p\rangle^{\mathfrak{k}} \mathbf{M}_{\mathfrak{c}}^{-\frac{\lambda}{2}} g\big\|_{L^{\infty}} \cdot \nu_{\mathfrak{c}}(p)\mathbf{M}_{\mathfrak{c}}^{\frac{\lambda}{2}}(p),
	\end{align*}
	and $0<\lambda_1<\lambda<1$,  one has that
	\begin{align*}
	\big\|\langle p\rangle^{\mathfrak{k}} \mathbf{M}_{\mathfrak{c}}^{-\frac{\lambda_1}{2}}\left[\partial_{\beta}^{\alpha}, \mathbf{L}_{\mathfrak{c}}\right]f\big\|_{L^{\infty}} \lesssim\big\|\langle p\rangle^{\mathfrak{k}} \mathbf{M}_{\mathfrak{c}}^{-\frac{\lambda}{2}} g \big\|_{L^{\infty}}.
	\end{align*}
	Consequently, by Lemma \ref{l2.10}, it holds
	\begin{align}\label{F.2}
	\left|\left(\mathbf{I}-\mathbf{P}_{\mathfrak{c}}\right) \partial_{\beta}^{\alpha} f(t, x, p)\right|
	\lesssim&\big\|\langle p\rangle^{\mathfrak{k}} \mathbf{M}_{\mathfrak{c}}^{-\frac{\lambda_1}{2}}\left(\partial_{\beta}^{\alpha} g-\left[\partial_{\beta}^{\alpha}, \mathbf{L}_{\mathfrak{c}}\right] f\right)\big\|_{L^{\infty}} \cdot \mathbf{M}_{\mathfrak{c}}^{\frac{\lambda_1}{2}}(p)\nonumber\\
	\lesssim&\Big\{\big\|\langle p\rangle^{\mathfrak{k}} \mathbf{M}_{\mathfrak{c}}^{-\frac{\lambda}{2}} g \big\|_{L^{\infty}}+\big\|\langle p\rangle^{\mathfrak{k}} \mathbf{M}_{\mathfrak{c}}^{-\frac{\lambda_1}{2}} \partial_{\beta}^{\alpha} g \big\|_{L^{\infty}}\Big\} \cdot \mathbf{M}_{\mathfrak{c}}^{\frac{\lambda_1}{2}}(p).
	\end{align}
	
	For the macroscopic part $\mathbf{P}_{\mathfrak{c}}(\partial_{\beta}^{\alpha} f(t,x,p))$, recalling the orthonormal basis of $\mathcal{N}_{\mathfrak{c}}$, $\chi_j^{\mathfrak{c}}(j=0,1, \cdots, 4)$ in \eqref{e2.36-0}, one gets
	\begin{align*}
		\mathbf{P}_{\mathfrak{c}} (\partial_{\beta}^{\alpha} f)=\sum_{j=0}^4\left\langle\partial_{\beta}^{\alpha} f, \chi_j^{\mathfrak{c}}\right\rangle \cdot \chi_j^{\mathfrak{c}}.
	\end{align*}
	Since $f \in \mathcal{N}_{\mathfrak{c}}^{\perp}$ and $\chi_j^{\mathfrak{c}} \in \mathcal{N}_{\mathfrak{c}}(j=0,1, \cdots, 4)$, by Leibniz's rule, it holds
	\begin{align*}
		\left\langle\partial^{\alpha} f, \chi_j^{\mathfrak{c}}\right\rangle=\partial^{\alpha}\left\langle f, \chi_j^{\mathfrak{c}}\right\rangle-\left\langle f, \partial^{\alpha} \chi_j^{\mathfrak{c}}\right\rangle=-\left\langle f, \partial^{\alpha} \chi_j^{\mathfrak{c}}\right\rangle
	\end{align*}
	when $|\alpha|=1$. While for $|\beta|=1$, it holds
	\begin{align*}
		\left\langle\partial_{\beta} f, \chi_j^{\mathfrak{c}}\right\rangle=-\left\langle f, \partial_{\beta} \chi_j^{\mathfrak{c}}\right\rangle.
	\end{align*}
	Then it follows that
	\begin{align*}
		\left|\left\langle\partial_{\beta}^{\alpha} f, \chi_j^{\mathfrak{c}}\right\rangle\right|  =\left|\left\langle f, \partial_{\beta}^{\alpha} \chi_j^{\mathfrak{c}}\right\rangle\right|   \lesssim \int_{\mathbb{R}^3}\big\|\langle p\rangle^{\mathfrak{k}} \mathbf{M}_{\mathfrak{c}}^{-\frac{\lambda}{2}} g \big\|_{L^{\infty}} \cdot\langle p\rangle^{3} \mathbf{M}_{\mathfrak{c}}^{\frac{\lambda}{2}}(p) \mathbf{M}_{\mathfrak{c}}^{\frac{1}{2}}(p) d p \lesssim\big\|\langle p\rangle^{\mathfrak{k}} \mathbf{M}_{\mathfrak{c}}^{-\frac{\lambda}{2}} g\big\|_{L^{\infty}},
	\end{align*}
	where we have used \eqref{F.1} and \eqref{e2.36} in the penultimate inequality.
	Noting $0<\lambda_1<\lambda<1$, we have
	\begin{align}\label{F.3}
		\left|\mathbf{P}_{\mathfrak{c}} \partial_{\beta}^{\alpha} f(t, x, p)\right| &\leq \sum_{j=0}^4\left|\left\langle\partial_{\beta}^{\alpha} f, \chi_j^{\mathfrak{c}}\right\rangle \chi_j^{\mathfrak{c}}\right| \lesssim\big\|\langle p\rangle^{\mathfrak{k}} \mathbf{M}_{\mathfrak{c}}^{-\frac{\lambda}{2}} g\big\|_{L^{\infty}} \cdot \langle p\rangle \mathbf{M}_{\mathfrak{c}}^{\frac{1}{2}}(p)\nonumber\\
		&\lesssim\big\|\langle p\rangle^{\mathfrak{k}} \mathbf{M}_{\mathfrak{c}}^{-\frac{\lambda}{2}} g\big\|_{L^{\infty}} \cdot \mathbf{M}_{\mathfrak{c}}^{\frac{\lambda_1}{2}}(p) ,
	\end{align}
	which, together with \eqref{F.2}, yields Lemma \ref{l2.11} for $l=1$.
	
	{\it Step 2.} We apply the induction argument.
	For any $|\alpha|+|\beta|=j$, $j \in\{1,2, \cdots, s-3\}$, assume that we have already obtained
	\begin{align}\label{F.6}
	\left|\partial_{\beta}^{\alpha} f(t, x, p)\right| \lesssim \sum_{i=0}^j \sum_{|\alpha^{\prime}|+|\beta^{\prime}|=i}\big\|\langle p\rangle^{\mathfrak{k}} \mathbf{M}_{\mathfrak{c}}^{-\frac{\lambda_i}{2}} \partial_{\beta^{\prime}}^{\alpha^{\prime}} g\big\|_{L^{\infty}} \cdot \mathbf{M}_{\mathfrak{c}}^{\frac{\lambda_j}{2}}(p),\quad p \in \mathbb{R}^3,
	\end{align}
	where the constant is independent of $\mathfrak{c}$. 
	
	Next, we proceed to prove \eqref{e2.33} for $l=j+1$.
	Noting
	\begin{align*}
	\mathbf{L}_{\mathfrak{c}}\left\{\left(\mathbf{I}-\mathbf{P}_{\mathfrak{c}}\right) \partial_{\beta}^{\alpha} f\right\}=\mathbf{L}_{\mathfrak{c}}\left\{\partial_{\beta}^{\alpha} f\right\}=\partial_{\beta}^{\alpha} g-\left[\partial_{\beta}^{\alpha}, \mathbf{L}_{\mathfrak{c}}\right] f,
	\end{align*}
	we have from Lemma \ref{l2.10} that
	\begin{align*}
	\left|\left(\mathbf{I}-\mathbf{P}_{\mathfrak{c}}\right) \partial_{\beta}^{\alpha} f(t, x, p)\right| \lesssim\Big\{\big\|\langle p\rangle^{\mathfrak{k}} \mathbf{M}_{\mathfrak{c}}^{-\frac{\lambda_{j+1}}{2}} \partial_{\beta}^{\alpha} g\big\|_{L^{\infty}}+\big\|\langle p\rangle^{\mathfrak{k}} \mathbf{M}_{\mathfrak{c}}^{-\frac{\lambda_{j+1}}{2}}\left[\partial_{\beta}^{\alpha}, \mathbf{L}_{\mathfrak{c}}\right] f\big\|_{L^{\infty}}\Big\} \cdot \mathbf{M}_{\mathfrak{c}}^{\frac{\lambda_{j+1}}{2}}(p).
	\end{align*}
	It follows from \eqref{F.6} and Lemma \ref{l2.5} that
	\begin{align}\label{F.7}
	 \left[\partial_{\beta}^{\alpha}, \mathbf{L}_{\mathfrak{c}}\right] f
	\lesssim \sum_{i=0}^{j} \sum_{|\alpha^{\prime}|+|\beta^{\prime}|=i}\big\|\langle p\rangle^{\mathfrak{k}} \mathbf{M}_{\mathfrak{c}}^{-\frac{\lambda_i}{2}} \partial_{\beta^{\prime}}^{\alpha^{\prime}} g\big\|_{L^{\infty}} \cdot\nu_{\mathfrak{c}}(p) \mathbf{M}_{\mathfrak{c}}^{\frac{\lambda_j}{2}}(p),
	\end{align}
	which yields
	\begin{align}\label{F.8}
	\left|\left(\mathbf{I}-\mathbf{P}_{\mathfrak{c}}\right) \partial_{\beta}^{\alpha} f(t, x, p)\right| \lesssim \sum_{i=0}^{j+1} \sum_{|\alpha^{\prime}|+|\beta^{\prime}|=i}\big\|\langle p\rangle^{\mathfrak{k}} \mathbf{M}_{\mathfrak{c}}^{-\frac{\lambda_i}{2}} \partial_{\beta^{\prime}}^{\alpha^{\prime}} g\big\|_{L^{\infty}} \cdot \mathbf{M}_{\mathfrak{c}}^{\frac{\lambda_{j+1}}{2}}(p)
	\end{align}
	for $|\alpha|+|\beta|=j+1$.
	
	Now we focus on $\mathbf{P}_{\mathfrak{c}} \partial_{\beta}^{\alpha} f(t, x, p)$ for $|\alpha|+|\beta|=j+1$. We divide the proof into two cases. 
	
	{\it Case 1:}
	$|\beta|= 0$. Splitting $\alpha=\tilde{\alpha}+\bar{\alpha}$ with $|\bar{\alpha}|=1$ and $|\tilde{\alpha}|=j$, one has $\partial^{\alpha}=\partial^{\tilde{\alpha}} \partial^{\bar{\alpha}}$.
	Then it holds $\langle\partial^{\bar{\alpha}} f, \chi_j^{\mathfrak{c}}\rangle=-\langle f, \partial^{\bar{\alpha}} \chi_j^{\mathfrak{c}}\rangle$, which implies that
	\begin{align*}
		\partial^{\tilde{\alpha}}\left\langle\partial^{\bar{\alpha}} f, \chi_j^{\mathfrak{c}}\right\rangle=-\partial^{\tilde{\alpha}}\left\langle f, \partial^{\bar{\alpha}} \chi_j^{\mathfrak{c}}\right\rangle,
	\end{align*}
	thus
	\begin{align*}
		\left\langle\partial^{\alpha} f, \chi_j^{\mathfrak{c}}\right\rangle=&\left\langle\partial^{\tilde{\alpha}+\bar{\alpha}} f, \chi_j^{\mathfrak{c}}\right\rangle-\partial^{\tilde{\alpha}}\left\langle\partial^{\bar{\alpha}} f, \chi_j^{\mathfrak{c}}\right\rangle-\partial^{\tilde{\alpha}}\langle f, \partial^{\bar{\alpha}} \chi_j^{\mathfrak{c}}\rangle \\
		=&-\sum_{0\neq\tilde{\alpha}_1 \leq \tilde{\alpha}} C_{\tilde{\alpha}}^{\tilde{\alpha}_1}\left\langle\partial^{\tilde{\alpha}+\bar{\alpha}-\widetilde{\alpha}_1} f, \partial^{\tilde{\alpha}_1} \chi_j^{\mathfrak{c}}\right\rangle -\sum_{\tilde{\alpha}_2 \leq \tilde{\alpha} } C_{\tilde{\alpha}}^{\tilde{\alpha}_2} \left\langle\partial^{\tilde{\alpha}-\tilde{\alpha}_2} f, \partial^{\tilde{\alpha}_2+\bar{\alpha}} \chi_j^{\mathfrak{c}}\right\rangle .\nonumber
	\end{align*}
	Noticing that $0\neq\tilde{\alpha}_1 \leq \tilde{\alpha}, \tilde{\alpha}_2 \leq \tilde{\alpha}, |\tilde{\alpha}|=j$ and $|\bar{\alpha}|=1$, one knows that
	\begin{align*}
		1\leq\left|\tilde{\alpha}+\bar{\alpha}-\tilde{\alpha}_1\right| \leq j,\quad 0\leq\left|\tilde{\alpha}-\tilde{\alpha}_2\right|\leq j .
	\end{align*}
	Then it follows from \eqref{e2.36} and \eqref{F.6} that
	\begin{align*}
	    \left|\left\langle\partial^{\alpha} f, \chi_j^{\mathfrak{c}}\right\rangle\right| \lesssim& \sum_{i=0}^j \sum_{|\alpha^{\prime}|=i}\big\|\langle p\rangle^{\mathfrak{k}} \mathbf{M}_{\mathfrak{c}}^{-\frac{\lambda_i}{2}} \partial^{\alpha^{\prime}} g\big\|_{L^{\infty}} \int_{\mathbb{R}^3}\langle p\rangle^{2|\alpha|+1} \mathbf{M}_{\mathfrak{c}}^{\frac{\lambda_j}{2}}(p) \mathbf{M}_{\mathfrak{c}}^{\frac{1}{2}}(p) d p\\
		\lesssim& \sum_{i=0}^j \sum_{|\alpha^{\prime}|=i}\big\|\langle p\rangle^{\mathfrak{k}} \mathbf{M}_{\mathfrak{c}}^{-\frac{\lambda_i}{2}} \partial^{\alpha^{\prime}} g\big\|_{L^{\infty}} .
	\end{align*}
	
	{\it Case 2:}
	$|\beta|\geq 1$. Splitting $\beta=\tilde{\beta}+\bar{\beta}$ with $|\bar{\beta}|=1$ and $|\alpha|+|\tilde{\beta}|=j$, one has $\partial_{\beta}^{\alpha}=\partial_{\bar{\beta}}\partial_{\tilde{\beta}}^\alpha$. 
	Applying integration by parts, we obtain
	\begin{align*}
	\langle\partial_{\bar{\beta}}\partial_{\tilde{\beta}}^\alpha f, \chi_j^{\mathfrak{c}}\rangle=-\langle \partial_{\tilde{\beta}}^\alpha f, \partial_{\bar{\beta}} \chi_j^{\mathfrak{c}}\rangle.
	\end{align*}
	Noting $|\alpha|+|\tilde{\beta}|=j$,
	 the induction hypothesis yields that
	\begin{align}\label{F.9}
	|\partial_{\tilde{\beta}}^\alpha f|\lesssim\sum_{i=0}^j \sum_{|\alpha^{\prime}|+|\beta^{\prime}|=i}\big\|\langle p\rangle^{\mathfrak{k}} \mathbf{M}_{\mathfrak{c}}^{-\frac{\lambda_i}{2}} \partial_{\beta^{\prime}}^{\alpha^{\prime}} g\big\|_{L^{\infty}} \cdot \mathbf{M}_{\mathfrak{c}}^{\frac{\lambda_j}{2}}(p).
	\end{align}
	With the help of estimate \eqref{e2.36}, we can easily obtain
	\begin{align}\label{F.10}
	\big| \big\langle \partial_{\tilde{\beta}}^\alpha f, \partial_{\bar{\beta}} \chi_j^{\mathfrak{c}}\big\rangle\big| \lesssim \sum_{i=0}^j \sum_{|\alpha^{\prime}|+|\beta^{\prime}|=i}\big\|\langle p\rangle^{\mathfrak{k}} \mathbf{M}_{\mathfrak{c}}^{-\frac{\lambda_i}{2}} \partial_{\beta^{\prime}}^{\alpha^{\prime}} g\big\|_{L^{\infty}}.
	\end{align}
	
	 Combining above two cases, one gets
	\begin{align}\label{F.11}
	\left|\mathbf{P}_{\mathfrak{c}} \partial_{\beta}^{\alpha} f(t, x, p)\right| \lesssim \sum_{i=0}^{j} \sum_{|\alpha^{\prime}|+|\beta^{\prime}|=i}\big\|\langle p\rangle^{\mathfrak{k}} \mathbf{M}_{\mathfrak{c}}^{-\frac{\lambda_i}{2}} \partial_{\beta^{\prime}}^{\alpha^{\prime}} g\big\|_{L^{\infty}} \cdot \mathbf{M}_{\mathfrak{c}}^{\frac{\lambda_{j+1}}{2}}(p)
	\end{align}	
	for $|\alpha|+|\beta|=j+1$. Using \eqref{F.8} and \eqref{F.11}, we obtain \eqref{e2.33} for  $l=j+1$. Therefore the proof of Lemma \ref{l2.11} is completed.
	\end{proof}
	
	\subsection{Uniform-in-$\mathfrak{c}$ estimates on $F_n^{\mathfrak{c}}$}
    \begin{proposition}\label{t5.3}
		Recall $F_0^{\mathfrak{c}}=\mathbf{M}_{\mathfrak{c}}\left(n_{\mathfrak{e}}, u_{\mathfrak{e}}, T_\mathfrak{e} ; p\right)$ in \eqref{e1.11}. Let $s$ be the sufficiently large integer defined in Lemma \ref{l3.2} and $T$ be the lifespan of Lemma \ref{t3.6}. Denote $s_0:=s-2$. Assume that for $n=1,2, \cdots,2k-1$, $\big(a_n, b_n, c_n, E_n^{\mathfrak{c}}, B_n^{\mathfrak{c}}\big)(0, x)$ be the initial data for systems \eqref{N.11}-\eqref{N.12} and satisfy
		\begin{align}\label{N.13-0}
			\big(a_n, b_n, c_n, E_n^{\mathfrak{c}}, B_n^{\mathfrak{c}}\big)(0, x) \in \mathcal H^{s_0-4n}(\mathbb{R}^3).
		\end{align}
		Then we can construct the smooth terms $F_n^{\mathfrak{c}}$ of the Hilbert expansion in $(t, x) \in[0, T] \times \mathbb{R}^3$ such that, for any $0<\lambda<1$, $n=1,2, \cdots, 2 k-1$, the following estimates hold
		\begin{align}\label{N.13}
			\sum_{l=0}^{s_0-4n} \sum_{|\alpha|+|\beta|=l}\left|\partial_{\beta}^{\alpha} F_{n}^{\mathfrak{c}}\right|\lesssim \mathbf{M}_{\mathfrak{c}}^{\frac{1+\lambda}{2}}\big(n_{\mathfrak{e}}(t,x),u_{\mathfrak{e}}(t,x),T_\mathfrak{e}(t,x);p\big), 
		\end{align}
		\begin{align}\label{N.14}
			\sup_{t\in[0,T]}\|\left(E_n^{\mathfrak{c}}, B_n^{\mathfrak{c}}\right)(t)\|_{\mathcal H^{s_0-4n}} \lesssim 1 .
		\end{align}
		 We emphasize that the constants in \eqref{N.13} and \eqref{N.14} are independent of $\mathfrak{c}$.
	\end{proposition}
	
	\begin{proof}
		For any $0<\lambda<1$ and the fixed integer $s_0$, we first define a set $\mathbf{S}_{s_0}(\lambda, 1)$ as
		\begin{align}\label{N.15}
			\mathbf{S}_{s_0}(\lambda, 1)=\{\lambda_{n l}|\lambda<\lambda_{n l}<1, \lambda_{n l}<\lambda_{n,l-1}, \lambda_{n+1,0}<\lambda_{n,l},n=1,2, \cdots, 2 k-1, 0\leq l \leq s_0-4n\}.
		\end{align}
		It is easy to see that $\mathbf{S}_{s_0}(\lambda, 1)\neq \emptyset$.
		Noting $F_n^{\mathfrak{c}}=\sqrt{\mathbf{M}_{\mathfrak{c}}}f_n^{\mathfrak{c}}$ and $\lambda_{n l}>\lambda$, one needs only to control 
		\begin{align}\label{N.16}
			\left|\partial_{\beta}^{\alpha} f_n^{\mathfrak{c}}(t, x, p)\right| \leq C(n_{\mathfrak{e}},u_{\mathfrak{e}}, T_\mathfrak{e})\mathbf{M}_{\mathfrak{c}}^{\frac{\lambda_{n l}}{2}}(t, x, p) , \quad n=1,2, \cdots, 2 k-1
		\end{align}
		for any $|\alpha|+|\beta|=l, 0 \leq l \leq s_0-4n$.
		
		We consider the following decomposition:
		\begin{align}\label{N.17}
		    f_n^{\mathfrak{c}}=\mathbf{P}_{\mathfrak{c}} f_n^{\mathfrak{c}}+\left(\mathbf{I}-\mathbf{P}_{\mathfrak{c}}\right) f_n^{\mathfrak{c}}, \quad \mathbf{P}_{\mathfrak{c}} f_n^{\mathfrak{c}}=\Big[a_n+b_n \cdot p+c_n \frac{p^0}{\mathfrak{c}}\Big] \sqrt{\mathbf{M}_{\mathfrak{c}}}.
		\end{align}
		The assertion \eqref{N.16} can be justified by induction on $n \geq 1$.
		
		{\it Step 1.}
		We begin by estimating  $f_1^{\mathfrak{c}}$.
		Denote
		\begin{align*}
		    g_1:=-\frac{1}{\sqrt{\mathbf{M}_{\mathfrak{c}}}}\Big[\left(\partial_t+\hat{p} \cdot \nabla_x\right) \mathbf{M}_{\mathfrak{c}}-\big(E_0^{\mathfrak{c}}+\frac{p}{p^0} \times B_0^{\mathfrak{c}}\big) \cdot \nabla_p \mathbf{M}_{\mathfrak{c}}\Big],
		\end{align*}
		then it follows from $\eqref{e1.8}_1$ that
		\begin{align}\label{N.19}
		\left(\mathbf{I}-\mathbf{P}_{\mathfrak{c}}\right) f_1^{\mathfrak{c}}=\mathbf{L}_{\mathfrak{c}}^{-1} g_1.
	\end{align}
	Recalling $\mathfrak{e}_l$ from \eqref{e2.34}, by Lemma \ref{l2.12} and Theorem \ref{c3.4}, a tedious computation gives
		\begin{align}\label{N.21}
			\sum_{i=0}^l \sum_{|\alpha^{\prime}|+|\beta^{\prime}|=i}\big|\partial_{\beta^{\prime}}^{\alpha^{\prime}} g_1\big| \lesssim \sum_{i=0}^l \sum_{|\alpha^{\prime}|+|\beta^{\prime}|=i}C(\mathfrak{e}_i)\langle p\rangle^{2\left|\alpha^{\prime}\right|+\left|\beta^{\prime}\right|+2} \mathbf{M}_{\mathfrak{c}}^{\frac{1}{2}} \leq C(\mathfrak{e}_l)\langle p\rangle^{2l+2} \mathbf{M}_{\mathfrak{c}}^{\frac{1}{2}}(t, x, p)
		\end{align}
		for $0 \leq l \leq s_0-1$.
		Then, the assumption $\lambda_{1i}<1$ and the Sobolev embedding inequality imply, via \eqref{N.21}, that 
		\begin{align*}
			\sum_{i=0}^l \sum_{|\alpha^{\prime}|+|\beta^{\prime}|=i}\big\|\langle p\rangle^{\mathfrak{k}} \mathbf{M}_{\mathfrak{c}}^{-\frac{\lambda_{1i}}{2}} \partial_{\beta^{\prime}}^{\alpha^{\prime}}g_1\big\|_{L^{\infty}}<\infty.
		\end{align*}
	Thus, using Lemmas \ref{l2.10}-\ref{l2.11}, one has
		\begin{align}\label{N.18}
			\left|\partial_{\beta}^{\alpha}\left(\mathbf{I}-\mathbf{P}_{\mathfrak{c}}\right) f_1^{\mathfrak{c}}\right| \lesssim \sum_{i=0}^l \sum_{|\alpha^{\prime}|+|\beta^{\prime}|=i}\big\|\langle p\rangle^{\mathfrak{k}} \mathbf{M}_{\mathfrak{c}}^{-\frac{\lambda_{1i}}{2}} \partial_{\beta^{\prime}}^{\alpha^{\prime}} g_1\big\|_{L^{\infty}} \cdot \mathbf{M}_{\mathfrak{c}}^{\frac{\lambda_{1l}}{2}}(t, x, p)\lesssim \mathbf{M}_{\mathfrak{c}}^{\frac{\lambda_{1l}}{2}}(t, x, p)
		\end{align}
	for $|\alpha|+|\beta|=l, 0 \leq l \leq s_0-1$.
		
		Now we turn to the macroscopic part $\mathbf{P}_{\mathfrak{c}} f_1^{\mathfrak{c}}$. Recall \eqref{N.11}-\eqref{N.12} for $\left(U_1, \bar{U}_1\right)$, by Theorem \ref{c3.4} and \eqref{N.19}, we get
		\begin{align}\label{N.24}
			\left\|\mathbf{S}_1\right\|_{H^{s_0-2}} \lesssim 1,\quad \left\|\bar{\mathbf{S}}_1\right\|_{H^{s_0-1}} \lesssim 1.
		\end{align}
		Note that $ \mathbf{A}_i(i=0,1,2,3)$, $\mathbf{B_1}$ and $\mathbf{B_2}$ in \eqref{N.11} depend only on the smooth functions $n_{\mathfrak{e}}(t, x)$, $u_{\mathfrak{e}}(t, x)$, $T_{\mathfrak{e}}(t, x)$, $E_0^{\mathfrak{c}}(t, x)$ and $B_0^{\mathfrak{c}}(t, x)$, then one obtains from Theorem \ref{c3.4} that
		\begin{align}\label{N.25}
			\left\|\partial_t \mathbf{A}_0\right\|_{H^{s_0-1}}+\sum_{i=0}^3\left\|\nabla_x \mathbf{A}_i\right\|_{H^{s_0-1}}+\|\mathbf{B}_1\|_{H^{s_0-1}}+\|\mathbf{B}_2\|_{H^{s_0}}+\|\mathbf{B}\|_{H^{s_0}} \lesssim 1.
		\end{align}
		Employing the classical energy estimate of symmetric hyperbolic equations to \eqref{N.11}-\eqref{N.12}, we obtain
		\begin{align}\label{N.26}
			\frac{d}{d t} \|\left(U_1, \bar{U}_1\right)(t)\|_{H^{s_0-2}}^2 \lesssim \|U_1\|_{H^{s_0-2}}^2+\|\bar{U}_1\|_{H^{s_0-2}}^2+1,
		\end{align}
		which, together with the Gronwall's inequality, yields
		$$
		\sup_{t\in[0,T]}\|\left(U_1, \bar{U}_1\right)(t)\|_{H^{s_0-2}} \lesssim 1 .
		$$
		Similarly, one also has
		\begin{align}\label{N.27}
		    \sup_{t\in[0,T]}\|\left(U_1, \bar{U}_1\right)(t)\|_{\mathcal H^{s_0-2}} \lesssim 1 .
		\end{align}
		Hence it follows from \eqref{N.17} and \eqref{N.27} that	for $|\alpha|+|\beta|=l$, $0 \leq l \leq s_0-4$
		\begin{align}\label{N.30}
			\left|\partial_{\beta}^{\alpha} (\mathbf{P}_{\mathfrak{c}} f_1^{\mathfrak{c}})\right| \leq& \sum_{\alpha^{\prime} \leq \alpha, \beta^{\prime}\leq\beta} C_\alpha^{\alpha^{\prime}} C_\beta^{\beta^{\prime}}\Big|\partial_{\beta-\beta^{\prime}}^{\alpha-\alpha^{\prime}}\Big(a_1+b_1 \cdot p+c_1 \frac{p^0}{\mathfrak{c}}\Big) \partial_{\beta^{\prime}}^{\alpha^{\prime}} \sqrt{\mathbf{M}_{\mathfrak{c}}}\Big| \nonumber\\
			\lesssim&\langle p\rangle^{2|\alpha|+|\beta|+1} \mathbf{M}_{\mathfrak{c}}^{\frac{1}{2}}(t, x, p)\lesssim \mathbf{M}_{\mathfrak{c}}^{\frac{\lambda_{1l}}{2}}(t, x, p),
		\end{align}
		where we have chosen $l \leq s_0-4$ to make the Sobolev embedding \(H^{m}(\mathbb{R}^3)\hookrightarrow L^\infty(\mathbb{R}^3)(m\ge 2)\) holds.
		Then we have from \eqref{N.18} and \eqref{N.30} that 
		\begin{align}\label{N.31}
			\left|\partial_{\beta}^{\alpha} f_1^{\mathfrak{c}}\right| \lesssim \mathbf{M}_{\mathfrak{c}}^{\frac{\lambda_{1l}}{2}}(t, x, p), \quad |\alpha|+|\beta|=l, \quad 0 \leq l \leq s_0-4.
		\end{align}
		
		{\it Step 2.}
		Assume the assertion \eqref{N.16} holds for all $n=j$, $j \in \{1,2, \cdots, \mathfrak{n}-1\}(\mathfrak{n}\geq 2)$. Namely, for any
		$|\alpha|+|\beta|=l$, $0 \leq l \leq s_0-4j$, $j=1,2, \cdots, \mathfrak{n}-1$, 
		\begin{align}\label{N.32}
			\left|\partial_{\beta}^{\alpha} f_j^{\mathfrak{c}}(t, x, p)\right| \lesssim \mathbf{M}_{\mathfrak{c}}^{\frac{\lambda_{j l}}{2}}(t, x, p),
		\end{align}
		and
		\begin{align}\label{N.28}
		    \sup_{t\in[0,T]}\|\left(U_j, \bar{U}_j\right)(t)\|_{\mathcal H^{s_0-4j}} \lesssim 1 .
		\end{align}
		
		Now we proceed to prove \eqref{N.32}-\eqref{N.28} for $n =\mathfrak{n}$.
		Denote
		\begin{gather*}
		    G_{\mathfrak{n}}:=-\Big\{\frac{1}{\sqrt{\mathbf{M}_{\mathfrak{c}}}}\Big[\left(\partial_t+\hat{p} \cdot \nabla_x\right)\left(f_{\mathfrak{n}-1}^{\mathfrak{c}} \sqrt{\mathbf{M}_{\mathfrak{c}}}\right)-\Big(E_{\mathfrak{n}-1}^{\mathfrak{c}}+\frac{p}{p^0} \times B_{\mathfrak{n}-1}^{\mathfrak{c}}\Big) \cdot \nabla_p\mathbf{M}_{\mathfrak{c}}\\
		    -\sum_{\substack{i+j=\mathfrak{n}-1 \\
				i\geq 0, j \geq 1}}\Big(E_i^{\mathfrak{c}}+\frac{p}{p^0} \times B_i^{\mathfrak{c}}\Big) \cdot \nabla_p\left(f_j^{\mathfrak{c}} \sqrt{\mathbf{M}_{\mathfrak{c}}}\right)\Big]\Big\}.
		\end{gather*}
		From $\eqref{e1.9}_1$, one has
		\begin{align}\label{N.33}
		\left(\mathbf{I}-\mathbf{P}_{\mathfrak{c}}\right) f_{\mathfrak{n}}^{\mathfrak{c}}=\mathbf{L}_{\mathfrak{c}}^{-1}\Big\{G_{\mathfrak{n}}+\sum_{\substack{i+j=\mathfrak{n} \\
				i, j\geq1}} \Gamma_{\mathfrak{c}}\left(f_i^{\mathfrak{c}}, f_j^{\mathfrak{c}}\right)\Big\} .
		\end{align}
	
		A straightforward calculation implies that
		\begin{align}\label{N.36}
			\big|\partial_{\beta^{\prime}}^{\alpha^{\prime}} G_{\mathfrak{n}}\big| \lesssim& \big|\partial_{\beta^{\prime}}^{\alpha^{\prime}} \left(\partial_t+\hat{p} \cdot \nabla_x\right)f_{\mathfrak{n}-1}^{\mathfrak{c}}\big|+\langle p\rangle^3 \sum_{\substack{\alpha^{\prime \prime} \leq \alpha^{\prime}\\
			\beta^{\prime \prime} \leq \beta^{\prime}}}\big|\partial^{\alpha^{\prime\prime}}_{\beta^{\prime \prime}} f_{\mathfrak{n}-1}^{\mathfrak{c}}\big|
			+ \sum_{\substack{\alpha^{\prime \prime} \leq \alpha^{\prime} \\
			\beta^{\prime \prime} \leq \beta^{\prime}}}\big|\partial^{\alpha^{\prime\prime}}_{\beta^{\prime \prime}} \nabla_p \mathbf{M}_{\mathfrak{c}}\big|\nonumber\\
			&+\sum_{\substack{i+j=\mathfrak{n}-1 \\
				i\geq 0, j \geq 1}}\sum_{\substack{\alpha^{\prime \prime} \leq \alpha^{\prime} \\
			\beta^{\prime \prime} \leq \beta^{\prime}}} \Big|\partial _{\beta^\prime-\beta^{\prime \prime}}^{\alpha^\prime-\alpha^{\prime\prime}}\Big(E_i^{\mathfrak{c}}+\frac{p}{p^0} \times B_i^{\mathfrak{c}}\Big)\Big|
			\Big\{\big|\partial_{\beta^{\prime \prime}}^{\alpha^{\prime\prime}}\nabla_p f_j^{\mathfrak{c}}\big|+\langle p\rangle
			 \big|\partial _{\beta^{\prime \prime}}^{\alpha^{\prime\prime}} f_j^{\mathfrak{c}} \big|\Big\}\nonumber\\
			\lesssim&\langle p\rangle^3\sum_{\substack{i+j=\mathfrak{n}-1 \\
				i\geq 0, j \geq 1}}\sum_{\substack{\alpha^{\prime \prime} \leq \alpha^{\prime} \\
			\beta^{\prime \prime} \leq \beta^{\prime}}} \Big|\partial _{\beta^\prime-\beta^{\prime \prime}}^{\alpha^\prime-\alpha^{\prime\prime}}\Big(E_i^{\mathfrak{c}}+\frac{p}{p^0} \times B_i^{\mathfrak{c}}\Big)\Big|\big|\partial_{\beta^{\prime \prime}}^{\alpha^{\prime \prime}} \partial_{t,x,p}f_j^{\mathfrak{c}}\big|.
		\end{align}
		Then it follows from \eqref{N.36} that for any $0 \leq \mathfrak{t} \leq l\leq s_0-4\mathfrak{n}+3$,
		\begin{align}\label{N.38}
			\sum_{|\alpha^{\prime}|+|\beta^{\prime}|=\mathfrak{t}}\big|\partial_{\beta^{\prime}}^{\alpha^{\prime}}G_{\mathfrak{n}}\big| \lesssim\langle p\rangle^3\sum_{1\leq j\leq \mathfrak{n}-1}\sum_{0\leq \mathfrak{r}\leq \mathfrak{t}+1}   \mathbf{M}_{\mathfrak{c}}^{\frac{\lambda_{j\mathfrak{r}}}{2}}\lesssim \mathbf{M}_{\mathfrak{c}}^{\frac{\lambda_{\mathfrak{n}-1,\mathfrak{t}+1}}{2}}(t, x, p).
		\end{align}
		Here the penultimate inequality is derived from \eqref{N.32}-\eqref{N.28}, and the last one is implied by the relation \eqref{N.15}.
		
		Next, we will control the quantities $\partial_{\beta^{\prime}}^{\alpha^{\prime}} \Gamma_{\mathfrak{c}}\left(f_i^{\mathfrak{c}}, f_j^{\mathfrak{c}}\right)$. For $0 \leq \mathfrak{t} \leq l\leq s_0-4\mathfrak{n}+4$ and $i+j=\mathfrak{n}$ with $i, j \geq 1$, the induction hypotheses and Lemma \ref{l2.5} imply that
		\begin{align}\label{N.39}
			\sum_{|\alpha^{\prime}|+|\beta^{\prime}|=\mathfrak{t}}\big|\partial_{\beta^{\prime}}^{\alpha^{\prime}} \Gamma_{\mathfrak{c}}\left(f_i^{\mathfrak{c}}, f_j^{\mathfrak{c}}\right)\big| \lesssim \nu_{\mathfrak{c}}(p) \mathbf{M}_{\mathfrak{c}}^{\frac{\lambda_{\mathfrak{n}-1,\mathfrak{t}}}{2}}\lesssim \mathbf{M}_{\mathfrak{c}}^{\frac{\lambda_{\mathfrak{n}-1,\mathfrak{t}+1}}{2}}(t,x,p)  .
		\end{align}
		It is thereby derived from  \eqref{N.38}-\eqref{N.39} and Lemmas \ref{l2.10}-\ref{l2.11} that for any $|\alpha|+|\beta|=l$, $0 \leq l\leq s_0-4\mathfrak{n}+3$,
		\begin{align}\label{N.40}
			\left|\partial_{\beta}^{\alpha}\left(\mathbf{I}-\mathbf{P}_{\mathfrak{c}}\right) f_{\mathfrak{n}}^{\mathfrak{c}}\right| \lesssim& \sum_{\mathfrak{t}=0}^l \sum_{|\alpha^{\prime}|+|\beta^{\prime}|=\mathfrak{t}}\big\|\langle p\rangle^{\mathfrak{k}} \mathbf{M}_{\mathfrak{c}}^{-\frac{\lambda_{\mathfrak{n}\mathfrak{t}}}{2}} \partial_{\beta^{\prime}}^{\alpha^{\prime}} \Big\{G_{\mathfrak{n}}+\sum_{\substack{i+j=\mathfrak{n} \\
				i, j\geq1}} \Gamma_{\mathfrak{c}}\left(f_i^{\mathfrak{c}}, f_j^{\mathfrak{c}}\right)\Big\}\big\|_{L^{\infty}}\cdot \mathbf{M}_{\mathfrak{c}}^{\frac{\lambda_{\mathfrak{n}l}}{2}} \nonumber\\
			\lesssim& \sum_{\mathfrak{t}=0}^l \sum_{|\alpha^{\prime}|+|\beta^{\prime}|=\mathfrak{t}}\Big\|\langle p\rangle^{\mathfrak{k}} \mathbf{M}_{\mathfrak{c}}^{\frac{\lambda_{\mathfrak{n}-1,\mathfrak{t}+1}- \lambda_{\mathfrak{n}\mathfrak{t}}}{2}} \Big\|_{L^{\infty}} \cdot \mathbf{M}_{\mathfrak{c}}^{\frac{\lambda_{\mathfrak{n}l}}{2}}\lesssim \mathbf{M}_{\mathfrak{c}}^\frac{\lambda_{\mathfrak{n}l}}{2}(t, x, p).
		\end{align}
		
    Next we prove the estimate \eqref{N.28} for the pair $(U_\mathfrak{n},\bar U_\mathfrak{n})$.
    Recalling the definitions of $\mathbf{S}_\mathfrak{n}$ and $\bar{\mathbf{S}}_\mathfrak{n}$ in \eqref{N.11}-\eqref{N.12}, every term in $\mathbf{S}_\mathfrak{n}$ and $\bar{\mathbf{S}}_\mathfrak{n}$ consists of either
    $\left(\mathbf{I}-\mathbf{P}_{\mathfrak{c}}\right) f_{\mathfrak{n}}^{\mathfrak{c}}$, which is controlled by \eqref{N.40}, or lower-order profiles
    $(f_j^{\mathfrak{c}},U_j,\bar U_j)$ with $1\le j\le \mathfrak{n}-1$, which are controlled by the induction
    assumptions \eqref{N.32}-\eqref{N.28}; hence we infer from the
    standard product estimate in Sobolev spaces that
    for all $0\le l\le s_0-4\mathfrak{n}+2$,
    \[
    \sup_{t\in[0,T]}
    \Big(
    \|\partial_t^l \mathbf{S}_\mathfrak{n}(t)\|_{H^{s_0-4\mathfrak{n}+2-l}}
    +\|\partial_t^l \bar{\mathbf{S}}_\mathfrak{n}(t)\|_{H^{s_0-4\mathfrak{n}+3-l}}
    \Big)\lesssim 1 .
    \]
    
Applying $\partial_t^l\partial_x^\alpha$ with $|\alpha|+l\le s_0-4\mathfrak{n}+2$ to
\eqref{N.11}-\eqref{N.12}, multiplying the resulting equations by
$\partial_t^l\partial_x^\alpha U_\mathfrak{n}$ and $\partial_t^l\partial_x^\alpha \bar U_\mathfrak{n}$,
respectively, and integrating over $\mathbb R^3$, we obtain, by the coefficient bounds in \eqref{N.25}, that
\[
\frac{d}{dt}
\| (U_\mathfrak{n},\bar U_\mathfrak{n})(t)\|_{\mathcal H^{s_0-4\mathfrak{n}+2}}^2
\lesssim 
\| (U_\mathfrak{n},\bar U_\mathfrak{n})(t)\|_{\mathcal H^{s_0-4\mathfrak{n}+2}}^2
+1.
\]
Combining this differential inequality with the corresponding well-posed initial data, and using Gronwall's inequality \eqref{N.13-0}, we conclude that
\begin{align}\label{N.41}
\sup_{t\in[0,T]}
\| (U_\mathfrak{n},\bar U_\mathfrak{n})(t)\|_{\mathcal H^{s_0-4\mathfrak{n}+2}}
\lesssim 1 ,
\end{align}
which thereby closes the induction hypothesis
\eqref{N.28}. Since $\bar U_\mathfrak{n}=(E_{\mathfrak{n}}^{\mathfrak{c}},B_{\mathfrak{n}}^{\mathfrak{c}})$, this also yields
\eqref{N.14}.

Finally, using \eqref{N.17}, the bound \eqref{N.41}, and arguing exactly as in
the proof of \eqref{N.30}, we obtain
\begin{align}\label{N.42}
\left|\partial_{\beta}^{\alpha} \mathbf{P}_{\mathfrak{c}} f_{\mathfrak{n}}^{\mathfrak{c}}\right| \lesssim \mathbf{M}_{\mathfrak{c}}^{\frac{\lambda_{\mathfrak{n}l}}{2}}(t, x, p),
\qquad |\alpha|+|\beta|=l,\ \ 0 \leq l \leq s_0-4\mathfrak{n}.
\end{align}
Similarly, here we choose  $l \leq s_0-4\mathfrak{n}$ so that the Sobolev embedding \(H^{m}(\mathbb{R}^3)\hookrightarrow L^\infty(\mathbb{R}^3)(m\ge 2)\) holds.
 Collecting the bounds \eqref{N.40} and \eqref{N.42}, one has
		\begin{align}\label{N.43}
			\left|\partial_{\beta}^{\alpha} f_{\mathfrak{n}}^{\mathfrak{c}}\right| \lesssim \mathbf{M}_{\mathfrak{c}}^{\frac{\lambda_{\mathfrak{n}l}}{2}}(t, x, p), \qquad |\alpha|+|\beta|=l,\ \ 0 \leq l \leq s_0-4\mathfrak{n}.
		\end{align}
		Therefore the proof of Proposition \ref{t5.3} is completed.
	\end{proof}


		\section{Characteristics and $L^{\infty}$-Estimates} \label{section 5}
	
	In this section, we first establish the estimates of characteristic under the  {\it a priori} assumptions \eqref{W.0}-\eqref{C.2}, then we shall derive the $L^{\infty}$ estimate of $h_R^{\varepsilon,\mathfrak{c}}$.
	
	\subsection{Characteristics estimates} 
	In this subsection, we will provide characteristics estimates for rVMB \eqref{e1.1}, which requires a uniform $W^{1,\infty}$ estimate for the electro-magnetic field $\left(E^{\varepsilon,\mathfrak{c}}, B^{\varepsilon,\mathfrak{c}}\right)(t,x)$. To this aim, we impose the following {\it a priori} assumptions:
	\begin{gather}
			\sup _{t \in[0, T]} \varepsilon^2\big\|h_R^{\varepsilon,\mathfrak{c}}(t)\big\|_{L^{\infty}} \leq \varepsilon^{\frac{1}{4}}, \quad \sup _{t \in[0, T]} \varepsilon^3\big\|\nabla_{x,p} h_R^{\varepsilon,\mathfrak{c}}(t)\big\|_{L^{\infty}}\leq \varepsilon^{\frac{1}{8}},\label{W.0}\\
	    \sup_{t \in[0, T]} \varepsilon^3\left\|\left(E_R^{\varepsilon,\mathfrak{c}}, B_R^{\varepsilon,\mathfrak{c}}\right)(t)\right\|_{L^{\infty}} \leq \varepsilon^{\frac{1}{4}}, \quad
		\sup_{t \in[0, T]} \varepsilon^4\left\|\big(\nabla_x  E_R^{\varepsilon,\mathfrak{c}}, \nabla_x B_R^{\varepsilon,\mathfrak{c}}\big)(t)\right\|_{L^{\infty}} \leq \varepsilon^{\frac{1}{8}}\label{C.2}.
	\end{gather}
	Thus, Theorem \ref{c3.4} and Proposition \ref{t5.3}, together with \eqref{C.2}, yield that for $k \geq 5$,
	\begin{align}\label{C.3}
		\sup _{t \in[0, T]}\left\|\left(E^{\varepsilon,\mathfrak{c}}, B^{\varepsilon,\mathfrak{c}}\right)(t)\right\|_{W^{1, \infty}}\leq \sum_{i=0}^{2 k-1} \varepsilon^i\left\|\left(E_i^{\mathfrak{c}}, B_i^{\mathfrak{c}}\right)(t)\right\|_{W^{1, \infty}}+\varepsilon^k\left\|\left(E_R^{\varepsilon,\mathfrak{c}}, B_R^{\varepsilon,\mathfrak{c}}\right)(t)\right\|_{W^{1, \infty}} \lesssim 1.
	\end{align}
	With the help of \eqref{C.3}, we shall study the curved trajectory of \eqref{e1.1}. We define the characteristics $\big(X(\tau; t, x, p), P(\tau; t, x, p)\big)$ passing through $(t, x, p)$ such that
	\begin{align} \label{C.4}
    \left\{\begin{aligned}
		&\frac{d X(\tau ; t, x, p)}{d \tau}=\hat{P}(\tau ; t, x, p),  \\
		&\frac{d P(\tau ; t, x, p)}{d \tau}=-E^{\varepsilon,\mathfrak{c}}(\tau, X(\tau ; t, x, p))-\frac{P}{P^0}(\tau ; t, x, p)\times B^{\varepsilon,\mathfrak{c}}(\tau, X(\tau ; t, x, p)), \\
        &X(t; t, x, p)=x, \quad P(t; t, x, p)=p.
        \end{aligned}\right.
	\end{align}
	
	For notational brevity, we denote
	$$
	X(\tau):=X(\tau; t, x, p) \quad \text{and} \quad P(\tau):=P(\tau; t, x, p).
	$$
	\begin{lemma}\label{lC.1}
		Let \eqref{C.2}-\eqref{C.3} hold. Then there exists a small constant $\bar{T} \in[0, T]$ such that for $0\leq\tau\leq t\leq \bar{T} $,
		\begin{equation}\label{C.6}
			\frac{\mathfrak{c}^5}{2\left(p^0\right)^5}|t-\tau|^3  \lesssim\left|\operatorname{det}\left(\frac{\partial X(\tau)}{\partial p}\right)\right| \lesssim \frac{2\mathfrak{c}^5}{\left(p^0\right)^5}|t-\tau|^3,
		\end{equation}
	where $\bar{T}>0$ is independent of $\mathfrak{c}$ and $\varepsilon$.
	\end{lemma}
	\begin{proof}
		We first prove that 
        $$
        \left|D_p X(\tau)\right| \lesssim \mathfrak{c}\frac{|t-\tau|}{p^0}.
        $$
        For $i, j=1,2,3$, applying $\partial_{p_i}$ to \eqref{C.4}, we have
		\begin{equation}\label{C.7}
			\frac{d \partial_{p_i} X_j(\tau ; t, x, p)}{d \tau}=\mathfrak{c}\left(\frac{\left(P^0\right)^2 \partial_{p_i} P_j-P_j\left(P \cdot \partial_{p_i} P\right)}{\left(P^0\right)^3}\right)(\tau ; t, x, p),
		\end{equation}
		and
		\begin{align}\label{C.8}
			&\frac{d \partial_{p_i} P_j(\tau ; t, x, p)}{d \tau} \nonumber\\
			=& -\left[\nabla_x E_j^{\varepsilon,\mathfrak{c}}(\tau, X(\tau ; t, x, p)) \cdot \partial_{p_i} X(\tau ; t, x, p)\right] \nonumber\\
			& -\left[\frac{P}{P^0}(\tau ; t, x, p) \times\left[\nabla_x B^{\varepsilon,\mathfrak{c}}(\tau, X(\tau ; t, x, p)) \cdot \partial_{p_i} X(\tau ; t, x, p)\right]\right]_j \\
			& -\left[\left(\frac{\left(P^0\right)^2 \partial_{p_i} P_j-P_j\left(P \cdot \partial_{p_i} P\right)}{\left(P^0\right)^3}\right)(\tau ; t, x, p) \times B^{\varepsilon,\mathfrak{c}}(\tau, X(\tau ; t, x, p))\right]_j .\nonumber
		\end{align}
		
		Using \eqref{C.3} and $\partial_{p_i} P_j(t ; t, x, p)=\delta_{i j}$, we integrate \eqref{C.8} over $[\tau, t]$ to get
		$$
		\left\|\partial_{p_i} P_j(\tau)\right\|_{L^{\infty}} \leq \delta_{i j}+C|\tau-t|\left(\sup _{\tau \in[0, \bar{T}]}\left\|\partial_{p_i} X(\tau)\right\|_{L^{\infty}}+\frac{\sup _{\tau \in[0, \bar{T}]}\left\|\partial_{p_i} P(\tau)\right\|_{L^{\infty}}}{\inf _{\tau \in[0, \bar{T}]}\left|P^0(\tau)\right|}\right) .
		$$
		Now, for sufficiently small $\bar{T}$, it holds that
		\begin{equation}\label{C.9}
			\left\|\partial_{p_i} P_j(\tau)\right\|_{L^{\infty}} \leq \frac{5}{4} \delta_{i j}+C|\tau-t| \sup _{\tau \in[0, \bar{T}]}\left\|\partial_{p_i} X(\tau)\right\|_{L^{\infty}} .
		\end{equation}
		
		 It's clear that
		\begin{align*}
			P^0(\tau)=p^0+(\tau-t)\frac{P}{P^0} \frac{d P}{d \tau}(\bar{\tau} ; t, x, p) ,
		\end{align*}
		which, together with $\eqref{C.4}_2$ and \eqref{C.3}, yields that
		\begin{align}\label{C.10}
		    \left|P^0(\tau)\right|\geq p^0-C|\tau-t|\geq\frac{p^0}{2}.
		\end{align}
		Noting $\partial_{p_i} X_j(t ; t, x, p)=0$, we integrate \eqref{C.7} over $[\tau, t]$ and use \eqref{C.9}, \eqref{C.10} to obtain
		\begin{equation}\label{C.11}
			\left|\partial_{p_i} X_j(\tau)\right| \leq \frac{2\mathfrak{c}|t-\tau|\sup _{\tau \in[0, \bar{T}]}\left\|\partial_{p_i} P_j(\tau)\right\|_{L^{\infty}}}{\inf _{\tau \in[0, \bar{T}]}\left|P^0(\tau)\right|} \lesssim \mathfrak{c}\frac{|t-\tau|}{p^0},
		\end{equation}
		provided $\bar{T}\ll 1$. Inserting \eqref{C.11} into \eqref{C.9}, one gets
		\begin{equation}\label{C.12}
			\left\|\partial_{p_i} P_j(\tau)\right\|_{L^{\infty}} \leq 2,
		\end{equation}
		provided $\bar{T}\ll 1$.
		
		From \eqref{C.4}, one has
		\begin{align*}
			\frac{d^2 \partial_{p_i} X_j(\tau)}{d \tau^2}= & \mathfrak{c}\partial_{p_i}\left(\frac{1}{P^0} \frac{d P_j}{d \tau}-\frac{P_j}{\left(P^0\right)^3} \frac{d P}{d \tau} \cdot P\right)(\tau) \\
			= & \mathfrak{c}\partial_{p_i}\left[\frac{-E_j^{\varepsilon,\mathfrak{c}}-\left[\frac{P}{P^0}\times B^{\varepsilon,\mathfrak{c}}\right]_j}{P^0}+\frac{P_j P \cdot E^{\varepsilon,\mathfrak{c}}}{\left(P^0\right)^3}\right](\tau) \\
			= & \mathfrak{c}\left\{\frac{-\nabla_x E_j^{\varepsilon,\mathfrak{c}} \cdot \partial_{p_i} X-\left[\frac{P}{P^0} \times\left(\nabla_x B^{\varepsilon,\mathfrak{c}} \cdot \partial_{p_i} X\right)\right]_j-\left[\partial_{p_i} \big(\frac{P}{P^0}\big) \times B^{\varepsilon,\mathfrak{c}}\right]_j}{P^0}\right. \\
			& +\frac{\partial_{p_i} P_j P \cdot E^{\varepsilon,\mathfrak{c}}+P_j \partial_{p_i} P \cdot E^{\varepsilon,\mathfrak{c}}+P_j P \cdot\left(\nabla_x E^{\varepsilon,\mathfrak{c}} \cdot \partial_{p_i} X\right)}{\left(P^0\right)^3} \\
			& \left.+\frac{E_j^{\varepsilon,\mathfrak{c}}+\left[\frac{P}{P^0} \times B^{\varepsilon,\mathfrak{c}}\right]_j}{\left(P^0\right)^3} P \cdot \partial_{p_i} P-\frac{3 P_j P \cdot E^{\varepsilon,\mathfrak{c}}\left(P \cdot \partial_{p_i} P\right)}{\left(P^0\right)^5}\right\}(\tau) ,
		\end{align*}
		which, together with \eqref{C.10}, \eqref{C.11} and \eqref{C.12}, yields that
		\begin{equation}\label{C.13}
			\frac{d^2 \partial_{p_i} X_j(\tau)}{d \tau^2} \lesssim \frac{\mathfrak{c}^2}{\left(p^0\right)^2}|t-\tau|+\frac{\mathfrak{c}}{\left(p^0\right)^2}. 
		\end{equation}
		
		It follows from the Taylor expansion and \eqref{C.7} that
		\begin{align*}
			\partial_{p_i} X_j(\tau)= & \partial_{p_i} X_j(t)+\left.\frac{d \partial_{p_i} X_j(\tau ; t, x, p)}{d \tau}\right|_{\tau=t}(\tau-t)+\frac{(\tau-t)^2}{2} \frac{d^2 \partial_{p_i} X_j(\bar{\tau})}{d \tau^2} \\
			= &\mathfrak{c} (\tau-t) \frac{\left(p^0\right)^2 \delta_{i j}-p_i p_j}{\left(p^0\right)^3}+\frac{(\tau-t)^2}{2}\frac{d^2 \partial_{p_i} X_j(\bar{\tau})}{d \tau^2},
		\end{align*}
		which implies that
		\begin{align}\label{C.13-1}
		\begin{aligned}
			\left|\operatorname{det}\left(\frac{\partial X(\tau)}{\partial p}\right)\right|= &\frac{\mathfrak{c}^3|\tau-t|^3}{\left(p^0\right)^3} \left|\operatorname{det}\left(\mathbf{I}-\frac{p \otimes p}{\left(p^0\right)^2}+\frac{(\tau-t)p^0}{2\mathfrak{c}} \frac{d^2 }{d \tau^2}\frac{\partial X(\bar{\tau})}{\partial p}\right)\right|,
		\end{aligned}
		\end{align}
		where $\mathbf{I}$ is the $3 \times 3$ identity matrix. A direct calculation shows
		\begin{align}\label{C.13-2}
		\left|\operatorname{det}\left(\mathbf{I}-\frac{p \otimes p}{\left(p^0\right)^2}\right)\right|=1-\frac{p^2}{(p^0)^2}=\frac{\mathfrak{c}^2}{(p^0)^2}.
		\end{align}
		It follows from \eqref{C.13} that
		\begin{align}\label{C.13-3}
		\left|\frac{(\tau-t)p^0}{2\mathfrak{c}} \frac{d^2 \partial_{p_i} X_j(\bar{\tau})}{d \tau^2}\right|\leq\frac{\mathfrak{c}|\tau-t|^2}{2p^0}+\frac{|\tau-t|}{2}\leq\frac{1}{4}.
		\end{align}
		Thus we have from \eqref{C.13-1}-\eqref{C.13-3} that
		$$
		\frac{\mathfrak{c}^5}{2\left(p^0\right)^5}|t-\tau|^3  \lesssim\left|\operatorname{det}\left(\frac{\partial X(\tau)}{\partial p}\right)\right| \lesssim \frac{2\mathfrak{c}^5}{\left(p^0\right)^5}|t-\tau|^3.
		$$
		Therefore the proof of Lemma \ref{lC.1} is completed.
	\end{proof} 
	
	\begin{lemma}\label{lC.2}
	    For $0 \leq s \leq t \leq \bar{T}$, there holds 
	    \begin{align}\label{C.14}
	        \int_0^t \exp \left\{-\frac{1}{\varepsilon}\int_s^t \nu_{\mathfrak{c}}(\tau, X(\tau), P(\tau)) d \tau \right\} \nu_{\mathfrak{c}}(s, X(s), P(s)) d s \lesssim \varepsilon.
	    \end{align}
	\end{lemma}
	\begin{proof}
    For $0 \leq s \leq t \leq \bar{T}$, it follows from $\eqref{C.4}_2$ that
	\begin{align}\label{C.0}
	    |P(s)-p|\leq |t-s|\sup _{s \in[0, T]}\left\|\left(E^{\varepsilon,\mathfrak{c}}(s), B^{\varepsilon,\mathfrak{c}}(s)\right)\right\|_{L^{\infty}}\leq \frac{1}{2} ,
	\end{align}
	which implies
	\begin{align}\label{C.00}
	    1+|P(s)|\cong 1+ |p|.
	\end{align}
	We divide the proof into three cases.
	
	{\it Case 1:} $|p|<\frac{1}{2}\mathfrak{c}$. Since $0 \leq s \leq \tau \leq t \leq \bar{T}$, it follows from \eqref{C.0} that
	\begin{align*}
	    |P(s)|\leq |P(s)-p|+|p|\leq \frac{1}{2}+\frac{1}{2}\mathfrak{c} <\mathfrak{c}.
	\end{align*}
	Then it follows from \eqref{C.00} and Lemma \ref{l2.1} that
	\begin{align*}
	    \nu_{\mathfrak{c}}(\tau, X(\tau), P(\tau))\cong 1+|P(\tau)| \cong 1+|p|, \quad 0\leq s \leq \tau \leq t \leq \bar{T},
	\end{align*}
	which yields immediately \eqref{C.14}.
	
	{\it Case 2:} $|p|>\frac{3}{2}\mathfrak{c}$. Similarly, it holds
    that
    \begin{align*}
	    |P(\tau)|\geq |p|-|P(\tau)-p|\geq \frac{3}{2}\mathfrak{c}-\frac{1}{2} >\mathfrak{c}.
	\end{align*}
	Then it follows from Lemma \ref{l2.1} that
	\begin{align*}
	    \nu_{\mathfrak{c}}(\tau, X(\tau), P(\tau))\cong \mathfrak{c}, \quad 0 \leq s \leq \tau \leq t \leq \bar{T},
	\end{align*}
	which yields \eqref{C.14} clearly.
	
	{\it Case 3:} $\frac{1}{2}\mathfrak{c}<|p|<\frac{3}{2}\mathfrak{c}$. It follows from $1+|P(\tau)| \cong 1+|p|\cong \mathfrak{c}$ and Lemma \ref{l2.1} that
	\begin{align*}
	    \nu_{\mathfrak{c}}(\tau, X(\tau), P(\tau))\cong 1+|p|\cong \mathfrak{c}, \quad 0 \leq s \leq \tau \leq t \leq \bar{T},
	\end{align*}
	which implies \eqref{C.14}.
	 
	We conclude \eqref{C.14} from above three cases. Therefore the proof of Lemma \ref{lC.2} is completed.
	\end{proof}
	
	\subsection{$L^{\infty}$-estimate of $h_R^{\varepsilon,\mathfrak{c}}$}
	
	\begin{lemma}\label{pC.2}
		Under the {\it a priori} assumptions \eqref{W.0}-\eqref{C.2}, for given $\bar{T}$ in Lemma \ref{lC.1}, there exists a $\varepsilon_0\ll1$ such that for $\varepsilon\in (0,\varepsilon_0],$
		\begin{align}\label{C.15}
			\sup_{s \in\left[0, \bar{T}\right]}\big\|\varepsilon^{\frac{3}{2}} h_R^{\varepsilon, \mathfrak{c}}(s)\big\|_{L^{\infty}} \leq& C\big\|\varepsilon^{\frac{3}{2}} h_R^{\varepsilon, \mathfrak{c}}(0)\big\|_{L^{\infty}}+C\sup_{s \in\left[0, \bar{T}\right]}\big\|\varepsilon^{\frac{5}{2}}\big(E_R^{\varepsilon,\mathfrak{c}}, B_R^{\varepsilon,\mathfrak{c}}\big)(s)\big\|_{L^{\infty}}\nonumber\\
			&+C \sup_{s \in\left[0, \bar{T}\right]}\big\|f_R^{\varepsilon, \mathfrak{c}}(s)\big\|+C \varepsilon^{k+\frac{5}{2}} ,
		\end{align}
		and
		\begin{align}\label{C.16}
			\big\|\varepsilon^{\frac{3}{2}} h_R^{\varepsilon, \mathfrak{c}}(\bar{T})\big\|_{L^{\infty}} \leq& \frac{1}{2}\big\|\varepsilon^{\frac{3}{2}} h_R^{\varepsilon, \mathfrak{c}}(0)\big\|_{L^{\infty}}+C\sup_{s \in\left[0, \bar{T}\right]}\big\|\varepsilon^{\frac{5}{2}}\big(E_R^{\varepsilon,\mathfrak{c}}, B_R^{\varepsilon,\mathfrak{c}}\big)(s)\big\|_{L^{\infty}}\nonumber\\
			&+C \sup_{s \in\left[0, \bar{T}\right]}\big\|f_R^{\varepsilon, \mathfrak{c}}(s)\big\|+C \varepsilon^{k+\frac{5}{2}} ,
		\end{align}
		where $C>0$ is independent of $\mathfrak{c}$.
	\end{lemma}
	\begin{proof}
		From $\eqref{e1.10}_1$ and \eqref{e1.14}, we know that
		{\small
		\begin{align}\label{C.18}
			\partial_t h_R^{\varepsilon,\mathfrak{c}}&+\hat{p} \cdot \nabla_x h_R^{\varepsilon,\mathfrak{c}}-\Big(E^{\varepsilon,\mathfrak{c}}+\frac{p}{p^0} \times B^{\varepsilon,\mathfrak{c}}\Big) \cdot \nabla_p h_R^{\varepsilon,\mathfrak{c}}+\frac{\nu_{\mathfrak{c}}}{\varepsilon} h_R^{\varepsilon,\mathfrak{c}}-\frac{1}{\varepsilon}\mathcal{K}_{\mathfrak{c},w}\left(h_R^{\varepsilon,\mathfrak{c}}\right) \nonumber\\
			=&\sum_{i=1}^{2 k-1} \varepsilon^{i-1} \frac{w_{\ell}}{\sqrt{J_{\mathbf{M}}}}\Big[Q_{\mathfrak{c}}\Big(F_i^{\mathfrak{c}}, \frac{\sqrt{J_{\mathbf{M}}}}{w_{\ell}} h_R^{\varepsilon,\mathfrak{c}}\Big)+Q_{\mathfrak{c}}\Big(\frac{\sqrt{J_{\mathbf{M}}}}{w_{\ell}} h_R^{\varepsilon,\mathfrak{c}}, F_i^{\mathfrak{c}}\Big)\Big] \nonumber\\
			& +\varepsilon^{k-1} \frac{w_{\ell}}{\sqrt{J_{\mathbf{M}}}} Q_{\mathfrak{c}}\left(\frac{\sqrt{J_{\mathbf{M}}}}{w_{\ell}} h_R^{\varepsilon,\mathfrak{c}}, \frac{\sqrt{J_{\mathbf{M}}}}{w_{\ell}} h_R^{\varepsilon,\mathfrak{c}}\right)+\Big(E^{\varepsilon,\mathfrak{c}}+\frac{p}{p^0} \times B^{\varepsilon,\mathfrak{c}}\Big) \cdot \nabla_p\left(\frac{\sqrt{J_{\mathbf{M}}}}{w_{\ell}}\right)\frac{w_{\ell}}{\sqrt{J_{\mathbf{M}}}}h_R^{\varepsilon,\mathfrak{c}} \nonumber\\
			& +\Big(E_R^{\varepsilon,\mathfrak{c}}+\frac{p}{p^0} \times B_R^{\varepsilon,\mathfrak{c}}\Big) \cdot \frac{w_{\ell}}{\sqrt{J_{\mathbf{M}}}}\Big[\nabla_p\Big(\mathbf{M}_{\mathfrak{c}}+\sum_{i=1}^{2k-1} \varepsilon^i F_i^{\mathfrak{c}}\Big)\Big]+\varepsilon^k \frac{w_{\ell}}{\sqrt{J_{\mathbf{M}}}} A , 
		\end{align}}
		where 
		$$
		\mathcal{K}_{\mathfrak{c},w}\left(h_R^{\varepsilon,\mathfrak{c}}\right)=w_{\ell}\mathcal{K}_{\mathfrak{c}}\left(\frac{h_R^{\varepsilon,\mathfrak{c}}}{w_{\ell}}\right) .
		$$ 
		For any $(t,x,p)$, integrating \eqref{C.18} along the backward trajectory \eqref{C.4}, one has that
		{\footnotesize
		\begin{align}\label{C.19}
			h_R^{\varepsilon,\mathfrak{c}}(t, x, p)=&\exp \left(-\frac{\tilde{\nu}_{\mathfrak{c}}(t, 0)}{\varepsilon}\right) h_R^{\varepsilon,\mathfrak{c}}(0, X(0), P(0)) \nonumber\\
			&+\frac{1}{\varepsilon} \int_0^t \exp \Big(-\frac{\tilde{\nu}_{\mathfrak{c}}(t, s)}{\varepsilon}\Big) \mathcal{K}_{\mathfrak{c},w}\left(h_R^{\varepsilon,\mathfrak{c}}\left(s, X(s), P(s)\right)\right) d s \nonumber\\
			&+\varepsilon^{k-1} \int_0^t \exp \Big(-\frac{\tilde{\nu}_{\mathfrak{c}}(t, s)}{\varepsilon}\Big) \left\{\frac{w_{\ell}}{\sqrt{J_{\mathbf{M}}}} Q_{\mathfrak{c}}\left(\frac{\sqrt{J_{\mathbf{M}}}}{w_{\ell}} h_R^{\varepsilon,\mathfrak{c}}, \frac{\sqrt{J_{\mathbf{M}}}}{w_{\ell}} h_R^{\varepsilon,\mathfrak{c}}\right)\right\}\left(s, X(s), P(s)\right) d s \nonumber\\
			&+\sum_{i=1}^{2 k-1} \varepsilon^{i-1} \int_0^t \exp \Big(-\frac{\tilde{\nu}_{\mathfrak{c}}(t, s)}{\varepsilon}\Big) \left\{\frac{w_{\ell}}{\sqrt{J_{\mathbf{M}}}}\left[Q_{\mathfrak{c}}\left(F_i^{\mathfrak{c}}, \frac{\sqrt{J_{\mathbf{M}}}}{w_{\ell}} h_R^{\varepsilon,\mathfrak{c}}\right)+Q_{\mathfrak{c}}\left(\frac{\sqrt{J_{\mathbf{M}}}}{w_{\ell}} h_R^{\varepsilon,\mathfrak{c}}, F_i^{\mathfrak{c}}\right)\right]\right\}\left(s, X(s), P(s)\right) d s \nonumber\\
			&+\int_0^t \exp \Big(-\frac{\tilde{\nu}_{\mathfrak{c}}(t, s)}{\varepsilon}\Big)\left\{\Big(E^{\varepsilon,\mathfrak{c}}+\frac{P(s)}{P^0(s)} \times B^{\varepsilon,\mathfrak{c}}\Big) \cdot \nabla_p\left(\frac{\sqrt{J_{\mathbf{M}}}}{w_{\ell}}\right)\frac{w_{\ell}}{\sqrt{J_{\mathbf{M}}}}h_R^{\varepsilon,\mathfrak{c}}\right\}\left(s, X(s), P(s)\right) d s \nonumber\\
			&+\int_0^t \exp \Big(-\frac{\tilde{\nu}_{\mathfrak{c}}(t, s)}{\varepsilon}\Big)\left\{\Big(E_R^{\varepsilon,\mathfrak{c}}+\frac{P(s)}{P^0(s)} \times B_R^{\varepsilon,\mathfrak{c}}\Big) \cdot \frac{w_{\ell}}{\sqrt{J_{\mathbf{M}}}}\big[\nabla_p \mathbf{M}_{\mathfrak{c}}+\sum_{i=1}^{2 k-1}\varepsilon^i \nabla_p F_i^{\mathfrak{c}}\big]\right\}\left(s, X(s), P(s)\right) d s \nonumber\\
			&+\varepsilon^k \int_0^t \exp \Big(-\frac{\tilde{\nu}_{\mathfrak{c}}(t, s)}{\varepsilon}\Big) \Big(\frac{w_{\ell}}{\sqrt{J_{\mathbf{M}}}} A\Big)\left(s, X(s), P(s)\right) d s 
			:=\sum_{i=1}^7 J_i,
		\end{align}}
		where we have denoted
		$$
		\tilde{\nu}_{\mathfrak{c}}(t, s):=\int_s^t \nu_{\mathfrak{c}}(\tau, X(\tau), P(\tau)) d \tau .
		$$
		
		In the following, we only prove \eqref{C.16}. Since \eqref{C.15} can be proved similarly.
		It follows from \eqref{C.0}-\eqref{C.00} and Lemma \ref{l2.1} that for $0 \leq s \leq \bar{T} \leq T$,
	    \begin{align*}
	         \nu_{\mathfrak{c}}(s):=\nu_{\mathfrak{c}}(s, X(s), P(s))\sim \nu_{\mathfrak{c}}(p)\geq \nu_0,
	    \end{align*}
        where $\nu_0$ is a positive constant which is independent of $\mathfrak{c}$.
		Then it is clear that
		\begin{align}\label{C.20-1}
		|\varepsilon^{\frac{3}{2}} J_1| \leq e^{-\frac{\nu_0 t}{\varepsilon}}\big\|\varepsilon^{\frac{3}{2}} h_R^{\varepsilon,\mathfrak{c}}(0)\big\|_{L^{\infty}}.
		\end{align}
		
		We have from \eqref{e2.4-1} that
		\begin{align}\label{C.001}
		    \Big|\frac{w_{\ell}}{\sqrt{J_{\mathbf{M}}}} Q_{\mathfrak{c}}\Big(\frac{\sqrt{J_{\mathbf{M}}}}{w_{\ell}} h_R^{\varepsilon,\mathfrak{c}}, \frac{\sqrt{J_{\mathbf{M}}}}{w_{\ell}} h_R^{\varepsilon,\mathfrak{c}}\Big)\left(s, X(s), P(s)\right)\Big|\lesssim\nu_{\mathfrak{c}}(s)\big\|h_R^{\varepsilon,\mathfrak{c}}(s)\big\|_{L^{\infty}}^2,
		\end{align}
        which, together with \eqref{C.14}, implies that
		\begin{align}\label{C.20-2}
			|\varepsilon^{\frac{3}{2}}J_3|\lesssim&\varepsilon^{k-\frac{5}{2}} \int_0^t \exp \Big(-\frac{\tilde{\nu}_{\mathfrak{c}}(t, s)}{\varepsilon}\Big) \nu_{\mathfrak{c}}(s)d s\sup_{s \in [0, t]} \big\|\varepsilon^{\frac{3}{2}} h_R^{\varepsilon,\mathfrak{c}}(s)\big\|_{L^{\infty}}^2 \nonumber\\
            \lesssim& \varepsilon^{k-\frac{3}{2}} e^{-\frac{\nu_0 t}{2\varepsilon}}\sup_{s \in [0, t]}\left\{e^{\frac{\nu_0 s}{2\varepsilon}}\big\|\varepsilon^{\frac{3}{2}} h_R^{\varepsilon,\mathfrak{c}}(s)\big\|_{L^{\infty}}^2\right\}.
		\end{align}
		Similarly, we have
		\begin{align}\label{C.20-3}
			|\varepsilon^{\frac{3}{2}} J_4|
			\lesssim& \sum_{i=1}^{2 k-1} \varepsilon^{i-1} \int_0^t \exp \Big(-\frac{\tilde{\nu}_{\mathfrak{c}}(t, s)}{\varepsilon}\Big) \nu_{\mathfrak{c}}(s) d s \cdot\sup_{s \in [0, t]}\Big\|\frac{w_{\ell} }{\sqrt{J_{\mathbf{M}}}}F_i^{\mathfrak{c}}\Big\|_{L^{\infty}} \cdot\sup_{s \in [0, t]}\big\|\varepsilon^{\frac{3}{2}} h_R^{\varepsilon,\mathfrak{c}}(s)\big\|_{L^{\infty}}  \nonumber\\
			\lesssim & \varepsilon e^{-\frac{\nu_0 t}{2\varepsilon}}\sup_{s \in\left[0, t\right]}\left\{e^{\frac{\nu_0 s}{2\varepsilon}}\big\|\varepsilon^{\frac{3}{2}} h_R^{\varepsilon,\mathfrak{c}}(s)\big\|_{L^{\infty}}\right\}. 
		\end{align}
		
		A direct calculation shows that  $\Big|\nabla_p\left(\frac{\sqrt{J_{\mathbf{M}}}}{w_{\ell}}\right)\frac{w_{\ell}}{\sqrt{J_{\mathbf{M}}}}\Big|\lesssim\nu_{\mathfrak{c}}(p)$, then we obtain
		\begin{align}\label{C.20-4}
			|\varepsilon^{\frac{3}{2}} J_5| &\lesssim \int_0^t \exp \Big(-\frac{\tilde{\nu}_{\mathfrak{c}}(t, s)}{\varepsilon}\Big)\nu_{\mathfrak{c}}(P(s))\Big|\Big(E^{\varepsilon,\mathfrak{c}}+\frac{P(s)}{P^0(s)} \times B^{\varepsilon,\mathfrak{c}}\Big) \varepsilon^{\frac{3}{2}} h_R^{\varepsilon,\mathfrak{c}}\left(s, X(s), P(s)\right)\Big| d s \nonumber\\
			& \lesssim \int_0^t \exp \Big(-\frac{\tilde{\nu}_{\mathfrak{c}}(t,s)}{\varepsilon}\Big)\nu_{\mathfrak{c}}(s) d s  \cdot\sup_{s \in [0, t]}\big\|\varepsilon^{\frac{3}{2}} h_R^{\varepsilon,\mathfrak{c}}(s)\big\|_{L^{\infty}}  \cdot\sup_{s \in [0, t]}\|\big(E^{\varepsilon,\mathfrak{c}}, B^{\varepsilon,\mathfrak{c}}\big)(s) \|_{L^{\infty}} \nonumber\\
			& \lesssim \varepsilon e^{-\frac{\nu_0 t}{2\varepsilon}} \sup_{s \in\left[0, t\right]}\left\{e^{\frac{\nu_0 s}{2\varepsilon}}\big\|\varepsilon^{\frac{3}{2}} h_R^{\varepsilon,\mathfrak{c}}(s)\big\|_{L^{\infty}}\right\}, 
		\end{align}
		where we have used the uniform estimate \eqref{C.3}.
		
		For $J_6$ and $J_7$, one has
		\begin{align}\label{C.20-5}
		    |\varepsilon^{\frac{3}{2}} J_6|+|\varepsilon^{\frac{3}{2}} J_7|\lesssim \sup_{s \in[0,t]}\varepsilon^{\frac{5}{2}}\big\|\big(E_R^{\varepsilon,\mathfrak{c}}, B_R^{\varepsilon,\mathfrak{c}}\big)(s)\big\|_{L^{\infty}}+\varepsilon^{k+\frac{5}{2}}.
		\end{align}
        
		Substituting \eqref{C.20-1} and \eqref{C.20-2}-\eqref{C.20-5} into \eqref{C.19}, one obtains
		\begin{align}\label{C.20}
			\big|\varepsilon^{\frac{3}{2}} h_R^{\varepsilon,\mathfrak{c}}(t,x,p)\big|\leq & e^{-\frac{\nu_0 t}{\varepsilon}}\big\|\varepsilon^{\frac{3}{2}} h_R^{\varepsilon,\mathfrak{c}}(0)\big\|_{L^{\infty}}+ C\varepsilon e^{-\frac{\nu_0 t}{2\varepsilon}}\sup_{s \in\left[0, t\right]}\left\{e^{\frac{\nu_0 s}{2\varepsilon}}\big\|\varepsilon^{\frac{3}{2}} h_R^{\varepsilon,\mathfrak{c}}(s)\big\|_{L^{\infty}}\right\}\nonumber\\
			&+C\varepsilon^{k-\frac{3}{2}} e^{-\frac{\nu_0 t}{2\varepsilon}}\sup_{s \in [0, t]}\left\{e^{\frac{\nu_0 s}{2\varepsilon}}\big\|\varepsilon^{\frac{3}{2}} h_R^{\varepsilon,\mathfrak{c}}(s)\big\|_{L^{\infty}}^2\right\}\nonumber\\
			& + C\sup_{s \in[0,t]}\varepsilon^{\frac{5}{2}}\big\|\big(E_R^{\varepsilon,\mathfrak{c}}, B_R^{\varepsilon,\mathfrak{c}}\big)(s)\big\|_{L^{\infty}}+C \varepsilon^{k+\frac{5}{2}}+ |\varepsilon^{\frac{3}{2}}J_2|,
		\end{align}
		where $C>0$ is independent of $\mathfrak{c}$.
        
        Now it remains to control the last term $J_2$. To this end, we iterate the mild formulation once more along the backward characteristics.
		Denote $$
		X^\prime(s^{\prime}):=X(s^\prime;s,X(s;t,x,p),q), \quad P^\prime(s^{\prime}):=P(s^\prime;s,X(s;t,x,p),q)
		$$ 
		and
		$$
		\tilde{\nu}_{\mathfrak{c}}^\prime\left(s, s^{\prime}\right):=\int_{s^{\prime}}^s \nu_{\mathfrak{c}}\left(\tau, X^\prime(\tau), P^\prime(\tau)\right) d \tau, \quad
		\nu_{\mathfrak{c}}^\prime(s^\prime):=\nu_{\mathfrak{c}}\left(s^\prime, X^\prime(s^\prime), P^\prime(s^\prime)\right).
		$$
		Similar to \eqref{C.0}-\eqref{C.00}, one knows that
        \begin{align*}
            1+|P^\prime(s^{\prime})|\cong 1+ |q|, \quad 0 \leq s' \leq s \leq t \leq \bar{T}.
        \end{align*}
		
			We substitute \eqref{C.19} into $J_2$ to obtain
		{\footnotesize
		\begin{align}
			|J_2| 
			\lesssim&  \frac{1}{\varepsilon} \int_0^t \exp \Big(-\frac{\tilde{\nu}_{\mathfrak{c}}(t, s)}{\varepsilon}\Big) \int_{\mathbb{R}^3} k_w(P(s), q)\Big|h_R^{\varepsilon,\mathfrak{c}}\left(s, X(s), q\right)\Big| d q d s \nonumber\\
			\lesssim & \frac{1}{\varepsilon} \int_0^t \exp \Big(-\frac{\tilde{\nu}_{\mathfrak{c}}(t, s)}{\varepsilon}\Big) d s \int_{\mathbb{R}^3} k_w(P(s), q) \exp \left(-\frac{\tilde{\nu}_{\mathfrak{c}}^{\prime}(s, 0)}{\varepsilon}\right)\Big|h_R^{\varepsilon,\mathfrak{c}}\left(0, X^\prime(0), P^\prime(0)\right)\Big| d q \nonumber\\
			& +\frac{1}{\varepsilon^2} \int_0^t \exp \Big(-\frac{\tilde{\nu}_{\mathfrak{c}}(t, s)}{\varepsilon}\Big) d s \int_0^s \exp \Big(-\frac{\tilde{\nu}_{\mathfrak{c}}^{\prime}\left(s, s^{\prime}\right)}{\varepsilon}\Big) d s^{\prime} \int_{\mathbb{R}^3} k_w(P(s), q)\Big|\mathcal{K}_{\mathfrak{c},w} h_R^{\varepsilon,\mathfrak{c}}\left(s^\prime, X^\prime(s^\prime), P^\prime(s^\prime)\right)\Big| d q\nonumber\\
			&+\varepsilon^{k-2}\int_0^t \exp \Big(-\frac{\tilde{\nu}_{\mathfrak{c}}(t, s)}{\varepsilon}\Big) d s \int_0^s \exp \Big(-\frac{\tilde{\nu}_{\mathfrak{c}}^{\prime}\left(s, s^{\prime}\right)}{\varepsilon}\Big) d s^{\prime} \nonumber\\
			&\quad\times\int_{\mathbb{R}^3} k_w(P(s), q)\frac{w_{\ell}}{\sqrt{J_{\mathbf{M}}}}\Big|Q_{\mathfrak{c}}\left(\frac{\sqrt{J_{\mathbf{M}}}}{w_{\ell}} h_R^{\varepsilon,\mathfrak{c}}, \frac{\sqrt{J_{\mathbf{M}}}}{w_{\ell}} h_R^{\varepsilon,\mathfrak{c}}\right)\left(s^\prime, X^\prime(s^\prime), P^\prime(s^\prime)\right)\Big| d q\nonumber\\
			& +\frac{1}{\varepsilon} \sum_{i=1}^{2 k-1} \varepsilon^{i-1} \int_0^t \exp \Big(-\frac{\tilde{\nu}_{\mathfrak{c}}(t, s)}{\varepsilon}\Big) d s \int_0^s \exp \Big(-\frac{\tilde{\nu}_{\mathfrak{c}}^{\prime}\left(s, s^{\prime}\right)}{\varepsilon}\Big) d s^{\prime} \nonumber\\ 
			&\quad \times\int_{\mathbb{R}^3} k_w(P(s), q) \frac{w_{\ell}}{\sqrt{J_{\mathbf{M}}}}\Big|\left\{Q_{\mathfrak{c}}\left(F_i^{\mathfrak{c}}, \frac{\sqrt{J_{\mathbf{M}}}}{w_{\ell}} h_R^{\varepsilon,\mathfrak{c}}\right)+\right.  \left.Q_{\mathfrak{c}}\left(\frac{\sqrt{J_{\mathbf{M}}}}{w_{\ell}} h_R^{\varepsilon,\mathfrak{c}}, F_i^{\mathfrak{c}}\right)\right\}\left(s^\prime, X^\prime(s^\prime), P^\prime(s^\prime)\right) \Big| d q \nonumber\\
			& +\frac{1}{\varepsilon} \int_0^t \exp \Big(-\frac{\tilde{\nu}_{\mathfrak{c}}(t, s)}{\varepsilon}\Big) d s \int_0^s \exp \Big(-\frac{\tilde{\nu}_{\mathfrak{c}}^{\prime}\left(s, s^{\prime}\right)}{\varepsilon}\Big) d s^{\prime}  \nonumber\\
			&\quad \times\int_{\mathbb{R}^3} k_w(P(s), q)\left(E^{\varepsilon,\mathfrak{c}}+\frac{P'(s')}{P'^0(s')} \times B^{\varepsilon,\mathfrak{c}}\right) \cdot \nabla_p\left(\frac{\sqrt{J_{\mathbf{M}}}}{w_{\ell}}\right)\frac{w_{\ell}}{\sqrt{J_{\mathbf{M}}}}h_R^{\varepsilon,\mathfrak{c}}\left(s^\prime, X^\prime(s^\prime), P^\prime(s^\prime)\right) d q \nonumber\\
			& +\frac{1}{\varepsilon} \int_0^t \exp \Big(-\frac{\tilde{\nu}_{\mathfrak{c}}(t, s)}{\varepsilon}\Big) d s \int_0^s \exp \Big(-\frac{\tilde{\nu}_{\mathfrak{c}}^{\prime}\left(s, s^{\prime}\right)}{\varepsilon}\Big) d s^{\prime}  \nonumber\\
			&\quad \times\int_{\mathbb{R}^3} k_w(P(s), q) \Big|\Big(E_R^{\varepsilon,\mathfrak{c}}+\frac{P'(s')}{P'^0(s')} \times B_R^{\varepsilon,\mathfrak{c}}\Big) \Big|\cdot \frac{w_{\ell}}{\sqrt{J_{\mathbf{M}}}}\big[\nabla_p \mathbf{M}_{\mathfrak{c}}+\sum_{i=1}^{2 k-1} \varepsilon^i \nabla_p F_i^{\mathfrak{c}}\big]\left(s^\prime, X^\prime(s^\prime), P^\prime(s^\prime)\right) d q \nonumber\\
			& +\varepsilon^{k-1} \int_0^t \exp \Big(-\frac{\tilde{\nu}_{\mathfrak{c}}(t, s)}{\varepsilon}\Big) d s \int_0^s \exp \Big(-\frac{\tilde{\nu}_{\mathfrak{c}}^{\prime}\left(s, s^{\prime}\right)}{\varepsilon}\Big) d s^{\prime} \int_{\mathbb{R}^3} k_w(P(s), q) \frac{w_{\ell}}{\sqrt{J_{\mathbf{M}}}}\Big|A\left(s^\prime, X^\prime(s^\prime), P^\prime(s^\prime)\right)\Big| d q \nonumber\\
			:=&\sum_{i=1}^7 J_{2 i} .\nonumber
		\end{align}}

         By similar arguments as in \eqref{C.20-1}-\eqref{C.20-5}, one obtains
		\begin{align}\label{C.25}
			\big|\varepsilon^{\frac{3}{2}} J_2\big|\lesssim & \frac{t}{\varepsilon} e^{-\frac{\nu_0 t}{\varepsilon}}\big\|\varepsilon^{\frac{3}{2}} h_R^{\varepsilon,\mathfrak{c}}(0)\big\|_{L^{\infty}}+ \varepsilon e^{-\frac{\nu_0 t}{2\varepsilon}}\sup_{s \in\left[0, t\right]}\left\{e^{\frac{\nu_0 s}{2\varepsilon}}\big\|\varepsilon^{\frac{3}{2}} h_R^{\varepsilon,\mathfrak{c}}(s)\big\|_{L^{\infty}}\right\}\nonumber\\
			&+\varepsilon^{k-\frac{3}{2}} e^{-\frac{\nu_0 t}{2\varepsilon}}\sup_{s \in [0, t]}\left\{e^{\frac{\nu_0 s}{2\varepsilon}}\big\|\varepsilon^{\frac{3}{2}} h_R^{\varepsilon,\mathfrak{c}}(s)\big\|_{L^{\infty}}^2\right\}\nonumber\\
			& + \varepsilon^{\frac{5}{2}}\sup_{s \in\left[0, t\right]}\big\|\big(E_R^{\varepsilon,\mathfrak{c}}, B_R^{\varepsilon,\mathfrak{c}}\big)(s)\big\|_{L^{\infty}}+ \varepsilon^{k+\frac{5}{2}}+ |\varepsilon^{\frac{3}{2}}J_{22}|.
		\end{align}
		Now we focus on the most difficult term $J_{22}$. It holds that
		\begin{align}\label{C.26}
			|\varepsilon^{\frac{3}{2}} J_{22}| \lesssim& \frac{1}{\varepsilon^2} \int_0^t \exp \left(-\frac{\tilde{\nu}_{\mathfrak{c}}\left(t, s\right)}{\varepsilon}\right) d s \int_0^s \exp \Big(-\frac{\tilde{\nu}_{\mathfrak{c}}^{\prime}\left(s, s^{\prime}\right)}{\varepsilon}\Big) d s^{\prime} \nonumber\\
			&\int_{\mathbb{R}^3} k_w(P(s), q) d q \int_{\mathbb{R}^3} k_w(P^\prime(s^\prime), q^\prime)\Big|\varepsilon^{\frac{3}{2}} h_R^{\varepsilon,\mathfrak{c}}\left(s^\prime, X^\prime(s^\prime), q^\prime\right)\Big| d q^\prime.
		\end{align}
		We divide the estimate of \eqref{C.26} into four cases.
		
		\noindent{\it Case 1. $|p| \geq N_0$.} By \eqref{C.0}, there holds $|P(s)|\geq \frac{N_0}{2}$. Then \eqref{e2.39} together with $\mathfrak{c}\gg N_0$ implies that 
		$$
		\int_{\mathbb{R}^3}\int_{\mathbb{R}^3} k_w(P(s), q)k_w(P^\prime(s^\prime), q^\prime) d q d q^\prime\lesssim \max \left\{\frac{1}{\mathfrak{c}}, \frac{1}{1+|P(s)|}\right\} \lesssim \frac{1}{N_0}.
		$$
		Thus we have
		$$
		|\varepsilon^{\frac{3}{2}} J_{22}| \lesssim \frac{1}{N_0} e^{-\frac{\nu_0 t}{2\varepsilon}} \sup_{s \in\left[0, t\right]}\left\{e^{\frac{\nu_0 s}{2\varepsilon}}\big\|\varepsilon^{\frac{3}{2}} h_R^{\varepsilon,\mathfrak{c}}(s)\big\|_{L^{\infty}}\right\}  .
		$$
		\noindent{\it Case 2. $|p| \leq N_0,|q| \geq 2 N_0$ or $|q| \leq 2 N_0,|q^\prime| \geq 3 N_0$.}  Observe that 
		$$
		\begin{gathered}
			|P(s)-q|\geq |q-p|-|P(s)-p|\geq |q|-|p|-|P(s)-p|,\\
			|P^\prime(s^\prime)-q^\prime|\geq |q^\prime-q|-|P^\prime(s^\prime)-q|\geq |q^\prime|-|q|-|P^\prime(s^\prime)-q|,
		\end{gathered}
		$$
		thus we have either $|P(s)-q|\geq \frac{N_0}{2}$ or $|P^\prime(s^\prime)-q^\prime|\geq \frac{N_0}{2}$. We use \eqref{e2.39} to get the following bound
		\begin{align*}
		    &\iint_{|p| \leq N_0,|q| \geq 2 N_0} k_w(P(s), q)k_w(P^\prime(s^\prime), q^\prime) d q d q^\prime\lesssim e^{-\frac{\delta}{8T_M} N_0}\max \left\{\frac{1}{\mathfrak{c}}, \frac{1}{1+|P(s)|}\right\} \lesssim e^{-\frac{\delta}{8T_M} N_0},\\
		    &\iint_{|q| \leq 2 N_0,|q^{\prime}| \geq 3 N_0}k_w(P(s), q)k_w(P^\prime(s^\prime), q^\prime) d q d q^\prime\lesssim e^{-\frac{\delta}{8T_M} N_0}\max \left\{\frac{1}{\mathfrak{c}}, \frac{1}{1+|P'(s')|}\right\} \lesssim e^{-\frac{\delta}{8T_M} N_0}.
		\end{align*}
		Consequently,
		\begin{align}
			& \frac{1}{\varepsilon^2} \int_0^t \exp \Big(-\frac{\tilde{\nu_{\mathfrak{c}}}(t, s)}{\varepsilon}\Big) d s \int_0^s \exp \Big(-\frac{\tilde{\nu}_{\mathfrak{c}}^{\prime}\left(s, s^{\prime}\right)}{\varepsilon}\Big) d s^{\prime} \nonumber\\
			& \quad \times\Big\{\iint_{|p| \leq N_0,|q| \geq 2 N_0}+\iint_{|q| \leq 2 N_0,|q^{\prime}| \geq 3 N_0}\Big\} \nonumber\\
			\lesssim& e^{-\frac{\delta}{8T_M} N_0} e^{-\frac{\nu_0 t}{2\varepsilon}} \sup_{s \in\left[0, t\right]}\left\{e^{\frac{\nu_0 s}{2\varepsilon}}\big\|\varepsilon^{\frac{3}{2}} h_R^{\varepsilon,\mathfrak{c}}(s)\big\|_{L^{\infty}}\right\}  \lesssim \frac{1}{N_0} e^{-\frac{\nu_0 t}{2\varepsilon}} \sup_{s \in\left[0, t\right]}\left\{e^{\frac{\nu_0 s}{2\varepsilon}}\big\|\varepsilon^{\frac{3}{2}} h_R^{\varepsilon,\mathfrak{c}}(s)\big\|_{L^{\infty}}\right\}. \nonumber
		\end{align}
		\noindent{\it Case 3.} For $s-s^{\prime} \leq \kappa \varepsilon$ and $|p| \leq N_0,|q| \leq 2 N_0,|q^\prime| \leq 3 N_0$, one has
		\begin{align}
			&\frac{1}{\varepsilon^2} \int_0^t \exp \Big(-\frac{\tilde{\nu_{\mathfrak{c}}}(t, s)}{\varepsilon}\Big) d s \int_{s-\kappa \varepsilon}^s \exp \Big(-\frac{\tilde{\nu}_{\mathfrak{c}}^{\prime}\left(s, s^{\prime}\right)}{\varepsilon}\Big) d s^{\prime}\nonumber\\
			& \quad\times \int_{|q| \leq 2 N_0} k_w(P(s), q) d q \int_{|q^{\prime}| \leq 3 N_0} k_w(P^\prime(s^\prime), q^\prime)\Big|\varepsilon^{\frac{3}{2}} h_R^{\varepsilon, \mathfrak{c}}\left(s^\prime, X^\prime(s^\prime), q^\prime\right)\Big| d q^\prime \nonumber\\
			\leq & C_{N_0} e^{-\frac{\nu_0 t}{2\varepsilon}} \sup_{s \in\left[0, t\right]}\left\{e^{\frac{\nu_0 s}{2\varepsilon}}\big\|\varepsilon^{\frac{3}{2}} h_R^{\varepsilon,\mathfrak{c}}(s)\big\|_{L^{\infty}}\right\} \Big(\frac{1}{\varepsilon} \int_0^t \exp \Big(-\frac{\tilde{\nu_{\mathfrak{c}}}(t, s)}{\varepsilon}\Big) d s\Big) \Big(\int_{s-\kappa \varepsilon}^s\frac{1}{\varepsilon} d s^\prime\Big)\nonumber\\
			\lesssim & \kappa e^{-\frac{\nu_0 t}{2\varepsilon}} \sup_{s \in\left[0, t\right]}\left\{e^{\frac{\nu_0 s}{2\varepsilon}}\big\|\varepsilon^{\frac{3}{2}} h_R^{\varepsilon,\mathfrak{c}}(s)\big\|_{L^{\infty}}\right\}. \nonumber
		\end{align}
		\noindent{\it Case 4.} For $s-s^{\prime} \geq \kappa \varepsilon$ and $|p| \leq N_0,|q| \leq 2 N_0,|q^\prime| \leq 3 N_0$, this is the last remaining case. By \eqref{e2.39}, one has
		\begin{align}
			& \int_{|q| \leq 2 N_0} \int_{|q^\prime| \leq 3 N_0} k_w(P(s), q) k_w\left(P^\prime(s^\prime), q^\prime\right)\Big| h_R^{\varepsilon, \mathfrak{c}}\left(s^\prime, X^\prime(s^\prime), q^\prime\right)\Big| d q^\prime d q\nonumber\\
			\leq& C_{N_0} \int_{|q| \leq 2 N_0} \int_{|q^\prime| \leq 3 N_0} k_w(P(s), q) k_w\left(P^\prime(s^\prime), q^\prime\right)\Big|f_R^{\varepsilon, \mathfrak{c}}\left(s^\prime, X^\prime(s^\prime), q^\prime\right)\Big| d q^\prime d q \nonumber\\
			\leq& C_{N_0}\Big(\int_{|q| \leq 2 N_0} \int_{|q^\prime| \leq 3 N_0} k_w^2(P(s), q) k_w^2\left(P^\prime(s^\prime), q^\prime\right) d q^\prime d q\Big)^{\frac{1}{2}} \nonumber\\
			& \quad \times\Big(\int_{|q| \leq 2 N_0} \int_{|q^\prime| \leq 3 N_0}\Big|f_R^{\varepsilon, \mathfrak{c}}\left(s^\prime, X^\prime(s^\prime), q^\prime\right)\Big|^2 d q^\prime d q\Big)^{\frac{1}{2}} \nonumber\\
			\leq& C_{N_0}\Big(\int_{\mathbb{R}^3} \int_{\mathbb{R}^3}\Big|f_R^{\varepsilon, \mathfrak{c}}\left(s^\prime, X^\prime(s^\prime), q^\prime\right)\Big|^2 \cdot \varepsilon^{-3} \kappa^{-3} d X^\prime(s^\prime) d q^\prime\Big)^{\frac{1}{2}} \nonumber\\
			\leq& \frac{C_{N_0, \kappa}}{\varepsilon^{\frac{3}{2}}} \sup_{s \in\left[0, \bar{T}\right]}\big\|f_R^{\varepsilon, \mathfrak{c}}(s)\big\|,\nonumber
		\end{align}
		where we have made a change of variables $q \mapsto X^\prime(s^\prime)$ with
		$$
		\Big|\operatorname{det}\Big(\frac{d X^\prime(s^\prime)}{d q}\Big)\Big|\geq\frac{\mathfrak{c}^5}{2\left(q^0\right)^5}\left(s-s^{\prime}\right)^3 \geq \frac{\kappa^3 \varepsilon^3}{2^6}
		$$
		for $1 \leq N_0 \ll \mathfrak{c}$. Thus it holds
		\begin{align}
			& \frac{1}{\varepsilon^2} \int_0^t \exp \Big(-\frac{\tilde{\nu}_{\mathfrak{c}}(t, s)}{\varepsilon}\Big) d s \int_0^{s-\kappa \varepsilon} \exp \Big(-\frac{\tilde{\nu}_{\mathfrak{c}}^\prime\left(s, s^{\prime}\right)}{\varepsilon}\Big) d s^{\prime}\nonumber\\
			& \quad \times \int_{|q| \leq 2 N_0} k_w(P(s), q) d q \int_{|q^\prime| \leq 3 N_0} k_w\left(P^\prime(s^\prime), q^\prime\right)\Big|\varepsilon^{\frac{3}{2}} h_R^{\varepsilon, \mathfrak{c}}\left(s^\prime, X^\prime(s^\prime), q^\prime\right)\Big| d q^\prime \nonumber\\
			& \leq C_{N_0, \kappa} \sup_{s \in\left[0, \bar{T}\right]}\big\|f_R^{\varepsilon, \mathfrak{c}}(s)\big\| .\nonumber
		\end{align}
        
		Combining the above four cases, we obtain
		\begin{equation}\label{C.21}
			|\varepsilon^{\frac{3}{2}} J_{22}| \leq C\Big(\kappa+\frac{1}{N_0}\Big) e^{-\frac{\nu_0 t}{2\varepsilon}} \sup_{s \in\left[0, t\right]}\left\{e^{\frac{\nu_0 s}{2\varepsilon}}\big\|\varepsilon^{\frac{3}{2}} h_R^{\varepsilon,\mathfrak{c}}(s)\big\|_{L^{\infty}}\right\} +C_{N_0, \kappa} \sup_{s \in\left[0, \bar{T}\right]}\big\|f_R^{\varepsilon, \mathfrak{c}}(s)\big\| .
		\end{equation}
		Therefore, combining \eqref{C.25} and \eqref{C.21}, one obtains
		\begin{align*}
			|\varepsilon^{\frac{3}{2}} J_{2}| \lesssim & \frac{t}{\varepsilon}e^{-\frac{\nu_0 t}{\varepsilon}}\big\|\varepsilon^{\frac{3}{2}} h_R^{\varepsilon,\mathfrak{c}}(0)\big\|_{L^{\infty}}+ \Big(\varepsilon+\kappa+\frac{1}{N_0}\Big) e^{-\frac{\nu_0 t}{2\varepsilon}}\sup_{s \in\left[0, t\right]}\left\{e^{\frac{\nu_0 s}{2\varepsilon}}\big\|\varepsilon^{\frac{3}{2}} h_R^{\varepsilon,\mathfrak{c}}(s)\big\|_{L^{\infty}}\right\} \nonumber\\
			&+\varepsilon^{k-\frac{3}{2}} e^{-\frac{\nu_0 t}{2\varepsilon}}\sup_{s \in [0, t]}\left\{e^{\frac{\nu_0 s}{2\varepsilon}}\big\|\varepsilon^{\frac{3}{2}} h_R^{\varepsilon,\mathfrak{c}}(s)\big\|_{L^{\infty}}^2\right\}+ \sup_{s \in\left[0, \bar{T}\right]}\big\|f_R^{\varepsilon, \mathfrak{c}}(s)\big\|\\
			&+ \sup_{s \in\left[0, \bar{T}\right]}\varepsilon^{\frac{5}{2}}\big\|\big(E_R^{\varepsilon,\mathfrak{c}}, B_R^{\varepsilon,\mathfrak{c}}\big)(s)\big\|_{L^{\infty}} + \varepsilon^{k+\frac{5}{2}},
		\end{align*}
		which, together with \eqref{C.20}, \eqref{W.0} and $k\geq 5$, yields that
		\begin{align}\label{C.22}
			\big\|\varepsilon^{\frac{3}{2}} h_R^{\varepsilon,\mathfrak{c}}(t)\big\|_{L^{\infty}} \lesssim & \Big(1+\frac{t}{\varepsilon}\Big)e^{-\frac{\nu_0 t}{\varepsilon}}\big\|\varepsilon^{\frac{3}{2}} h_R^{\varepsilon,\mathfrak{c}}(0)\big\|_{L^{\infty}}+ \Big(\varepsilon+\kappa+\frac{1}{N_0}\Big) e^{-\frac{\nu_0 t}{2\varepsilon}}\sup_{s \in\left[0, t\right]}\left\{e^{\frac{\nu_0 s}{2\varepsilon}}\big\|\varepsilon^{\frac{3}{2}} h_R^{\varepsilon,\mathfrak{c}}(s)\big\|_{L^{\infty}}\right\} \nonumber\\
			&+ \sup_{s \in\left[0, \bar{T}\right]}\varepsilon^{\frac{5}{2}}\big\|\big(E_R^{\varepsilon,\mathfrak{c}}, B_R^{\varepsilon,\mathfrak{c}}\big)(s)\big\|_{L^{\infty}}+  \sup_{s \in\left[0, \bar{T}\right]}\big\|f_R^{\varepsilon, \mathfrak{c}}(s)\big\| + \varepsilon^{k+\frac{5}{2}}.
		\end{align}
		Multiplying \eqref{C.22} by $e^{\frac{\nu_0 t}{2\varepsilon}}$ and take the supremum of $t\in [0, \bar{T}]$, we obtain
		\begin{align}\label{C.23}
			&\sup_{s \in\left[0, \bar{T}\right]}\left\{e^{\frac{\nu_0 s}{2\varepsilon}}\big\|\varepsilon^{\frac{3}{2}} h_R^{\varepsilon, \mathfrak{c}}(s)\big\|_{L^{\infty}}\right\} \nonumber\\
			\lesssim& \sup_{s \in\left[0, \bar{T}\right]}\left\{\left(1+\frac{s}{\varepsilon}\right)e^{-\frac{\nu_0 s}{2\varepsilon}}\big\|\varepsilon^{\frac{3}{2}} h_R^{\varepsilon,\mathfrak{c}}(0)\big\|_{L^{\infty}}\right\}  +e^{\frac{\nu_0 \bar{T}}{2\varepsilon}}\sup_{s \in\left[0, \bar{T}\right]}\big\|\varepsilon^{\frac{5}{2}}\big(E_R^{\varepsilon,\mathfrak{c}}, B_R^{\varepsilon,\mathfrak{c}}\big)(s)\big\|_{L^{\infty}}\nonumber\\
			&+\Big(\varepsilon+\kappa+\frac{1}{N_0}\Big) \sup_{s \in\left[0, \bar{T}\right]}\left\{e^{\frac{\nu_0 s}{2\varepsilon}}\big\|\varepsilon^{\frac{3}{2}} h_R^{\varepsilon, \mathfrak{c}}(s)\big\|_{L^{\infty}}\right\}+ e^{\frac{\nu_0 \bar{T}}{2\varepsilon}}\sup_{s \in\left[0, \bar{T}\right]}\big\|f_R^{\varepsilon, \mathfrak{c}}(s)\big\|+e^{\frac{\nu_0 \bar{T}}{2\varepsilon}}\varepsilon^{k+\frac{5}{2}} .
		\end{align}
		Note that $\left(1+\frac{s}{\varepsilon}\right)e^{-\frac{\nu_0 s}{2\varepsilon}}$ is uniformly bounded in $s$ and $\varepsilon$. Choosing $N_0$ suitably large and $\kappa, \varepsilon$ suitably small that $\varepsilon+\kappa+\frac{1}{N_0}<\frac{1}{2}$, one gets from \eqref{C.23} that
		\begin{align}
			\sup_{s \in\left[0, \bar{T}\right]}\left\{e^{\frac{\nu_0 s}{2\varepsilon}}\big\|\varepsilon^{\frac{3}{2}} h_R^{\varepsilon, \mathfrak{c}}(s)\big\|_{L^{\infty}}\right\} \leq& C\big\|\varepsilon^{\frac{3}{2}} h_R^{\varepsilon,\mathfrak{c}}(0)\big\|_{L^{\infty}}+e^{\frac{\nu_0 \bar{T}}{2\varepsilon}}\sup_{s \in\left[0, \bar{T}\right]}\varepsilon^{\frac{5}{2}}\big\|\big(E_R^{\varepsilon,\mathfrak{c}}, B_R^{\varepsilon,\mathfrak{c}}\big)(s)\big\|_{L^{\infty}}\nonumber\\
			&+ e^{\frac{\nu_0 \bar{T}}{2\varepsilon}}\sup_{s \in\left[0, \bar{T}\right]}\big\|f_R^{\varepsilon, \mathfrak{c}}(s)\big\|+e^{\frac{\nu_0 \bar{T}}{2\varepsilon}}\varepsilon^{k+\frac{5}{2}},\nonumber
		\end{align}
		which yields immediately that
		\begin{align}\label{C.24}
			\big\|\varepsilon^{\frac{3}{2}} h_R^{\varepsilon, \mathfrak{c}}(\bar{T})\big\|_{L^{\infty}} \leq& \frac{1}{2}\big\|\varepsilon^{\frac{3}{2}} h_R^{\varepsilon, \mathfrak{c}}(0)\big\|_{L^{\infty}}+C\sup_{s \in\left[0, \bar{T}\right]}\varepsilon^{\frac{5}{2}}\big\|\big(E_R^{\varepsilon,\mathfrak{c}}, B_R^{\varepsilon,\mathfrak{c}}\big)(s)\big\|_{L^{\infty}}\nonumber\\
			&+C \sup_{s \in\left[0, \bar{T}\right]}\big\|f_R^{\varepsilon, \mathfrak{c}}(s)\big\|+C \varepsilon^{k+\frac{5}{2}},
		\end{align}
		where we have taken $\varepsilon$ sufficiently small that  $Ce^{-\frac{\nu_0 \bar{T}}{2\varepsilon}}\leq \frac{1}{2}$.
		Therefore the proof of Lemma  \ref{pC.2} is completed.
	\end{proof}
	
	\begin{lemma}\label{pC.3}
		 Assume the {\it a priori} assumptions \eqref{W.0}-\eqref{C.2}. For any given $T>0$, taking $\varepsilon\ll 1$, we have that for $t \in [0,T]$,
			\begin{align}\label{C.17}
				\big\|\varepsilon^{\frac{3}{2}} h_R^{\varepsilon, \mathfrak{c}}(t)\big\|_{L^{\infty}} \leq& C\big\|\varepsilon^{\frac{3}{2}} h_R^{\varepsilon, \mathfrak{c}}(0)\big\|_{L^{\infty}}+C\sup_{s \in\left[0, t\right]}\big\|\varepsilon^{\frac{5}{2}}\big(E_R^{\varepsilon,\mathfrak{c}}, B_R^{\varepsilon,\mathfrak{c}}\big)(s)\big\|_{L^{\infty}}\nonumber\\
				&+C \sup_{s \in\left[0, t\right]}\big\|f_R^{\varepsilon, \mathfrak{c}}(s)\big\|+C \varepsilon^{k+\frac{5}{2}},
			\end{align}
		where $C>0$ is independent of $\mathfrak{c}$.
	\end{lemma}
	\begin{proof}
		For $t \leq \bar{T}$, \eqref{C.17} follows directly from \eqref{C.15}. 
        
        For $\bar{T} \leq t \leq T$, there exists a positive integer $n$ so that $t=n \bar{T}+\tau$, where $0 \leq \tau \leq \bar{T}$. Applying \eqref{C.16} repeatedly, we get
		{\small
		\begin{align*}
			&\varepsilon^{3 / 2}\|h_R^{\varepsilon, \mathfrak{c}}(t)\|_{L^{\infty}} \\
			\leq& \frac{1}{2}\big\|\varepsilon^{3 / 2} h_R^{\varepsilon, \mathfrak{c}}\left(\{n-1\} \bar{T}+\tau\right)\big\|_{L^{\infty}}+C\Big\{\sup_{s \in\left[0, t\right]}\big\|\varepsilon^{\frac{5}{2}}\big(E_R^{\varepsilon,\mathfrak{c}}, B_R^{\varepsilon,\mathfrak{c}}\big)(s)\big\|_{L^{\infty}}+\sup_{s \in\left[0, t\right]}\big\|f_R^{\varepsilon, \mathfrak{c}}(s)\big\|+\varepsilon^{k+\frac{5}{2}} \Big\} \\
			\leq& \frac{1}{4}\big\|\varepsilon^{3 / 2} h_R^{\varepsilon, \mathfrak{c}}\left(\{n-2\} \bar{T}+\tau\right)\big\|_{L^{\infty}}+\Big\{\frac{C}{2}+C\Big\}\Big\{\sup_{s \in\left[0, t\right]}\big\|\varepsilon^{\frac{5}{2}}\big(E_R^{\varepsilon,\mathfrak{c}}, B_R^{\varepsilon,\mathfrak{c}}\big)(s)\big\|_{L^{\infty}}+\sup_{s \in\left[0, t\right]}\big\|f_R^{\varepsilon, \mathfrak{c}}(s)\big\|+\varepsilon^{k+\frac{5}{2}}\Big\} \\
			\leq& \cdots 
			\leq \frac{1}{2^n}\big\|\varepsilon^{3 / 2} h_R^{\varepsilon, \mathfrak{c}}(\tau)\big\|_{L^{\infty}}+2 C\Big\{\sup_{s \in\left[0, t\right]}\big\|\varepsilon^{\frac{5}{2}}\big(E_R^{\varepsilon,\mathfrak{c}}, B_R^{\varepsilon,\mathfrak{c}}\big)(s)\big\|_{L^{\infty}}+\sup_{s \in\left[0, t\right]}\big\|f_R^{\varepsilon, \mathfrak{c}}(s)\big\|+\varepsilon^{k+\frac{5}{2}}\Big\},
		\end{align*}}
		which, together with $0 \leq \tau \leq \bar{T}$ and \eqref{C.15}, concludes \eqref{C.17}. Therefore the proof of Lemma \ref{pC.3} is completed. 
	\end{proof}
	
		
	\section{$L^{ \infty}$-Estimate for $\nabla_{x,p} h_R^{\varepsilon,\mathfrak{c}}$} \label{section 6}
	In this section, we present the $L^{ \infty}$ estimates for $\nabla_{x,p} h_R^{\varepsilon,\mathfrak{c}}$.
	
    	\begin{lemma}\label{pD.2}
		For $\bar{T}$ obtained in Lemma \ref{lC.1}, there exists a sufficiently small $\varepsilon_0>0$ such that for $\varepsilon\in (0,\varepsilon_0]$,
		\begin{align}\label{D.2}
			&\sup_{s \in\left[0, \bar{T}\right]}\big\|\varepsilon^{\frac{3}{2}} \nabla_{x,p}h_R^{\varepsilon,\mathfrak{c}}(s)\big\|_{L^{\infty}} \leq C\big\|\varepsilon^{\frac{3}{2}} h_R^{\varepsilon,\mathfrak{c}}(0)\big\|_{L^{\infty}}+C\big\|\varepsilon^{\frac{3}{2}} \nabla_{x,p}h_R^{\varepsilon,\mathfrak{c}}(0)\big\|_{L^{\infty}}\nonumber\\
			&\quad\quad+C\varepsilon^{\frac{5}{2}} \sup_{s \in\left[0, \bar{T}\right]}\big\|\big(E_R^{\varepsilon,\mathfrak{c}}, B_R^{\varepsilon,\mathfrak{c}}\big)(s)  \big\|_{W^{1,\infty}}+C\sup_{s \in\left[0, \bar{T}\right]}\big\|f_R^{\varepsilon, \mathfrak{c}}(s)\big\|_{H^1}+C\varepsilon^{k+\frac{5}{2}},
		\end{align}
		and
		\begin{align}\label{D.3}
			&\big\|\varepsilon^{\frac{3}{2}} \nabla_{x,p}h_R^{\varepsilon,\mathfrak{c}}(\bar{T})\big\|_{L^{\infty}} \leq \frac{1}{2}\left\{\big\|\varepsilon^{\frac{3}{2}} h_R^{\varepsilon,\mathfrak{c}}(0)\big\|_{L^{\infty}}+\big\|\varepsilon^{\frac{3}{2}} \nabla_{x,p}h_R^{\varepsilon,\mathfrak{c}}(0)\big\|_{L^{\infty}}\right\}\nonumber\\
			&\quad\quad+C\varepsilon^{\frac{5}{2}} \sup_{s \in\left[0, \bar{T}\right]}\big\|\big(E_R^{\varepsilon,\mathfrak{c}}, B_R^{\varepsilon,\mathfrak{c}}\big)(s)  \big\|_{W^{1,\infty}}+C\sup_{s \in\left[0, \bar{T}\right]}\big\|f_R^{\varepsilon, \mathfrak{c}}(s)\big\|_{H^1}+C\varepsilon^{k+\frac{5}{2}},
		\end{align}
		where $C>0$ is independent of $\mathfrak{c}$.
	\end{lemma}
	
	\begin{proof}
		We only prove \eqref{D.3}, since the estimate \eqref{D.2} can be done in the same way. 
		
		\noindent{\it Step 1. Estimate on $\big\| \nabla_x h_R^{\varepsilon,\mathfrak{c}}(s)\big\|_{L^{\infty}}$.}
		Let $D_x$ be any $x$ derivative. We apply $D_x$ to \eqref{C.18} to have
		{\small
		\begin{align}\label{D.4}
			\partial_t&\left(D_x h_R^{\varepsilon,\mathfrak{c}}\right)+\hat{p} \cdot \nabla_x\left(D_x h_R^{\varepsilon,\mathfrak{c}}\right)-\Big(E^{\varepsilon,\mathfrak{c}}+\frac{p}{p^0} \times B^{\varepsilon,\mathfrak{c}}\Big) \cdot \nabla_p\left(D_x h_R^{\varepsilon,\mathfrak{c}}\right)+\frac{1}{\varepsilon}\nu_{\mathfrak{c}}(p)\left(D_x h_R^{\varepsilon,\mathfrak{c}}\right) \nonumber\\
			=&D_x\Big(E^{\varepsilon,\mathfrak{c}}+\frac{p}{p^0} \times B^{\varepsilon,\mathfrak{c}}\Big) \cdot \nabla_p h_R^{\varepsilon,\mathfrak{c}}-\frac{1}{\varepsilon} D_x \big(\nu_{\mathfrak{c}}(p)\big) h_R^{\varepsilon,\mathfrak{c}}+\frac{1}{\varepsilon}D_x\big(\mathcal{K}_{\mathfrak{c},w}\left(h_R^{\varepsilon,\mathfrak{c}}\right)\big) \nonumber\\
			& +\sum_{i=1}^{2k-1} \varepsilon^{i-1} \frac{w_{\ell}}{\sqrt{J_{\mathbf{M}}}} D_x\Big[Q_{\mathfrak{c}}\left(F_i^{\mathfrak{c}}, \frac{\sqrt{J_{\mathbf{M}}}}{w_{\ell}} h_R^{\varepsilon,\mathfrak{c}}\right)+Q_{\mathfrak{c}}\left(\frac{\sqrt{J_{\mathbf{M}}}}{w_{\ell}} h_R^{\varepsilon,\mathfrak{c}}, F_i^{\mathfrak{c}}\right)\Big] \nonumber\\
			&	+\varepsilon^{k-1} \frac{w_{\ell}}{\sqrt{J_{\mathbf{M}}}} D_x \Big[Q_{\mathfrak{c}}\left(\frac{\sqrt{J_{\mathbf{M}}}}{w_{\ell}} h_R^{\varepsilon,\mathfrak{c}}, \frac{\sqrt{J_{\mathbf{M}}}}{w_{\ell}} h_R^{\varepsilon,\mathfrak{c}}\right)\Big]+D_x\Big[\Big(E^{\varepsilon,\mathfrak{c}}+\frac{p}{p^0} \times B^{\varepsilon,\mathfrak{c}}\Big) \cdot \nabla_p\left(\frac{\sqrt{J_{\mathbf{M}}}}{w_{\ell}}\right)\frac{w_{\ell}}{\sqrt{J_{\mathbf{M}}}}h_R^{\varepsilon,\mathfrak{c}}\Big]\nonumber\\
			& +D_x\Big[\frac{w_{\ell}}{\sqrt{J_{\mathbf{M}}}}\Big(E_R^{\varepsilon,\mathfrak{c}}+\frac{p}{p^0} \times B_R^{\varepsilon,\mathfrak{c}}\Big) \cdot \nabla_p\Big(\mathbf{M}_{\mathfrak{c}}+\sum_{i=1}^{2k-1} \varepsilon^i F_i^{\mathfrak{c}}\Big)\Big]+\varepsilon^k \frac{w_{\ell}}{\sqrt{J_{\mathbf{M}}}} D_x A .
		\end{align}}
		Then $D_x h_R^{\varepsilon,\mathfrak{c}}$ can be expressed as follows:
		{\footnotesize
		\begin{align}\label{D.5}
			D_x h_R^{\varepsilon,\mathfrak{c}}(t, x, p)=&\exp \Big(-\frac{\tilde{\nu}_{\mathfrak{c}}(t, 0)}{\varepsilon}\Big) D_x h_R^{\varepsilon,\mathfrak{c}}(0, X(0 ; t, x, p), P(0 ; t, x, p)) \nonumber\\
			&\quad+\int_0^t \exp \Big(-\frac{\tilde{\nu}_{\mathfrak{c}}(t, s)}{\varepsilon}\Big) \left\{D_x\Big(E^{\varepsilon,\mathfrak{c}}+\frac{P(s)}{P^0(s)} \times B^{\varepsilon,\mathfrak{c}}\Big) \cdot \nabla_p h_R^{\varepsilon,\mathfrak{c}}\right\}(s, X(s), P(s)) d s \nonumber\\
			&\quad-\int_0^t \exp \Big(-\frac{\tilde{\nu}_{\mathfrak{c}}(t, s)}{\varepsilon}\Big)\left(\frac{1}{\varepsilon}D_x \big(\nu_{\mathfrak{c}}(p)\big) h_R^{\varepsilon,\mathfrak{c}}\right)(s, X(s), P(s)) d s\nonumber \\
			&\quad+\int_0^t \exp \Big(-\frac{\tilde{\nu}_{\mathfrak{c}}(t, s)}{\varepsilon}\Big)\left(\frac{1}{\varepsilon} D_x \big(\mathcal{K}_{\mathfrak{c},w}\left(h_R^{\varepsilon,\mathfrak{c}}\right)\big)\right)(s, X(s), P(s)) d s \nonumber\\
			&\quad+\varepsilon^{k-1} \int_0^t \exp \Big(-\frac{\tilde{\nu}_{\mathfrak{c}}(t, s)}{\varepsilon}\Big) \left\{\frac{w_{\ell}}{\sqrt{J_{\mathbf{M}}}} D_x \Big[Q_{\mathfrak{c}}\left(\frac{\sqrt{J_{\mathbf{M}}}}{w_{\ell}} h_R^{\varepsilon,\mathfrak{c}}, \frac{\sqrt{J_{\mathbf{M}}}}{w_{\ell}} h_R^{\varepsilon,\mathfrak{c}}\right)\Big]\right\}(s, X(s), P(s)) d s \nonumber\\
			&\quad+\sum_{i=1}^{2k-1} \varepsilon^{i-1} \int_0^t \exp \Big(-\frac{\tilde{\nu}_{\mathfrak{c}}(t, s)}{\varepsilon}\Big) \left\{\frac{w_{\ell}}{\sqrt{J_{\mathbf{M}}}} D_x\Big[Q_{\mathfrak{c}}\left(F_i^{\mathfrak{c}}, \frac{\sqrt{J_{\mathbf{M}}}}{w_{\ell}} h_R^{\varepsilon,\mathfrak{c}}\right)+Q_{\mathfrak{c}}\left(\frac{\sqrt{J_{\mathbf{M}}}}{w_{\ell}} h_R^{\varepsilon,\mathfrak{c}}, F_i^{\mathfrak{c}}\right)\Big]\right\}(s, X(s), P(s)) d s \nonumber\\
			&\quad+\int_0^t \exp \Big(-\frac{\tilde{\nu}_{\mathfrak{c}}(t, s)}{\varepsilon}\Big) D_x\Big[\Big(E^{\varepsilon,\mathfrak{c}}+\frac{P(s)}{P^0(s)} \times B^{\varepsilon,\mathfrak{c}}\Big) \cdot \nabla_p\left(\frac{\sqrt{J_{\mathbf{M}}}}{w_{\ell}}\right)\frac{w_{\ell}}{\sqrt{J_{\mathbf{M}}}}h_R^{\varepsilon,\mathfrak{c}}\Big](s, X(s), P(s)) d s \nonumber\\
			&\quad+\int_0^t \exp \Big(-\frac{\tilde{\nu}_{\mathfrak{c}}(t, s)}{\varepsilon}\Big) D_x\Big[\frac{w_{\ell}}{\sqrt{J_{\mathbf{M}}}}\Big(E_R^{\varepsilon,\mathfrak{c}}+\frac{P(s)}{P^0(s)} \times B_R^{\varepsilon,\mathfrak{c}}\Big) \cdot \nabla_p\Big(\mathbf{M}_{\mathfrak{c}}+\sum_{i=1}^{2k-1} \varepsilon^i F_i^{\mathfrak{c}}\Big)\Big](s, X(s), P(s)) d s \nonumber\\
			&\quad+\int_0^t \exp \Big(-\frac{\tilde{\nu}_{\mathfrak{c}}(t, s)}{\varepsilon}\Big) \left(\varepsilon^k \frac{w_{\ell}}{\sqrt{J_{\mathbf{M}}}} D_x A\right)(s, X(s), P(s)) d s
			:=\sum_{i=1}^9 G_i .
		\end{align}}
		Next we estimate the terms on RHS of \eqref{D.5} one by one.
		
		It is clear that
		$$
		\big|\varepsilon^{\frac{3}{2}} G_1\big| \leq e^{-\frac{\nu_0 t}{\varepsilon}}\big\|\varepsilon^{\frac{3}{2}} D_x h_R^{\varepsilon,\mathfrak{c}}(0)\big\|_{L^{\infty}}.
		$$
		Noting $\eqref{C.3}$, one has
		$$
		\big|\varepsilon^{\frac{3}{2}}G_2\big|\lesssim\varepsilon e^{-\frac{\nu_0 t}{2\varepsilon}} \sup_{s \in\left[0, t\right]}\left\{e^{\frac{\nu_0 s}{2\varepsilon}}\big\|\varepsilon^{\frac{3}{2}} \nabla_{p}h_R^{\varepsilon,\mathfrak{c}}(s)\big\|_{L^{\infty}}\right\}.
		$$
		By \eqref{e2.1}, together with the parameter-derivative estimate on 
    $\mathbf{M}_{\mathfrak{c}}$ in Lemma \ref{l2.12}, one has $|D_x \nu_{\mathfrak{c}}(p)|\lesssim\nu_{\mathfrak{c}}(p)$. Thus it holds
		$$
		\big|\varepsilon^{\frac{3}{2}}G_3\big|\lesssim e^{-\frac{\nu_0 t}{2\varepsilon}} \sup_{s \in\left[0, t\right]}\left\{e^{\frac{\nu_0 s}{2\varepsilon}}\big\|\varepsilon^{\frac{3}{2}} h_R^{\varepsilon,\mathfrak{c}}(s)\big\|_{L^{\infty}}\right\}.
		$$
		
		For $G_5$, we note
		$$
		D_x \Big[Q_{\mathfrak{c}}\Big(\frac{\sqrt{J_{\mathbf{M}}}}{w_{\ell}} h_R^{\varepsilon,\mathfrak{c}}, \frac{\sqrt{J_{\mathbf{M}}}}{w_{\ell}} h_R^{\varepsilon,\mathfrak{c}}\Big)\Big]=Q_{\mathfrak{c}}\Big(\frac{\sqrt{J_{\mathbf{M}}}}{w_{\ell}} D_x h_R^{\varepsilon,\mathfrak{c}}, \frac{\sqrt{J_{\mathbf{M}}}}{w_{\ell}} h_R^{\varepsilon,\mathfrak{c}}\Big)+Q_{\mathfrak{c}}\Big(\frac{\sqrt{J_{\mathbf{M}}}}{w_{\ell}} h_R^{\varepsilon,\mathfrak{c}}, \frac{\sqrt{J_{\mathbf{M}}}}{w_{\ell}} D_x h_R^{\varepsilon,\mathfrak{c}}\Big),
		$$
		which, together with Lemma \ref{l2.2}, yields that
		\begin{align}\label{D.5-1}
		\Big|\frac{w_{\ell}}{\sqrt{J_{\mathbf{M}}}} D_x \Big[Q_{\mathfrak{c}}\Big(\frac{\sqrt{J_{\mathbf{M}}}}{w_{\ell}} h_R^{\varepsilon,\mathfrak{c}}, \frac{\sqrt{J_{\mathbf{M}}}}{w_{\ell}} h_R^{\varepsilon,\mathfrak{c}}\Big)\Big]\Big|\lesssim  \nu_{\mathfrak{c}}(p) \big\|h_R^{\varepsilon,\mathfrak{c}}\big\|_{L^{\infty}}\big\| D_x h_R^{\varepsilon,\mathfrak{c}}\big\|_{L^{\infty}}.
		\end{align}
		Thus we have from \eqref{D.5-1} that
		\begin{align}
			\big|\varepsilon^{\frac{3}{2}}G_5\big|&\lesssim \varepsilon^{k-\frac{3}{2}} e^{-\frac{\nu_0 t}{\varepsilon}} \sup_{s \in\left[0, t\right]}\left\{\left(e^{\frac{\nu_0 s}{2\varepsilon}}\big\|\varepsilon^{\frac{3}{2}} h_R^{\varepsilon,\mathfrak{c}}(s)\big\|_{L^{\infty}}\right)\left(e^{\frac{\nu_0 s}{2\varepsilon}}\big\|\varepsilon^{\frac{3}{2}} \nabla_{x}h_R^{\varepsilon,\mathfrak{c}}(s)\big\|_{L^{\infty}}\right)\right\} \nonumber\\
			&\lesssim \varepsilon e^{-\frac{\nu_0 t}{2\varepsilon}} \sup_{s \in\left[0, t\right]}\left\{e^{\frac{\nu_0 s}{2\varepsilon}}\big\|\varepsilon^{\frac{3}{2}}\nabla_{x}h_R^{\varepsilon,\mathfrak{c}}(s)\big\|_{L^{\infty}}\right\},\nonumber
		\end{align}
		where we have used $\eqref{W.0}$ and $k\geq 5$.
		
		Noting
		$\Big|\nabla_p\left(\frac{\sqrt{J_{\mathbf{M}}}}{w_{\ell}}\right)\frac{w_{\ell}}{\sqrt{J_{\mathbf{M}}}}\Big|\lesssim\nu_{\mathfrak{c}}(p)$,
		we get from \eqref{C.3}
		\begin{align}
			\big|\varepsilon^{\frac{3}{2}}G_7\big|& \lesssim \int_0^t \exp \Big(-\frac{\tilde{\nu}_{\mathfrak{c}}(t, s)}{\varepsilon}\Big)\nu_{\mathfrak{c}}(s)
			\cdot\sup_{s \in [0, t]}\left\{\big\|\left(E^{\varepsilon,\mathfrak{c}}, B^{\varepsilon,\mathfrak{c}}\right)(s)\big\|_{W^{1,\infty}}\big\|\varepsilon^{\frac{3}{2}}h_R^{\varepsilon,\mathfrak{c}}(s)\big\|_{W^{1, \infty}}\right\} d s \nonumber\\
			& \lesssim \varepsilon e^{-\frac{\nu_0 t}{2\varepsilon}} \sup_{s \in\left[0, t\right]}\left\{e^{\frac{\nu_0 s}{2\varepsilon}}\big\|\varepsilon^{\frac{3}{2}}h_R^{\varepsilon,\mathfrak{c}}(s)\big\|_{W^{1, \infty}}\right\} .\nonumber
		\end{align}
		Similarly, there holds 
		$$
		\big|\varepsilon^{\frac{3}{2}}\big(G_6,G_8,G_9\big)\big| \lesssim \varepsilon e^{-\frac{\nu_0 t}{2\varepsilon}} \sup_{s \in\left[0, t\right]}\left\{e^{\frac{\nu_0 s}{2\varepsilon}}\big\|\varepsilon^{\frac{3}{2}}h_R^{\varepsilon,\mathfrak{c}}(s)\big\|_{W^{1, \infty}}\right\}+\varepsilon^{\frac{5}{2}} \sup_{s \in\left[0, \bar{T}\right]}\big\|\big(E_R^{\varepsilon,\mathfrak{c}}, B_R^{\varepsilon,\mathfrak{c}}\big)(s)  \big\|_{W^{1,\infty}}+\varepsilon^{k+\frac{5}{2}}.
		$$
		
	We shall concentrate on the fourth term on the right-hand side of \eqref{D.5}. It follows from \eqref{e2.41} that
		\begin{align}
			G_4=&\frac{1}{\varepsilon} \int_0^t \exp \Big(-\frac{\tilde{\nu}_{\mathfrak{c}}(t, s)}{\varepsilon}\Big) \big(D_x \left(\mathcal{K}_{\mathfrak{c},w}\left(h_R^{\varepsilon,\mathfrak{c}}\right)\right)\big)(s, X(s), P(s)) d s \nonumber\\
			\leq&\frac{1}{\varepsilon} \int_0^t \exp \Big(-\frac{\tilde{\nu}_{\mathfrak{c}}(t, s)}{\varepsilon}\Big)\int_{\mathbb{R}^3} \hat{k}_w(P(s), q)\big|h_R^{\varepsilon,\mathfrak{c}}(s, X(s), q)\big| d q d s\nonumber\\
			&+\frac{1}{\varepsilon} \int_0^t \exp \Big(-\frac{\tilde{\nu}_{\mathfrak{c}}(t, s)}{\varepsilon}\Big)\int_{\mathbb{R}^3} k_w(P(s), q)\big|D_x h_R^{\varepsilon,\mathfrak{c}}(s, X(s), q)\big| d q d s\nonumber\\
			:=&\uppercase\expandafter{\romannumeral1}+\uppercase\expandafter{\romannumeral2}. \nonumber
		\end{align}
		By \eqref{e2.41-1}, one obtains
		\begin{align}
			\varepsilon^{\frac{3}{2}}\uppercase\expandafter{\romannumeral1}\lesssim& \frac{t}{\varepsilon}e^{-\frac{\nu_0 t}{\varepsilon}}\big\|\varepsilon^{\frac{3}{2}} h_R^{\varepsilon,\mathfrak{c}}(0)\big\|_{L^{\infty}}+ \Big(\varepsilon+\kappa+\frac{1}{N_0}\Big) e^{-\frac{\nu_0 t}{2\varepsilon}}\sup_{s \in\left[0, t\right]}\left\{e^{\frac{\nu_0 s}{2\varepsilon}}\big\|\varepsilon^{\frac{3}{2}}h_R^{\varepsilon,\mathfrak{c}}(s)\big\|_{L^{\infty}}\right\} \nonumber\\
			&+ \sup_{s \in\left[0, \bar{T}\right]}\varepsilon^{\frac{5}{2}}\big\|\big(E_R^{\varepsilon,\mathfrak{c}}, B_R^{\varepsilon,\mathfrak{c}}\big)(s)\big\|_{L^{\infty}}+ \sup_{s \in\left[0, \bar{T}\right]}\big\|f_R^{\varepsilon, \mathfrak{c}}(s)\big\| +\varepsilon^{k+\frac{5}{2}}.\nonumber
		\end{align}
		
		We now use \eqref{D.5} again to evaluate \uppercase\expandafter{\romannumeral2}. The estimate is similar to that in Lemma \ref{pC.2}, except for the final part, where \(\big\|f_R^{\varepsilon, \mathfrak{c}}\big\|_{H^1}\) appears in place of \(\big\|f_R^{\varepsilon, \mathfrak{c}}\big\|\). Here we will only present this last case: recall the
		Case 4. As before we denote 
		$$
		X^\prime(s^{\prime}):=X(s^\prime;s,X(s;t,x,p),q), \quad P^\prime(s^\prime):=P(s^\prime;s,X(s;t,x,p),q)
		$$ 
		and
		$$
		\tilde{\nu}_{\mathfrak{c}}^{\prime}\left(s, s^{\prime}\right):=\int_{s^{\prime}}^s \nu_{\mathfrak{c}}\left(\tau, X^\prime(\tau), P^\prime(\tau)\right) d \tau .
		$$
		For $s-s^{\prime} \geq \kappa \varepsilon$ and $|p| \leq N_0,|q| \leq 2 N_0,|q^\prime| \leq 3 N_0$. By \eqref{e2.41-1}, Cauchy's inequality, and the change of variables, one has
		\begin{align}
			& \int_{|q| \leq 2 N_0} \int_{|q^\prime| \leq 3 N_0} \hat{k}_w(P(s), q) \hat{k}_w\left(P^\prime(s^\prime), q^\prime\right)\big| D_x h_R^{\varepsilon, \mathfrak{c}}\left(s^\prime, X^\prime(s^\prime), q^\prime\right)\big| d q^\prime d q \nonumber\\
			\leq& C_{N_0} \int_{|q| \leq 2 N_0} \int_{|q^\prime| \leq 3 N_0} \hat{k}_w(P(s), q) \hat{k}_w\left(P^\prime(s^\prime), q^\prime\right)\left(\big|f_R^{\varepsilon, \mathfrak{c}}\big|+\big|D_x f_R^{\varepsilon, \mathfrak{c}}\big|\right)\left(s^\prime, X^\prime(s^\prime), q^\prime\right) d q^\prime d q \nonumber\\
			\leq& C_{N_0}\Big(\int_{|q| \leq 2 N_0} \int_{|q^\prime| \leq 3 N_0} \hat{k}_w^2(P(s), q) \hat{k}_w^2\left(P^\prime(s^\prime), q^\prime\right) d q^\prime d q\Big)^{\frac{1}{2}} \nonumber\\
			& \quad \times\Big(\int_{|q| \leq 2 N_0} \int_{|q^\prime| \leq 3 N_0}\big|f_R^{\varepsilon, \mathfrak{c}}\left(s^\prime, X^\prime(s^\prime), q^\prime\right)\big|^2 d q^\prime d q\Big)^{\frac{1}{2}} \nonumber\\
			&+ C_{N_0}\Big(\int_{|q| \leq 2 N_0} \int_{|q^\prime| \leq 3 N_0} \hat{k}_w^2(P(s), q) \hat{k}_w^2\left(P^\prime(s^\prime), q^\prime\right) d q^\prime d q\Big)^{\frac{1}{2}} \nonumber\\
			& \quad \times\Big(\int_{|q| \leq 2 N_0} \int_{|q^\prime| \leq 3 N_0}\big|D_x f_R^{\varepsilon, \mathfrak{c}}\left(s^\prime, X^\prime(s^\prime), q^\prime\right)\big|^2 d q^\prime d q\Big)^{\frac{1}{2}} \nonumber\\
			\leq& C_{N_0}\Big(\int_{\mathbb{R}^3} \int_{\mathbb{R}^3}\big|f_R^{\varepsilon, \mathfrak{c}}\left(s^\prime, X^\prime(s^\prime), q^\prime\right)\big|^2 \cdot \varepsilon^{-3} \kappa^{-3} d X^\prime(s^\prime) d q^\prime\Big)^{\frac{1}{2}} \nonumber\\
			&+ C_{N_0}\Big(\int_{\mathbb{R}^3} \int_{\mathbb{R}^3}\big|D_x f_R^{\varepsilon, \mathfrak{c}}\left(s^\prime, X^\prime(s^\prime), q^\prime\right)\big|^2 \cdot \varepsilon^{-3} \kappa^{-3} d X^\prime(s^\prime) d q^\prime\Big)^{\frac{1}{2}} \nonumber\\
			\leq& \frac{C_{N_0, \kappa}}{\varepsilon^{\frac{3}{2}}} \sup_{s \in\left[0, \bar{T}\right]}\big\|f_R^{\varepsilon, \mathfrak{c}}(s)\big\|_{H^1},\nonumber
		\end{align}
		where we have made a change of variables $q \mapsto X^\prime(s^\prime)$ with
		$$
		\Big|\operatorname{det}\Big(\frac{d X^\prime(s^\prime)}{d q}\Big)\Big|\geq\frac{\mathfrak{c}^5}{2\left(q^0\right)^5}\left(s-s^{\prime}\right)^3 \geq \frac{\kappa^3 \varepsilon^3}{2^6} 
		$$
		for $1 \leq N_0 \ll \mathfrak{c}$. Thus we have
		\begin{align}
			& \frac{1}{\varepsilon^2} \int_0^t \exp \Big(-\frac{\tilde{\nu}_{\mathfrak{c}}(t, s)}{\varepsilon}\Big) d s \int_0^{s-\kappa \varepsilon} \exp \Big(-\frac{\tilde{\nu}_{\mathfrak{c}}^\prime\left(s, s^{\prime}\right)}{\varepsilon}\Big) d s^{\prime}\nonumber\\
			& \quad \times \int_{|q| \leq 2 N_0} \hat{k}_w(P(s), q) d q \int_{|q^\prime| \leq 3 N_0} \hat{k}_w\left(P^\prime(s^\prime), q^\prime\right)\Big|\varepsilon^{\frac{3}{2}} D_x h_R^{\varepsilon, \mathfrak{c}}\left(s^\prime, X^\prime(s^\prime), q^\prime\right)\Big| d q^\prime \nonumber\\
			& \leq C_{N_0, \kappa} \sup_{s \in\left[0, \bar{T}\right]}\big\|f_R^{\varepsilon, \mathfrak{c}}(s)\big\|_{H^1} .\nonumber
		\end{align}	
		Similarly we have
		\begin{align}
			&\varepsilon^{\frac{3}{2}}\uppercase\expandafter{\romannumeral2}\lesssim\frac{t}{\varepsilon}e^{-\frac{\nu_0 t}{\varepsilon}}\big\|\varepsilon^{\frac{3}{2}} D_x h_R^{\varepsilon,\mathfrak{c}}(0)\big\|_{L^{\infty}}+\varepsilon e^{-\frac{\nu_0 t}{2\varepsilon}} \sup_{s \in\left[0, t\right]}\left\{e^{\frac{\nu_0 s}{2\varepsilon}}\big\|\varepsilon^{\frac{3}{2}} \nabla_{p}h_R^{\varepsilon,\mathfrak{c}}(s)\big\|_{L^{\infty}}\right\}\nonumber\\
			&+(1+\varepsilon)e^{-\frac{\nu_0 t}{2\varepsilon}} \sup_{s \in\left[0, t\right]}\left\{e^{\frac{\nu_0 s}{2\varepsilon}}\big\|\varepsilon^{\frac{3}{2}} h_R^{\varepsilon,\mathfrak{c}}(s)\big\|_{L^{\infty}}\right\}+\Big(\varepsilon+\kappa+\frac{1}{N_0}\Big) e^{-\frac{\nu_0 t}{2\varepsilon}} \sup_{s \in\left[0, t\right]}\left\{e^{\frac{\nu_0 s}{2\varepsilon}}\big\|\varepsilon^{\frac{3}{2}}\nabla_x h_R^{\varepsilon,\mathfrak{c}}(s)\big\|_{L^{\infty}}\right\}\nonumber\\
			&+\varepsilon^{\frac{5}{2}} \sup_{s \in\left[0, \bar{T}\right]}\big\|\big(E_R^{\varepsilon,\mathfrak{c}}, B_R^{\varepsilon,\mathfrak{c}}\big)(s)  \big\|_{W^{1,\infty}}+\sup_{s \in\left[0, \bar{T}\right]}\big\|f_R^{\varepsilon, \mathfrak{c}}(s)\big\|_{H^1}+\varepsilon^{k+\frac{5}{2}}.\nonumber
		\end{align}
		In summary, for any $x$ derivative $D_x$, we have established
		\begin{align}\label{D.6}
			\varepsilon^{\frac{3}{2}}\big\| D_x h_R^{\varepsilon,\mathfrak{c}}(t)\big\|_{L^{\infty}}\lesssim&\frac{t}{\varepsilon}e^{-\frac{\nu_0 t}{\varepsilon}}\big\|\varepsilon^{\frac{3}{2}} h_R^{\varepsilon,\mathfrak{c}}(0)\big\|_{L^{\infty}}+\Big(1+\frac{t}{\varepsilon}\Big)e^{-\frac{\nu_0 t}{\varepsilon}}\big\|\varepsilon^{\frac{3}{2}} D_x h_R^{\varepsilon,\mathfrak{c}}(0)\big\|_{L^{\infty}}\nonumber\\
			&+\varepsilon e^{-\frac{\nu_0 t}{2\varepsilon}} \sup_{s \in\left[0, t\right]}\left\{e^{\frac{\nu_0 s}{2\varepsilon}}\big\|\varepsilon^{\frac{3}{2}} \nabla_{p}h_R^{\varepsilon,\mathfrak{c}}(s)\big\|_{L^{\infty}}\right\}+e^{-\frac{\nu_0 t}{2\varepsilon}} \sup_{s \in\left[0, t\right]}\left\{e^{\frac{\nu_0 s}{2\varepsilon}}\big\|\varepsilon^{\frac{3}{2}} h_R^{\varepsilon,\mathfrak{c}}(s)\big\|_{L^{\infty}}\right\}\nonumber\\
			&+\Big(\varepsilon+\kappa+\frac{1}{N_0}\Big) e^{-\frac{\nu_0 t}{2\varepsilon}} \sup_{s \in\left[0, t\right]}\left\{e^{\frac{\nu_0 s}{2\varepsilon}}\big\|\varepsilon^{\frac{3}{2}}h_R^{\varepsilon,\mathfrak{c}}(s)\big\|_{W^{1,\infty}}\right\}\nonumber\\
			&+\varepsilon^{\frac{5}{2}} \sup_{s \in\left[0, \bar{T}\right]}\big\|\big(E_R^{\varepsilon,\mathfrak{c}}, B_R^{\varepsilon,\mathfrak{c}}\big)(s)  \big\|_{W^{1,\infty}}+\sup_{s \in\left[0, \bar{T}\right]}\big\|f_R^{\varepsilon, \mathfrak{c}}(s)\big\|_{H^1}+\varepsilon^{k+\frac{5}{2}}.
		\end{align}
		
		\noindent{\it Step 2. Estimate on $\big\| \nabla_p h_R^{\varepsilon,\mathfrak{c}}(s)\big\|_{L^{\infty}}$.}
		 Let $D_p$ be any $p$ derivative. Applying $D_p$ to \eqref{C.18}, one has
		{\small
		\begin{align}\label{D.7}
			&\partial_t\left(D_p h_R^{\varepsilon,\mathfrak{c}}\right)+\hat{p} \cdot \nabla_x\left(D_p h_R^{\varepsilon,\mathfrak{c}}\right)-\Big(E^{\varepsilon,\mathfrak{c}}+\frac{p}{p^0} \times B^{\varepsilon,\mathfrak{c}}\Big) \cdot \nabla_p \left(D_p h_R^{\varepsilon,\mathfrak{c}}\right)+\frac{1}{\varepsilon}\nu_{\mathfrak{c}}(p) D_p h_R^{\varepsilon,\mathfrak{c}} \nonumber\\
			&=  \frac{1}{\varepsilon} D_p\left(\mathcal{K}_{\mathfrak{c},w}\left(h_R^{\varepsilon,\mathfrak{c}}\right)\right) -D_p(\hat{p}) \cdot \nabla_x h_R^{\varepsilon,\mathfrak{c}}+\left(D_p\big(\frac{p}{p^0}\big) \times B^{\varepsilon,\mathfrak{c}}\right) \cdot \nabla_p h_R^{\varepsilon,\mathfrak{c}}-\frac{1}{\varepsilon}D_p \nu_{\mathfrak{c}}(p) h_R^{\varepsilon,\mathfrak{c}}\nonumber\\
			&\quad+\sum_{i=1}^{2 k-1} \varepsilon^{i-1} D_p\left\{\frac{w_{\ell}}{\sqrt{J_{\mathbf{M}}}}\Big[Q_{\mathfrak{c}}\left(F_i^{\mathfrak{c}}, \frac{\sqrt{J_{\mathbf{M}}}}{w_{\ell}} h_R^{\varepsilon,\mathfrak{c}}\right)+Q_{\mathfrak{c}}\left(\frac{\sqrt{J_{\mathbf{M}}}}{w_{\ell}} h_R^{\varepsilon,\mathfrak{c}}, F_i^{\mathfrak{c}}\right)\Big]\right\}\nonumber\\
			&\quad+\varepsilon^{k-1} D_p\left\{\frac{w_{\ell}}{\sqrt{J_{\mathbf{M}}}}Q_{\mathfrak{c}}\left(\frac{\sqrt{J_{\mathbf{M}}}}{w_{\ell}} h_R^{\varepsilon,\mathfrak{c}}, \frac{\sqrt{J_{\mathbf{M}}}}{w_{\ell}} h_R^{\varepsilon,\mathfrak{c}}\right)\right\}+D_p\left\{\Big(E^{\varepsilon,\mathfrak{c}}+\frac{p}{p^0} \times B^{\varepsilon,\mathfrak{c}}\Big)\cdot\nabla_p(\frac{\sqrt{J_{\mathbf{M}}}}{w_{\ell}})\frac{w_{\ell}}{\sqrt{J_{\mathbf{M}}}}h_R^{\varepsilon,\mathfrak{c}}\right\}\nonumber\\
			&\quad+D_p\Big\{\Big(E_R^{\varepsilon,\mathfrak{c}}+\frac{p}{p^0} \times B_R^{\varepsilon,\mathfrak{c}}\Big) \cdot \frac{w_{\ell}}{\sqrt{J_{\mathbf{M}}}} \nabla_p\Big(\mathbf{M}_{\mathfrak{c}}+\sum_{i=1}^{2 k-1} \varepsilon^i F_i^{\mathfrak{c}}\Big)\Big\}+\varepsilon^k D_p\left(\frac{w_{\ell}}{\sqrt{J_{\mathbf{M}}}} A\right).
		\end{align}}
		Similarly, $D_p h_R^{\varepsilon,\mathfrak{c}}$ can be expressed as follows:
		{\footnotesize
		\begin{align}\label{D.8}
			& D_p h_R^{\varepsilon,\mathfrak{c}}(t, x, p)\nonumber\\
            =&\exp \left(-\frac{\tilde{\nu}_{\mathfrak{c}}(t, 0)}{\varepsilon}\right) D_p h_R^{\varepsilon,\mathfrak{c}}(0, X(0 ; t, x, p), P(0 ; t, x, p)) \nonumber\\
			& \quad-\int_0^t \exp \Big(-\frac{\tilde{\nu}_{\mathfrak{c}}(t, s)}{\varepsilon}\Big) \big(D_p(\hat{p}) \cdot \nabla_x h_R^{\varepsilon,\mathfrak{c}}\big)(s, X(s), P(s)) d s \nonumber\\
			&\quad +\int_0^t \exp \Big(-\frac{\tilde{\nu}_{\mathfrak{c}}(t, s)}{\varepsilon}\Big)\Big(D_p\big(\frac{p}{p^0}\big) \times B^{\varepsilon,\mathfrak{c}} \cdot \nabla_p h_R^{\varepsilon,\mathfrak{c}}\Big)(s, X(s), P(s)) d s \nonumber\\
			&\quad-\int_0^t \exp \Big(-\frac{\tilde{\nu}_{\mathfrak{c}}(t, s)}{\varepsilon}\Big)\frac{1}{\varepsilon}\left(D_p \nu_{\mathfrak{c}}(p) h_R^{\varepsilon,\mathfrak{c}}\right)(s, X(s), P(s)) d s \nonumber\\
			& \quad+\int_0^t \exp \Big(-\frac{\tilde{\nu}_{\mathfrak{c}}(t, s)}{\varepsilon}\Big) \frac{1}{\varepsilon} D_p\left(\mathcal{K}_{\mathfrak{c},w}\left(h_R^{\varepsilon,\mathfrak{c}}\right)\right)(s, X(s), P(s)) d s \nonumber\\
			&\quad +\varepsilon^{k-1} \int_0^t \exp \Big(-\frac{\tilde{\nu}_{\mathfrak{c}}(t, s)}{\varepsilon}\Big) D_p\left\{\frac{w_{\ell}}{\sqrt{J_{\mathbf{M}}}}Q_{\mathfrak{c}}\left(\frac{\sqrt{J_{\mathbf{M}}}}{w_{\ell}} h_R^{\varepsilon,\mathfrak{c}}, \frac{\sqrt{J_{\mathbf{M}}}}{w_{\ell}} h_R^{\varepsilon,\mathfrak{c}}\right)\right\}(s, X(s), P(s)) d s \nonumber\\
			&\quad +\sum_{i=1}^{2 k-1} \varepsilon^{i-1} \int_0^t \exp \Big(-\frac{\tilde{\nu}_{\mathfrak{c}}(t, s)}{\varepsilon}\Big) D_p\left\{\frac{w_{\ell}}{\sqrt{J_{\mathbf{M}}}}\left[Q_{\mathfrak{c}}\left(F_i^{\mathfrak{c}}, \frac{\sqrt{J_{\mathbf{M}}}}{w_{\ell}} h_R^{\varepsilon,\mathfrak{c}}\right)+Q_{\mathfrak{c}}\left(\frac{\sqrt{J_{\mathbf{M}}}}{w_{\ell}} h_R^{\varepsilon,\mathfrak{c}}, F_i^{\mathfrak{c}}\right)\right]\right\}(s, X(s), P(s)) d s \nonumber\\
			&\quad +\int_0^t \exp \Big(-\frac{\tilde{\nu}_{\mathfrak{c}}(t, s)}{\varepsilon}\Big) D_p\left\{\Big(E^{\varepsilon,\mathfrak{c}}+\frac{P(s)}{P^0(s)}\times B^{\varepsilon,\mathfrak{c}}\Big) \cdot \nabla_p\left(\frac{\sqrt{J_{\mathbf{M}}}}{w_{\ell}}\right)\frac{w_{\ell}}{\sqrt{J_{\mathbf{M}}}}h_R^{\varepsilon,\mathfrak{c}}\right\}(s, X(s), P(s)) d s \nonumber\\
			& \quad+\int_0^t \exp \Big(-\frac{\tilde{\nu}_{\mathfrak{c}}(t, s)}{\varepsilon}\Big) D_p\left\{\Big(E_R^{\varepsilon,\mathfrak{c}}+\frac{P(s)}{P^0(s)} \times B_R^{\varepsilon,\mathfrak{c}}\Big) \cdot \frac{w_{\ell}}{\sqrt{J_{\mathbf{M}}}} \nabla_p\left(\mathbf{M}_{\mathfrak{c}}+\sum_{i=1}^{2 k-1} \varepsilon^i F_i^{\mathfrak{c}}\right)\right\}(s, X(s), P(s)) d s \nonumber\\
			&\quad +\varepsilon^k \int_0^t \exp \Big(-\frac{\tilde{\nu}_{\mathfrak{c}}(t, s)}{\varepsilon}\Big) D_p\left(\frac{w_{\ell}}{\sqrt{J_{\mathbf{M}}}} A\right)(s, X(s), P(s)) d s 
			:= \sum_{i=1}^{10}Z_i.
		\end{align}}
		Noting \eqref{e2.41}-\eqref{e2.41-1}, by similar arguments as Step 1, one gets from \eqref{D.8}
		\begin{align}\label{D.9}
			&\varepsilon^{\frac{3}{2}}\big\| D_p h_R^{\varepsilon,\mathfrak{c}}(t)\big\|_{L^{\infty}}\lesssim\frac{t}{\varepsilon}e^{-\frac{\nu_0 t}{\varepsilon}}\big\|\varepsilon^{\frac{3}{2}} h_R^{\varepsilon,\mathfrak{c}}(0)\big\|_{L^{\infty}}+\Big(1+\frac{t}{\varepsilon}\Big)e^{-\frac{\nu_0 t}{\varepsilon}}\big\|\varepsilon^{\frac{3}{2}} D_p h_R^{\varepsilon,\mathfrak{c}}(0)\big\|_{L^{\infty}}\nonumber\\
			&\quad+\varepsilon e^{-\frac{\nu_0 t}{2\varepsilon}} \sup_{s \in\left[0, t\right]}\left\{e^{\frac{\nu_0 s}{2\varepsilon}}\big\|\varepsilon^{\frac{3}{2}} \nabla_{x}h_R^{\varepsilon,\mathfrak{c}}(s)\big\|_{L^{\infty}}\right\}+e^{-\frac{\nu_0 t}{2\varepsilon}} \sup_{s \in\left[0, t\right]}\left\{e^{\frac{\nu_0 s}{2\varepsilon}}\big\|\varepsilon^{\frac{3}{2}} h_R^{\varepsilon,\mathfrak{c}}(s)\big\|_{L^{\infty}}\right\}\nonumber\\
			&\quad+\Big(\varepsilon+\kappa+\frac{1}{N_0}\Big) e^{-\frac{\nu_0 t}{2\varepsilon}} \sup_{s \in\left[0, t\right]}\left\{e^{\frac{\nu_0 s}{2\varepsilon}}\big\|\varepsilon^{\frac{3}{2}}h_R^{\varepsilon,\mathfrak{c}}(s)\big\|_{W^{1,\infty}}\right\}\\
			&\quad+\varepsilon^{\frac{5}{2}} \sup_{s \in\left[0, \bar{T}\right]}\big\|\big(E_R^{\varepsilon,\mathfrak{c}}, B_R^{\varepsilon,\mathfrak{c}}\big)(s)  \big\|_{W^{1,\infty}}+\sup_{s \in\left[0, \bar{T}\right]}\big\|f_R^{\varepsilon, \mathfrak{c}}(s)\big\|_{H^1}+\varepsilon^{k+\frac{5}{2}}.\nonumber
		\end{align}
		Combining \eqref{D.6} and \eqref{D.9}, one obtains
		\begin{align}\label{D.10}
			&\sup_{s \in\left[0, \bar{T}\right]}\left\{e^{\frac{\nu_0 s}{2\varepsilon}}\big\|\varepsilon^{\frac{3}{2}} \nabla_{x,p}h_R^{\varepsilon,\mathfrak{c}}(s)\big\|_{L^{\infty}}\right\} \leq C\big\|\varepsilon^{\frac{3}{2}} h_R^{\varepsilon,\mathfrak{c}}(0)\big\|_{L^{\infty}}+C\big\|\varepsilon^{\frac{3}{2}} \nabla_{x,p}h_R^{\varepsilon,\mathfrak{c}}(0)\big\|_{L^{\infty}}\nonumber\\
			&+ \sup_{s \in\left[0, \bar{T}\right]}\left\{e^{\frac{\nu_0 s}{2\varepsilon}}\big\|\varepsilon^{\frac{3}{2}} h_R^{\varepsilon,\mathfrak{c}}(s)\big\|_{L^{\infty}}\right\}+\Big(\varepsilon+\kappa+\frac{1}{N_0}\Big)  \sup_{s \in\left[0, \bar{T}\right]}\left\{e^{\frac{\nu_0 s}{2\varepsilon}}\big\|\varepsilon^{\frac{3}{2}}h_R^{\varepsilon,\mathfrak{c}}(s)\big\|_{W^{1,\infty}}\right\}\nonumber\\
			&+\varepsilon^{\frac{5}{2}} e^{\frac{\nu_0 \bar{T}}{2\varepsilon}}\sup_{s \in\left[0, \bar{T}\right]}\big\|\big(E_R^{\varepsilon,\mathfrak{c}}, B_R^{\varepsilon,\mathfrak{c}}\big)(s)  \big\|_{W^{1,\infty}}+e^{\frac{\nu_0 \bar{T}}{2\varepsilon}}\sup_{s \in\left[0, \bar{T}\right]}\big\|f_R^{\varepsilon, \mathfrak{c}}(s)\big\|_{H^1}+\varepsilon^{k+\frac{5}{2}}e^{\frac{\nu_0 \bar{T}}{2\varepsilon}}.
		\end{align}
		From Lemma  \ref{pC.3}, we have the estimate of $\sup_{s \in\left[0, \bar{T}\right]}\left\{\big\|\varepsilon^{\frac{3}{2}} h_R^{\varepsilon,\mathfrak{c}}(s)\big\|_{L^{\infty}}\right\}$. Using the weighted $L^\infty$ estimate obtained in the proof of Lemma \ref{pC.2}
        (on $[0,\bar T]$), together with Lemma \ref{pC.3},  and choose $\kappa$ small, $N_0$ large to deduce that for sufficiently small $\varepsilon$,
		\begin{align}\label{D.11}
			&\big\|\varepsilon^{\frac{3}{2}} \nabla_{x,p}h_R^{\varepsilon,\mathfrak{c}}(\bar{T})\big\|_{L^{\infty}} \leq \frac{1}{2}\left[\big\|\varepsilon^{\frac{3}{2}} h_R^{\varepsilon,\mathfrak{c}}(0)\big\|_{L^{\infty}}+C\big\|\varepsilon^{\frac{3}{2}} \nabla_{x,p}h_R^{\varepsilon,\mathfrak{c}}(0)\big\|_{L^{\infty}}\right]\nonumber\\
			&\quad\quad+C\varepsilon^{\frac{5}{2}} \sup_{s \in\left[0, \bar{T}\right]}\big\|\big(E_R^{\varepsilon,\mathfrak{c}}, B_R^{\varepsilon,\mathfrak{c}}\big)(s)  \big\|_{W^{1,\infty}}+C\sup_{s \in\left[0, \bar{T}\right]}\big\|f_R^{\varepsilon, \mathfrak{c}}(s)\big\|_{H^1}+C\varepsilon^{k+\frac{5}{2}}.
		\end{align}
		We thus conclude the proof of \eqref{D.3}, and the proof of \eqref{D.2} can be carried out similarly without using the factor $e^{\frac{\nu_0 s}{2 \varepsilon}}$. Therefore the proof of Lemma \ref{pD.2} is completed.
	\end{proof}

	With the help of Lemma \ref{pD.2}, we can easily obtain the following Lemma \ref{pD.1}. Since the proof is similar to Lemma \ref{pC.3}, we omit it.
	\begin{lemma}\label{pD.1}
		Assume the {\it a priori} assumptions \eqref{W.0}-\eqref{C.2}. For given $T>0$, taking $\varepsilon\ll 1$, we have that for $t \in [0,T]$,
		\begin{equation}\label{D.1}
			\begin{aligned}
				\big\|\varepsilon^{\frac{3}{2}} \nabla_{x,p}h_R^{\varepsilon, \mathfrak{c}}(t)\big\|_{L^{\infty}} \lesssim& \big\|\varepsilon^{\frac{3}{2}} h_R^{\varepsilon, \mathfrak{c}}(0)\big\|_{W^{1,\infty}}+\varepsilon^{\frac{5}{2}}\sup_{s \in\left[0, t\right]}\big\|\big(E_R^{\varepsilon,\mathfrak{c}}, B_R^{\varepsilon,\mathfrak{c}}\big)(s)\big\|_{W^{1,\infty}}\\
				&+ \sup_{s \in\left[0, t\right]}\big\|f_R^{\varepsilon, \mathfrak{c}}(s)\big\|_{H^1}+ \varepsilon^{k+\frac{5}{2}},
			\end{aligned}
		\end{equation}
		where the constant is independent of $\mathfrak{c}$.
	\end{lemma}
	
	\section{$W^{1, \infty}$-Estimates for $\left(E_R^{\varepsilon,\mathfrak{c}}, B_R^{\varepsilon,\mathfrak{c}}\right)$} \label{section 7}
	
	In this section, we use the well-known Glassey--Strauss representation in \cite{Glassey-ARMA-1986} to perform $W^{1, \infty}$-estimates for $\left(E_R^{\varepsilon,\mathfrak{c}}, B_R^{\varepsilon,\mathfrak{c}}\right)$.

	\subsection{$L^{\infty}$-estimate of $\left(E_R^{\varepsilon,\mathfrak{c}}, B_R^{\varepsilon,\mathfrak{c}}\right)$ } 
	
	It follows from \eqref{e1.10} that $\left(E_R^{\varepsilon, \mathfrak{c}}, B_R^{\varepsilon, \mathfrak{c}}, F_R^{\varepsilon, \mathfrak{c}}\right)$ satisfy the following wave equation:
	\begin{align}\label{W0.1}
		\left\{\begin{aligned}
			&\partial_t^2E_R^{\varepsilon, \mathfrak{c}}-\mathfrak{c}^{2}\Delta E_R^{\varepsilon, \mathfrak{c}}=4\pi\mathfrak{c}^{2}  \int_{\mathbb{R}^3}\nabla_x F_R^{\varepsilon, \mathfrak{c}}d p+4\pi\int_{\mathbb{R}^3}\hat{p}\partial_t F_R^{\varepsilon, \mathfrak{c}}d p,\\
			&\partial_t^2B_R^{\varepsilon, \mathfrak{c}}-\mathfrak{c}^{2}\Delta B_R^{\varepsilon, \mathfrak{c}}=-4\pi\mathfrak{c}\nabla\times\int_{\mathbb{R}^3}\hat{p}  F_R^{\varepsilon, \mathfrak{c}} d p,
		\end{aligned}\right.
	\end{align}
    with $t\geq 0, x\in \mathbb{R}^3$. The initial data are 
    $$
    \left(E_R^{\varepsilon, \mathfrak{c}}, \partial_t E_R^{\varepsilon, \mathfrak{c}},B_R^{\varepsilon, \mathfrak{c}},\partial_t B_R^{\varepsilon, \mathfrak{c}}\right)(0,x).
    $$
	Recall the definition of $\mathcal{I}_0$ in \eqref{e1.16-0}.
    \begin{align}\label{W0.2}
    \mathcal{I}_0 := \sup_{(x,p) \in \mathbb{R}^3 \times \mathbb{R}^3} \Bigg\{ (1+|x|) \Bigg( &\sum_{|\alpha|+|\beta|\leq 1} \big|\partial_x^\alpha \partial_p^\beta h_R^{\varepsilon, \mathfrak{c}}(0,x,p)\big| \nonumber\\
    &+ \sum_{1\leq |\alpha|\leq 2} \big( \big|\partial_x^\alpha E_R^{\varepsilon, \mathfrak{c}}(0,x)\big| + \big|\partial_x^\alpha B_R^{\varepsilon, \mathfrak{c}}(0,x)\big| \big) \Bigg) \Bigg\} < \infty.
    \end{align}

    \begin{remark}
    Assumption \eqref{W0.2} is needed to derive uniform-in-$\mathfrak{c}$ $W^{1,\infty}$ bounds: the factor $\mathfrak{c}t$ in the Glassey--Strauss initial-data terms is canceled by the $O(|x|^{-1})$ decay via $\int_{\mathbb{S}^2}|x+\mathfrak{c}t\omega|^{-1}d\omega \le 4\pi(\mathfrak{c}t)^{-1}$. 
    The first momentum derivative of $h_R^{\varepsilon,\mathfrak{c}}(0)$ is included because it is used in the differentiated characteristic estimates of Section \ref{section 6}. The second-order spatial derivatives in \eqref{W0.2} are required for estimating $\nabla_x(E_R^{\varepsilon,\mathfrak{c}},B_R^{\varepsilon,\mathfrak{c}})$.
    \end{remark}
    
	For later use, we define the transport operator
	$$
	S:=\partial_t+\hat{p} \cdot \nabla_x.
	$$
	By Kirchhoff's formula for the wave equation in $\mathbb{R}^3$ and a similar calculation to \cite{Glassey-ARMA-1986}, 
	 the solution of \eqref{W0.1} can be explicitly expressed in terms of $F_R^{\varepsilon, \mathfrak{c}}$ and $S F_R^{\varepsilon, \mathfrak{c}}$ as follows:
	\begin{lemma}\label{lW.2}
		For $x \in \mathbb{R}^3$, $i=1,2,3$, we have
		\begin{align}
			& E_{R,i}^{\varepsilon,\mathfrak{c}}(t,x)=\mathbf{I}_{E_{R,i}^{\varepsilon,\mathfrak{c}}}(t,x)+E_{R_T,i}^{\varepsilon,\mathfrak{c}}(t, x)+E_{R_S,i}^{\varepsilon,\mathfrak{c}}(t, x), \label{W.2-1}\\
			& B_{R,i}^{\varepsilon,\mathfrak{c}}(t,x)=\mathbf{I}_{B_{R,i}^{\varepsilon,\mathfrak{c}}}(t,x)+B_{R_T,i}^{\varepsilon,\mathfrak{c}}(t, x)+B_{R_S,i}^{\varepsilon,\mathfrak{c}}(t, x),\label{W.2-2}
		\end{align}
		where 
		\begin{align}
		 &\mathbf{I}_{E_{R,i}^{\varepsilon,\mathfrak{c}}}(t,x)=\frac{1}{\mathfrak{c}^2}\partial_t\Big(\frac{1}{t} \int_{|x-y|=\mathfrak{c} t}  E_{R,i}^{\varepsilon,\mathfrak{c}}(0, y) d S_y \Big)+\frac{1}{\mathfrak{c}^2 t} \int_{|x-y|=\mathfrak{c} t} \partial_t E_{R,i}^{\varepsilon,\mathfrak{c}}(0,y) d S_y \nonumber\\
			& \qquad \qquad\quad+ \iint_{|x-y|=\mathfrak{c} t} \frac{d S_y d p}{|x-y|} \frac{\omega_i-\frac{\hat{p_i}}{\mathfrak{c}} \frac{\hat{p}}{\mathfrak{c}} \cdot \omega}{1+\frac{\hat{p}}{\mathfrak{c}} \cdot \omega} F_R^{\varepsilon, \mathfrak{c}}(0, y, p), \nonumber\\
			 &E_{R_T,i}^{\varepsilon,\mathfrak{c}}(t, x)= \iint_{|x-y|<\mathfrak{c}t} \frac{d y d p}{|x-y|^2}\frac{\big(\omega_i+\frac{\hat{p}_i}{\mathfrak{c}}\big)\big(1-\frac{|\hat{p}|^2}{\mathfrak{c}^2} \big)}{\big(1+\frac{\hat{p}}{\mathfrak{c}} \cdot \omega\big)^2} F_R^{\varepsilon,\mathfrak{c}}\Big(t-\frac{|x-y|}{\mathfrak{c}}, y, p\Big), \nonumber\\
			 &E_{R_S,i}^{\varepsilon,\mathfrak{c}}(t, x)=\frac{1}{\mathfrak{c}^2} \iint_{|x-y|<\mathfrak{c}t}\frac{d y d p}{|x-y|} \frac{\mathfrak{c}\omega_i+\hat{p}_i}{1+\frac{\hat{p}}{\mathfrak{c}} \cdot \omega}\left(S F_R^{\varepsilon,\mathfrak{c}}\right)\Big(t-\frac{|x-y|}{\mathfrak{c}}, y, p\Big), \nonumber\\
			&\nonumber \\
			 &\mathbf{I}_{B_{R,i}^{\varepsilon,\mathfrak{c}}}(t,x)=\frac{1}{\mathfrak{c}^2}\partial_t\Big(\frac{1}{t} \int_{|x-y|=\mathfrak{c} t}  B_{R,i}^{\varepsilon,\mathfrak{c}}(0, y) d S_y \Big)+\frac{1}{\mathfrak{c}^2 t} \int_{|x-y|=\mathfrak{c} t} \partial_t B_{R,i}^{\varepsilon,\mathfrak{c}}(0,y) d S_y \nonumber\\
			&\qquad \qquad\quad - \iint_{|x-y|=\mathfrak{c} t} \frac{d S_y d p}{|x-y|} \frac{\big(\omega\times \frac{\hat{p}}{\mathfrak{c}}\big)_i}{1+\frac{\hat{p}}{\mathfrak{c}} \cdot \omega} F_R^{\varepsilon, \mathfrak{c}}(0, y, p), \nonumber\\
			 &B_{R_T,i}^{\varepsilon,\mathfrak{c}}(t, x)=- \iint_{|x-y|<\mathfrak{c} t} \frac{d y d p}{|x-y|^2} \frac{\big(\omega\times\frac{\hat{p}}{\mathfrak{c}}\big)_i\big(1-\frac{|\hat{p}|^2}{\mathfrak{c}^2}\big)}{\big(1+\frac{\hat{p}}{\mathfrak{c}} \cdot \omega\big)^2} F_R^{\varepsilon,\mathfrak{c}}\Big(t-\frac{|x-y|}{\mathfrak{c}}, y, p\Big), \nonumber\\
			&B_{R_S,i}^{\varepsilon,\mathfrak{c}}(t, x)=-\frac{1}{\mathfrak{c}^2} \iint_{|x-y|<\mathfrak{c} t}\frac{d y d p}{|x-y|} \frac{(\omega\times \hat{p})_i}{1+\frac{\hat{p}}{\mathfrak{c}} \cdot \omega}(S F_R^{\varepsilon,\mathfrak{c}})\Big(t-\frac{|x-y|}{\mathfrak{c}}, y, p\Big) \nonumber
		\end{align}
		with $\omega_i=\frac{y_i-x_i}{|x-y|}$.
	\end{lemma}
	\begin{remark}
	    We point out that Lemma \ref{lW.2} has been proved in \cite{Glassey-ARMA-1986} when $\mathfrak{c}=1$. For the general case, the
        proof is very similar to the one in \cite{Glassey-ARMA-1986} and we relegate the details to Appendix \ref{appendix C.1}.
	\end{remark}

    Throughout the rest of Section \ref{section 7}, when no confusion can arise, we suppress the superscript $(\varepsilon,\mathfrak{c})$ and write
    \[
    F_R=F_R^{\varepsilon,\mathfrak{c}},\quad E_R=E_R^{\varepsilon,\mathfrak{c}},\quad B_R=B_R^{\varepsilon,\mathfrak{c}},\quad h_R=h_R^{\varepsilon,\mathfrak{c}},\quad f_R=f_R^{\varepsilon,\mathfrak{c}}.
    \]
    We also set
    \begin{equation}\label{W.alpha}
        \alpha_{\mathfrak{c}}(p):=\min\left\{\frac{|p|}{\mathfrak{c}},1\right\},\qquad p^0=(\mathfrak{c}^2+|p|^2)^{\frac12},\qquad \frac{\hat p}{\mathfrak{c}}=\frac{p}{p^0}.
    \end{equation}
    Define the charge density and current density associated with the kinetic remainder by
    \begin{equation}\label{W.rhoj}
        \rho_R(t,x):=\int_{\mathbb R^3}F_R(t,x,p)\,dp,\qquad
        j_R(t,x):=\int_{\mathbb R^3}\hat p\,F_R(t,x,p)\,dp.
    \end{equation}
    Taking the divergence in the Maxwell equation for $E_R$ in \eqref{e1.10} and using the constraint equation for $\operatorname{div}E_R$, one obtains
    \begin{equation}\label{W.cont}
        \partial_t\rho_R+\nabla_x\cdot j_R=0.
    \end{equation}
    Moreover, by \eqref{e1.13}-\eqref{e1.14}, Lemma \ref{l2.1}, and the elementary bound $|\hat p|\le |p|$, for $0\le t\le T$,
    \begin{align}\label{W.rhoj.basic}
        &\|\rho_R(t)\|_{L_x^\infty}+\|j_R(t)\|_{L_x^\infty}\lesssim \|h_R(t)\|_{L_{x,p}^\infty},\nonumber\\
        &\|\nabla_x\rho_R(t)\|_{L_x^\infty}+\|\nabla_xj_R(t)\|_{L_x^\infty}\lesssim \|\nabla_xh_R(t)\|_{L_{x,p}^\infty},\nonumber\\
        &\|\rho_R(t)\|_{L_x^2}\lesssim \|f_R(t)\|_{L_{x,p}^2}.
    \end{align}
    All constants in \eqref{W.rhoj.basic} are independent of $\mathfrak{c}$ and $\varepsilon$.

    We shall also use the following elementary bounds for the Glassey--Strauss kernels. For every fixed integer $m\geq0$, $\omega\in\mathbb S^2$, and $i=1,2,3$,
    \begin{gather}\label{W.5}
        \mathfrak{c}^2-|\hat{p}|^2=\frac{\mathfrak{c}^4}{\mathfrak{c}^2+|p|^2}\leq \mathfrak{c}^2,\qquad
        \left|\mathfrak{c}\omega_i+\hat{p}_i\right|\lesssim \mathfrak{c},\nonumber\\
        \Big(1+\frac{\hat{p}}{\mathfrak{c}} \cdot \omega\Big)^{-m}
        \lesssim (1+|p|)^{2m},\qquad
        \bigg|\nabla_p\bigg(\frac{\mathfrak{c}\omega_i+\hat{p}_i}{1+\frac{\hat{p}}{\mathfrak{c}} \cdot \omega}\bigg)\bigg|\lesssim (1+|p|)^4 .
    \end{gather}

    \begin{lemma}\label{lW.3}
    For each fixed integer $m\ge 8$, there exists a constant $C_m>0$, independent of $\mathfrak{c}$, such that
    \begin{equation}\label{W.kernel.1}
        \int_{\mathbb R^3}\alpha_{\mathfrak{c}}(p)\langle p\rangle^m\frac{\sqrt{J_{\mathbf M}}(p)}{w_\ell(p)}\,dp\le \frac{C_m}{\mathfrak{c}}.
    \end{equation}
    Let
    \begin{align*}
        K_i^E(\omega,p)&:=\frac{\big(\omega_i+\frac{\hat p_i}{\mathfrak{c}}\big)\big(1-\frac{|\hat p|^2}{\mathfrak{c}^2}\big)}{\big(1+\frac{\hat p}{\mathfrak{c}}\cdot\omega\big)^2},\\
        K_i^B(\omega,p)&:=-\frac{\big(\omega\times\frac{\hat p}{\mathfrak{c}}\big)_i\big(1-\frac{|\hat p|^2}{\mathfrak{c}^2}\big)}{\big(1+\frac{\hat p}{\mathfrak{c}}\cdot\omega\big)^2}
    \end{align*}
    be the kernels in the $T$-parts of the Glassey--Strauss representation. Then
    \begin{equation}\label{W.kernel.2}
        |K_i^E(\omega,p)-\omega_i|+|K_i^B(\omega,p)|\le C_m\alpha_{\mathfrak{c}}(p)\langle p\rangle^m.
    \end{equation}
    \end{lemma}
    \begin{proof}
    We first prove \eqref{W.kernel.1}. On $\{|p|\le \mathfrak{c}\}$, $\alpha_{\mathfrak{c}}(p)=|p|/\mathfrak{c}$, and hence
    \[
        \int_{|p|\le\mathfrak{c}}\alpha_{\mathfrak{c}}(p)\langle p\rangle^m\frac{\sqrt{J_{\mathbf M}}(p)}{w_\ell(p)}\,dp
        \le \frac1{\mathfrak{c}}\int_{\mathbb R^3}|p|\langle p\rangle^m\frac{\sqrt{J_{\mathbf M}}(p)}{w_\ell(p)}\,dp\lesssim\frac1{\mathfrak{c}}.
    \]
    On $\{|p|>\mathfrak{c}\}$, the uniform exponential tail following from \eqref{e1.15} and Lemma \ref{l2.1} gives, for any $N>0$,
    \[
        \int_{|p|>\mathfrak{c}}\langle p\rangle^m\frac{\sqrt{J_{\mathbf M}}(p)}{w_\ell(p)}\,dp\lesssim \mathfrak{c}^{-N}.
    \]
    Taking $N=1$ proves \eqref{W.kernel.1}. To prove \eqref{W.kernel.2}, observe that when $|p|\le\mathfrak{c}$,
    \[
        \left|\frac{\hat p}{\mathfrak{c}}\right|=\frac{|p|}{p^0}\le \frac1{\sqrt2}.
    \]
    Thus the kernels are smooth functions of $\hat p/\mathfrak{c}$ in a fixed compact subset of the unit ball. Since $K_i^E(\omega,0)=\omega_i$ and $K_i^B(\omega,0)=0$, the mean-value theorem gives
    \[
        |K_i^E(\omega,p)-\omega_i|+|K_i^B(\omega,p)|\lesssim \frac{|p|}{\mathfrak{c}}=\alpha_{\mathfrak{c}}(p).
    \]
    If $|p|>\mathfrak{c}$, then $\alpha_{\mathfrak{c}}(p)=1$, and the polynomial kernel bound already used in \eqref{W.5} yields \eqref{W.kernel.2} after increasing $m$ if necessary.
    \end{proof}

    \begin{lemma}\label{lW.4}
    For $R>0$, define
    \begin{align}\label{W.PH.def}
        \mathcal P_R[\rho]_i(x)&:=\int_{|z|<R}\frac{z_i}{|z|^3}\rho(x+z)\,dz,\nonumber\\
        \mathcal H_R[\rho]_{ij}(x)&:=\operatorname{p.v.}\int_{|z|<R}\frac{3\omega_i\omega_j-\delta_{ij}}{|z|^3}\rho(x+z)\,dz,\qquad \omega=\frac{z}{|z|}.
    \end{align}
    Then
    \begin{equation}\label{W.poisson.1}
        \sup_{R>0}\|\mathcal P_R[\rho]\|_{L_x^\infty}\lesssim \|\rho\|_{L_x^\infty}^{\frac13}\|\rho\|_{L_x^2}^{\frac23},
    \end{equation}
    and
    \begin{equation}\label{W.poisson.2}
        \sup_{R>0}\|\mathcal H_R[\rho]\|_{L_x^\infty}\lesssim \|\nabla_x\rho\|_{L_x^\infty}^{\frac35}\|\rho\|_{L_x^2}^{\frac25}.
    \end{equation}
    \end{lemma}
    \begin{proof}
    We give the proof since the uniformity in the truncation radius will be used below. Let $\eta>0$. For \eqref{W.poisson.1}, the near-field part is bounded by
    \[
        \int_{|z|<\eta}\frac{|\rho(x+z)|}{|z|^2}\,dz\lesssim \eta\|\rho\|_{L_x^\infty},
    \]
    whereas the far-field part is estimated by Cauchy's inequality as
    \[
        \int_{\eta<|z|<R}\frac{|\rho(x+z)|}{|z|^2}\,dz
        \le \left(\int_{|z|>\eta}|z|^{-4}\,dz\right)^{\frac12}\|\rho\|_{L_x^2}
        \lesssim \eta^{-\frac12}\|\rho\|_{L_x^2}.
    \]
    Optimizing in $\eta$ gives \eqref{W.poisson.1}. For \eqref{W.poisson.2}, we use
    \[
        \int_{\mathbb S^2}(3\omega_i\omega_j-\delta_{ij})\,d\omega=0.
    \]
    Hence, in the near field, the principal value can be written with the difference $\rho(x+z)-\rho(x)$, and
    \[
        \int_{|z|<\eta}\frac{|\rho(x+z)-\rho(x)|}{|z|^3}\,dz\lesssim \eta\|\nabla_x\rho\|_{L_x^\infty}.
    \]
    For the far field,
    \[
        \int_{\eta<|z|<R}\frac{|\rho(x+z)|}{|z|^3}\,dz
        \le \left(\int_{|z|>\eta}|z|^{-6}\,dz\right)^{\frac12}\|\rho\|_{L_x^2}
        \lesssim \eta^{-\frac32}\|\rho\|_{L_x^2}.
    \]
    Optimizing in $\eta$ proves \eqref{W.poisson.2}.
    \end{proof}

    \begin{lemma}\label{lW.5}
    For $t\in[0,T]$, define the retarded Coulomb and Hessian-Coulomb principal parts by
    \begin{align}\label{W.ret.def}
        \mathscr C_i(t,x)&:=\int_{|z|<\mathfrak{c} t}\frac{z_i}{|z|^3}\rho_R\left(t-\frac{|z|}{\mathfrak{c}},x+z\right)\,dz,\nonumber\\
        \mathscr H_{ij}(t,x)&:=\operatorname{p.v.}\int_{|z|<\mathfrak{c} t}\frac{3\omega_i\omega_j-\delta_{ij}}{|z|^3}\rho_R\left(t-\frac{|z|}{\mathfrak{c}},x+z\right)\,dz,
    \end{align}
    where $\omega=z/|z|$. Then
    \begin{align}\label{W.ret.1}
        \|\mathscr C(t)-\mathcal P_{\mathfrak{c} t}[\rho_R(t)]\|_{L_x^\infty}&\le C_T\sup_{0\le\tau\le t}\|j_R(\tau)\|_{L_x^\infty},\\
        \|\mathscr H(t)-\mathcal H_{\mathfrak{c} t}[\rho_R(t)]\|_{L_x^\infty}&\le C_T\sup_{0\le\tau\le t}\|\nabla_xj_R(\tau)\|_{L_x^\infty}.\label{W.ret.2}
    \end{align}
    The constant $C_T$ is independent of $\mathfrak{c}$ and $\varepsilon$.
    \end{lemma}
    \begin{proof}
    The comparison is made with the instantaneous operator truncated at the same radius $\mathfrak{c} t$; no tail outside this light cone is introduced. By \eqref{W.cont}, for $r=|z|$,
    \[
        \rho_R\left(t-\frac r{\mathfrak{c}},x+z\right)-\rho_R(t,x+z)
        =\int_0^{r/\mathfrak{c}}\nabla_x\cdot j_R(t-s,x+z)\,ds.
    \]
    Therefore, after exchanging the order of integration,
    \begin{equation}\label{W.ret.proof1}
        \mathscr C_i(t,x)-\mathcal P_{\mathfrak{c} t}[\rho_R(t)]_i(x)
        =\int_0^t\int_{\mathfrak{c} s<|z|<\mathfrak{c} t}\frac{z_i}{|z|^3}\nabla_x\cdot j_R(t-s,x+z)\,dz\,ds.
    \end{equation}
    On the annulus $A_s:=\{z:\mathfrak{c} s<|z|<\mathfrak{c} t\}$, use $\nabla_x\cdot j_R(t-s,x+z)=\nabla_z\cdot j_R(t-s,x+z)$ and integrate by parts. The boundary contribution is bounded by $C\|j_R(t-s)\|_{L_x^\infty}$, while
    \[
        \left|\int_{A_s}\nabla_z\left(\frac{z_i}{|z|^3}\right)\cdot j_R(t-s,x+z)\,dz\right|
        \lesssim \left(1+\log\frac{t}{s}\right)\|j_R(t-s)\|_{L_x^\infty}.
    \]
    Since $\int_0^t\log(t/s)\,ds=t$, \eqref{W.ret.1} follows from \eqref{W.ret.proof1}. For the Hessian-Coulomb part, first work with the principal value truncated to $\delta<|z|<\mathfrak{c} t$ and then let $\delta\downarrow0$. Again using \eqref{W.cont},
    \[
        \left|\rho_R\left(t-\frac{|z|}{\mathfrak{c}},x+z\right)-\rho_R(t,x+z)\right|
        \le \frac{|z|}{\mathfrak{c}}\sup_{0\le\tau\le t}\|\nabla_xj_R(\tau)\|_{L_x^\infty}.
    \]
    Hence
    \[
        \int_{|z|<\mathfrak{c} t}\frac1{|z|^3}\frac{|z|}{\mathfrak{c}}\,dz
        \lesssim \frac1{\mathfrak{c}}\int_0^{\mathfrak{c} t}dr=t,
    \]
    which proves \eqref{W.ret.2}.
    \end{proof}

		\begin{proposition}\label{pW.5}
		Assume the {\it a priori} assumptions \eqref{W.0}-\eqref{C.2}. For every $T>0$, there exist constants $C_T>0$ and $\varepsilon_0=\varepsilon_0(T)>0$, independent of $\mathfrak{c}$ and $\varepsilon$, such that, for $0<\varepsilon\leq\varepsilon_0$, $\mathfrak{c}\geq\mathfrak{c}_0$, and $t\in[0,T]$, one has
		\begin{align}\label{W.4}
			\sup_{0\leq s\leq t}\varepsilon^{\frac{5}{2}}&\big\|\big(E_R^{\varepsilon,\mathfrak{c}}, B_R^{\varepsilon,\mathfrak{c}}\big)(s)\big\|_{L^{\infty}} \nonumber\\
			&\leq C_T\Big[ \varepsilon^{\frac{5}{2}}\big\|\big(E_R^{\varepsilon,\mathfrak{c}}, B_R^{\varepsilon,\mathfrak{c}}\big)(0)\big\|_{L^{\infty}}+\varepsilon^{\frac{3}{2}}\mathcal{I}_0+\sup_{0\leq s\leq t}\big\| f_R^{\varepsilon,\mathfrak{c}}(s)\big\|+\varepsilon^{k+\frac{5}{2}}\Big].
		\end{align}
		The constant $C_T$ is independent of $\mathfrak{c}$ and $\varepsilon$.
	\end{proposition}
	\begin{proof}
		We use the shortened notation introduced before Lemma \ref{lW.3}. The elementary bounds \eqref{W.5} will be used only for the source terms; the transport terms are treated below by exploiting the cancellations isolated in Lemmas \ref{lW.3}-\ref{lW.5}.

        We divide the proof into several steps.
        
        {\it Step 1.1. }
        Performing the change of variables $y=x+\mathfrak{c} t\omega$ in the boundary terms of Lemma \ref{lW.2}, using the Maxwell equations at $t=0$, and applying the decay assumption \eqref{W0.2}, one obtains
		\begin{align}\label{W.6-1}
		\big\|\big(\mathbf I_{E_R}^{\varepsilon,\mathfrak{c}},\mathbf I_{B_R}^{\varepsilon,\mathfrak{c}}\big)(t)\big\|_{L^\infty}
		\le C\big\|\big(E_R^{\varepsilon,\mathfrak{c}},B_R^{\varepsilon,\mathfrak{c}}\big)(0)\big\|_{L^\infty}+C\mathcal I_0 .
		\end{align}
		Indeed, the only point which is not immediate is the factor $\mathfrak{c} t$ produced by differentiating the spherical means. It is canceled by
		\[
		\mathfrak{c} t\int_{\mathbb S^2}\frac{d\omega}{1+|x+\mathfrak{c} t\omega|}
		\leq \mathfrak{c} t\int_{\mathbb S^2}\frac{d\omega}{|x+\mathfrak{c} t\omega|}
		=\frac{4\pi\mathfrak{c} t}{\max\{|x|,\mathfrak{c} t\}}\leq 4\pi,
		\]
		and the momentum integrations are finite uniformly in $\mathfrak{c}$ by Lemma \ref{l2.1}.

        {\it Step 1.2. }
        We decompose the kernel in $E_{R_T}$ as
        \[
        K_i^E(\omega,p,\mathfrak{c})
        :=\frac{\big(\omega_i+\frac{\hat p_i}{\mathfrak{c}}\big)\big(1-\frac{|\hat p|^2}{\mathfrak{c}^2}\big)}{\big(1+\frac{\hat p}{\mathfrak{c}}\cdot\omega\big)^2}
        =\omega_i+\widetilde K_i^E(\omega,p,\mathfrak{c}).
        \]
        Hence, with $z=y-x$,
        \[
        E_{R_T,i}(t,x)=\mathscr C_i(t,x)+\widetilde E_{R_T,i}(t,x),
        \]
        where $\mathscr C_i$ is defined in \eqref{W.ret.def}. By Lemma \ref{lW.3},
        \begin{align}\label{W.6-2}
        |\widetilde E_{R_T,i}(t,x)|
        &\lesssim \int_0^{\mathfrak{c} t}\int_{\mathbb S^2}\int_{\mathbb R^3}
        \alpha_{\mathfrak{c}}(p)\langle p\rangle^m\frac{\sqrt{J_{\mathbf M}}}{w_\ell}
        |h_R(t-r/\mathfrak{c},x+r\omega,p)|\,dp\,d\omega\,dr\nonumber\\
        &\lesssim \int_0^{\mathfrak{c} t}\frac{dr}{\mathfrak{c}}\sup_{0\leq\tau\leq t}\|h_R(\tau)\|_{L^\infty}
        \leq C_T\sup_{0\leq\tau\leq t}\|h_R(\tau)\|_{L^\infty}.
        \end{align}
        For the principal Coulomb part, Lemmas \ref{lW.4} and \ref{lW.5}, together with \eqref{W.rhoj.basic}, give
        \begin{align*}
        \|\mathscr C(t)\|_{L_x^\infty}
        &\leq \|\mathscr C(t)-\mathcal P_{\mathfrak{c} t}[\rho_R(t)]\|_{L_x^\infty}
        +\|\mathcal P_{\mathfrak{c} t}[\rho_R(t)]\|_{L_x^\infty}\\
        &\leq C_T\sup_{0\leq\tau\leq t}\|j_R(\tau)\|_{L_x^\infty}
        +C\|\rho_R(t)\|_{L_x^\infty}^{\frac13}\|\rho_R(t)\|_{L_x^2}^{\frac23}\\
        &\leq C_T\sup_{0\leq\tau\leq t}\|h_R(\tau)\|_{L^\infty}
        +C\|h_R(t)\|_{L^\infty}^{\frac13}\|f_R(t)\|^{\frac23}.
        \end{align*}
        Combining this estimate with \eqref{W.6-2}, we obtain the improved bound
		\begin{align}\label{W.6-3}
		\|E_{R_T}^{\varepsilon,\mathfrak{c}}(t)\|_{L^\infty}
		\leq C_T\sup_{0\leq\tau\leq t}\|h_R^{\varepsilon,\mathfrak{c}}(\tau)\|_{L^\infty}
		+C\|h_R^{\varepsilon,\mathfrak{c}}(t)\|_{L^\infty}^{\frac13}\|f_R^{\varepsilon,\mathfrak{c}}(t)\|^{\frac23}.
		\end{align}
		This estimate cancels the factor $\mathfrak{c}$ coming from the size of the light cone.
		
        {\it Step 1.3. }
        We retain the notation
		\begin{align}\label{W.7}
			E_{R_S,i}^{\varepsilon,\mathfrak{c}}(t, x)=&\frac{1}{\mathfrak{c}^2} \iint_{|x-y|<\mathfrak{c}t}\frac{d y d p}{|x-y|} \frac{\mathfrak{c}\omega_i+\hat{p}_i}{1+\frac{\hat{p}}{\mathfrak{c}} \cdot \omega}\left(S F_R^{\varepsilon,\mathfrak{c}}\right)\Big(t-\frac{|x-y|}{\mathfrak{c}}, y, p\Big).
		\end{align}
		Substituting $S F_R^{\varepsilon,\mathfrak{c}}$ by the first equation of \eqref{e1.10}, using \eqref{W.5}, and integrating by parts in $p$ for the force terms, we also use the pointwise collision bound
		\begin{align}\label{W.18}
		\Big|\frac{1}{\sqrt{\mathbf{M}_{\mathfrak{c}}}} Q_{\mathfrak{c}}\left(F_R^{\varepsilon, \mathfrak{c}}, F^{\varepsilon, \mathfrak{c}}\right)\Big|
		\lesssim \nu_{\mathfrak{c}}(p)\Big\|\frac{F_R^{\varepsilon, \mathfrak{c}}}{\sqrt{\mathbf{M}_{\mathfrak{c}}}}\Big\|_{L^{\infty}}\Big\|\frac{F^{\varepsilon, \mathfrak{c}}}{\sqrt{\mathbf{M}_{\mathfrak{c}}}}\Big\|_{L^{\infty}}
		\lesssim \nu_{\mathfrak{c}}(p)\big\|h_R^{\varepsilon, \mathfrak{c}}\big\|_{L^{\infty}},
		\end{align}
		which follows from Lemma \ref{l2.2} and \eqref{W.0}; the same bound holds with the two arguments of $Q_{\mathfrak{c}}$ interchanged. The possible leading factor in the collision part is removed before estimating absolute values. Indeed,
        \[
        \frac{\mathfrak{c}\omega_i+\hat p_i}{1+\frac{\hat p}{\mathfrak{c}}\cdot\omega}
        =\mathfrak{c}\omega_i+
        \frac{\hat p_i-\omega_i(\hat p\cdot\omega)}{1+\frac{\hat p}{\mathfrak{c}}\cdot\omega},
        \]
        and the contribution of $\mathfrak{c}\omega_i$ vanishes after integration in $p$ by the conservation law
        \[
        \int_{\mathbb R^3}
        \{Q_{\mathfrak{c}}(g_1,g_2)+Q_{\mathfrak{c}}(g_2,g_1)\}(p)\,dp=0.
        \]
        The remaining factor is polynomially bounded by \eqref{W.5}, and therefore \eqref{W.18} gives the displayed $\varepsilon^{-1}$ contribution. These estimates give
		\begin{align}\label{W.6-4}
            \|E_{R_S,i}^{\varepsilon,\mathfrak{c}}(t)\|_{L^\infty}
            &\leq C_T\int_0^t(t-\tau)\|\big(E_R^{\varepsilon,\mathfrak{c}},B_R^{\varepsilon,\mathfrak{c}}\big)(\tau)\|_{L^\infty}\,d\tau\nonumber\\
            &\quad +\frac{C_T}{\varepsilon}\int_0^t(t-\tau)\|h_R^{\varepsilon,\mathfrak{c}}(\tau)\|_{L^\infty}\,d\tau+C_T\varepsilon^k.
        \end{align}
        We recall the origin of the two terms in \eqref{W.6-4}. The factor $\mathfrak{c}^{-2}|x-y|^{-1}(\mathfrak{c}\omega_i+\hat p_i)$ gives, after the change $y=x+\mathfrak{c}(t-\tau)\omega$, the harmless weight $t-\tau$. The linear collision and transport remainders in \eqref{e1.10} produce the $\varepsilon^{-1}\|h_R\|_{L^\infty}$ contribution, while the force terms involving $\nabla_pF_R$ are moved onto the kernel and the known coefficients; the derivative bound in \eqref{W.5}, the coefficient estimates from Section \ref{section 4}, and the bootstrap assumption \eqref{W.0} then give the displayed right-hand side. The purely known source terms are bounded by $C_T\varepsilon^k$.

        {\it Step 2. }
        The initial part satisfies the analogue of \eqref{W.6-1}. For the transport part, the magnetic kernel vanishes at $p=0$; hence Lemma \ref{lW.3} directly yields
        \begin{align}\label{W.8-1}
        \|B_{R_T}^{\varepsilon,\mathfrak{c}}(t)\|_{L^\infty}
        \leq C_T\sup_{0\leq\tau\leq t}\|h_R^{\varepsilon,\mathfrak{c}}(\tau)\|_{L^\infty}.
        \end{align}
        The source term $B_{R_S}$ is estimated in the same way as \eqref{W.6-4}. Consequently, from Lemma \ref{lW.2}, \eqref{W.6-1}, \eqref{W.6-3}, \eqref{W.6-4}, and \eqref{W.8-1},
		\begin{align}\label{W.10}
			\|\big(E_R^{\varepsilon,\mathfrak{c}},B_R^{\varepsilon,\mathfrak{c}}\big)(t)\|_{L^\infty}
			&\leq C_T\big\|\big(E_R^{\varepsilon,\mathfrak{c}},B_R^{\varepsilon,\mathfrak{c}}\big)(0)\big\|_{L^\infty}+C_T\mathcal I_0\nonumber\\
			&\quad+C_T\sup_{0\leq\tau\leq t}\|h_R^{\varepsilon,\mathfrak{c}}(\tau)\|_{L^\infty}
			+C\|h_R^{\varepsilon,\mathfrak{c}}(t)\|_{L^\infty}^{\frac13}\|f_R^{\varepsilon,\mathfrak{c}}(t)\|^{\frac23}\nonumber\\
			&\quad+C_T\int_0^t(t-\tau)\|\big(E_R^{\varepsilon,\mathfrak{c}},B_R^{\varepsilon,\mathfrak{c}}\big)(\tau)\|_{L^\infty}\,d\tau\nonumber\\
			&\quad+\frac{C_T}{\varepsilon}\int_0^t(t-\tau)\|h_R^{\varepsilon,\mathfrak{c}}(\tau)\|_{L^\infty}\,d\tau+C_T\varepsilon^k .
		\end{align}

        It remains to close this inequality with Lemma \ref{pC.3}. Set
        \[
        \Phi_0(t):=\sup_{0\leq s\leq t}\varepsilon^{\frac52}\|\big(E_R^{\varepsilon,\mathfrak{c}},B_R^{\varepsilon,\mathfrak{c}}\big)(s)\|_{L^\infty},\quad
        H_0(t):=\sup_{0\leq s\leq t}\varepsilon^{\frac32}\|h_R^{\varepsilon,\mathfrak{c}}(s)\|_{L^\infty},
        \]
        and
        \[
        F_0(t):=\sup_{0\leq s\leq t}\|f_R^{\varepsilon,\mathfrak{c}}(s)\|.
        \]
        Multiplying \eqref{W.10} by $\varepsilon^{5/2}$ and taking the supremum over $0\leq s\leq t$ gives
        \begin{align}\label{W.10a}
        \Phi_0(t)&\leq C_T\varepsilon^{\frac52}\|\big(E_R^{\varepsilon,\mathfrak{c}},B_R^{\varepsilon,\mathfrak{c}}\big)(0)\|_{L^\infty}+C_T\varepsilon^{\frac52}\mathcal I_0
        +C_T\varepsilon H_0(t)\nonumber\\
        &\quad+C_T\varepsilon^2 H_0(t)^{\frac13}F_0(t)^{\frac23}
        +C_T\int_0^t(t-\tau)\Phi_0(\tau)\,d\tau
        +C_T\int_0^t(t-\tau)H_0(\tau)\,d\tau+C_T\varepsilon^{k+\frac52}.
        \end{align}
        By Young's inequality,
        \[
        \varepsilon^2H_0(t)^{\frac13}F_0(t)^{\frac23}\leq C\varepsilon^2H_0(t)+CF_0(t).
        \]
        Moreover, Lemma \ref{pC.3} implies
        \begin{align}\label{W.10b}
        H_0(t)\leq C\varepsilon^{\frac32}\|h_R^{\varepsilon,\mathfrak{c}}(0)\|_{L^\infty}+C\Phi_0(t)+CF_0(t)+C\varepsilon^{k+\frac52}.
        \end{align}
        Since $\|h_R^{\varepsilon,\mathfrak{c}}(0)\|_{L^\infty}\leq\mathcal I_0$ and $\varepsilon\leq1$, \eqref{W.10a}-\eqref{W.10b} yield
        \begin{align}\label{W.10c}
        \Phi_0(t)&\leq C_T\Big[\varepsilon^{\frac52}\|\big(E_R^{\varepsilon,\mathfrak{c}},B_R^{\varepsilon,\mathfrak{c}}\big)(0)\|_{L^\infty}+\varepsilon^{\frac32}\mathcal I_0+F_0(t)+\varepsilon^{k+\frac52}\Big]\nonumber\\
        &\quad+C_T\varepsilon\Phi_0(t)+C_T\int_0^t(t-\tau)\Phi_0(\tau)\,d\tau.
        \end{align}
        Choosing $\varepsilon_0(T)>0$ so small that $C_T\varepsilon_0\leq1/2$, we absorb the term $C_T\varepsilon\Phi_0(t)$ into the left-hand side. The Volterra form of Gronwall's inequality then gives \eqref{W.4}. This completes the proof of Proposition \ref{pW.5}.
	\end{proof}
	\subsection[$L^{\infty}$-estimate of $\nabla_x(E_R,B_R)$]{$L^{\infty}$-estimate of $\nabla_x(E_R^{\varepsilon,\mathfrak{c}},B_R^{\varepsilon,\mathfrak{c}})$}
	The following representation gives the spatial gradient of the field remainder explicitly:
	\begin{lemma}\label{lW.6}
		Denote $\partial_j:=\partial_{x_j}$. For any indices $i, j=1,2,3$, we have
		\begin{align}
			\partial_j  E_{R,i}^{\varepsilon,\mathfrak{c}}(t, x)= & \mathbf{I}_{E_{R,i j}^{\varepsilon,\mathfrak{c}}}(t,x)+\iint_{|x-y|<\mathfrak{c}t} \frac{a_A(\omega, \hat{p})}{|x-y|^3} F_R^{\varepsilon,\mathfrak{c}} \Big(t-\frac{|x-y|}{\mathfrak{c}}, y, p\Big)d y d p\nonumber\\
			&+\frac{1}{\mathfrak{c}} \iint_{|x-y|<\mathfrak{c}t} \frac{b_A(\omega, \hat{p})}{|x-y|^2} (S F_R^{\varepsilon,\mathfrak{c}})\Big(t-\frac{|x-y|}{\mathfrak{c}}, y, p\Big)d y d p \nonumber\\
			&+\frac{1}{\mathfrak{c}^2} \iint_{|x-y|<\mathfrak{c}t} \frac{c_A(\omega, \hat{p})}{|x-y|} (S^2 F_R^{\varepsilon,\mathfrak{c}})\Big(t-\frac{|x-y|}{\mathfrak{c}}, y, p\Big)d y d p,\label{W.11}
		\end{align}
		where 
		\begin{align}
			\mathbf{I}_{E_{R,i j}^{\varepsilon,\mathfrak{c}}}(t,x)=&\partial_j \left(\mathbf{I}_{E_{R,i}^{\varepsilon,\mathfrak{c}}}(t,x)\right)+t \iint_{|\omega|=1} \frac{\omega_j\big(\omega_i+\frac{\hat{p}_i}{\mathfrak{c}}\big)}{\big(1+\frac{\hat{p}}{\mathfrak{c}} \cdot \omega\big)^2} (S F_R^{\varepsilon,\mathfrak{c}})(0, x+\mathfrak{c} t \omega, p) d \omega d p\nonumber\\
			&+\frac{1}{\mathfrak{c}^2} \iint_{|\omega|=1} \frac{\omega_j\big(\omega_i+\frac{\hat{p}_i}{\mathfrak{c}}\big)\left(\mathfrak{c}^2-|\hat{p}|^2\right)}{\big(1+\frac{\hat{p}}{\mathfrak{c}} \cdot \omega\big)^3} F_R^{\varepsilon,\mathfrak{c}}(0, x+\mathfrak{c} t \omega, p) d \omega d p, \label{W.11-1}\\
        a_A(\omega, \hat{p}) =& \frac{3(\omega_i+\frac{\hat{p_i}}{\mathfrak{c}})\left[\omega_j(1-\frac{|\hat{p}|^2}{\mathfrak{c}^2})+\frac{\hat{p}_j}{\mathfrak{c}}(1+\frac{\hat{p}}{\mathfrak{c}} \cdot \omega)\right]-(1+\frac{\hat{p}}{\mathfrak{c}} \cdot \omega)^2 \delta_{i j}}{(1+\frac{p^2}{\mathfrak{c}^2})(1+\frac{\hat{p}}{\mathfrak{c}} \cdot \omega)^4},\label{W.11-3}\\
			b_A(\omega, \hat{p})=&-\frac{(\omega_i+\frac{\hat{p_i}}{\mathfrak{c}})\left[3\omega_j(\frac{|\hat{p}|^2}{\mathfrak{c}^2}-1)-2\frac{\hat{p}_j}{\mathfrak{c}}(1+\frac{\hat{p}}{\mathfrak{c}} \cdot \omega)\right]+(1+\frac{\hat{p}}{\mathfrak{c}} \cdot \omega)^2 \delta_{i j}}{(1+\frac{\hat{p}}{\mathfrak{c}} \cdot \omega)^3},\label{W.11-4}\\
			c_A(\omega, \hat{p})=&\frac{(\omega_i+\frac{\hat{p}_i}{\mathfrak{c}})\omega_j}{(1+\frac{\hat{p}}{\mathfrak{c}} \cdot \omega)^2},\label{W.11-5}
		\end{align}
		and
		\begin{align}
		    \partial_j  B_{R,i}^{\varepsilon,\mathfrak{c}}(t, x)= & \mathbf{I}_{B_{R,i j}^{\varepsilon,\mathfrak{c}}}(t,x)+\iint_{|x-y|<\mathfrak{c}t} \frac{a_B(\omega, \hat{p})}{|x-y|^3} F_R^{\varepsilon,\mathfrak{c}} \Big(t-\frac{|x-y|}{\mathfrak{c}}, y, p\Big)d y d p\nonumber\\
			&+\frac{1}{\mathfrak{c}} \iint_{|x-y|<\mathfrak{c}t} \frac{b_B(\omega, \hat{p})}{|x-y|^2} (S F_R^{\varepsilon,\mathfrak{c}})\Big(t-\frac{|x-y|}{\mathfrak{c}}, y, p\Big)d y d p \nonumber\\
			&+\frac{1}{\mathfrak{c}^2} \iint_{|x-y|<\mathfrak{c}t} \frac{c_B(\omega, \hat{p})}{|x-y|} (S^2 F_R^{\varepsilon,\mathfrak{c}})\Big(t-\frac{|x-y|}{\mathfrak{c}}, y, p\Big)d y d p,\label{W.120}
		\end{align}
		where
		\begin{align}
		    \mathbf{I}_{B_{R,i j}^{\varepsilon,\mathfrak{c}}}(t,x)=&\partial_j \left(\mathbf{I}_{B_{R,i}^{\varepsilon,\mathfrak{c}}}(t,x)\right)-t \iint_{|\omega|=1} \frac{\omega_j(\omega\times\frac{\hat{p}}{\mathfrak{c}})_i}{1+\frac{\hat{p}}{\mathfrak{c}} \cdot \omega} (S F_R^{\varepsilon,\mathfrak{c}})(0, x+\mathfrak{c} t \omega, p) d \omega d p\nonumber\\
			&-\frac{1}{\mathfrak{c}^2} \iint_{|\omega|=1} \frac{\omega_j(\omega\times\frac{\hat{p}}{\mathfrak{c}})_i\left(\mathfrak{c}^2-|\hat{p}|^2\right)}{\big(1+\frac{\hat{p}}{\mathfrak{c}} \cdot \omega\big)^2} F_R^{\varepsilon,\mathfrak{c}}(0, x+\mathfrak{c} t \omega, p) d \omega d p, \label{W.12-1}\\
			a_B(\omega, \hat{p}) =& \frac{-3(\omega\times\frac{\hat{p}}{\mathfrak{c}})_i\left[\omega_j(1-\frac{|\hat{p}|^2}{\mathfrak{c}^2})+\frac{\hat{p}_j}{\mathfrak{c}}(1+\frac{\hat{p}}{\mathfrak{c}} \cdot \omega)\right]+(1+\frac{\hat{p}}{\mathfrak{c}} \cdot \omega)^2 (\delta_{ j}\times\frac{\hat{p}}{\mathfrak{c}})_i}{(1+\frac{p^2}{\mathfrak{c}^2})(1+\frac{\hat{p}}{\mathfrak{c}} \cdot \omega)^4} ,\label{W.12-3}\\
			b_B(\omega, \hat{p})=&\frac{(\omega\times\frac{\hat{p}}{\mathfrak{c}})_i\left[3\omega_j(\frac{|\hat{p}|^2}{\mathfrak{c}^2}-1)-2\frac{\hat{p}_j}{\mathfrak{c}}(1+\frac{\hat{p}}{\mathfrak{c}} \cdot \omega)\right]+(1+\frac{\hat{p}}{\mathfrak{c}} \cdot \omega)^2 (\delta_{ j}\times\frac{\hat{p}}{\mathfrak{c}})_i}{(1+\frac{\hat{p}}{\mathfrak{c}} \cdot \omega)^3},\label{W.12-4}\\
			c_B(\omega, \hat{p})=&-\frac{(\omega\times\frac{\hat{p}}{\mathfrak{c}})_i\omega_j}{(1+\frac{\hat{p}}{\mathfrak{c}} \cdot \omega)^2}.\label{W.12-5}
		\end{align}
		Note that the kernels $a_A$, $b_A$, $c_A$, $a_B$, $b_B$, and $c_B$ are smooth functions of $(\omega,p)$, since their denominators are uniformly bounded away from zero, and they satisfy the inequalities
		\begin{align}\label{W.12}
			\sum_{i=0}^1\left|\nabla_p^i (a_A,b_A,c_A,a_B,b_B,c_B)(\omega, \hat{p})\right| \lesssim (1+|p|)^8.
		\end{align}
		Furthermore,
		\begin{align}\label{W.13}
			\int_{|\omega|=1} a_A(\omega, \hat{p}) d \omega=0 \quad \text{and} \quad  \int_{|\omega|=1} a_B(\omega, \hat{p}) d \omega=0.
		\end{align}
	 The specific calculations are provided in the Appendix \ref{Appendix C.2}.
	\end{lemma}

    \begin{lemma}\label{lW.7}
    Set
    \begin{equation}\label{W.sing.a0}
        a^0_{ij}(\omega):=3\omega_i\omega_j-\delta_{ij}.
    \end{equation}
    For each fixed integer $m\ge 8$, there exists a constant $C_m>0$, independent of $\mathfrak{c}$, such that
    \begin{equation}\label{W.sing.1}
        |a_A(\omega,\hat p)-a^0_{ij}(\omega)|+|a_B(\omega,\hat p)|\le C_m\alpha_{\mathfrak{c}}(p)\langle p\rangle^m.
    \end{equation}
    Moreover,
    \begin{equation}\label{W.sing.2}
        \int_{\mathbb S^2}\big(a_A(\omega,\hat p)-a^0_{ij}(\omega)\big)\,d\omega=0,\qquad
        \int_{\mathbb S^2}a_B(\omega,\hat p)\,d\omega=0.
    \end{equation}
    \end{lemma}
    \begin{proof}
    Substituting $\hat p=0$ in \eqref{W.11-3} and \eqref{W.12-3} gives
    \[
        a_A(\omega,0)=3\omega_i\omega_j-\delta_{ij}=a^0_{ij}(\omega),\qquad a_B(\omega,0)=0.
    \]
    If $|p|\le \mathfrak{c}$, then $|\hat p|/\mathfrak{c}=|p|/p^0\le 1/\sqrt2$. Hence the kernels in \eqref{W.11-3} and \eqref{W.12-3} are smooth functions of $\hat p/\mathfrak{c}$ in a fixed compact subset of the unit ball, and the mean-value theorem yields
    \[
        |a_A(\omega,\hat p)-a^0_{ij}(\omega)|+|a_B(\omega,\hat p)|\lesssim \frac{|p|}{\mathfrak{c}}.
    \]
    If $|p|>\mathfrak{c}$, then $\alpha_{\mathfrak{c}}(p)=1$, and the polynomial bound \eqref{W.12} gives \eqref{W.sing.1} after increasing $m$ if necessary. Finally, \eqref{W.sing.2} follows from \eqref{W.13} and from $\int_{\mathbb S^2}a^0_{ij}(\omega)\,d\omega=0$.
    \end{proof}

	\begin{proposition}\label{pW.7}
Assume the {\it a priori} assumptions \eqref{W.0}-\eqref{C.2}. For every $T>0$, there exist constants $C_T>0$ and $\varepsilon_0=\varepsilon_0(T)>0$, independent of $\mathfrak{c}$ and $\varepsilon$, such that, for $0<\varepsilon\leq\varepsilon_0$, $\mathfrak{c}\geq\mathfrak{c}_0$, and $t\in[0,T]$, one has
\begin{align}\label{W.14}
\sup_{0\leq s\leq t}\varepsilon^{\frac{7}{2}}&\big\|\big(\nabla_x E_R^{\varepsilon,\mathfrak{c}},\nabla_x B_R^{\varepsilon,\mathfrak{c}}\big)(s)\big\|_{L^{\infty}} \nonumber\\
&\leq C_T\Big[\varepsilon^{\frac{5}{2}}\big\|\big(E_R^{\varepsilon,\mathfrak{c}},B_R^{\varepsilon,\mathfrak{c}}\big)(0)\big\|_{L^\infty}+\varepsilon^{\frac{3}{2}}\mathcal I_0
+\sup_{0\leq s\leq t}\varepsilon^{\frac{5}{2}}\big\|\big(E_R^{\varepsilon,\mathfrak{c}},B_R^{\varepsilon,\mathfrak{c}}\big)(s)\big\|_{L^\infty}\nonumber\\
&\qquad\qquad+\sup_{0\leq s\leq t}\big\|f_R^{\varepsilon,\mathfrak{c}}(s)\big\|+\varepsilon\sup_{0\leq s\leq t}\big\|f_R^{\varepsilon,\mathfrak{c}}(s)\big\|_{H^1}+\varepsilon^{k+\frac52}\Big].
\end{align}
The constant $C_T$ is independent of $\mathfrak{c}$ and $\varepsilon$.
\end{proposition}
\begin{proof}
We use the shortened notation introduced before Lemma \ref{lW.3}. We divide the proof into three steps.

{\it Step 1. }
The estimate of the initial terms is the same as in the proof of Proposition \ref{pW.5}, except that one differentiates the spherical means once more. Using the Maxwell equations at $t=0$, the decay assumption \eqref{W0.2}, and
\[
\mathfrak{c} t\int_{\mathbb S^2}\frac{d\omega}{1+|x+\mathfrak{c} t\omega|}\leq C,
\]
one obtains, for $i,j=1,2,3$,
\begin{align}\label{W.14-1}
\big|\mathbf I_{E_{R,ij}^{\varepsilon,\mathfrak{c}}}(t,x)\big|+\big|\mathbf I_{B_{R,ij}^{\varepsilon,\mathfrak{c}}}(t,x)\big|
\leq C_T\big\|\big(E_R^{\varepsilon,\mathfrak{c}},B_R^{\varepsilon,\mathfrak{c}}\big)(0)\big\|_{L^\infty}+C\Big(1+\frac1\varepsilon\Big)\mathcal I_0+C_T\varepsilon^k.
\end{align}
The factor $\varepsilon^{-1}$ comes only from replacing $SF_R(0)$ by the remainder equation \eqref{e1.10}. It is harmless after multiplication by $\varepsilon^{7/2}$.

{\it Step 2. }
We first treat the singular term in \eqref{W.11}. Put $z=y-x$, $\omega=z/|z|$, and decompose
\[
a_A(\omega,\hat p)=a^0_{ij}(\omega)+\widetilde a_A(\omega,
\hat p),\qquad a^0_{ij}(\omega):=3\omega_i\omega_j-\delta_{ij}.
\]
Let
\begin{align*}
\mathcal A^E_{ij}(t,x)&:=\iint_{|z|<\mathfrak{c} t}\frac{a_A(\omega,\hat p)}{|z|^3}F_R\Big(t-\frac{|z|}{\mathfrak{c}},x+z,p\Big)\,dz\,dp\\
&=\mathcal A^{E,0}_{ij}(t,x)+\widetilde{\mathcal A}^{E}_{ij}(t,x)
\end{align*}
be the corresponding decomposition. Since $\int_{\mathbb S^2}\widetilde a_A(\omega,\hat p)\,d\omega=0$ by \eqref{W.sing.2}, the near-field part of $\widetilde{\mathcal A}^{E}$ is written with a spatial difference:
\begin{align*}
&\left|\int_{0}^{1\wedge \mathfrak{c} t}\int_{\mathbb S^2}\int_{\mathbb R^3}\frac{\widetilde a_A(\omega,\hat p)}{r}
\Big[F_R\Big(t-\frac r{\mathfrak{c}},x+r\omega,p\Big)-F_R\Big(t-\frac r{\mathfrak{c}},x,p\Big)\Big]dp\,d\omega\,dr\right|\\
&\qquad\lesssim \sup_{0\leq \tau\leq t}\|\nabla_x h_R(\tau)\|_{L^\infty}
\int_{\mathbb R^3}\alpha_{\mathfrak{c}}(p)\langle p\rangle^m\frac{\sqrt{J_{\mathbf M}}}{w_\ell}\,dp
\lesssim \sup_{0\leq \tau\leq t}\|\nabla_x h_R(\tau)\|_{L^\infty}.
\end{align*}
For the far field, which is absent if $\mathfrak{c} t\leq1$, Lemma \ref{lW.7} gives
\begin{align*}
&\left|\int_{1<|z|<\mathfrak{c} t}\int_{\mathbb R^3}\frac{\widetilde a_A(\omega,\hat p)}{|z|^3}F_R\Big(t-\frac{|z|}{\mathfrak{c}},x+z,p\Big)\,dp\,dz\right|\\
&\qquad\lesssim \sup_{0\leq \tau\leq t}\|h_R(\tau)\|_{L^\infty}
\int_1^{\mathfrak{c} t}\frac{dr}{r}\int_{\mathbb R^3}\alpha_{\mathfrak{c}}(p)\langle p\rangle^m\frac{\sqrt{J_{\mathbf M}}}{w_\ell}\,dp\\
&\qquad\lesssim \frac{\log(1+\mathfrak{c} t)}{\mathfrak{c}}\sup_{0\leq \tau\leq t}\|h_R(\tau)\|_{L^\infty}
\leq C_T\sup_{0\leq \tau\leq t}\|h_R(\tau)\|_{L^\infty}.
\end{align*}
Therefore
\begin{equation}\label{W.14-2a}
\|\widetilde{\mathcal A}^{E}(t)\|_{L_x^\infty}
\leq C_T\sup_{0\leq \tau\leq t}\Big(\|h_R(\tau)\|_{L^\infty}+\|\nabla_x h_R(\tau)\|_{L^\infty}\Big).
\end{equation}
The principal part is exactly the retarded Hessian-Coulomb operator:
\[
\mathcal A^{E,0}_{ij}(t,x)=\mathscr H_{ij}(t,x).
\]
Hence Lemmas \ref{lW.5} and \ref{lW.4}, together with \eqref{W.rhoj.basic}, imply
\begin{align*}
\|\mathcal A^{E,0}(t)\|_{L_x^\infty}
&\leq \|\mathscr H(t)-\mathcal H_{\mathfrak{c} t}[\rho_R(t)]\|_{L_x^\infty}
+\|\mathcal H_{\mathfrak{c} t}[\rho_R(t)]\|_{L_x^\infty}\\
&\leq C_T\sup_{0\leq \tau\leq t}\|\nabla_x j_R(\tau)\|_{L_x^\infty}
+C\|\nabla_x\rho_R(t)\|_{L_x^\infty}^{\frac35}\|\rho_R(t)\|_{L_x^2}^{\frac25}\\
&\leq C_T\sup_{0\leq \tau\leq t}\|\nabla_x h_R(\tau)\|_{L^\infty}
+C\|\nabla_x h_R(t)\|_{L^\infty}^{\frac35}\|f_R(t)\|^{\frac25}.
\end{align*}
Combining this with \eqref{W.14-2a}, we obtain
\begin{align}\label{W.14-2}
\|\mathcal A^E(t)\|_{L_x^\infty}
&\leq C_T\sup_{0\leq \tau\leq t}\Big(\|h_R(\tau)\|_{L^\infty}+\|\nabla_x h_R(\tau)\|_{L^\infty}\Big)
+C\|\nabla_x h_R(t)\|_{L^\infty}^{\frac35}\|f_R(t)\|^{\frac25}.
\end{align}
The magnetic singular term is easier. Since $a_B(\omega,0)=0$, there is no Hessian-Coulomb principal part. Using Lemma \ref{lW.7} in the same near-field and far-field decomposition gives
\begin{align}\label{W.14-3}
\left\|\iint_{|z|<\mathfrak{c} t}\frac{a_B(\omega,\hat p)}{|z|^3}F_R\Big(t-\frac{|z|}{\mathfrak{c}},x+z,p\Big)\,dz\,dp\right\|_{L_x^\infty}
\leq C_T\sup_{0\leq \tau\leq t}\Big(\|h_R(\tau)\|_{L^\infty}+\|\nabla_xh_R(\tau)\|_{L^\infty}\Big).
\end{align}
Equations \eqref{W.14-2}-\eqref{W.14-3} give the singular-kernel bound without any logarithmic growth in $\mathfrak{c}$.

{\it Step 3. }
The terms in \eqref{W.11} and \eqref{W.120} containing $b_A,b_B,c_A,c_B$ do not contain the singular Hessian-Coulomb kernel. In the following two estimates, $b_{A,B}$ denotes either $b_A$ or $b_B$, and $c_{A,B}$ denotes either $c_A$ or $c_B$. Substituting $SF_R$ and $S^2F_R$ by the remainder equation \eqref{e1.10}, using the kernel bound \eqref{W.12}, and integrating by parts in $p$ in the force terms as in the proof of Proposition \ref{pW.5}, one obtains
\begin{align}\label{W.15-0}
&\frac1{\mathfrak{c}}\left\|\iint_{|z|<\mathfrak{c} t}\frac{b_{A,B}(\omega,\hat p)}{|z|^2}(SF_R)\Big(t-\frac{|z|}{\mathfrak{c}},x+z,p\Big)\,dz\,dp\right\|_{L_x^\infty}\nonumber\\
&\qquad\leq C_T\sup_{0\leq s\leq t}\big\|\big(E_R,B_R\big)(s)\big\|_{L^\infty}+\frac{C_T}{\varepsilon}\sup_{0\leq s\leq t}\|h_R(s)\|_{L^\infty}+C_T\varepsilon^k,
\end{align}
and
\begin{align}\label{W.19}
&\frac1{\mathfrak{c}^2}\left\|\iint_{|z|<\mathfrak{c} t}\frac{c_{A,B}(\omega,\hat p)}{|z|}(S^2F_R)\Big(t-\frac{|z|}{\mathfrak{c}},x+z,p\Big)\,dz\,dp\right\|_{L_x^\infty}\nonumber\\
&\qquad\leq C_T\int_0^t(t-\tau)\big\|\big(\nabla_xE_R,\nabla_xB_R\big)(\tau)\big\|_{L^\infty}\,d\tau
+\frac{C_T}{\varepsilon}\int_0^t(t-\tau)\|\nabla_{x,p}h_R(\tau)\|_{L^\infty}\,d\tau\nonumber\\
&\qquad\quad+\frac{C_T}{\varepsilon}\sup_{0\leq s\leq t}\big\|\big(E_R,B_R\big)(s)\big\|_{L^\infty}
+\frac{C_T}{\varepsilon^2}\sup_{0\leq s\leq t}\|h_R(s)\|_{L^\infty}+C_T\varepsilon^{k-1}.
\end{align}
We briefly recall the origin of the two singular powers of $\varepsilon$ in \eqref{W.15-0}-\eqref{W.19}. They come from the linearized collision part in $SF_R$ and from applying $S$ once more in $S^2F_R$; no factor of $\mathfrak{c}$ is lost, since the volume factor of the light cone is balanced by $\mathfrak{c}^{-1}$ or $\mathfrak{c}^{-2}$ in front of the kernels. Combining \eqref{W.14-1}-\eqref{W.19} gives
\begin{align}\label{W.20}
&\big\|\big(\nabla_xE_R,\nabla_xB_R\big)(t)\big\|_{L^\infty}
\leq C_T\Big[\big\|\big(E_R,B_R\big)(0)\big\|_{L^\infty}+\frac1\varepsilon\mathcal I_0+\frac1\varepsilon\sup_{s\leq t}\big\|\big(E_R,B_R\big)(s)\big\|_{L^\infty}\nonumber\\
&\quad+\int_0^t(t-\tau)\big\|\big(\nabla_xE_R,\nabla_xB_R\big)(\tau)\big\|_{L^\infty}\,d\tau
+\frac1\varepsilon\int_0^t(t-\tau)\|\nabla_{x,p}h_R(\tau)\|_{L^\infty}\,d\tau\nonumber\\
&\quad+\sup_{s\leq t}\big(\|h_R(s)\|_{L^\infty}+\|\nabla_xh_R(s)\|_{L^\infty}\big)
+\|\nabla_x h_R(t)\|_{L^\infty}^{\frac35}\|f_R(t)\|^{\frac25}
+\frac1{\varepsilon^2}\sup_{s\leq t}\|h_R(s)\|_{L^\infty}+\varepsilon^{k-1}\Big].
\end{align}
Here and below $C_T$ absorbs harmless powers of $T$.
Define
\begin{align*}
\Phi_0(t)&:=\sup_{0\leq s\leq t}\varepsilon^{\frac52}\big\|\big(E_R,B_R\big)(s)\big\|_{L^\infty},\qquad
\Phi_1(t):=\sup_{0\leq s\leq t}\varepsilon^{\frac72}\big\|\big(\nabla_xE_R,\nabla_xB_R\big)(s)\big\|_{L^\infty},\\
H_0(t)&:=\sup_{0\leq s\leq t}\varepsilon^{\frac32}\|h_R(s)\|_{L^\infty},
\qquad
G_1(t):=\sup_{0\leq s\leq t}\varepsilon^{\frac32}\|\nabla_{x,p}h_R(s)\|_{L^\infty},\\
F_0(t)&:=\sup_{0\leq s\leq t}\|f_R(s)\|,
\qquad
F_1(t):=\sup_{0\leq s\leq t}\|f_R(s)\|_{H^1}.
\end{align*}
Multiplying \eqref{W.20} by $\varepsilon^{7/2}$ and taking the supremum over $[0,t]$ yields
\begin{align}\label{W.22}
\Phi_1(t)&\leq C_T\Big[\varepsilon^{\frac52}\big\|\big(E_R,B_R\big)(0)\big\|_{L^\infty}+\varepsilon^{\frac32}\mathcal I_0+\Phi_0(t)+H_0(t)+\varepsilon^2 G_1(t)\nonumber\\
&\qquad+\varepsilon\int_0^t(t-\tau)G_1(\tau)\,d\tau+
\int_0^t(t-\tau)\Phi_1(\tau)\,d\tau
+\varepsilon^{\frac{13}{5}}G_1(t)^{\frac35}F_0(t)^{\frac25}+\varepsilon^{k+\frac52}\Big].
\end{align}
Indeed,
\[
\varepsilon^{\frac72}\|\nabla_x h_R(t)\|_{L^\infty}^{\frac35}\|f_R(t)\|^{\frac25}
=\varepsilon^{\frac{13}{5}}G_1(t)^{\frac35}F_0(t)^{\frac25}.
\]
By Lemmas \ref{pC.3} and \ref{pD.1}, together with the definition of $\mathcal I_0$,
\begin{align}\label{W.23}
H_0(t)&\leq C\varepsilon^{\frac32}\mathcal I_0+C\Phi_0(t)+CF_0(t)+C\varepsilon^{k+\frac52},\nonumber\\
G_1(t)&\leq C\varepsilon^{\frac32}\mathcal I_0+C\Phi_0(t)+C\varepsilon^{-1}\Phi_1(t)+CF_1(t)+C\varepsilon^{k+\frac52}.
\end{align}
Moreover,
\[
\varepsilon^{\frac{13}{5}}G_1(t)^{\frac35}F_0(t)^{\frac25}
\leq C\varepsilon^{\frac{13}{5}}\big(G_1(t)+F_0(t)\big).
\]
The terms $\varepsilon^2G_1(t)$ and $\varepsilon\int_0^t(t-\tau)G_1(\tau)\,d\tau$ are estimated by \eqref{W.23}; the former produces only $C_T\varepsilon\Phi_1(t)$, while the latter contributes another Volterra term $C_T\int_0^t(t-\tau)\Phi_1(\tau)\,d\tau$. The part of the interpolation term containing $G_1$ gives at most $C_T\varepsilon^{8/5}\Phi_1(t)$ plus lower-order terms of the same type as in \eqref{W.23}. Substitution of \eqref{W.23} into \eqref{W.22} therefore gives
\begin{align}\label{W.24}
\Phi_1(t)&\leq C_T\Big[\varepsilon^{\frac52}\big\|\big(E_R,B_R\big)(0)\big\|_{L^\infty}+\varepsilon^{\frac32}\mathcal I_0+\Phi_0(t)+F_0(t)+\varepsilon F_1(t)+\varepsilon^{k+\frac52}\Big]\nonumber\\
&\quad+C_T\varepsilon\Phi_1(t)+C_T\int_0^t(t-\tau)\Phi_1(\tau)\,d\tau.
\end{align}
Choose $\varepsilon_0(T)>0$ so small that $C_T\varepsilon_0\leq1/2$. The term $C_T\varepsilon\Phi_1(t)$ is then absorbed into the left-hand side, and the Volterra form of Gronwall's inequality proves \eqref{W.14}. This completes the proof of Proposition \ref{pW.7}.
\end{proof}

	\section{$H^1$-Estimate for $f_R^{\varepsilon, \mathfrak{c}}$}\label{section 8}
	To close Lemmas \ref{pC.3} and \ref{pD.1}, and Propositions \ref{pW.5} and \ref{pW.7}, we still need the $H^1$ estimate of $f_R^{\varepsilon, \mathfrak{c}}$.
    
	\subsection{$L^2$-estimate of $f_R^{\varepsilon,\mathfrak{c}}$}
	To perform the $L^2$ energy estimate, we first use \eqref{e1.13} to rewrite \eqref{e1.10} as
	{\small
	\begin{align}\label{L.1}
		&\partial_t f_R^{\varepsilon,\mathfrak{c}}+\hat{p} \cdot \nabla_x f_R^{\varepsilon,\mathfrak{c}}+\frac{1}{\varepsilon}\mathbf{L}_{\mathfrak{c}} f_R^{\varepsilon,\mathfrak{c}}+(E_R^{\varepsilon,\mathfrak{c}}+\frac{p}{p^0}\times B_R^{\varepsilon,\mathfrak{c}})\cdot\Big(\frac{p}{p^0}\frac{u_{\mathfrak{e}}^0}{T_\mathfrak{e}}-\frac{u_{\mathfrak{e}}}{T_\mathfrak{e}}\Big)\sqrt{\mathbf{M}_{\mathfrak{c}}}
		-\Big(E^{\varepsilon,\mathfrak{c}}+\frac{p}{p^0} \times B^{\varepsilon,\mathfrak{c}}\Big) \cdot \nabla_p f_R^{\varepsilon,\mathfrak{c}}\nonumber \\
		&\quad=-\frac{f_R^{\varepsilon,\mathfrak{c}}}{\sqrt{\mathbf{M}_{\mathfrak{c}}}}\Big\{\partial_t+\hat{p} \cdot \nabla_x-\big(E_0^{\mathfrak{c}}+\frac{p}{p^0} \times B_0^{\mathfrak{c}}\big) \cdot \nabla_p\Big\} \sqrt{\mathbf{M}_{\mathfrak{c}}}+\varepsilon^{k-1} \Gamma_{\mathfrak{c}}\left(f_R^{\varepsilon,\mathfrak{c}}, f_R^{\varepsilon,\mathfrak{c}}\right) \nonumber\\
		&\qquad +\sum_{i=1}^{2 k-1} \varepsilon^{i-1}\Big\{\Gamma_{\mathfrak{c}}\Big(\frac{F_i^{\mathfrak{c}}}{\sqrt{\mathbf{M}_{\mathfrak{c}}}}, f_R^{\varepsilon,\mathfrak{c}}\Big)+\Gamma_{\mathfrak{c}}\Big(f_R^{\varepsilon,\mathfrak{c}}, \frac{F_i^{\mathfrak{c}}}{\sqrt{\mathbf{M}_{\mathfrak{c}}}}\Big)\Big\} -\varepsilon^k \Big(\frac{p}{p^0}\frac{u_{\mathfrak{e}}^0}{2T_\mathfrak{e}}-\frac{u_{\mathfrak{e}}}{2T_\mathfrak{e}}\Big) \cdot\Big(E_R^{\varepsilon,\mathfrak{c}}+\frac{p}{p^0} \times B_R^{\varepsilon,\mathfrak{c}}\Big) f_R^{\varepsilon,\mathfrak{c}} \nonumber\\
		&\qquad +\sum_{i=1}^{2 k-1} \varepsilon^i\Big(E_R^{\varepsilon,\mathfrak{c}}+\frac{p}{p^0} \times B_R^{\varepsilon,\mathfrak{c}}\Big) \cdot \frac{\nabla_p F_i^{\mathfrak{c}}}{\sqrt{\mathbf{M}_{\mathfrak{c}}}} -\sum_{i=1}^{2 k-1} \varepsilon^i\Big\{\Big(E_i^{\mathfrak{c}}+\frac{p}{p^0} \times B_i^{\mathfrak{c}}\Big) \cdot \Big(\frac{p}{p^0}\frac{u_{\mathfrak{e}}^0}{2T_\mathfrak{e}}-\frac{u_{\mathfrak{e}}}{2T_\mathfrak{e}}\Big) f_R^{\varepsilon,\mathfrak{c}}\Big\}+\varepsilon^k \bar{A}, 
	\end{align}}
	and
	\begin{align}\label{L.2}
    \left\{\begin{aligned}
		& \partial_t E_R^{\varepsilon,\mathfrak{c}}-\mathfrak{c} \nabla_x \times B_R^{\varepsilon,\mathfrak{c}}=4\pi\int_{\mathbb{R}^3} \hat{p} \sqrt{\mathbf{M}_{\mathfrak{c}}} f_R^{\varepsilon,\mathfrak{c}} d p, \\
		& \partial_t B_R^{\varepsilon, \mathfrak{c}}+\mathfrak{c} \nabla_x \times E_R^{\varepsilon, \mathfrak{c}}=0 \text {, } \\
		& \operatorname{div} E_R^{\varepsilon,\mathfrak{c}}=-4\pi\int_{\mathbb{R}^3} \sqrt{\mathbf{M}_{\mathfrak{c}}} f_R^{\varepsilon,\mathfrak{c}} d p \text {, } \\
		& \operatorname{div} B_R^{\varepsilon,\mathfrak{c}}=0 \text {, }
    \end{aligned}\right.
	\end{align}
	where $\bar{A}=\frac{A}{\sqrt{\mathbf{M}_{\mathfrak{c}}}}$.
	\begin{lemma}\label{p4.1}
		Recall $\mathbf{M}_{\mathfrak{c}}\left(n_{\mathfrak{e}}, u_{\mathfrak{e}}, T_\mathfrak{e} ; p\right), f_R^{\varepsilon,\mathfrak{c}}, h_R^{\varepsilon,\mathfrak{c}}$ defined in \eqref{e1.11}, \eqref{e1.13} and \eqref{e1.14}, respectively, and recall $\zeta_0 > 0$ in Lemma \ref{l2.1}. Then there exists constants $\varepsilon_0>0$ and $C>0$, such that for all $\varepsilon \in (0, \varepsilon_0]$, 
		\begin{align}\label{L.3}
		&\frac{1}{2}\frac{d}{d t}\Big\|\sqrt{\frac{2\mathfrak{c}T_\mathfrak{e}}{u_{\mathfrak{e}}^0}}f_R^{\varepsilon,\mathfrak{c}}\Big\|^2+\frac{d}{d t}\left\|\big(E_R^{\varepsilon,\mathfrak{c}}, B_R^{\varepsilon,\mathfrak{c}}\big)\right\|^2+\frac{\zeta_0c_0}{\varepsilon}\|\{\mathbf{I}-\mathbf{P}_{\mathfrak{c}}\}f_R^{\varepsilon,\mathfrak{c}}\|_{\nu_{\mathfrak{c}}}^2\nonumber\\
		&\quad\leq C\Big\{1+\varepsilon^{k-1}\left\|h_R^{\varepsilon,\mathfrak{c}}\right\|_{L^{\infty}}\Big\} \Big\{\left\|\big(E_R^{\varepsilon,\mathfrak{c}}, B_R^{\varepsilon,\mathfrak{c}}\big)\right\|^2+\left\|f_R^{\varepsilon,\mathfrak{c}}\right\|^2\Big\}\nonumber\\
		&\qquad +C\Big\{1+\varepsilon^{\frac{19}{6}}\left\|h_R^{\varepsilon,\mathfrak{c}}\right\|_{L^{\infty}}\Big\}\left\|f_R^{\varepsilon,\mathfrak{c}}\right\|,
		\end{align}
		where the constant $C>0$ depends on $ F_i^{\mathfrak{c}}$, $E_i^{\mathfrak{c}}$, $B_i^{\mathfrak{c}}$$(i=0,1,\cdots,2k-1)$, and is independent of $\mathfrak{c}$.
	\end{lemma}
	\begin{proof}
		Multiplying \eqref{L.1} by $\frac{2 T_\mathfrak{e}}{u_{\mathfrak{e}}^0} f_R^{\varepsilon, \mathfrak{c}}$ and integrating over $\mathbb{R}^3\times \mathbb{R}^3$, one obtains that
		{\footnotesize
		\begin{align}\label{L.4}
			&\frac{1}{2}\frac{d}{d t}\Big\|\sqrt{\frac{2T_\mathfrak{e}}{u_{\mathfrak{e}}^0}}f_R^{\varepsilon,\mathfrak{c}}\Big\|^2+ \frac{1}{\varepsilon}\Big\langle \mathbf{L}_{\mathfrak{c}} f_R^{\varepsilon,\mathfrak{c}}, \frac{2T_\mathfrak{e}}{u_{\mathfrak{e}}^0}f_R^{\varepsilon,\mathfrak{c}} \Big\rangle +\Big\langle\big(E_R^{\varepsilon,\mathfrak{c}}+\frac{p}{p^0}\times B_R^{\varepsilon,\mathfrak{c}}\big)\cdot2\Big(\frac{p}{p^0}-\frac{u_{\mathfrak{e}}}{u_{\mathfrak{e}}^0}\Big),\sqrt{\mathbf{M}_{\mathfrak{c}}}f_R^{\varepsilon,\mathfrak{c}}\Big\rangle
			\nonumber\\
            &=\Big\langle\Big\{\left(\partial_t+\hat{p} \cdot \nabla_x\right)\left(\frac{ T_\mathfrak{e}}{u_{\mathfrak{e}}^0}\right)\Big\} f_R^{\varepsilon,\mathfrak{c}}, f_R^{\varepsilon,\mathfrak{c}}\Big\rangle+\Big\langle\Big(E^{\varepsilon,\mathfrak{c}}+\frac{p}{p^0} \times B^{\varepsilon,\mathfrak{c}}\Big) \cdot \nabla_p f_R^{\varepsilon,\mathfrak{c}},\frac{2T_\mathfrak{e}}{u_{\mathfrak{e}}^0}f_R^{\varepsilon,\mathfrak{c}}\Big\rangle\nonumber\\
			& \quad-\Big\langle\frac{f_R^{\varepsilon,\mathfrak{c}}}{\sqrt{\mathbf{M}_{\mathfrak{c}}}}\Big\{\partial_t+\hat{p} \cdot \nabla_x-\big(E_0^{\mathfrak{c}}+\frac{p}{p^0} \times B_0^{\mathfrak{c}}\big) \cdot \nabla_p\Big\} \sqrt{\mathbf{M}_{\mathfrak{c}}},\frac{2T_\mathfrak{e}}{u_{\mathfrak{e}}^0}f_R^{\varepsilon,\mathfrak{c}}\Big\rangle\nonumber\\
			&\quad+\varepsilon^{k-1}\Big\langle \Gamma_{\mathfrak{c}}\left(f_R^{\varepsilon,\mathfrak{c}}, f_R^{\varepsilon,\mathfrak{c}}\right),\frac{2T_\mathfrak{e}}{u_{\mathfrak{e}}^0}f_R^{\varepsilon,\mathfrak{c}}\Big\rangle +\sum_{i=1}^{2 k-1} \varepsilon^{i-1}\Big\langle\Big\{\Gamma_{\mathfrak{c}}\Big(\frac{F_i^{\mathfrak{c}}}{\sqrt{\mathbf{M}_{\mathfrak{c}}}}, f_R^{\varepsilon,\mathfrak{c}}\Big)+\Gamma_{\mathfrak{c}}\Big(f_R^{\varepsilon,\mathfrak{c}}, \frac{F_i^{\mathfrak{c}}}{\sqrt{\mathbf{M}_{\mathfrak{c}}}}\Big)\Big\},\frac{2T_\mathfrak{e}}{u_{\mathfrak{e}}^0}f_R^{\varepsilon,\mathfrak{c}}\Big\rangle\nonumber\\
			&\quad-\varepsilon^k \Big\langle \Big(E_R^{\varepsilon,\mathfrak{c}}+\frac{p}{p^0} \times B_R^{\varepsilon,\mathfrak{c}}\Big)\cdot\left( \frac{p}{p^0}-\frac{u_{\mathfrak{e}}}{u_{\mathfrak{e}}^0}\right) f_R^{\varepsilon,\mathfrak{c}} ,f_R^{\varepsilon,\mathfrak{c}}\Big\rangle+\sum_{i=1}^{2 k-1} \varepsilon^i\Big\langle\Big(E_R^{\varepsilon, \mathfrak{c}}+\frac{p}{p^0} \times B_R^{\varepsilon, \mathfrak{c}}\Big) \cdot \frac{\nabla_p F_i^{\mathfrak{c}}}{\sqrt{\mathbf{M}_{\mathfrak{c}}}}, \frac{2 T_\mathfrak{e}}{u_{\mathfrak{e}}^0} f_R^{\varepsilon, \mathfrak{c}}\Big\rangle\nonumber\\
			& \quad-\sum_{i=1}^{2 k-1} \varepsilon^i\Big\langle\Big(E_i^{\mathfrak{c}}+\frac{p}{p^0} \times B_i^{\mathfrak{c}}\Big) \cdot\left(\frac{p}{p^0}-\frac{u_{\mathfrak{e}}}{u_{\mathfrak{e}}^0}\right) f_R^{\varepsilon,\mathfrak{c}}, f_R^{\varepsilon, \mathfrak{c}}\Big\rangle+\varepsilon^k\Big\langle\bar{A}, \frac{2T_\mathfrak{e}}{u_{\mathfrak{e}}^0} f_R^{\varepsilon, \mathfrak{c}}\Big\rangle.
		\end{align}}
		
		According to Lemma \ref{l2.1}, it's clear that
		\begin{align}\label{L.7}
			\frac{1}{\varepsilon}\Big\langle \mathbf{L}_{\mathfrak{c}} f_R^{\varepsilon,\mathfrak{c}}, \frac{2T_\mathfrak{e}}{u_{\mathfrak{e}}^0}f_R^{\varepsilon,\mathfrak{c}} \Big\rangle \geq \frac{4\zeta_0c_0}{\varepsilon\mathfrak{c}}\|\{\mathbf{I}-\mathbf{P}_{\mathfrak{c}}\}f_R^{\varepsilon,\mathfrak{c}}\|_{\nu_{\mathfrak{c}}}^2,
		\end{align}
		where we have used $u_{\mathfrak{e}}^0<2\mathfrak{c}$ and  \eqref{e.30}.
		
		Multiplying $\eqref{L.2}_1$ by $E_R^{\varepsilon,\mathfrak{c}}$ and $\eqref{L.2}_2$ by $B_R^{\varepsilon,\mathfrak{c}}$, respectively, we obtain that
		\begin{align}\label{L.5}
			\frac{1}{2}\frac{d}{d t}\|\big(E_R^{\varepsilon,\mathfrak{c}},B_R^{\varepsilon,\mathfrak{c}}\big)\|^2=4\pi\Big\langle\hat{p}\cdot E_R^{\varepsilon,\mathfrak{c}}\sqrt{\mathbf{M}}_{\mathfrak{c}},f_R^{\varepsilon,\mathfrak{c}} \Big\rangle,
		\end{align}
		where we have used the following fact:
		$$
	\int_{\mathbb{R}^3}-\mathfrak{c} \nabla_x \times B_R^{\varepsilon,\mathfrak{c}}\cdot E_R^{\varepsilon,\mathfrak{c}} +\mathfrak{c} \nabla_x \times E_R^{\varepsilon,\mathfrak{c}}\cdot B_R^{\varepsilon,\mathfrak{c}}d x=\mathfrak{c} \int_{\mathbb{R}^3} \operatorname{div}\left(E_R^{\varepsilon,\mathfrak{c}} \times B_R^{\varepsilon,\mathfrak{c}}\right) d x=0.
	$$
	Then it follows from \eqref{L.5} that
	\begin{align}\label{L.6}
		&\Big\langle\big(E_R^{\varepsilon,\mathfrak{c}}+\frac{p}{p^0}\times B_R^{\varepsilon,\mathfrak{c}}\big)\cdot2\Big(\frac{p}{p^0}-\frac{u_{\mathfrak{e}}}{u_{\mathfrak{e}}^0}\Big),\sqrt{\mathbf{M}_{\mathfrak{c}}}f_R^{\varepsilon,\mathfrak{c}}\Big\rangle\nonumber\\
		=&\Big\langle 2E_R^{\varepsilon,\mathfrak{c}}\cdot \frac{p}{p^0},\sqrt{\mathbf{M}_{\mathfrak{c}}}f_R^{\varepsilon,\mathfrak{c}}\Big\rangle -\Big\langle2 E_R^{\varepsilon,\mathfrak{c}}\cdot \frac{u_{\mathfrak{e}}}{u_{\mathfrak{e}}^0},\sqrt{\mathbf{M}}_{\mathfrak{c}}f_R^{\varepsilon,\mathfrak{c}} \Big\rangle-\Big\langle 2\Big(\frac{p}{p^0}\times B_R^{\varepsilon,\mathfrak{c}}\Big)\cdot\frac{u_{\mathfrak{e}}}{u_{\mathfrak{e}}^0},\sqrt{\mathbf{M}}_{\mathfrak{c}}f_R^{\varepsilon,\mathfrak{c}} \Big\rangle \nonumber\\
		\gtrsim&\frac{1}{\mathfrak{c}}\left\{\frac{d}{d t}\|\big(E_R^{\varepsilon,\mathfrak{c}},B_R^{\varepsilon,\mathfrak{c}}\big)\|^2-\|\big(E_R^{\varepsilon,\mathfrak{c}},B_R^{\varepsilon,\mathfrak{c}}\big)\|^2-\|f_R^{\varepsilon,\mathfrak{c}}\|^2\right\}.
		\end{align}
		
		Now we turn to the right-hand side of \eqref{L.4}. Integrating by parts with respect to $p$, we obtain
		$$
		\Big\langle\Big(E^{\varepsilon,\mathfrak{c}}+\frac{p}{p^0} \times B^{\varepsilon,\mathfrak{c}}\Big) \cdot \nabla_p f_R^{\varepsilon, \mathfrak{c}}, \frac{2 T_\mathfrak{e}}{u_{\mathfrak{e}}^0} f_R^{\varepsilon, \mathfrak{c}}\Big\rangle=0. 
		$$
		
		For $\partial=\partial_t$ or $\partial=\partial_{x_i}$, it holds that
		\begin{align}\label{L.8}
			& \frac{\partial \mathbf{M}_{\mathfrak{c}}}{\mathbf{M}_{\mathfrak{c}}}=\frac{\partial n_\mathfrak{e}}{n_\mathfrak{e}}-3 \frac{\partial T_\mathfrak{e}}{T_\mathfrak{e}}+\frac{\partial T_\mathfrak{e}}{T_\mathfrak{e}^2}\left(u_{\mathfrak{e}}^0 p^0-\mathfrak{c}^2 \frac{K_1(\gamma)}{K_2(\gamma)}\right)-\frac{\partial T_\mathfrak{e}}{T_\mathfrak{e}^2} \sum_{i=1}^3 u_{\mathfrak{e},i} p_i+\frac{1}{T_\mathfrak{e}}\left(\sum_{i=1}^3 p_i \partial u_{\mathfrak{e},i}-\frac{\partial u_\mathfrak{e} \cdot u_\mathfrak{e}}{u_{\mathfrak{e}}^0} p^0\right), \nonumber\\
			& \frac{\partial p_i \mathbf{M}_{\mathfrak{c}}}{\mathbf{M}_{\mathfrak{c}}}=\frac{u_{\mathfrak{e},i}-u_{\mathfrak{e}}^0 \frac{p_i}{p^0}}{T_\mathfrak{e}}.
		\end{align}	
		A direct calculation shows that
		$$
		\left|u_{\mathfrak{e}}^0 p^0-\mathfrak{c}^2 \frac{K_1(\gamma)}{K_2(\gamma)}\right| \approx\left(1+|p|^2\right) C\left(n_\mathfrak{e}, u_\mathfrak{e}, T_\mathfrak{e}\right)\text {, }
		$$
		which, together with \eqref{L.8}, yields that
		$$
		\Big|\frac{1}{\sqrt{\mathbf{M}_{\mathfrak{c}}}}\left\{\partial_t+\hat{p} \cdot \nabla_x-\big(E_0^{\mathfrak{c}}+\frac{p}{p^0} \times B_0^{\mathfrak{c}}\big) \cdot \nabla_p\right\} \sqrt{\mathbf{M}_{\mathfrak{c}}}\Big| \leq C\left( \partial n_\mathfrak{e}, \partial u_\mathfrak{e}, \partial T_\mathfrak{e}, E_0^{\mathfrak{c}}, B_0^{\mathfrak{c}}\right)(1+|p|)^3.
		$$
		Noting from \eqref{e1.16}-\eqref{e1.17}, it holds $(1+|p|)^3|f_R^{\varepsilon, \mathfrak{c}}| \leq(1+|p|)^{-11}|h_R^{\varepsilon, \mathfrak{c}}|$, then we have
		$$
		\Big(\int_{1+|p| \geq \sqrt[3]{\frac{\kappa}{\varepsilon}}}(1+|p|)^{-11 \times 2} d p\Big)^{1 / 2} \lesssim\left(\frac{\varepsilon}{\kappa}\right)^{\frac{19}{6}},
		$$
		where $\kappa$ is a sufficiently small positive constant to be determined later. Thus one gets 
		\begin{align}\label{L.9}
			&\Big|\Big\langle \frac{1}{\sqrt{\mathbf{M}_{\mathfrak{c}}}}\left\{\partial_t+\hat{p} \cdot \nabla_x-\big(E_0^{\mathfrak{c}}+\frac{p}{p^0} \times B_0^{\mathfrak{c}}\big) \cdot \nabla_p\right\} \sqrt{\mathbf{M}_{\mathfrak{c}}}f_R^{\varepsilon, \mathfrak{c}},\frac{2 T_\mathfrak{e}}{u_{\mathfrak{e}}^0} f_R^{\varepsilon, \mathfrak{c}}\Big\rangle\Big|\nonumber\\
			\leq& \Big|\iint_{1+|p|\geq\sqrt[3]{\frac{\kappa}{\varepsilon}}}d x d p\Big|+\Big|\iint_{1+|p|\leq\sqrt[3]{\frac{\kappa}{\varepsilon}}}d x d p\Big|\nonumber\\
			\leq& \frac{C_{\kappa}}{\mathfrak{c}}\varepsilon^{\frac{19}{6}}\|h_R^{\varepsilon,\mathfrak{c}}\|_{L^{\infty}}\cdot\|\left(\partial n_\mathfrak{e}, \partial u_\mathfrak{e}, \partial T_\mathfrak{e}, E_0^{\mathfrak{c}}, B_0^{\mathfrak{c}}\right)\|\cdot\|f_R^{\varepsilon,\mathfrak{c}}\|\nonumber\\
			&+\frac{C}{\mathfrak{c}}\|\left(\partial n_\mathfrak{e}, \partial u_\mathfrak{e}, \partial T_\mathfrak{e}, E_0^{\mathfrak{c}}, B_0^{\mathfrak{c}}\right)\|_{L^{\infty}}\cdot\left\|(1+|p|)^{\frac{3}{2}}f_R^{\varepsilon,\mathfrak{c}}\mathbf{I}_{1+|p|\leq\sqrt[3]{\frac{\kappa}{\varepsilon}}}\right\|^2\nonumber\\
			\leq&\frac{C_{\kappa}}{\mathfrak{c}} \varepsilon^{\frac{19}{6}}\|h_R^{\varepsilon,\mathfrak{c}}\|_{L^{\infty}}\cdot\|f_R^{\varepsilon,\mathfrak{c}}\|+\frac{C}{\mathfrak{c}}\left\|(1+|p|)^{\frac{3}{2}}\mathbf{P}_{\mathfrak{c}} f_R^{\varepsilon, \mathfrak{c}} \mathbf{I}_{1+|p| \leq \sqrt[3]{\frac{\kappa}{\varepsilon}}}\right\|^2\nonumber\\
			& +\frac{C}{\mathfrak{c}}\left\|(1+|p|)^{\frac{3}{2}}\left\{\mathbf{I}-\mathbf{P}_{\mathfrak{c}}\right\} f_R^{\varepsilon, \mathfrak{c}} \mathbf{I}_{1+|p| \leq \sqrt[3]{\frac{\kappa}{\varepsilon}}}\right\|^2 \nonumber\\
			\leq& \frac{C_{\kappa}}{\mathfrak{c}} \varepsilon^{\frac{19}{6}}\left\|h_R^{\varepsilon, \mathfrak{c}}\right\|_{L^{\infty}} \cdot\left\|f_R^{\varepsilon, \mathfrak{c}}\right\|+\frac{C}{\mathfrak{c}}\left\|f_R^{\varepsilon, \mathfrak{c}}\right\|^2+\frac{C \kappa}{\varepsilon\mathfrak{c}} \|\left\{\mathbf{I}-\mathbf{P}_{\mathfrak{c}}\right\} f_R^{\varepsilon, \mathfrak{c}} \|_{\nu_{\mathfrak{c}}}^2,
		\end{align}
		where $C$ and $C_{\kappa}$ are independent of $\mathfrak{c}$. 
		
		Similar as \eqref{L.9}, there holds
		\begin{align}
			& \sum_{i=1}^{2 k-1} \varepsilon^i\Big\langle\Big(E_i^{\mathfrak{c}}+\frac{p}{p^0} \times B_i^{\mathfrak{c}}\Big) \cdot\Big(\frac{p}{p^0}-\frac{u_{\mathfrak{e}}}{u_{\mathfrak{e}}^0}\Big) f_R^{\varepsilon, \mathfrak{c}}, f_R^{\varepsilon, \mathfrak{c}}\Big\rangle \nonumber\\
			\lesssim& \frac{C_{\kappa}}{\mathfrak{c}} \varepsilon^{\frac{25}{2}}\left\|h_R^{\varepsilon, \mathfrak{c}}\right\|_{L^{\infty}} \cdot\left\|f_R^{\varepsilon, \mathfrak{c}}\right\|+\frac{\varepsilon C}{\mathfrak{c}}\left\|f_R^{\varepsilon, \mathfrak{c}}\right\|^2+\frac{C \kappa}{\mathfrak{c}}\|\left\{\mathbf{I}-\mathbf{P}_{\mathfrak{c}}\right\} f_R^{\varepsilon, \mathfrak{c}} \|_{\nu_{\mathfrak{c}}}^2,\nonumber
		\end{align}
		and
		\begin{align*}
		    &\Big\langle\Big\{\left(\partial_t+\hat{p} \cdot \nabla_x\right)\Big(\frac{ T_\mathfrak{e}}{u_{\mathfrak{e}}^0}\Big)\Big\} f_R^{\varepsilon,\mathfrak{c}}, f_R^{\varepsilon,\mathfrak{c}}\Big\rangle\\
		 \lesssim&\frac{C_{\kappa}}{\mathfrak{c}} \varepsilon^{\frac{23}{2}}\left\|h_R^{\varepsilon, \mathfrak{c}}\right\|_{L^{\infty}} \cdot\left\|f_R^{\varepsilon, \mathfrak{c}}\right\|+\frac{C}{\mathfrak{c}}\left\|f_R^{\varepsilon, \mathfrak{c}}\right\|^2+\frac{C \kappa}{\varepsilon\mathfrak{c}} \|\left\{\mathbf{I}-\mathbf{P}_{\mathfrak{c}}\right\} f_R^{\varepsilon, \mathfrak{c}} \|_{\nu_{\mathfrak{c}}}^2. \nonumber
		\end{align*}
		
		For the fourth and fifth terms on RHS of \eqref{L.4}, it follows from Lemma \ref{l2.3} that
		$$
		\varepsilon^{k-1}\Big\langle\Gamma_{\mathfrak{c}}\left(f_R^{\varepsilon, \mathfrak{c}}, f_R^{\varepsilon, \mathfrak{c}}\right), \frac{2 T_\mathfrak{e}}{u_{\mathfrak{e}}^0} f_R^{\varepsilon, \mathfrak{c}}\Big\rangle \lesssim \frac{\varepsilon^{k-1}}{\mathfrak{c}}\left\|h_R^{\varepsilon,\mathfrak{c}}\right\|_{L^{\infty}}\left\|f_R^{\varepsilon,\mathfrak{c}}\right\|^2,
		$$
		and
		\begin{align}
			& \sum_{i=1}^{2 k-1} \varepsilon^{i-1}\Big\langle\Big\{\Gamma_{\mathfrak{c}}\Big(\frac{F_i^{\mathfrak{c}}}{\sqrt{\mathbf{M}_{\mathfrak{c}}}}, f_R^{\varepsilon,\mathfrak{c}}\Big)+\Gamma_{\mathfrak{c}}\Big(f_R^{\varepsilon,\mathfrak{c}}, \frac{F_i^{\mathfrak{c}}}{\sqrt{\mathbf{M}_{\mathfrak{c}}}}\Big)\Big\}, \frac{2 T_\mathfrak{e}}{u_{\mathfrak{e}}^0} f_R^{\varepsilon,\mathfrak{c}}\Big\rangle \nonumber\\
			\lesssim &\sum_{i=1}^{2 k-1} \frac{\varepsilon^{i-1}}{\mathfrak{c}}\left\| f_R^{\varepsilon,\mathfrak{c}}\right\|_{\nu_{\mathfrak{c}}}^2 
			\lesssim \frac{1}{\mathfrak{c}}\left\|\{\mathbf{I}-\mathbf{P}_{\mathfrak{c}}\} f_R^{\varepsilon,\mathfrak{c}}\right\|^2_{\nu_{\mathfrak{c}}}+\frac{1}{\mathfrak{c}}\left\|\mathbf{P}_{\mathfrak{c}} f_R^{\varepsilon,\mathfrak{c}}\right\|^2_{\nu_{\mathfrak{c}}} \nonumber\\
			\lesssim& \frac{1}{\mathfrak{c}}\left\|\{\mathbf{I}-\mathbf{P}_{\mathfrak{c}}\}f_R^{\varepsilon,\mathfrak{c}}\right\|^2_{\nu_{\mathfrak{c}}}+\frac{1}{\mathfrak{c}}\left\|f_R^{\varepsilon,\mathfrak{c}}\right\|^2 .\nonumber
		\end{align}
		
		A direct calculation yields that the sixth and seventh terms on RHS of \eqref{L.4} can be bounded by 
		$$
		\Big(\frac{\varepsilon^k}{\mathfrak{c}}\left\|h_R^{\varepsilon,\mathfrak{c}}\right\|_{L^{\infty}}+\frac{\varepsilon}{\mathfrak{c}}\Big)\Big\{\left\|\big(E_R^{\varepsilon,\mathfrak{c}},B_R^{\varepsilon,\mathfrak{c}}\big)\right\|^2+\left\|f_R^{\varepsilon,\mathfrak{c}}\right\|^2\Big\}.
		$$
		For the last term, it holds
		\begin{align}
			\Big|\Big\langle\varepsilon^k \bar{A},f_R^{\varepsilon,\mathfrak{c}}\Big\rangle\Big|
			\lesssim& \varepsilon^k \sum_{\substack{i+j \geq 2k+1 \\ 2 \leq i, j \leq 2k-1}}\varepsilon^{i+j-2 k-1}\Big|\Big\langle\Gamma_{\mathfrak{c}}\Big(\frac{F_i^{\mathfrak{c}}}{\sqrt{\mathbf{M}_{\mathfrak{c}}}}, \frac{F_j^{\mathfrak{c}}}{\sqrt{\mathbf{M}_{\mathfrak{c}}}}\Big),\frac{2T_\mathfrak{e}}{u_{\mathfrak{e}}^0}f_R^{\varepsilon,\mathfrak{c}}\Big\rangle\Big|\nonumber\\
			& +\varepsilon^k \sum_{\substack{i+j\geq2 k \\
					1 \leq i, j \leq 2 k-1}} \varepsilon^{i+j-2 k} \Big|\Big\langle\Big(E_i^{\mathfrak{c}}+\frac{p}{p^0} \times B_i^{\mathfrak{c}}\Big) \cdot \frac{\nabla_p F_j^{\mathfrak{c}}}{\sqrt{\mathbf{M}_{\mathfrak{c}}}}, \frac{2T_\mathfrak{e}}{u_{\mathfrak{e}}^0} f_R^{\varepsilon, \mathfrak{c}}\Big\rangle\Big|\nonumber\\
			\lesssim&\frac{\varepsilon^k}{\mathfrak{c}}\left\|f_R^{\varepsilon,\mathfrak{c}}\right\|. \nonumber
		\end{align}
		
		Collecting all the above estimates, one has
		\begin{align}
			&\frac{1}{2}\frac{d}{d t}\Big\|\sqrt{\frac{2T_\mathfrak{e}}{u_{\mathfrak{e}}^0}}f_R^{\varepsilon,\mathfrak{c}}\Big\|^2+\frac{1}{\mathfrak{c}}\frac{d}{d t}\left\|\big(E_R^{\varepsilon,\mathfrak{c}}, B_R^{\varepsilon,\mathfrak{c}}\big)\right\|^2+\frac{4\zeta_0c_0}{\mathfrak{c}\varepsilon}\|\{\mathbf{I}-\mathbf{P}_{\mathfrak{c}}\}f_R^{\varepsilon,\mathfrak{c}}\|_{\nu_{\mathfrak{c}}}^2\nonumber\\
			\leq& C\Big(\frac{1}{\mathfrak{c}}+\frac{\varepsilon^{k-1}}{\mathfrak{c}}\left\|h_R^{\varepsilon,\mathfrak{c}}\right\|_{L^{\infty}}\Big) \Big\{\left\|\big(E_R^{\varepsilon,\mathfrak{c}}, B_R^{\varepsilon,\mathfrak{c}}\big)\right\|^2+\left\|f_R^{\varepsilon,\mathfrak{c}}\right\|^2\Big\}\nonumber\\
			& +C\Big(\frac{1}{\mathfrak{c}}+\frac{\varepsilon^{\frac{19}{6}}}{\mathfrak{c}}\left\|h_R^{\varepsilon,\mathfrak{c}}\right\|_{L^{\infty}}\Big)\left\|f_R^{\varepsilon,\mathfrak{c}}\right\|+C\Big(\frac{ \kappa}{\mathfrak{c}\varepsilon}+\frac{ \kappa}{\mathfrak{c}}+\frac{1}{\mathfrak{c}}\Big)\left\|\left\{\mathbf{I}-\mathbf{P}_{\mathfrak{c}}\right\} f_R^{\varepsilon, \mathfrak{c}}\right\|_{\nu_{\mathfrak{c}}}^2.\nonumber
		\end{align}
		First, by taking $\kappa$ sufficiently small and then taking $\varepsilon_0$ sufficiently small, one concludes \eqref{L.3}.
		Therefore the proof of Lemma \ref{p4.1} is completed.
	\end{proof}
	\subsection{$L^2$-estimate of $\nabla_{x,p}f_R^{\varepsilon,\mathfrak{c}}$}
	We first establish the spatial derivative estimate.
	\begin{lemma}\label{pH.1}
		There exists constants $\varepsilon_0>0$ and $C>0$, such that for all $\varepsilon \in (0, \varepsilon_0]$,
		\begin{align}\label{H.1}
			\frac{1}{2}\frac{d}{d t}&\Big\|  \sqrt{\frac{2\mathfrak{c}T_\mathfrak{e}}{u_{\mathfrak{e}}^0}} \nabla_x f_R^{\varepsilon, \mathfrak{c}}\Big\|^2+\frac{d}{d t}\left\|\big(\nabla_x E_R^{\varepsilon, \mathfrak{c}}, \nabla_x B_R^{\varepsilon, \mathfrak{c}}\big) \right\|^2+\frac{\zeta_0 c_0}{\varepsilon}\left\|\{\mathbf{I}-\mathbf{P}_{\mathfrak{c}}\}\nabla_x f_R^{\varepsilon, \mathfrak{c}}\right\|_{\nu_{\mathfrak{c}}}^2 \nonumber\\
			& \leq C\Big\{1+\varepsilon^{k-1}\left\|h_R^{\varepsilon, \mathfrak{c}}\right\|_{W^{1,\infty}}\Big\}\Big\{\left\|f_R^{\varepsilon, \mathfrak{c}}\right\|_{H^1}^2+\left\|\big(E_R^{\varepsilon, \mathfrak{c}}, B_R^{\varepsilon, \mathfrak{c}}\big)\right\|_{H^1}^2\Big\}+\frac{C}{\varepsilon^2}\left\| f_R^{\varepsilon, \mathfrak{c}}\right\|^2 \nonumber\\
			& \quad+C\varepsilon^{\frac{19}{6}}\left\|h_R^{\varepsilon, \mathfrak{c}}\right\|_{L^{\infty}}\left\|f_R^{\varepsilon, \mathfrak{c}}\right\|+C\Big\{\varepsilon^{\frac{19}{6}}\left\|h_R^{\varepsilon, \mathfrak{c}}\right\|_{W^{1,\infty}}+\varepsilon^{k}\Big\}\left\|\nabla_x f_R^{\varepsilon, \mathfrak{c}}\right\| \nonumber\\
			&\quad+\frac{C}{\varepsilon^2}\left\|\{\mathbf{I}-\mathbf{P}_{\mathfrak{c}}\}f_R^{\varepsilon, \mathfrak{c}}\right\|_{\nu_{\mathfrak{c}}}^2+C\left\|\nabla_x f_R^{\varepsilon, \mathfrak{c}}\right\|\left\|\nabla_pf_R^{\varepsilon, \mathfrak{c}}\right\|,
		\end{align}
	where the constant $C>0$ depends on $ F_i^{\mathfrak{c}}$, $E_i^{\mathfrak{c}}$, $B_i^{\mathfrak{c}}(i=0,1,\cdots,2k-1)$, and is independent of $\mathfrak{c}$.
	\end{lemma}
	
	\begin{proof}
		Applying $D_x$ to \eqref{L.1}, then multiplying the resultant equation by $\frac{2 T_\mathfrak{e}}{u_{\mathfrak{e}}^0} D_x f_R^{\varepsilon, \mathfrak{c}}$ and integrating over $\mathbb{R}^3\times \mathbb{R}^3$, one has
		{\footnotesize
		\begin{align}\label{H.2}
		    \frac{1}{2} \frac{d}{d t}&\Big\|\sqrt{\frac{2T_\mathfrak{e}}{u_{\mathfrak{e}}^0}} D_x f_R^{\varepsilon, \mathfrak{c}}\Big\|^2+\frac{1}{\varepsilon}\Big\langle D_x \big(\mathbf{L}_{\mathfrak{c}}f_R^{\varepsilon, \mathfrak{c}}\big), \frac{2T_\mathfrak{e}}{u_{\mathfrak{e}}^0}D_x f_R^{\varepsilon, \mathfrak{c}}\Big\rangle \nonumber\\
			&+\Big\langle D_x\Big\{\Big(E_R^{\varepsilon, \mathfrak{c}}+\frac{p}{p^0} \times B_R^{\varepsilon, \mathfrak{c}}\Big) \Big(\frac{p}{p^0}\frac{u_{\mathfrak{e}}^0}{T_\mathfrak{e}}-\frac{u_{\mathfrak{e}}}{T_\mathfrak{e}}\Big)\sqrt{\mathbf{M}_{\mathfrak{c}}}\Big\}, \frac{2T_\mathfrak{e}}{u_{\mathfrak{e}}^0} D_x f_R^{\varepsilon, \mathfrak{c}}\Big\rangle \nonumber\\
			=&-\Big\langle\frac{D_x f_R^{\varepsilon, \mathfrak{c}}}{\sqrt{\mathbf{M}_{\mathfrak{c}}}}\Big\{\partial_t+\hat{p} \cdot \nabla_x-\big(E_0^{\mathfrak{c}}+\frac{p}{p^0} \times B_0^{\mathfrak{c}}\big) \cdot \nabla_p\Big\} \sqrt{\mathbf{M}_{\mathfrak{c}}}, \frac{2T_\mathfrak{e}}{u_{\mathfrak{e}}^0} D_x f_R^{\varepsilon, \mathfrak{c}}\Big\rangle  \nonumber\\
			& -\Big\langle f_R^{\varepsilon, \mathfrak{c}} D_x\left\{\frac{1}{\sqrt{\mathbf{M}_{\mathfrak{c}}}}\Big[\partial_t+\frac{p}{p^0} \cdot \nabla_x-\big(E_0^{\mathfrak{c}}+\frac{p}{p^0} \times B_0^{\mathfrak{c}}\big) \cdot \nabla_p\Big] \sqrt{\mathbf{M}_{\mathfrak{c}}}\right\}, \frac{2T_\mathfrak{e}}{u_{\mathfrak{e}}^0} D_x f_R^{\varepsilon, \mathfrak{c}}\Big\rangle  \nonumber\\
			&+\Big\langle\Big\{\left(\partial_t+\hat{p} \cdot \nabla_x\right)\big(\frac{T_\mathfrak{e}}{u_{\mathfrak{e}}^0}\big)\Big\} D_x f_R^{\varepsilon, \mathfrak{c}}, D_x f_R^{\varepsilon, \mathfrak{c}}\Big\rangle+\varepsilon^{k-1}\Big\langle D_x \big(\Gamma_{\mathfrak{c}}\left(f_R^{\varepsilon, \mathfrak{c}}, f_R^{\varepsilon, \mathfrak{c}}\right)\big), \frac{2T_\mathfrak{e}}{u_{\mathfrak{e}}^0} D_x f_R^{\varepsilon, \mathfrak{c}}\Big\rangle  \nonumber\\
			& +\sum_{i=1}^{2 k-1} \varepsilon^{i-1}\Big\langle D_x\Big\{\Gamma_{\mathfrak{c}}\Big(\frac{F_i^{\mathfrak{c}}}{\sqrt{\mathbf{M}_{\mathfrak{c}}}}, f_R^{\varepsilon, \mathfrak{c}}\Big)+\Gamma_{\mathfrak{c}}\Big(f_R^{\varepsilon, \mathfrak{c}}, \frac{F_i^{\mathfrak{c}}}{\sqrt{\mathbf{M}_{\mathfrak{c}}}}\Big)\Big\}, \frac{2T_\mathfrak{e}}{u_{\mathfrak{e}}^0} D_x f_R^{\varepsilon, \mathfrak{c}}\Big\rangle  \nonumber\\
			& -\varepsilon^k\Big\langle D_x\Big\{\Big(E_R^{\varepsilon, \mathfrak{c}}+\frac{p}{p^0} \times B_R^{\varepsilon, \mathfrak{c}}\Big) \cdot \Big(\frac{p}{p^0}\frac{u_{\mathfrak{e}}^0}{2T_\mathfrak{e}}-\frac{u_{\mathfrak{e}}}{2T_\mathfrak{e}}\Big) f_R^{\varepsilon, \mathfrak{c}}\Big\}, \frac{2T_\mathfrak{e}}{u_{\mathfrak{e}}^0} D_x f_R^{\varepsilon, \mathfrak{c}}\Big\rangle  \nonumber\\
			& -\sum_{i=1}^{2 k-1} \varepsilon^i\Big\langle D_x\Big\{\Big(E_i^{\mathfrak{c}}+\frac{p}{p^0} \times B_i^{\mathfrak{c}}\Big) \cdot \Big(\frac{p}{p^0}\frac{u_{\mathfrak{e}}^0}{2T_\mathfrak{e}}-\frac{u_{\mathfrak{e}}}{2T_\mathfrak{e}}\Big) f_R^{\varepsilon, \mathfrak{c}}\Big\}, \frac{2T_\mathfrak{e}}{u_{\mathfrak{e}}^0} D_x f_R^{\varepsilon, \mathfrak{c}}\Big\rangle \nonumber\\
			& +\Big\langle D_x\Big\{\Big(E^{\varepsilon, \mathfrak{c}}+\frac{p}{p^0}\times B^{\varepsilon, \mathfrak{c}}\Big) \cdot \nabla_p f_R^{\varepsilon, \mathfrak{c}}\Big\}, \frac{2T_\mathfrak{e}}{u_{\mathfrak{e}}^0} D_x f_R^{\varepsilon, \mathfrak{c}}\Big\rangle  \nonumber\\
			& +\sum_{i=1}^{2 k-1} \varepsilon^i\Big\langle D_x\Big\{\Big(E_R^{\varepsilon, \mathfrak{c}}+\frac{p}{p^0} \times B_R^{\varepsilon, \mathfrak{c}}\Big) \cdot \frac{\nabla_p F_i^{\mathfrak{c}}}{\sqrt{\mathbf{M}_{\mathfrak{c}}}}\Big\}, \frac{2T_\mathfrak{e}}{u_{\mathfrak{e}}^0} D_x f_R^{\varepsilon, \mathfrak{c}}\Big\rangle+\varepsilon^k\Big\langle D_x \bar{A}, \frac{2T_\mathfrak{e}}{u_{\mathfrak{e}}^0} D_x f_R^{\varepsilon, \mathfrak{c}}\Big\rangle, 
		\end{align}}
		where $\bar{A}=\frac{A}{\sqrt{\mathbf{M}_{\mathfrak{c}}}}$.
		
		We begin by examining the second term on LHS of \eqref{H.2}. By the expression of the operator $\mathbf{L}_{\mathfrak{c}}$, there holds
		\begin{align}
			D_x \left(\mathbf{L}_{\mathfrak{c}}f_R^{\varepsilon, \mathfrak{c}}\right) 
			=&\mathbf{L}_{\mathfrak{c}}\left(D_x f_R^{\varepsilon, \mathfrak{c}}\right)-\Big\{\Gamma_{\mathfrak{c}}\left(f_R^{\varepsilon, \mathfrak{c}},D_x \sqrt{\mathbf{M}_{\mathfrak{c}}}\right)+\Gamma_{\mathfrak{c}}\left(D_x \sqrt{\mathbf{M}_{\mathfrak{c}}},f_R^{\varepsilon, \mathfrak{c}}\right)\Big\} \nonumber\\
			&-\int_{\mathbb{R}^3} d q \int_{\mathbb{S}^2} d \omega v_\phi(p,q) D_x \sqrt{\mathbf{M}_{\mathfrak{c}}}(q)\Big\{ f_R^{\varepsilon, \mathfrak{c}}\left(p^{\prime}\right) \sqrt{\mathbf{M}_{\mathfrak{c}}}\left(q^{\prime}\right)\nonumber\\
			&+ f_R^{\varepsilon, \mathfrak{c}}\left(q^{\prime}\right) \sqrt{\mathbf{M}_{\mathfrak{c}}}\left(p^{\prime}\right)-f_R^{\varepsilon, \mathfrak{c}}(p) \sqrt{\mathbf{M}_{\mathfrak{c}}}(q)- f_R^{\varepsilon, \mathfrak{c}}(q) \sqrt{\mathbf{M}_{\mathfrak{c}}}(p)\Big\}.\nonumber
		\end{align}
		Noting $|D_x \sqrt{\mathbf{M}_{\mathfrak{c}}}(q)|\leq C(1+|q|)^2 \sqrt{\mathbf{M}_{\mathfrak{c}}}$, 
		by Lemmas \ref{l2.1} and \ref{l2.3}, we can obtain
		\begin{align}
			&\frac{1}{\varepsilon}\Big\langle D_x \left(\mathbf{L}_{\mathfrak{c}}f_R^{\varepsilon, \mathfrak{c}}\right), \frac{2T_\mathfrak{e}}{u_{\mathfrak{e}}^0} D_x f_R^{\varepsilon, \mathfrak{c}}\Big\rangle\nonumber\\
			\geq&\frac{4\zeta_0c_0}{\varepsilon\mathfrak{c}}\left\|\{\mathbf{I}-\mathbf{P}_{\mathfrak{c}}\}D_x f_R^{\varepsilon, \mathfrak{c}}\right\|_{\nu_{\mathfrak{c}}}^2-\frac{1}{2\varepsilon\mathfrak{c} c_0}\left\|D_x f_R^{\varepsilon, \mathfrak{c}}\right\|_{\nu_{\mathfrak{c}}}\left\| f_R^{\varepsilon, \mathfrak{c}}\right\|_{\nu_{\mathfrak{c}}}\nonumber\\
			\geq&\frac{3\zeta_0c_0}{\varepsilon\mathfrak{c}}\left\|\{\mathbf{I}-\mathbf{P}_{\mathfrak{c}}\}D_x f_R^{\varepsilon, \mathfrak{c}}\right\|_{\nu_{\mathfrak{c}}}^2-\frac{C}{\mathfrak{c}}\left\| D_x f_R^{\varepsilon, \mathfrak{c}}\right\|^2-\frac{C}{\varepsilon^2\mathfrak{c}}\left\| f_R^{\varepsilon, \mathfrak{c}}\right\|^2-\frac{C}{\varepsilon^2\mathfrak{c}}\left\|\{\mathbf{I}-\mathbf{P}_{\mathfrak{c}}\}f_R^{\varepsilon, \mathfrak{c}}\right\|_{\nu_{\mathfrak{c}}}^2.\nonumber
		\end{align}
	
		Applying $D_x$ to $\eqref{L.2}_{1,2}$ and taking the $L^2$ inner product with $D_x E_R^{\varepsilon,\mathfrak{c}}$, $D_x B_R^{\varepsilon,\mathfrak{c}}$ respectively, we can obtain 
		\begin{align*}
		\frac{1}{2}\frac{d}{d t}\big\|\big(D_x E_R^{\varepsilon,\mathfrak{c}},D_x B_R^{\varepsilon,\mathfrak{c}}\big)\big\|^2=4\pi\Big\langle D_x E_R^{\varepsilon,\mathfrak{c}},\hat{p}\cdot D_x\left(\sqrt{\mathbf{M}}_{\mathfrak{c}}f_R^{\varepsilon,\mathfrak{c}}\right) \Big\rangle.
		\end{align*}
		Now we estimate the third term on LHS of \eqref{H.2},
		\begin{align}
			&\Big\langle D_x\Big\{\Big(E_R^{\varepsilon, \mathfrak{c}}+\frac{p}{p^0} \times B_R^{\varepsilon, \mathfrak{c}}\Big) \Big(\frac{p}{p^0}\frac{u_{\mathfrak{e}}^0}{T_\mathfrak{e}}-\frac{u_{\mathfrak{e}}}{T_\mathfrak{e}}\Big)\sqrt{\mathbf{M}_{\mathfrak{c}}}\Big\}, \frac{2T_\mathfrak{e}}{u_{\mathfrak{e}}^0} D_x f_R^{\varepsilon, \mathfrak{c}}\Big\rangle\nonumber\\
			=& \frac{1}{4\pi\mathfrak{c}}\frac{d}{d t}\left\|\big(D_x E_R^{\varepsilon, \mathfrak{c}}, D_x B_R^{\varepsilon, \mathfrak{c}}\big)\right\|^2-2\Big\langle D_x E_R^{\varepsilon, \mathfrak{c}},\frac{p}{p^0}D_x\sqrt{\mathbf{M}}_{\mathfrak{c}}f_R^{\varepsilon,\mathfrak{c}}\Big\rangle \nonumber\\
			&+\Big\langle E_R^{\varepsilon, \mathfrak{c}}D_x\Big(\frac{p}{p^0}\frac{u_{\mathfrak{e}}^0}{T_\mathfrak{e}}\sqrt{\mathbf{M}}_{\mathfrak{c}}\Big)-D_x\Big\{\Big(E_R^{\varepsilon, \mathfrak{c}}+\frac{p}{p^0}\times B_R^{\varepsilon, \mathfrak{c}}\Big)\frac{u_{\mathfrak{e}}}{T_\mathfrak{e}}\sqrt{\mathbf{M}_{\mathfrak{c}}}\Big\}, \frac{2T_\mathfrak{e}}{u_{\mathfrak{e}}^0} D_x f_R^{\varepsilon, \mathfrak{c}}\Big\rangle\nonumber\\
			\gtrsim& \frac{1}{\mathfrak{c}}\frac{d}{d t}\left\|\big(D_x E_R^{\varepsilon, \mathfrak{c}}, D_x B_R^{\varepsilon, \mathfrak{c}}\big)\right\|^2-\frac{1}{\mathfrak{c}}\left\|f_R^{\varepsilon, \mathfrak{c}}\right\|_{H^1}^2-\frac{1}{\mathfrak{c}}\left\|\left(E_R^{\varepsilon, \mathfrak{c}}, B_R^{\varepsilon, \mathfrak{c}}\right)\right\|_{H^1}^2.\nonumber
		\end{align}
		
		For the second term on RHS of \eqref{H.2}, we treat it via integration by parts with respect to $x$. A direct calculation yields that
		$$
		D_x^2\left\{\frac{1}{\sqrt{\mathbf{M}_{\mathfrak{c}}}}\Big[\partial_t+\hat{p} \cdot \nabla_x-\big(E_0^{\mathfrak{c}}+\frac{p}{p^0} \times B_0^{\mathfrak{c}}\big) \cdot \nabla_p\Big] \sqrt{\mathbf{M}_{\mathfrak{c}}}\right\}\lesssim (1+|p|)^3,
		$$ 
		where we have used \eqref{L.8}.
		 Then it follows from $(1+|p|)^3\left|f_R^{\varepsilon, \mathfrak{c}}\right| \leq(1+|p|)^{-11}\left|h_R^{\varepsilon, \mathfrak{c}}\right|$ that
		$$
		\Big(\int_{|p| \geq \sqrt[3]{\frac{\kappa}{\varepsilon}}}(1+|p|)^{-11 \times 2} d p\Big)^{1 / 2} \lesssim\left(\frac{\varepsilon}{\kappa}\right)^{\frac{19}{6}}.
		$$
		Consequently, similar to the estimation of \eqref{L.9}, there holds
		\begin{align}
		&\left|\Big\langle f_R^{\varepsilon, \mathfrak{c}} D_x\left\{\frac{1}{\sqrt{\mathbf{M}_{\mathfrak{c}}}}\Big[\partial_t+\frac{p}{p^0} \cdot \nabla_x-\big(E_0^{\mathfrak{c}}+\frac{p}{p^0} \times B_0^{\mathfrak{c}}\big) \cdot \nabla_p\Big] \sqrt{\mathbf{M}_{\mathfrak{c}}}\right\}, \frac{2T_\mathfrak{e}}{u_{\mathfrak{e}}^0} D_x f_R^{\varepsilon, \mathfrak{c}}\Big\rangle\right|\nonumber \\
		\lesssim&\frac{C}{\mathfrak{c}}\left\|f_R^{\varepsilon, \mathfrak{c}}\right\|^2+\frac{C\kappa}{\varepsilon\mathfrak{c}}\left\|\{\mathbf{I}-\mathbf{P}_{\mathfrak{c}}\} f_R^{\varepsilon, \mathfrak{c}}\right\|_{\nu_{\mathfrak{c}}}^2+\frac{C_\kappa}{\mathfrak{c}}\varepsilon^{\frac{19}{6}}\left\|h_R^{\varepsilon, \mathfrak{c}}\right\|_{L^{\infty}}\left\|f_R^{\varepsilon, \mathfrak{c}}\right\|.
		\end{align}
		
		Similarly,
    \begin{align*}
    & \Big|\Big\langle\frac{D_x f_R^{\varepsilon, \mathfrak{c}}}{\sqrt{\mathbf{M}_{\mathfrak{c}}}}\Big\{\partial_t+\hat{p} \cdot \nabla_x-\big(E_0^{\mathfrak{c}}+\frac{p}{p^0} \times B_0^{\mathfrak{c}}\big) \cdot \nabla_p\Big\} \sqrt{\mathbf{M}_{\mathfrak{c}}}, \frac{2T_\mathfrak{e}}{u_{\mathfrak{e}}^0} D_x f_R^{\varepsilon, \mathfrak{c}}\Big\rangle \Big| \\
    &\quad+\Big\langle\Big\{\left(\partial_t+\hat{p} \cdot \nabla_x\right)\big(\frac{T_\mathfrak{e}}{u_{\mathfrak{e}}^0}\big)\Big\} D_x f_R^{\varepsilon, \mathfrak{c}}, D_x f_R^{\varepsilon, \mathfrak{c}}\Big\rangle\\
    \lesssim& \frac{C_{\kappa}}{\mathfrak{c}} \varepsilon^{\frac{19}{6}}\left\|h_R^{\varepsilon, \mathfrak{c}}\right\|_{W^{1, \infty}} \cdot\left\|D_x f_R^{\varepsilon, \mathfrak{c}}\right\|+\frac{C}{\mathfrak{c}}\left\|D_x f_R^{\varepsilon, \mathfrak{c}}\right\|^2+\frac{C \kappa}{\varepsilon\mathfrak{c}} \|\left\{\mathbf{I}-\mathbf{P}_{\mathfrak{c}}\right\}D_x f_R^{\varepsilon, \mathfrak{c}} \|_{\nu_{\mathfrak{c}}}^2,\nonumber
	\end{align*}
	and
    \begin{align*}
        &\Big|\sum_{i=1}^{2 k-1} \varepsilon^i\Big\langle D_x\Big\{\Big(E_i^{\mathfrak{c}}+\frac{p}{p^0} \times B_i^{\mathfrak{c}}\Big) \cdot \Big(\frac{p}{p^0}\frac{u_{\mathfrak{e}}^0}{2T_\mathfrak{e}}-\frac{u_{\mathfrak{e}}}{2T_\mathfrak{e}}\Big) f_R^{\varepsilon, \mathfrak{c}}\Big\}, \frac{2T_\mathfrak{e}}{u_{\mathfrak{e}}^0} D_x f_R^{\varepsilon, \mathfrak{c}}\Big\rangle\Big|\\
        \lesssim& \frac{C_\kappa}{\mathfrak{c}}\varepsilon^{\frac{25}{6}}\left\|h_R^{\varepsilon, \mathfrak{c}}\right\|_{L^{\infty}}\left\|f_R^{\varepsilon, \mathfrak{c}}\right\|+\frac{C \varepsilon}{\mathfrak{c}}\left\|f_R^{\varepsilon, \mathfrak{c}}\right\|^2+\frac{C\kappa}{\mathfrak{c}}\left\|\{\mathbf{I}-\mathbf{P}_{\mathfrak{c}}\} f_R^{\varepsilon, \mathfrak{c}}\right\|_{\nu_{\mathfrak{c}}}^2\\
        &+\frac{C_{\kappa}}{\mathfrak{c}} \varepsilon^{\frac{25}{6}}\left\|h_R^{\varepsilon, \mathfrak{c}}\right\|_{W^{1, \infty}} \cdot\left\|D_x f_R^{\varepsilon, \mathfrak{c}}\right\|+\frac{C \varepsilon}{\mathfrak{c}}\left\|D_x f_R^{\varepsilon, \mathfrak{c}}\right\|^2+\frac{C \kappa}{\mathfrak{c}} \|\left\{\mathbf{I}-\mathbf{P}_{\mathfrak{c}}\right\}D_x f_R^{\varepsilon, \mathfrak{c}} \|_{\nu_{\mathfrak{c}}}^2.
    \end{align*}
		Using Lemma \ref{l2.3}, we can obtain
        \begin{align*}
        \varepsilon^{k-1}\Big\langle D_x \big(\Gamma_{\mathfrak{c}}\left(f_R^{\varepsilon, \mathfrak{c}}, f_R^{\varepsilon, \mathfrak{c}}\right)\big), \frac{2 T_{\mathfrak{e}}}{u_{\mathfrak{e}}^{0, \mathfrak{c}}} D_x f_R^{\varepsilon, \mathfrak{c}}\Big\rangle \lesssim \frac{\varepsilon^{k-1}}{\mathfrak{c}}\left\|h_R^{\varepsilon, \mathfrak{c}}\right\|_{W^{1, \infty}}\Big\{\left\|f_R^{\varepsilon, \mathfrak{c}}\right\|\left\|D_x f_R^{\varepsilon, \mathfrak{c}}\right\|+\left\|f_R^{\varepsilon, \mathfrak{c}}\right\|^2\Big\},
         \end{align*}
    and
    \begin{align*}
    & \sum_{i=1}^{2 k-1} \varepsilon^{i-1}\Big\langle D_x\Big\{\Gamma_{\mathfrak{c}}\Big(\frac{F_i^{\mathfrak{c}}}{\sqrt{\mathbf{M}_{\mathfrak{c}}}}, f_R^{\varepsilon, \mathfrak{c}}\Big)+\Gamma_{\mathfrak{c}}\Big(f_R^{\varepsilon, \mathfrak{c}}, \frac{F_i^{\mathfrak{c}}}{\sqrt{\mathbf{M}_{\mathfrak{c}}}}\Big)\Big\}, \frac{2 T_{\mathfrak{e}}}{u_{\mathfrak{e}}^{0, \mathfrak{c}}} D_x f_R^{\varepsilon, \mathfrak{c}}\Big\rangle \\
    & \lesssim \sum_{i=1}^{2 k-1} \frac{\varepsilon^{i-1}}{\mathfrak{c}}\Big\{\left\|f_R^{\varepsilon, \mathfrak{c}}\right\|_{\nu_{\mathfrak{c}}}\left\|D_x f_R^{\varepsilon, \mathfrak{c}}\right\|_{\nu_{\mathfrak{c}}}+\left\|D_x f_R^{\varepsilon, \mathfrak{c}}\right\|_{\nu_{\mathfrak{c}}}^2\Big\} \\
    & \lesssim \frac{1}{\mathfrak{c}}\Big\{\left\|f_R^{\varepsilon, \mathfrak{c}}\right\|_{H^1}^2+\left\|\left\{\mathbf{I}-\mathbf{P}_{\mathfrak{c}}\right\} f_R^{\varepsilon, \mathfrak{c}}\right\|_{\nu_{\mathfrak{c}}}^2+\left\|\left\{\mathbf{I}-\mathbf{P}_{\mathfrak{c}}\right\} D_x f_R^{\varepsilon, \mathfrak{c}}\right\|_{\nu_{\mathfrak{c}}}^2\Big\}.
\end{align*}

    The 6th, 8th, 9th and 10th terms on the right-hand side of \eqref{H.2} can be estimated similarly by using the boundedness of the coefficients and Cauchy--Schwarz's inequality, so that they are bounded by
    \begin{align*}
     \frac{1}{\mathfrak{c}}\left\|f_R^{\varepsilon, \mathfrak{c}}\right\|_{H^1}^2+\frac{1}{\mathfrak{c}}\left\|\left(E_R^{\varepsilon, \mathfrak{c}}, B_R^{\varepsilon, \mathfrak{c}}\right)\right\|_{H^1}^2+\frac{1}{\mathfrak{c}}\left\|\nabla_x f_R^{\varepsilon, \mathfrak{c}}\right\|\left\|\nabla_pf_R^{\varepsilon, \mathfrak{c}}\right\|+\frac{\varepsilon^k}{\mathfrak{c}} \left\|\nabla_x f_R^{\varepsilon, \mathfrak{c}}\right\|.
    \end{align*}
    
	Collecting all the above estimates, one has
		\begin{align*}
			\frac{1}{2}&\frac{d}{d t}\Big\|  \sqrt{\frac{2\mathfrak{c}T_\mathfrak{e}}{u_{\mathfrak{e}}^0}} \nabla_x f_R^{\varepsilon, \mathfrak{c}}\Big\|^2+\frac{d}{d t}\left\|\big(\nabla_x E_R^{\varepsilon, \mathfrak{c}}, \nabla_x B_R^{\varepsilon, \mathfrak{c}}\big) \right\|^2+\Big(\frac{3\zeta_0 c_0}{\varepsilon}-\frac{C \kappa}{\varepsilon}-C\Big)\left\|\{\mathbf{I}-\mathbf{P}_{\mathfrak{c}}\}\nabla_x f_R^{\varepsilon, \mathfrak{c}}\right\|_{\nu_{\mathfrak{c}}}^2 \nonumber\\
			& \lesssim\Big\{1+\varepsilon^{k-1}\left\|h_R^{\varepsilon, \mathfrak{c}}\right\|_{W^{1,\infty}}\Big\}\Big\{\left\|f_R^{\varepsilon, \mathfrak{c}}\right\|_{H^1}^2+\left\|\big(E_R^{\varepsilon, \mathfrak{c}}, B_R^{\varepsilon, \mathfrak{c}}\big)\right\|_{H^1}^2\Big\}+\frac{C}{\varepsilon^2}\left\| f_R^{\varepsilon, \mathfrak{c}}\right\|^2 \nonumber\\
			& \quad+C_\kappa\varepsilon^{\frac{19}{6}}\left\|h_R^{\varepsilon, \mathfrak{c}}\right\|_{L^{\infty}}\left\|f_R^{\varepsilon, \mathfrak{c}}\right\|+\Big\{C_\kappa\varepsilon^{\frac{19}{6}}\left\|h_R^{\varepsilon, \mathfrak{c}}\right\|_{W^{1,\infty}}+\varepsilon^{k}\Big\}\left\|\nabla_x f_R^{\varepsilon, \mathfrak{c}}\right\| \nonumber\\
			&\quad+C\Big(\frac{1}{\varepsilon^2}+\frac{\kappa}{\varepsilon}+1\Big)\left\|\{\mathbf{I}-\mathbf{P}_{\mathfrak{c}}\}f_R^{\varepsilon, \mathfrak{c}}\right\|_{\nu_{\mathfrak{c}}}^2+\left\|\nabla_x f_R^{\varepsilon, \mathfrak{c}}\right\|\left\|\nabla_pf_R^{\varepsilon, \mathfrak{c}}\right\|.
		\end{align*}
		First taking $\kappa$ sufficiently small and then taking $\varepsilon_0$ sufficiently small, thus one concludes \eqref{H.1}.
		Therefore the proof of Lemma \ref{pH.1} is completed.
	\end{proof}
	Next we deal with the momentum derivative estimate for $ f_R^{\varepsilon, \mathfrak{c}}$. 
	\begin{lemma}\label{pH.2}
		There exists constants $\varepsilon_0>0$ and $C>0$, such that for all $\varepsilon \in (0, \varepsilon_0]$,
	\begin{align}\label{H.5}
    \frac{1}{2} \frac{d}{d t}&\left\|\nabla_p f_R^{\varepsilon, \mathfrak{c}}\right\|^2+\frac{\zeta_0}{4 \varepsilon}\left\|\left\{\mathbf{I}-\mathbf{P}_{\mathfrak{c}}\right\} \nabla_p f_R^{\varepsilon, \mathfrak{c}}\right\|_{\nu_{\mathfrak{c}}}^2 \nonumber\\
    \leq & C\Big\{1+\varepsilon^{k-1}\left\|h_R^{\varepsilon, \mathfrak{c}}\right\|_{W^{1, \infty}}\Big\}\left\|f_R^{\varepsilon, \mathfrak{c}}\right\|_{H^1}^2+C\left\|\big(E_R^{\varepsilon, \mathfrak{c}}, B_R^{\varepsilon, \mathfrak{c}}\big)\right\|^2 \nonumber\\
    & +\frac{C}{\varepsilon^2}\left\|\left\{\mathbf{I}-\mathbf{P}_{\mathfrak{c}}\right\} f_R^{\varepsilon, \mathfrak{c}}\right\|_{\nu_{\mathfrak{c}}}^2+\frac{C}{\varepsilon^2}\left\|f_R^{\varepsilon, \mathfrak{c}}\right\|^2+C\left\|\nabla_x f_R^{\varepsilon, \mathfrak{c}}\right\|\left\|\nabla_p f_R^{\varepsilon, \mathfrak{c}}\right\|\nonumber\\
    &+C \varepsilon^{\frac{23}{2}}\left\|h_R^{\varepsilon, \mathfrak{c}}\right\|_{L^{\infty}}\left\|f_R^{\varepsilon, \mathfrak{c}}\right\|+C\left( \varepsilon^{\frac{19}{6}}\left\|h_R^{\varepsilon, \mathfrak{c}}\right\|_{W^{1, \infty}}+\varepsilon^k\right)\left\|\nabla_p f_R^{\varepsilon, \mathfrak{c}}\right\|,
	\end{align}
	where the constant $C>0$ depends on $ F_i^{\mathfrak{c}}$, $E_i^{\mathfrak{c}}$, $B_i^{\mathfrak{c}}(i=0,1,\cdots,2k-1)$, and is independent of $\mathfrak{c}$.
	\end{lemma}
	
	\begin{proof}
		Applying $D_p$ to \eqref{L.1}, then multiplying the resultant equation by $ D_p f_R^{\varepsilon, \mathfrak{c}}$ and integrating over $\mathbb{R}^3\times \mathbb{R}^3$, one has
		{\footnotesize
		\begin{align}\label{H.6}
		 	&\frac{1}{2} \frac{d}{d t}\left\|D_p f_R^{\varepsilon, \mathfrak{c}}\right\|^2+\frac{1}{\varepsilon}\Big\langle D_p \left(\mathbf{L}_{\mathfrak{c}}f_R^{\varepsilon, \mathfrak{c}}\right), D_p f_R^{\varepsilon, \mathfrak{c}}\Big\rangle \nonumber\\
			=&-\Big\langle D_p(\hat{p})\cdot\nabla_{x}f_R^{\varepsilon, \mathfrak{c}}, D_p f_R^{\varepsilon, \mathfrak{c}}\Big\rangle -\Big\langle D_p\Big\{\Big(E_R^{\varepsilon, \mathfrak{c}}+\frac{p}{p^0} \times B_R^{\varepsilon, \mathfrak{c}}\Big) \Big(\frac{p}{p^0}\frac{u_{\mathfrak{e}}^0}{T_\mathfrak{e}}-\frac{u_{\mathfrak{e}}}{T_\mathfrak{e}}\Big)\sqrt{\mathbf{M}_{\mathfrak{c}}}\Big\}, D_p f_R^{\varepsilon, \mathfrak{c}}\Big\rangle \nonumber\\
			&-\Big\langle\frac{D_p f_R^{\varepsilon, \mathfrak{c}}}{\sqrt{\mathbf{M}_{\mathfrak{c}}}}\Big\{\partial_t+\hat{p} \cdot \nabla_x-\big(E_0^{\mathfrak{c}}+\frac{p}{p^0} \times B_0^{\mathfrak{c}}\big) \cdot \nabla_p\Big\} \sqrt{\mathbf{M}_{\mathfrak{c}}},  D_p f_R^{\varepsilon, \mathfrak{c}}\Big\rangle  \nonumber\\
			& -\Big\langle f_R^{\varepsilon, \mathfrak{c}} D_p\left\{\frac{1}{\sqrt{\mathbf{M}_{\mathfrak{c}}}}\Big\{\partial_t+\hat{p}\cdot \nabla_x-\big(E_0^{\mathfrak{c}}+\frac{p}{p^0} \times B_0^{\mathfrak{c}}\big) \cdot \nabla_p\Big\} \sqrt{\mathbf{M}_{\mathfrak{c}}}\right\}, D_p f_R^{\varepsilon, \mathfrak{c}}\Big\rangle  \nonumber\\
			&+\varepsilon^{k-1}\Big\langle D_p \Gamma_{\mathfrak{c}}\left(f_R^{\varepsilon, \mathfrak{c}}, f_R^{\varepsilon, \mathfrak{c}}\right), D_p f_R^{\varepsilon, \mathfrak{c}}\Big\rangle +\sum_{i=1}^{2 k-1} \varepsilon^{i-1}\Big\langle D_p\Big\{\Gamma_{\mathfrak{c}}\Big(\frac{F_i^{\mathfrak{c}}}{\sqrt{\mathbf{M}_{\mathfrak{c}}}}, f_R^{\varepsilon, \mathfrak{c}}\Big)+\Gamma_{\mathfrak{c}}\Big(f_R^{\varepsilon, \mathfrak{c}}, \frac{F_i^{\mathfrak{c}}}{\sqrt{\mathbf{M}_{\mathfrak{c}}}}\Big)\Big\}, D_p f_R^{\varepsilon, \mathfrak{c}}\Big\rangle \nonumber\\
			& -\varepsilon^k\Big\langle D_p\Big\{\Big(E_R^{\varepsilon, \mathfrak{c}}+\frac{p}{p^0} \times B_R^{\varepsilon, \mathfrak{c}}\Big) \cdot \Big(\frac{p}{p^0}\frac{u_{\mathfrak{e}}^0}{2T_\mathfrak{e}}-\frac{u_{\mathfrak{e}}}{2T_\mathfrak{e}}\Big) f_R^{\varepsilon, \mathfrak{c}}\Big\}, D_p f_R^{\varepsilon, \mathfrak{c}}\Big\rangle \nonumber\\
			& -\sum_{i=1}^{2 k-1} \varepsilon^i\Big\langle D_p\Big\{\Big(E_i^{\mathfrak{c}}+\frac{p}{p^0} \times B_i^{\mathfrak{c}}\Big) \cdot \Big(\frac{p}{p^0}\frac{u_{\mathfrak{e}}^0}{2T_\mathfrak{e}}-\frac{u_{\mathfrak{e}}}{2T_\mathfrak{e}}\Big) f_R^{\varepsilon, \mathfrak{c}}\Big\}, D_p f_R^{\varepsilon, \mathfrak{c}}\Big\rangle\nonumber\\
			& +\Big\langle D_p\Big\{\Big(E^{\varepsilon, \mathfrak{c}}+\frac{p}{p^0}\times B^{\varepsilon, \mathfrak{c}}\Big) \cdot \nabla_p f_R^{\varepsilon, \mathfrak{c}}\Big\}, D_p f_R^{\varepsilon, \mathfrak{c}}\Big\rangle  \nonumber\\
			& +\sum_{i=1}^{2 k-1} \varepsilon^i\Big\langle D_p\Big\{\Big(E_R^{\varepsilon, \mathfrak{c}}+\frac{p}{p^0} \times B_R^{\varepsilon, \mathfrak{c}}\Big) \cdot \frac{\nabla_p F_i^{\mathfrak{c}}}{\sqrt{\mathbf{M}_{\mathfrak{c}}}}\Big\}, D_p f_R^{\varepsilon, \mathfrak{c}}\Big\rangle+\varepsilon^k\Big\langle D_p \bar{A}, D_p f_R^{\varepsilon, \mathfrak{c}}\Big\rangle,
		\end{align}}
		where $\bar{A}=\frac{A}{\sqrt{\mathbf{M}_{\mathfrak{c}}}}$.
		
		By the expression of the operator $\mathbf{L}_{\mathfrak{c}}$, there holds
    	\begin{align}
			D_p \left(\mathbf{L}_{\mathfrak{c}}f_R^{\varepsilon, \mathfrak{c}}\right) 
			=&\mathbf{L}_{\mathfrak{c}}\left(D_p f_R^{\varepsilon, \mathfrak{c}}\right)-\Big\{\Gamma_{\mathfrak{c}}\left(f_R^{\varepsilon, \mathfrak{c}},D_p \sqrt{\mathbf{M}_{\mathfrak{c}}}\right)+\Gamma_{\mathfrak{c}}\left(D_p \sqrt{\mathbf{M}_{\mathfrak{c}}},f_R^{\varepsilon, \mathfrak{c}}\right)\Big\} \nonumber\\
			&-\int_{\mathbb{R}^3} d q \int_{\mathbb{S}^2} d \omega v_\phi(p,q) D_p \sqrt{\mathbf{M}_{\mathfrak{c}}}(q)\Big\{ f_R^{\varepsilon, \mathfrak{c}}\left(p^{\prime}\right) \sqrt{\mathbf{M}_{\mathfrak{c}}}\left(q^{\prime}\right)\nonumber\\
			&+ f_R^{\varepsilon, \mathfrak{c}}\left(q^{\prime}\right) \sqrt{\mathbf{M}_{\mathfrak{c}}}\left(p^{\prime}\right)-f_R^{\varepsilon, \mathfrak{c}}(p) \sqrt{\mathbf{M}_{\mathfrak{c}}}(q)- f_R^{\varepsilon, \mathfrak{c}}(q) \sqrt{\mathbf{M}_{\mathfrak{c}}}(p)\Big\}.\nonumber
		\end{align}
    By Lemmas \ref{l2.1} and \ref{l2.4}, we get
		\begin{align}
			&\frac{1}{\varepsilon}\Big\langle D_p \left(\mathbf{L}_{\mathfrak{c}}f_R^{\varepsilon, \mathfrak{c}}\right), D_p f_R^{\varepsilon, \mathfrak{c}}\Big\rangle\nonumber\\
			\geq& \frac{ \zeta_0 }{\varepsilon}\left\|\{\mathbf{I}-\mathbf{P}_{\mathfrak{c}}\}D_p f_R^{\varepsilon, \mathfrak{c}}\right\|_{\nu_{\mathfrak{c}}}^2-\frac{C}{\varepsilon}\left\|D_pf_R^{\varepsilon, \mathfrak{c}}\right\|_{\nu_{\mathfrak{c}}}\left\| f_R^{\varepsilon, \mathfrak{c}}\right\|_{\nu_{\mathfrak{c}}}\nonumber\\
			\geq&\frac{ \zeta_0}{2\varepsilon}\left\|\{\mathbf{I}-\mathbf{P}_{\mathfrak{c}}\}D_p f_R^{\varepsilon, \mathfrak{c}}\right\|_{\nu_{\mathfrak{c}}}^2-\frac{C}{2\varepsilon^2}\left\|\{\mathbf{I}-\mathbf{P}_{\mathfrak{c}}\}f_R^{\varepsilon, \mathfrak{c}}\right\|_{\nu_{\mathfrak{c}}}^2-\frac{C}{2}\left\|D_pf_R^{\varepsilon, \mathfrak{c}}\right\|^2-\frac{C}{2\varepsilon^2}\left\| f_R^{\varepsilon, \mathfrak{c}}\right\|^2.\nonumber
		\end{align}
		
	Now we turn to the first term on RHS of \eqref{H.6}. It's clear that
		$$
		\Big\langle D_p(\hat{p})\cdot\nabla_{x}f_R^{\varepsilon, \mathfrak{c}}, D_p f_R^{\varepsilon, \mathfrak{c}}\Big\rangle \leq\left\|\nabla_{x}f_R^{\varepsilon, \mathfrak{c}}\right\|\left\|\nabla_p f_R^{\varepsilon, \mathfrak{c}}\right\|.
		$$
	For the fourth term on RHS of \eqref{H.6}, we treat it via integration by parts with respect to $p$. Note that 
	$$
	\Big|D_p^2\Big\{\frac{1}{\sqrt{\mathbf{M}_{\mathfrak{c}}}}\Big[\partial_t+\hat{p} \cdot \nabla_x-\big(E_0^{\mathfrak{c}}+\frac{p}{p^0} \times B_0^{\mathfrak{c}}\big) \cdot \nabla_p\Big] \sqrt{\mathbf{M}_{\mathfrak{c}}}\Big\}\Big| \lesssim 1+|p|,
	$$
	then it follows from $(1+|p|)\left|f_R^{\varepsilon, \mathfrak{c}}\right| \leq(1+|p|)^{-13}\left|h_R^{\varepsilon, \mathfrak{c}}\right|$ that
    $$
    \Big(\int_{1+|p| \geq \frac{\kappa}{\varepsilon}}(1+|p|)^{-13 \times 2} d p\Big)^{1 / 2} \lesssim\left(\frac{\varepsilon}{\kappa}\right)^{\frac{23}{2}}.
    $$
    Consequently, similar to the estimation of \eqref{L.9}, 
    \begin{align*}
    &\Big|\Big\langle f_R^{\varepsilon, \mathfrak{c}} D_p\left\{\frac{1}{\sqrt{\mathbf{M}_{\mathfrak{c}}}}\Big\{\partial_t+\hat{p}\cdot \nabla_x-\big(E_0^{\mathfrak{c}}+\frac{p}{p^0} \times B_0^{\mathfrak{c}}\big) \cdot \nabla_p\Big\} \sqrt{\mathbf{M}_{\mathfrak{c}}}\right\}, D_p f_R^{\varepsilon, \mathfrak{c}}\Big\rangle\Big| \\
    \leq &C\left\|f_R^{\varepsilon, \mathfrak{c}}\right\|^2+\frac{C\kappa}{\varepsilon}\left\|\left\{\mathbf{I}-\mathbf{P}_{\mathfrak{c}}\right\} f_R^{\varepsilon, \mathfrak{c}}\right\|_{\nu_{\mathfrak{c}}}^2+C_\kappa \varepsilon^{\frac{23}{2}}\left\|h_R^{\varepsilon, \mathfrak{c}}\right\|_{L^{\infty}}\left\|f_R^{\varepsilon, \mathfrak{c}}\right\|,
    \end{align*}
    where $C$ and $C_{\kappa}$ are independent of $\mathfrak{c}$. 
    
    Similarly,
    \begin{align*}
    & \Big|\Big\langle\frac{D_p f_R^{\varepsilon, \mathfrak{c}}}{\sqrt{\mathbf{M}_{\mathfrak{c}}}}\Big\{\partial_t+\hat{p} \cdot \nabla_x-\big(E_0^{\mathfrak{c}}+\frac{p}{p^0} \times B_0^{\mathfrak{c}}\big) \cdot \nabla_p\Big\} \sqrt{\mathbf{M}_{\mathfrak{c}}},  D_p f_R^{\varepsilon, \mathfrak{c}}\Big\rangle \Big| \\
    \leq& C_{\kappa} \varepsilon^{\frac{19}{6}}\left\|h_R^{\varepsilon, \mathfrak{c}}\right\|_{W^{1, \infty}} \cdot\left\|D_p f_R^{\varepsilon, \mathfrak{c}}\right\|+C\left\|D_p f_R^{\varepsilon, \mathfrak{c}}\right\|^2+\frac{C \kappa}{\varepsilon} \|\left\{\mathbf{I}-\mathbf{P}_{\mathfrak{c}}\right\}D_p f_R^{\varepsilon, \mathfrak{c}} \|_{\nu_{\mathfrak{c}}}^2,\nonumber
	\end{align*}
		and
    \begin{align*}
        &\Big|\sum_{i=1}^{2 k-1} \varepsilon^i\Big\langle D_p\Big\{\Big(E_i^{\mathfrak{c}}+\frac{p}{p^0} \times B_i^{\mathfrak{c}}\Big) \cdot \Big(\frac{p}{p^0}\frac{u_{\mathfrak{e}}^0}{2T_\mathfrak{e}}-\frac{u_{\mathfrak{e}}}{2T_\mathfrak{e}}\Big) f_R^{\varepsilon, \mathfrak{c}}\Big\}, D_p f_R^{\varepsilon, \mathfrak{c}}\Big\rangle\Big|\\
        \lesssim& C_\kappa \varepsilon^{\frac{25}{6}}\left\|h_R^{\varepsilon, \mathfrak{c}}\right\|_{L^{\infty}}\left\|f_R^{\varepsilon, \mathfrak{c}}\right\|+C \varepsilon \left\|f_R^{\varepsilon, \mathfrak{c}}\right\|^2+C\kappa\left\|\{\mathbf{I}-\mathbf{P}_{\mathfrak{c}}\} f_R^{\varepsilon, \mathfrak{c}}\right\|_{\nu_{\mathfrak{c}}}^2\\
        &+C_{\kappa} \varepsilon^{\frac{25}{6}}\left\|h_R^{\varepsilon, \mathfrak{c}}\right\|_{W^{1, \infty}} \cdot\left\|D_p f_R^{\varepsilon, \mathfrak{c}}\right\|+C \varepsilon \left\|D_p f_R^{\varepsilon, \mathfrak{c}}\right\|^2+C \kappa\|\left\{\mathbf{I}-\mathbf{P}_{\mathfrak{c}}\right\}D_p f_R^{\varepsilon, \mathfrak{c}} \|_{\nu_{\mathfrak{c}}}^2.
    \end{align*}
    
    By \eqref{e2.7} in Lemma \ref{l2.4}, we can obtain
    \begin{align*}
    \varepsilon^{k-1}\Big\langle D_p \Gamma_{\mathfrak{c}}\left(f_R^{\varepsilon, \mathfrak{c}}, f_R^{\varepsilon, \mathfrak{c}}\right), D_p f_R^{\varepsilon, \mathfrak{c}}\Big\rangle \lesssim \varepsilon^{k-1}\left\|h_R^{\varepsilon, \mathfrak{c}}\right\|_{W^{1, \infty}}\Big\{\left\|f_R^{\varepsilon, \mathfrak{c}}\right\|\left\|D_p f_R^{\varepsilon, \mathfrak{c}}\right\|+\left\|f_R^{\varepsilon, \mathfrak{c}}\right\|^2\Big\},
    \end{align*}
    and
    \begin{align*}
    & \sum_{i=1}^{2 k-1} \varepsilon^{i-1}\Big\langle D_p\Big\{\Gamma_{\mathfrak{c}}\Big(\frac{F_i^{\mathfrak{c}}}{\sqrt{\mathbf{M}_{\mathfrak{c}}}}, f_R^{\varepsilon, \mathfrak{c}}\Big)+\Gamma_{\mathfrak{c}}\Big(f_R^{\varepsilon, \mathfrak{c}}, \frac{F_i^{\mathfrak{c}}}{\sqrt{\mathbf{M}_{\mathfrak{c}}}}\Big)\Big\}, D_p f_R^{\varepsilon, \mathfrak{c}}\Big\rangle \\
    & \lesssim \sum_{i=1}^{2 k-1} \varepsilon^{i-1}\Big\{\left\|f_R^{\varepsilon, \mathfrak{c}}\right\|_{\nu_{\mathfrak{c}}}\left\|D_p f_R^{\varepsilon, \mathfrak{c}}\right\|_{\nu_{\mathfrak{c}}}+\left\|D_p f_R^{\varepsilon, \mathfrak{c}}\right\|_{\nu_{\mathfrak{c}}}^2\Big\} . \\
    & \lesssim\left\|f_R^{\varepsilon, \mathfrak{c}}\right\|_{H^1}^2+\left\|\left\{\mathbf{I}-\mathbf{P}_{\mathfrak{c}}\right\} f_R^{\varepsilon, \mathfrak{c}}\right\|_{\nu_{\mathfrak{c}}}^2+\left\|\left\{\mathbf{I}-\mathbf{P}_{\mathfrak{c}}\right\} D_p f_R^{\varepsilon, \mathfrak{c}}\right\|_{\nu_{\mathfrak{c}}}^2 .
    \end{align*}
    
    The 2nd, 7th, 8th, 9th, 10th and 11th terms on the right-hand side of \eqref{H.6} can be estimated similarly by using the boundedness of the coefficients and Cauchy--Schwarz's inequality, so that they are bounded by
    \begin{align*}
     \left\|\big(E_R^{\varepsilon, \mathfrak{c}}, B_R^{\varepsilon, \mathfrak{c}}\big)\right\|^2+\left\|f_R^{\varepsilon, \mathfrak{c}}\right\|^2+\left\|\nabla_p f_R^{\varepsilon, \mathfrak{c}}\right\|^2+\varepsilon^k \left\|\nabla_p f_R^{\varepsilon, \mathfrak{c}}\right\|.
    \end{align*}
    
    Collecting all the above estimates, one has
	\begin{align*}
    \frac{1}{2} \frac{d}{d t}&\left\|\nabla_p f_R^{\varepsilon, \mathfrak{c}}\right\|^2+\left(\frac{\zeta_0}{2 \varepsilon}-\frac{C \kappa}{\varepsilon}-C\right)\left\|\left\{\mathbf{I}-\mathbf{P}_{\mathfrak{c}}\right\} \nabla_p f_R^{\varepsilon, \mathfrak{c}}\right\|_{\nu_{\mathfrak{c}}}^2 \nonumber\\
    \lesssim & \Big\{1+\varepsilon^{k-1}\left\|h_R^{\varepsilon, \mathfrak{c}}\right\|_{W^{1, \infty}}\Big\}\left\|f_R^{\varepsilon, \mathfrak{c}}\right\|_{H^1}^2+\left\|\big(E_R^{\varepsilon, \mathfrak{c}}, B_R^{\varepsilon, \mathfrak{c}}\big)\right\|^2 \nonumber\\
    & +\left(\frac{C}{\varepsilon^2}+\frac{C \kappa}{\varepsilon}+C\right)\left\|\left\{\mathbf{I}-\mathbf{P}_{\mathfrak{c}}\right\} f_R^{\varepsilon, \mathfrak{c}}\right\|_{\nu_{\mathfrak{c}}}^2+\frac{C}{\varepsilon^2}\left\|f_R^{\varepsilon, \mathfrak{c}}\right\|^2+\left\|\nabla_x f_R^{\varepsilon, \mathfrak{c}}\right\|\left\|\nabla_p f_R^{\varepsilon, \mathfrak{c}}\right\|\nonumber\\
    &+C_\kappa \varepsilon^{\frac{23}{2}}\left\|h_R^{\varepsilon, \mathfrak{c}}\right\|_{L^{\infty}}\left\|f_R^{\varepsilon, \mathfrak{c}}\right\|+\left(C_\kappa \varepsilon^{\frac{19}{6}}\left\|h_R^{\varepsilon, \mathfrak{c}}\right\|_{W^{1, \infty}}+\varepsilon^k\right)\left\|\nabla_p f_R^{\varepsilon, \mathfrak{c}}\right\|,
	\end{align*}
	First taking $\kappa$ sufficiently small and then taking $\varepsilon_0$ sufficiently small, thus one gets \eqref{H.5}. Therefore the proof of Lemma \ref{pH.2} is completed.
	\end{proof}
	


	\section{Proof of Theorems \ref{t1.1} and \ref{t1.5}}\label{section 9}
	
	This section is devoted to the proof of Theorem \ref{t1.1} and Theorem \ref{t1.5}.
	\begin{proof}[Proof of Theorem \ref{t1.1}]
		Using \eqref{W.0} and Young's inequality, one has from \eqref{L.3} that
		\begin{align}\label{U.1}
			\frac{d}{d t}&\Big\{\frac{1}{2}\Big\|\sqrt{\frac{2\mathfrak{c}T_\mathfrak{e}}{u_{\mathfrak{e}}^0}}f_R^{\varepsilon,\mathfrak{c}}\Big\|^2+\left\|\big(E_R^{\varepsilon,\mathfrak{c}}, B_R^{\varepsilon,\mathfrak{c}}\big)\right\|^2\Big\}+\frac{\zeta_0c_0}{\varepsilon}\|\{\mathbf{I}-\mathbf{P}_{\mathfrak{c}}\}f_R^{\varepsilon,\mathfrak{c}}\|_{\nu_{\mathfrak{c}}}^2\nonumber\\
			&\lesssim \left\{\left\|f_R^{\varepsilon,\mathfrak{c}}\right\|^2+\left\|\big(E_R^{\varepsilon,\mathfrak{c}}, B_R^{\varepsilon,\mathfrak{c}}\big)\right\|^2+1\right\}.
		\end{align}
		By \eqref{e.30}, we have $c_0<\frac{2\mathfrak{c}T_\mathfrak{e}}{u_{\mathfrak{e}}^0}< \frac{2}{c_0}$, which implies $\big|\sqrt{\frac{2\mathfrak{c}T_\mathfrak{e}}{u_{\mathfrak{e}}^0}}f_R^{\varepsilon,\mathfrak{c}}\big|^2\cong |f_R^{\varepsilon,\mathfrak{c}}|^2$. Then, applying Gronwall's inequality, one derives
		\begin{align}\label{U.2}
			\sup_{t \in[0,T]}\left\|f_R^{\varepsilon,\mathfrak{c}}(t)\right\|+\sup_{t \in[0,T]}\left\|\left(E_R^{\varepsilon,\mathfrak{c}}, B_R^{\varepsilon,\mathfrak{c}}\right)(t)\right\|\leq C_T\left\{\left\|f_R^{\varepsilon,\mathfrak{c}}(0)\right\|+\left\|\big(E_R^{\varepsilon,\mathfrak{c}},B_R^{\varepsilon,\mathfrak{c}}\big)(0)\right\|+1\right\},
		\end{align}
		where the constant $C_T>0$ depends on $T$, but is independent of $\mathfrak{c}$ and $\varepsilon$.
		
		Combining \eqref{W.0}, \eqref{H.1} and \eqref{H.5}, and noting $k\geq 5$, one obtains
		\begin{align*}
			\frac{d}{d t} & \Big\{\Big\|\sqrt{\frac{2\mathfrak{c}T_\mathfrak{e}}{u_{\mathfrak{e}}^0}}\nabla_x f_R^{\varepsilon,\mathfrak{c}}\Big\|^2+\left\| \nabla_p f_R^{\varepsilon, \mathfrak{c}}\right\|^2+\left\|\big(\nabla_x E_R^{\varepsilon,\mathfrak{c}}, \nabla_x B_R^{\varepsilon,\mathfrak{c}}\big)\right\|^2\Big\}+\frac{\zeta_0}{\varepsilon}\|\{\mathbf{I}-\mathbf{P}_{\mathfrak{c}}\}\nabla_{x,p}f_R^{\varepsilon,\mathfrak{c}}\|_{\nu_{\mathfrak{c}}}^2\nonumber\\
			&\lesssim\left\{\left\|f_R^{\varepsilon,\mathfrak{c}}\right\|_{H^1}^2+\left\|\big(E_R^{\varepsilon,\mathfrak{c}}, B_R^{\varepsilon,\mathfrak{c}}\big)\right\|_{H^1}^2+\frac{1}{\varepsilon^2}\left\|f_R^{\varepsilon,\mathfrak{c}}\right\|^2+\frac{1}{\varepsilon^2}\left\|\{\mathbf{I}-\mathbf{P}_{\mathfrak{c}}\}f_R^{\varepsilon, \mathfrak{c}}\right\|_{\nu_{\mathfrak{c}}}^2+1\right\},
		\end{align*}
		which, together with \eqref{U.1}, yields
		\begin{align*}
			\frac{d}{d t} & \Big\{\Big\|\sqrt{\frac{2\mathfrak{c}T_\mathfrak{e}}{u_{\mathfrak{e}}^0}}\nabla_x f_R^{\varepsilon,\mathfrak{c}}\Big\|^2+\left\| \nabla_p f_R^{\varepsilon, \mathfrak{c}}\right\|^2+\left\|\big(\nabla_x E_R^{\varepsilon,\mathfrak{c}}, \nabla_x B_R^{\varepsilon,\mathfrak{c}}\big)\right\|^2\Big\}+\frac{\zeta_0}{\varepsilon}\|\{\mathbf{I}-\mathbf{P}_{\mathfrak{c}}\}\nabla_{x,p}f_R^{\varepsilon,\mathfrak{c}}\|_{\nu_{\mathfrak{c}}}^2\nonumber\\
			&\qquad+\frac{2}{\varepsilon\zeta_0c_0 }\frac{d}{d t}\Big\{\frac{1}{2}\Big\|\sqrt{\frac{2\mathfrak{c}T_\mathfrak{e}}{u_{\mathfrak{e}}^0}}f_R^{\varepsilon,\mathfrak{c}}\Big\|^2+\left\|\big(E_R^{\varepsilon,\mathfrak{c}}, B_R^{\varepsilon,\mathfrak{c}}\big)\right\|^2\Big\}\nonumber\\
			&\lesssim\left\{\left\|f_R^{\varepsilon,\mathfrak{c}}\right\|_{H^1}^2+\left\|\big(E_R^{\varepsilon,\mathfrak{c}}, B_R^{\varepsilon,\mathfrak{c}}\big)\right\|_{H^1}^2+\frac{1}{\varepsilon^2}\left\|f_R^{\varepsilon,\mathfrak{c}}\right\|^2+1\right\}\nonumber\\
			&\qquad+\frac{2}{\varepsilon\zeta_0c_0 } \left\{\left\|f_R^{\varepsilon,\mathfrak{c}}\right\|^2+\left\|\big(E_R^{\varepsilon,\mathfrak{c}}, B_R^{\varepsilon,\mathfrak{c}}\big)\right\|^2+1\right\},
		\end{align*}
		Noting $\Big\|\sqrt{\frac{2\mathfrak{c}T_\mathfrak{e}}{u_{\mathfrak{e}}^0}}\nabla_x f_R^{\varepsilon,\mathfrak{c}}\Big\|^2+\left\|\sqrt{\frac{2\mathfrak{c}T_\mathfrak{e}}{u_{\mathfrak{e}}^0}}f_R^{\varepsilon, \mathfrak{c}}\right\|^2\cong \left\|\nabla_x f_R^{\varepsilon,\mathfrak{c}}\right\|^2+\left\|f_R^{\varepsilon, \mathfrak{c}}\right\|^2$, we further apply Gronwall's inequality and \eqref{U.2} to obtain that
		\begin{align}\label{U.3}
			&\sup_{t \in[0,T]}\left\|\nabla_{x,p}f_R^{\varepsilon,\mathfrak{c}}(t)\right\|+\sup_{t \in[0,T]}\left\|\big(\nabla_x E_R^{\varepsilon,\mathfrak{c}}, \nabla_x B_R^{\varepsilon,\mathfrak{c}}\big)(t)\right\|\nonumber\\ \leq&C_T\Big\{\left\|\nabla_{x,p}f_R^{\varepsilon,\mathfrak{c}}(0)\right\|+\left\|\big(\nabla_x E_R^{\varepsilon,\mathfrak{c}},\nabla_x B_R^{\varepsilon,\mathfrak{c}}\big)(0)\right\|+\frac{1}{\varepsilon}\left\{\left\|f_R^{\varepsilon,\mathfrak{c}}(0)\right\|+\left\|\big(E_R^{\varepsilon,\mathfrak{c}},B_R^{\varepsilon,\mathfrak{c}}\big)(0)\right\|+1\right\}\Big\} .
		\end{align}
		
		Now we turn to the $W^{1,\infty}$-estimate of $\big(E_R^{\varepsilon,\mathfrak{c}},B_R^{\varepsilon,\mathfrak{c}},h_R^{\varepsilon,\mathfrak{c}}\big)$. For the rest of the proof, denote by $\mathcal D_0$ the size of the initial data appearing on the right-hand side of \eqref{e1.16}, namely
		\begin{align}\label{U.D0}
		\mathcal D_0:=&\,\varepsilon^{\frac52}\big\|\big(E_R^{\varepsilon,\mathfrak{c}},B_R^{\varepsilon,\mathfrak{c}}\big)(0)\big\|_{L^\infty}
		+\varepsilon^{\frac32}\mathcal I_0+\big\|\nabla_{x,p}f_R^{\varepsilon,\mathfrak{c}}(0)\big\|\nonumber\\
		&+\big\|\big(\nabla_xE_R^{\varepsilon,\mathfrak{c}},\nabla_xB_R^{\varepsilon,\mathfrak{c}}\big)(0)\big\|+\big\|f_R^{\varepsilon,\mathfrak{c}}(0)\big\|
		+\big\|\big(E_R^{\varepsilon,\mathfrak{c}},B_R^{\varepsilon,\mathfrak{c}}\big)(0)\big\|+1.
		\end{align}
		In particular, \eqref{e.16-1} gives $\mathcal D_0\leq C$. From \eqref{U.2} and \eqref{U.3}, and using $0<\varepsilon\leq1$, we first record
		\begin{align}\label{U.3a}
		&\sup_{0\leq s\leq T}\big\|f_R^{\varepsilon,\mathfrak{c}}(s)\big\|
		+\sup_{0\leq s\leq T}\big\|\big(E_R^{\varepsilon,\mathfrak{c}},B_R^{\varepsilon,\mathfrak{c}}\big)(s)\big\|
		+\varepsilon\sup_{0\leq s\leq T}\big\|f_R^{\varepsilon,\mathfrak{c}}(s)\big\|_{H^1}
		\leq C_T\mathcal D_0 .
		\end{align}
		The factor $\varepsilon$ in front of the $H^1$ norm is important: it compensates the only $\varepsilon^{-1}$ contribution in \eqref{U.3} and reduces it to the lower-order initial size already contained in $\mathcal D_0$.

		Combining the improved field estimate \eqref{W.4} with \eqref{U.3a}, and using $\varepsilon^{k+5/2}\leq1$, gives
		\begin{align}\label{U.4}
		\sup_{0\leq s\leq T}\varepsilon^{\frac52}\big\|\big(E_R^{\varepsilon,\mathfrak{c}},B_R^{\varepsilon,\mathfrak{c}}\big)(s)\big\|_{L^\infty}
		\leq C_T\mathcal D_0 .
		\end{align}
		The estimate \eqref{W.4} has a constant uniform for all $\mathfrak{c}\geq\mathfrak{c}_0$ and contains no exponential loss involving the product of the two singular parameters.

		Next, applying \eqref{C.17}, and recalling that $\varepsilon^{3/2}\|h_R^{\varepsilon,\mathfrak{c}}(0)\|_{L^\infty}\leq\varepsilon^{3/2}\mathcal I_0$, we obtain from \eqref{U.3a} and \eqref{U.4} that
		\begin{align}\label{U.5}
		\sup_{0\leq s\leq T}\varepsilon^{\frac32}\big\|h_R^{\varepsilon,\mathfrak{c}}(s)\big\|_{L^\infty}
		\leq C_T\mathcal D_0 .
		\end{align}

		We now use the improved gradient-field estimate \eqref{W.14}. Together with \eqref{U.3a} and \eqref{U.4}, it yields
		\begin{align}\label{U.6}
		\sup_{0\leq s\leq T}\varepsilon^{\frac72}
		\big\|\big(\nabla_xE_R^{\varepsilon,\mathfrak{c}},\nabla_xB_R^{\varepsilon,\mathfrak{c}}\big)(s)\big\|_{L^\infty}
		\leq C_T\mathcal D_0 .
		\end{align}
		No upper bound on $\mathfrak{c}$ is used in this step. The gradient-field estimate \eqref{W.14} is uniform in $\mathfrak{c}$; the only smallness restriction is the choice of $\varepsilon_0(T)$ used there to absorb $C_T\varepsilon\Phi_1$.

		Finally, multiplying \eqref{D.1} by $\varepsilon$ gives
		\begin{align*}
		\varepsilon^{\frac52}\big\|\nabla_{x,p}h_R^{\varepsilon,\mathfrak{c}}(t)\big\|_{L^\infty}
		\lesssim&\ \varepsilon^{\frac52}\big\|h_R^{\varepsilon,\mathfrak{c}}(0)\big\|_{W^{1,\infty}}
		+\varepsilon^{\frac72}\sup_{0\leq s\leq t}\big\|\big(E_R^{\varepsilon,\mathfrak{c}},B_R^{\varepsilon,\mathfrak{c}}\big)(s)\big\|_{W^{1,\infty}}\\
		&+\varepsilon\sup_{0\leq s\leq t}\big\|f_R^{\varepsilon,\mathfrak{c}}(s)\big\|_{H^1}+\varepsilon^{k+\frac72}.
		\end{align*}
		The first term is bounded by $\varepsilon^{3/2}\mathcal I_0$, by the strengthened definition \eqref{e1.16-0}. Moreover,
		\[
		\varepsilon^{\frac72}\sup_{0\leq s\leq T}\big\|\big(E_R^{\varepsilon,\mathfrak{c}},B_R^{\varepsilon,\mathfrak{c}}\big)(s)\big\|_{W^{1,\infty}}
		\lesssim
		\sup_{0\leq s\leq T}\varepsilon^{\frac52}\big\|\big(E_R^{\varepsilon,\mathfrak{c}},B_R^{\varepsilon,\mathfrak{c}}\big)(s)\big\|_{L^\infty}
		+\sup_{0\leq s\leq T}\varepsilon^{\frac72}\big\|\big(\nabla_xE_R^{\varepsilon,\mathfrak{c}},\nabla_xB_R^{\varepsilon,\mathfrak{c}}\big)(s)\big\|_{L^\infty},
		\]
		since $0<\varepsilon\leq1$. Hence \eqref{U.3a}, \eqref{U.4}, and \eqref{U.6} imply
		\begin{align}\label{U.13}
		\sup_{0\leq s\leq T}\varepsilon^{\frac52}
		\big\|\nabla_{x,p}h_R^{\varepsilon,\mathfrak{c}}(s)\big\|_{L^\infty}
		\leq C_T\mathcal D_0 .
		\end{align}

		From \eqref{U.2}--\eqref{U.6} and \eqref{U.13}, we conclude
		\begin{align*}
		 \sup _{t \in [0,T]}&\Big\{\varepsilon^{\frac{5}{2}}\left\|\nabla_{x,p}h_R^{\varepsilon,\mathfrak{c}}(t)\right\|_{L^{\infty}}+\varepsilon^{\frac{3}{2}}\left\|h_R^{\varepsilon,\mathfrak{c}}(t)\right\|_{L^{\infty}}+\varepsilon^{\frac{7}{2}}\big\|\big(\nabla_x E_R^{\varepsilon,\mathfrak{c}}, \nabla_x B_R^{\varepsilon,\mathfrak{c}}\big)(t)\big\|_{L^{\infty}}\nonumber\\
	   &+\varepsilon^{\frac{5}{2}}\big\|\big(E_R^{\varepsilon,\mathfrak{c}}, B_R^{\varepsilon,\mathfrak{c}}\big)(t)\big\|_{L^{\infty}}
		+\left\|\nabla_{x,p}f_R^{\varepsilon,\mathfrak{c}}(t)\right\|+\left\|\big(\nabla_{x}E_R^{\varepsilon,\mathfrak{c}},\nabla_{x}B_R^{\varepsilon,\mathfrak{c}}\big)(t)\right\|\nonumber\\
	   &+\left\|f_R^{\varepsilon,\mathfrak{c}}(t)\right\|+\left\|\big(E_R^{\varepsilon,\mathfrak{c}}, B_R^{\varepsilon,\mathfrak{c}}\big)(t)\right\|\Big\}
		\leq C_T\mathcal D_0,
		\end{align*}
		which yields \eqref{e1.16}. The derivation of \eqref{U.4}--\eqref{U.13} uses only the smallness of $\varepsilon\leq\varepsilon_0(T)$ and never invokes an upper bound on $\mathfrak{c}$ in terms of $\varepsilon$.
 
		Noting \eqref{e.16-1}, we actually obtain
         \begin{align}\label{U.14}
         \begin{gathered}
        \sup_{t\in[0,T]}\varepsilon^{2}\|h_R^{\varepsilon,\mathfrak{c}}(t)\|_{L^\infty}\leq C_T \varepsilon^{\frac{1}{2}}, \quad
         \sup_{t\in[0,T]}\varepsilon^{3}\|\nabla_{x,p}h_R^{\varepsilon,\mathfrak{c}}(t)\|_{L^\infty}
        \leq C_T \varepsilon^{\frac{1}{2}},\\
        \sup_{t \in[0,T]}\varepsilon^3\big\|\big(E_R^{\varepsilon,\mathfrak{c}}, B_R^{\varepsilon,\mathfrak{c}}\big)(t)\big\|_{L^{\infty}}\leq C_T \varepsilon^{\frac{1}{2}}, \quad
       \sup_{t \in [0,T]} \varepsilon^4\big\|\big(\nabla_x E_R^{\varepsilon,\mathfrak{c}}, \nabla_x B_R^{\varepsilon,\mathfrak{c}}\big)(t)\big\|_{L^{\infty}}\leq C_T \varepsilon^{\frac{1}{2}}.
        \end{gathered}
        \end{align}
        For any fixed finite $T>0$, choosing  $\varepsilon_0 \in (0,1)$ sufficiently small such that $C_T \varepsilon_0^{\frac{1}{4}}<\frac{1}{2}$, then for all $0<\varepsilon\leq \varepsilon_0$, \eqref{U.14} improves the bootstrap assumptions \eqref{W.0}-\eqref{C.2}, and therefore the bootstrap argument is closed.
		
		Finally, we conclude \eqref{e1.17} from \eqref{e.16-2}-\eqref{e1.16} and Proposition \ref{t5.3} for $k\geq 5$. Therefore the proof of Theorem \ref{t1.1} is completed. 
	\end{proof}
	\begin{proof}[Proof of Theorem \ref{t1.5}]
		Recall the local Maxwellian \eqref{e1.18}. It follows from \eqref{e1.17} that
		\begin{equation}\label{U.7}
			\left|F^{\varepsilon, \mathfrak{c}}(t, x, p)-\mathbf{M}_{\mathfrak{c}}(t, x, p)\right| \lesssim \varepsilon \sqrt{J_{\mathbf{M}}} \lesssim \varepsilon e^{-\frac{|p|}{2 T_M}} .
		\end{equation}
		A direct calculation shows that
		\begin{align}\label{U.8}
				& \mu(t, x, p)-\mathbf{M}_{\mathfrak{c}}(t, x, p) \nonumber\\
				=& \frac{\rho}{(2 \pi \theta)^{\frac{3}{2}}} \exp \left\{-\frac{|p-\mathfrak{u}|^2}{2 \theta}\right\}-\frac{n_{\mathfrak{e}} \gamma}{4 \pi \mathfrak{c}^3 K_2(\gamma)} \exp \left\{\frac{u_{\mathfrak{e}}^\mu p_\mu}{T_\mathfrak{e}}\right\}\nonumber \\
				=& \frac{\rho}{(2 \pi \theta)^{\frac{3}{2}}} \exp \left\{-\frac{|p-\mathfrak{u}|^2}{2 \theta}\right\}-\frac{n_{\mathfrak{e}}}{\left(2 \pi T_\mathfrak{e}\right)^{\frac{3}{2}}} \exp \left\{\frac{\mathfrak{c}^2+u_{\mathfrak{e}}^\mu p_\mu}{T_\mathfrak{e}}\right\}\left(1+O\left(\gamma^{-1}\right)\right) \nonumber\\
				=& O\left(\gamma^{-1}\right) \frac{n_{\mathfrak{e}}}{\left(2 \pi T_\mathfrak{e}\right)^{\frac{3}{2}}} \exp \left\{\frac{\mathfrak{c}^2+u_{\mathfrak{e}}^\mu p_\mu}{T_\mathfrak{e}}\right\}+\left(\frac{\rho}{(2 \pi \theta)^{\frac{3}{2}}}-\frac{n_{\mathfrak{e}}}{\left(2 \pi T_\mathfrak{e}\right)^{\frac{3}{2}}}\right) \exp \left\{-\frac{|p-\mathfrak{u}|^2}{2 \theta}\right\} \nonumber\\
				& \quad+\frac{n_{\mathfrak{e}}}{\left(2 \pi T_\mathfrak{e}\right)^{\frac{3}{2}}}\left(\exp \left\{-\frac{|p-\mathfrak{u}|^2}{2 \theta}\right\}-\exp \left\{\frac{\mathfrak{c}^2+u_{\mathfrak{e}}^\mu p_\mu}{T_\mathfrak{e}}\right\}\right) \nonumber\\
				:=& \mathcal{A}_1+\mathcal{A}_2+\mathcal{A}_3 .
		\end{align}
		
		By the Lorentz transformation, together with \cite[Lemmas 4.2--4.3]{Wang-Xiao-JLMS-2026}, we obtain
    \[
    \mathfrak{c}^2+u_{\mathfrak e}^\mu p_\mu
    \le -\frac12|p|+C,
    \]
    which, combined with \eqref{e.30}, yields
    \begin{align}\label{U.9}
    \left|\mathcal{A}_1\right|
    \lesssim \frac{1}{\mathfrak{c}^2} e^{- \frac{c_0}{2}|p|}.
    \end{align}
		
		For $\mathcal{A}_2$, it follows from Theorem \ref{c3.4} that
		$$
		\Big|\frac{\rho}{(2 \pi \theta)^{\frac{3}{2}}}-\frac{n_{\mathfrak{e}}}{\left(2 \pi T_\mathfrak{e}\right)^{\frac{3}{2}}}\Big|\lesssim|n_{\mathfrak{e}}-\rho|+|T_\mathfrak{e}-\theta|\lesssim \mathfrak{c}^{-1},
		$$
		which, together with $|p-\mathfrak u|^2
        \ge |p|^2-2|p||\mathfrak u|
        \ge \frac12|p|^2-2|\mathfrak u|^2\ge \frac12|p|^2-C$ and \eqref{e.30}, yields
		\begin{align}\label{U.10}
			\quad\left|\mathcal{A}_2\right| \lesssim \frac{1}{\mathfrak{c}} e^{-\frac{c_0 |p|}{4}} .
		\end{align}

    For $\mathcal{A}_3$, by the Mean Value Theorem applied to the function $f(z) = \exp\{z\}$, there exists a value $C(p)$ lying between the exponents $A(p):=-\frac{|p-\mathfrak u|^2}{2\theta}$ and $B(p):=\frac{\mathfrak{c}^2 + u_\mathfrak{e}^\mu p_\mu}{T_\mathfrak{e}}$ such that
    \begin{equation*}
    \left| \exp \left\{A(p) \right\} - \exp \left\{ B(p) \right\} \right| = \left|A(p)-B(p)\right|e^{C(p)} = \left| \frac{|p-\mathfrak u|^2}{2\theta} + \frac{\mathfrak{c}^2 + u_\mathfrak{e}^\mu p_\mu}{T_\mathfrak{e}} \right| e^{C(p)}.
    \end{equation*}
    A direct calculation shows that
		\begin{align*}
			\frac{|p-\mathfrak{u}|^2}{2\theta}+\frac{\mathfrak{c}^2+u_{\mathfrak{e}}^\mu p_\mu}{T_\mathfrak{e}}&=\frac{|p-\mathfrak{u}|^2}{2\theta T_\mathfrak{e}}(T_\mathfrak{e}-\theta)+\frac{1}{2T_\mathfrak{e}}(|p-\mathfrak{u}|^2+2\mathfrak{c}^2+2u_{\mathfrak{e}}^\mu p_\mu)\nonumber\\
			&=\frac{|p-\mathfrak{u}|^2}{2\theta T_\mathfrak{e}}(T_\mathfrak{e}-\theta)+\frac{1}{2T_\mathfrak{e}}\Big[\frac{|p|^4}{(p^0+u_{\mathfrak{e}}^0)^2}+2p\cdot (u_{\mathfrak{e}}-\mathfrak{u})+(|\mathfrak{u}|^2-|u_{\mathfrak{e}}|^2)\Big]\nonumber\\
			&\qquad +\frac{1}{2T_\mathfrak{e}}\Big[\frac{|u_{\mathfrak{e}}|^4}{(u_{\mathfrak{e}}^0+p^0)^2}-2\frac{|p|^2|u_{\mathfrak{e}}|^2}{(u_{\mathfrak{e}}^0+p^0)^2}\Big],
		\end{align*}
		which implies that 
		\begin{align*}
			\Big|\frac{|p-\mathfrak{u}|^2}{2\theta}+\frac{\mathfrak{c}^2+u_{\mathfrak{e}}^\mu p_\mu}{T_\mathfrak{e}}\Big|\lesssim \frac{1+|p|^2}{\mathfrak{c}}+\frac{|p|^4}{\mathfrak{c}^2}+\frac{1+|p|}{\mathfrak{c}^3}.
		\end{align*}

    Noting $\exp \left\{A(p) \right\}\lesssim e^{-\frac{c_0 |p|}{4}} $ and $\exp \left\{B(p) \right\}\lesssim e^{- \frac{c_0}{2}|p|}$ from \eqref{U.9}-\eqref{U.10}, it follows that
    $$
    e^{C(p)} \leq \max \left\{\exp \left\{A(p) \right\}, \exp \left\{B(p) \right\}\right\} \lesssim e^{-\frac{c_0 |p|}{4}}.
    $$
    Consequently, the term $e^{C(p)}$ provides sufficient decay to bound any polynomial of $p$. Then it holds that
    \begin{align}\label{U.11}
    |\mathcal{A}_3| &\le \frac{n_\mathfrak{e}}{(2\pi T_\mathfrak{e})^{\frac{3}{2}}} \left( \frac{1+|p|^2}{\mathfrak{c}}+\frac{|p|^4}{\mathfrak{c}^2}+\frac{1+|p|}{\mathfrak{c}^3} \right) e^{-\frac{c_0 |p|}{4}} \nonumber\\
    &\lesssim \mathfrak{c}^{-1} \left( (1+|p|)^4 e^{-\frac{c_0 |p|}{4}} \right) \nonumber\\
    &\lesssim \mathfrak{c}^{-1} e^{-\frac{c_0}{8}|p|}, 
    \end{align}

		Combining \eqref{U.8}-\eqref{U.11}, one has
		\begin{align}\label{U.12}
			\left|\mu(t, x, p)-\mathbf{M}_{\mathfrak{c}}(t, x, p)\right| \lesssim \mathfrak{c}^{-1}e^{-\frac{c_0 |p|}{8}} .
		\end{align}
		
		Using \eqref{U.7}, \eqref{U.12} and taking
		$$
		\delta_0:=\min \left(\frac{1}{2 T_M}, \frac{c_0}{8}\right)>0,
		$$
		one has
		$$
		\left|F^{\varepsilon, \mathfrak{c}}(t)-\mu(t)\right| \lesssim \varepsilon e^{-\frac{|p|}{2 T_M}}+\mathfrak{c}^{-1}e^{-\frac{c_0 |p|}{8}}\lesssim\left(\varepsilon+\mathfrak{c}^{-1}\right) e^{-\delta_0|p|},
		$$
		which implies that
		$$
		\sup _{0 \leq t \leq T}\big\|\left(F^{\varepsilon, \mathfrak{c}}-\mu\right)(t) e^{\delta_0|p|}\big\|_{L^{\infty}} \lesssim \varepsilon +\mathfrak{c}^{-1} .
		$$
		For \eqref{e1.19*}, it follows directly from \eqref{e1.17} and \eqref{e.29}-\eqref{e.26}.
		Therefore the proof of Theorem \ref{t1.5} is completed.
	\end{proof}

	\appendix
	\section{High-order Momentum Derivative Estimates for the Collision Operator}\label{Appendix A}

	Recalling the decomposition of $A$ and $A^c$ in \eqref{A.1}, consider the smooth test function $\chi \in C_0^{\infty}([0, \infty))$ such that $0 \leq \chi \leq 1$, and $\chi(r)=1$ for $r \in[0,1]$ with $\chi(r)=0$ for $r>2$. We use the splitting $1=\chi_A(p, q)+\chi_{A^c}(p, q)$ with
	\begin{align*}
		& \chi_A(p, q) := \chi\Big(\frac{p^0}{\mathfrak{c}}\Big)+\left(1-\chi\Big(\frac{p^0}{\mathfrak{c}}\Big)\right) \chi\Big(\frac{2}{3}\frac{|p|}{q^0}\Big), \\
		& \chi_{A^c}(p, q) :=\left(1-\chi\Big(\frac{p^0}{\mathfrak{c}}\Big)\right)\left(1-\chi\Big(\frac{2}{3}\frac{|p|}{q^0}\Big)\right).
	\end{align*}
	Here $\chi_A+\chi_{A^c}=1$. We emphasize that $\chi_A$ and $\chi_{A^c}$ are smooth cutoffs adapted to the regions $A$ and $A^c$, rather than exact cutoffs equal to $1$ on $A$ and $A^c$, respectively. In particular,
    \[
    \operatorname{supp}\chi_A\subset\{|p|\le 3q^0\},
    \qquad
    \operatorname{supp}\chi_{A^c}\subset\Big\{|p|\ge \frac32 q^0\Big\}.
    \]
	We split $\Gamma_{\mathfrak{c}}\left(f_1, f_2\right)=\Gamma_{\mathfrak{c}}^A+\Gamma_{\mathfrak{c}}^{A^c}$ using \eqref{04}, \eqref{01} and \eqref{07} as
	\begin{align}
		\Gamma_{\mathfrak{c}}^A & =\int_{\mathbb{R}^3 \times \mathbb{S}^2} d \omega d q \frac{\mathfrak{s} \mathfrak{B}}{p^0 q^0} \sqrt{\mathbf{M}_{\mathfrak{c}}(q)}\left[f_1\left(p^{\prime\prime}\right) f_2\left(q^{\prime\prime}\right)-f_1(p) f_2(q)\right] \chi_A(p, q), \label{A.2}\\
		\Gamma_{\mathfrak{c}}^{A^c} & =\int_{\mathbb{R}^3 \times \mathbb{S}^2} d \omega d q v_\phi \sqrt{\mathbf{M}_{\mathfrak{c}}(q)}\left[f_1\left(p^{\prime}\right) f_2\left(q^{\prime}\right)-f_1(p) f_2(q)\right] \chi_{A^c}(p, q).\label{A.3}
	\end{align}
	Here without loss of generality, we have taken $f_1$ and $f_2$ to be scalar functions.
	
	\subsection{The Glassey--Strauss frame}
	To avoid taking derivatives for the singular factor of $\left|\omega \cdot\Big(\frac{p}{p^0}-\frac{q}{q^0}\Big)\right|$ inside $\mathfrak{B}(p, q, \omega)$ for $\partial_{\beta} \Gamma_{\mathfrak{c}}^A$ in \eqref{A.2}), we introduce the following change of variables $q \rightarrow \bar{u}$ (for fixed $p$ ) as:
	\begin{align}\label{A.4}
		\bar{u}=p^0 q-q^0 p	
	\end{align}
	By \eqref{A.4}, we have that $q=\frac{q^0}{p^0} p+\frac{\bar{u}}{p^0}$ and taking norms on both sides yields
	$$
	q^0=\left(\frac{\bar{u}}{\mathfrak{c}^2}\cdot p\right)+\sqrt{\left(\frac{\bar{u}}{\mathfrak{c}^2}\cdot p\right)^2+\frac{\bar{u}^2}{\mathfrak{c}^2}+\left(p^0\right)^2} .
	$$
	Such a transformation \eqref{A.4} therefore defines an invertible mapping with
	\begin{align*}
		\frac{\partial \bar{u}_i}{\partial q_j} & =p^0 \delta_{i j}-\frac{q_j p_i}{q^0}, \quad(i, j=1,2,3), \\
		\Big|\frac{\partial \bar{u}}{\partial q}\Big| & =\operatorname{det}\left(\frac{\partial \bar{u}_i}{\partial q_j}\right)=\frac{\left(p^0\right)^2}{q^0}\left(p^0 q^0-p \cdot q\right) \geq \mathfrak{c}^2\frac{\left(p^0\right)^2}{q^0}.
	\end{align*}
	Since $\left|\omega \cdot\Big(\frac{p}{p^0}-\frac{q}{q^0}\Big)\right|=\frac{|\omega \cdot \bar{u}|}{p^0 q^0}$, by \eqref{A.4} we can express $\Gamma_{\mathfrak{c}}^A(f_1, f_2)$ as
	\begin{align}\label{A.5}
		\Gamma_{\mathfrak{c}}^A(f_1, f_2)=\int_{\mathbb{R}^3 \times \mathbb{S}^2} d \omega d \bar{u}\Big|\frac{\partial q}{\partial \bar{u}}\Big| \frac{\mathfrak{s} \tilde{\mathfrak{B}}|\omega \cdot \bar{u}|}{p^0 q^0}\sqrt{\mathbf{M}_{\mathfrak{c}}(q)}\left\{f_1\left(p^{\prime\prime}\right) f_2\left(q^{\prime\prime}\right)-f_1(p) f_2(q)\right\} \chi_A(p, q),
	\end{align}
	where now
	$$
	\tilde{\mathfrak{B}}:= \frac{\mathfrak{c}\left(p^0+q^0\right)^2}{\left[\left(p^0+q^0\right)^2-(\omega \cdot[p+q])^2\right]^2}.
	$$
    
	We take a high order derivative $\partial_{\beta}$ of \eqref{A.5} to obtain
	\begin{align}\label{A.6}
		\left|\partial_{\beta} \Gamma_{\mathfrak{c}}^A\right| \lesssim & \sum \int_{\mathbb{R}^3 \times \mathbb{S}^2} \mathbf{1}_{|p| \leq 3 q^0} d \omega d \bar{u} K_{\beta_0}^A \sqrt{\mathbf{M}_{\mathfrak{c}}(q)}\left|\left(\partial_{\beta_1} f_1\right)\left(p^{\prime\prime}\right)\left(\partial_{\beta_2} f_2\right)\left(q^{\prime\prime}\right) \iota_{q^{\prime\prime}}^{p^{\prime\prime}}\right| \nonumber\\
		& +\sum \int_{\mathbb{R}^3 \times \mathbb{S}^2} \mathbf{1}_{|p| \leq 3 q^0} d \omega d\bar{u} K_{\beta_0}^A\sqrt{\mathbf{M}_{\mathfrak{c}}(q)} \left|\partial_{\beta_1} f_1(p) \partial_{\beta_2} f_2(q)\iota_q\right|,
	\end{align}
	where the sum is over $\beta_0+\beta_1+\beta_2 \leq \beta$ and 
	\begin{align}\label{A.7}
		K_{\beta_0}^A=K_{\beta_0}^A(\bar{u}, p, \omega) :=|\omega \cdot \bar{u}|\left|\partial_{\beta_0}\Big(\Big|\frac{\partial q}{\partial \bar{u}}\Big| \frac{\mathfrak{s} \tilde{\mathfrak{B}}}{p^0 q^0}\mathbf{M}_{\mathfrak{c}}^{\frac{1}{2}}(q) \chi_A(p, q)\Big)\right|\mathbf{M}_{\mathfrak{c}}^{-\frac{1}{2}}(q).
	\end{align}
	Also $\iota_{q^{\prime\prime}}^{p^{\prime\prime}}=\iota_{q^{\prime\prime}}^{p^{\prime\prime}}(\bar{u}, p, \omega)$ and $\iota_q=\iota_q(\bar{u}, p, \omega)$ are the terms which result from applying the chain rule to $p^{\prime\prime}$, $q^{\prime\prime}$ and $q$. 
	\begin{lemma}\label{lA.1}
		On the set $A$, we have the following estimates for $|\beta|> 0$:
		\begin{align}\label{A.8}
			\Big|\frac{\partial \left(p^0 q-q^0 p\right)}{\partial q}\Big|\left|K_{\beta_0}^A\left(p^0 q-q^0 p, p, \omega\right)\right| \lesssim \frac{\mathfrak{s} \mathfrak{B}}{p^0 q^0}\langle q\rangle^n.
		\end{align}
		Similarly, we also have the upper bound of
		\begin{align}\label{A.9}
			\left|\iota_{q^{\prime\prime}}^{p^{\prime\prime}}\left(p^0 q-q^0 p, p, \omega\right)\right| \lesssim\langle q\rangle^n,\quad \left|\iota_q\left(p^0 q-q^0 p, p, \omega\right)\right| \lesssim\langle q\rangle^n.
		\end{align}
		
		Above $n>1$ is a fixed large integer which depends upon $\beta, \beta_0, \beta_1$, and $\beta_2$.
	\end{lemma}

\begin{proof}
\smallskip
\noindent
{\it Step 1. Estimates for $\iota_{q^{\prime\prime}}^{p^{\prime\prime}}$ and $\iota_q$.}
Recall that after the change of variables
\[
\bar u=p^0q-q^0p,
\]
we may write
\[
q_i=\frac{q^0}{p^0}p_i+\frac{\bar u_i}{p^0},
\qquad
q^0=\frac{\bar u\cdot p}{\mathfrak{c}^2}+\Theta,
\]
where
\[
\Theta:=\Bigg[\Big(\frac{\bar u}{\mathfrak{c}^2}\cdot p\Big)^2+\frac{|\bar u|^2}{\mathfrak{c}^2}+(p^0)^2\Bigg]^{1/2}.
\]
Moreover,
\[
p_i^{\prime\prime}=p_i+\tilde a\,\omega_i,
\qquad
q_i^{\prime\prime}=q_i-\tilde a\,\omega_i,
\qquad
\tilde a=\frac{X}{D},
\]
with
\[
X:=2(p^0+q^0)\,\omega\cdot\bar u,
\qquad
D:=(p^0+q^0)^2-\{\omega\cdot(p+q)\}^2.
\]

On the region $\{|p|\le 3q^0\}$, one has
\[
\frac{|p|}{\mathfrak{c}}+\frac{p^0}{\mathfrak{c}}+\frac{q^0}{\mathfrak{c}}\lesssim \langle q\rangle,
\qquad
\frac{|\bar u|}{\mathfrak{c}^2}\lesssim \langle q\rangle^2.
\]
In addition,
\[
p^0\ge \mathfrak{c},\qquad
\Theta\ge \mathfrak{c},\qquad
D\ge 4\mathfrak{c}^2,\qquad
Y:=p^0q^0-p\cdot q\ge \mathfrak{c}^2.
\]
Thus no singular denominator appears in the above quantities.

To keep track of the $\mathfrak{c}$-scaling systematically, we now introduce the dimensionless variables
\[
\widetilde p:=\frac{p}{\mathfrak{c}},\qquad
\widetilde q:=\frac{q}{\mathfrak{c}},\qquad
\widetilde p^0:=\frac{p^0}{\mathfrak{c}},\qquad
\widetilde q^0:=\frac{q^0}{\mathfrak{c}},\qquad
\widetilde{\bar u}:=\frac{\bar u}{\mathfrak{c}^2},
\]
together with
\[
\widetilde\Theta:=\frac{\Theta}{\mathfrak{c}},
\qquad
\widetilde D:=\frac{D}{\mathfrak{c}^2},
\qquad
\widetilde Y:=\frac{Y}{\mathfrak{c}^2}.
\]
Then
\[
\widetilde q^0=\widetilde{\bar u}\cdot \widetilde p+\widetilde\Theta,
\qquad
\frac{q_i}{\mathfrak{c}}
=
\frac{\widetilde q^0}{\widetilde p^0}\widetilde p_i+\frac{\widetilde{\bar u}_i}{\widetilde p^0},
\]
and
\[
\frac{p_i^{\prime\prime}}{\mathfrak{c}}
=
\widetilde p_i+\widetilde a\,\omega_i,
\qquad
\frac{q_i^{\prime\prime}}{\mathfrak{c}}
=
\frac{q_i}{\mathfrak{c}}-\widetilde a\,\omega_i,
\qquad
\widetilde a
=
\frac{2(\widetilde p^0+\widetilde q^0)\,\omega\cdot\widetilde{\bar u}}
{\widetilde D}.
\]
Furthermore, on $\{|p|\le 3q^0\}$,
\[
|\widetilde p|+\widetilde p^0+\widetilde q^0+|\widetilde{\bar u}|
\lesssim \langle q\rangle^2,
\]
while
\[
\widetilde p^0\ge 1,\qquad
\widetilde\Theta\ge 1,\qquad
\widetilde D\ge 4,\qquad
\widetilde Y\ge 1.
\]

We now introduce the class $\mathfrak R$ of rational functions of the form
\[
R(\widetilde p,\widetilde{\bar u},\omega)
=
\frac{
P(\widetilde p,\widetilde{\bar u},\widetilde p^0,\widetilde q^0,\omega)
}{
(\widetilde p^0)^a
(\widetilde\Theta)^b
(\widetilde D)^c
(\widetilde Y)^d
},
\]
where $a,b,c,d\ge0$ and $P$ is a polynomial in its arguments.

We claim that $\mathfrak R$ is stable under $\partial_{\widetilde p}$. Indeed,
\[
\partial_{\widetilde p_j}\widetilde p^0=\frac{\widetilde p_j}{\widetilde p^0},
\]
and
\[
\partial_{\widetilde p_j}\widetilde q^0
=
\widetilde{\bar u}_j
+
\frac{
(\widetilde{\bar u}\cdot \widetilde p)\widetilde{\bar u}_j+\widetilde p_j
}{
\widetilde\Theta
}.
\]
Similarly, $\partial_{\widetilde p_j}\widetilde D$ and $\partial_{\widetilde p_j}\widetilde Y$ also belong to $\mathfrak R$. Hence, by the product and quotient rules,
\[
\partial_{\widetilde p_j}\mathfrak R\subset\mathfrak R.
\]
Since
\[
\frac{q_i}{\mathfrak{c}},
\qquad
\widetilde a,
\qquad
\frac{p_i^{\prime\prime}}{\mathfrak{c}},
\qquad
\frac{q_i^{\prime\prime}}{\mathfrak{c}}
\]
all belong to $\mathfrak R$, it follows by induction on $|\beta|$ that
\[
\partial_{\widetilde p}^{\beta}\!\left(\frac{q_i}{\mathfrak{c}}\right),\qquad
\partial_{\widetilde p}^{\beta}\!\left(\frac{p_i^{\prime\prime}}{\mathfrak{c}}\right),\qquad
\partial_{\widetilde p}^{\beta}\!\left(\frac{q_i^{\prime\prime}}{\mathfrak{c}}\right)
\in\mathfrak R.
\]
Consequently, for each multi-index $\beta$,
\[
\left|
\partial_{\widetilde p}^{\beta}\!\left(\frac{q_i}{\mathfrak{c}}\right)
\right|
+
\left|
\partial_{\widetilde p}^{\beta}\!\left(\frac{p_i^{\prime\prime}}{\mathfrak{c}}\right)
\right|
+
\left|
\partial_{\widetilde p}^{\beta}\!\left(\frac{q_i^{\prime\prime}}{\mathfrak{c}}\right)
\right|
\le C_\beta \langle q\rangle^{n_\beta}
\]
for some $n_\beta>0$.

Noting $\partial_{p_j}=\mathfrak{c}^{-1}\partial_{\widetilde p_j}$, we obtain
\[
|\partial_\beta q_i|
+
|\partial_\beta p_i^{\prime\prime}|
+
|\partial_\beta q_i^{\prime\prime}|
\lesssim
\langle q\rangle^{n_\beta},
\]
uniformly in $\mathfrak{c}$. Finally, $\iota_{q^{\prime\prime}}^{p^{\prime\prime}}$ and $\iota_q$ are finite sums of products of such derivatives arising from the chain rule. Therefore, there exists some positive integer $n$ such that
\[
\left|\iota_{q^{\prime\prime}}^{p^{\prime\prime}}\right|\lesssim \langle q\rangle^n,
\qquad
\left|\iota_q\right|\lesssim \langle q\rangle^n,
\]
which proves \eqref{A.9}.

\smallskip
\noindent
{\it Step 2. Estimate for $K_{\beta_0}^A$.}
Recall that
\[
K_{\beta_0}^A
=
|\omega\cdot\bar u|
\left|
\partial_{\beta_0}
\Big(
\Big|\frac{\partial q}{\partial\bar u}\Big|
\frac{\mathfrak s\widetilde{\mathfrak B}}{p^0q^0}
\mathbf M_{\mathfrak{c}}^{1/2}(q)\chi_A(p,q)
\Big)
\right|
\mathbf M_{\mathfrak{c}}^{-1/2}(q).
\]
We first note that
\[
\Big|\frac{\partial q}{\partial\bar u}\Big|
=
\frac{q^0}{(p^0)^2Y},
\qquad
\widetilde{\mathfrak B}
=
\frac{\mathfrak{c} (p^0+q^0)^2}{D^2}.
\]
Set
\[
F(\bar u,p,\omega)
:=
\Big|\frac{\partial q}{\partial\bar u}\Big|
\frac{\mathfrak s\widetilde{\mathfrak B}}{p^0q^0}.
\]
In terms of the dimensionless variables introduced in Step 1,
\[
F
=
\frac{1}{\mathfrak{c}^4}
\frac{\widetilde{\mathfrak s}\,(\widetilde p^0+\widetilde q^0)^2}
{(\widetilde p^0)^3\,\widetilde Y\,(\widetilde D)^2},
\qquad
\widetilde{\mathfrak s}:=\frac{\mathfrak s}{\mathfrak{c}^2},
\]
where $\widetilde{\mathfrak s}=2(1+\widetilde Y)$.
On $\{|p|\le 3q^0\}$, we have
\[
\widetilde p^0\ge 1,\qquad
\widetilde q^0\ge 1,\qquad
\widetilde Y\ge 1,\qquad
\widetilde D\ge 4,
\]
and hence
\[
\widetilde p^0+\widetilde q^0\ge 2,
\qquad
\widetilde{\mathfrak s}\ge 4.
\]
Moreover, since
\[
\partial_{p_j}=\mathfrak{c}^{-1}\partial_{\widetilde p_j},
\]
the derivative bounds established in Step 1 imply that, for every multi-index $\beta$,
\[
|\partial_\beta \widetilde{\mathfrak s}|
+
|\partial_\beta \widetilde p^0|
+
|\partial_\beta \widetilde q^0|
+
|\partial_\beta \widetilde Y|
+
|\partial_\beta \widetilde D|
\lesssim
\langle q\rangle^{N_\beta}.
\]

For $|\beta|=1$, the product and quotient rules give
\[
F^{-1}\partial_j F
=
\frac{\partial_j\widetilde{\mathfrak s}}{\widetilde{\mathfrak s}}
+2\frac{\partial_j(\widetilde p^0+\widetilde q^0)}{\widetilde p^0+\widetilde q^0}
-3\frac{\partial_j\widetilde p^0}{\widetilde p^0}
-\frac{\partial_j\widetilde Y}{\widetilde Y}
-2\frac{\partial_j\widetilde D}{\widetilde D},
\]
and therefore
\[
|F^{-1}\partial_j F|
\lesssim
\langle q\rangle^{N_1}.
\]
For general $\beta$, repeated differentiation shows that $F^{-1}\partial_\beta F$ is a finite sum of products of factors of the form
\[
\frac{\partial_{\beta'} \widetilde{\mathfrak s}}{\widetilde{\mathfrak s}},
\qquad
\frac{\partial_{\beta'} (\widetilde p^0+\widetilde q^0)}{\widetilde p^0+\widetilde q^0},
\qquad
\frac{\partial_{\beta'} \widetilde p^0}{\widetilde p^0},
\qquad
\frac{\partial_{\beta'} \widetilde Y}{\widetilde Y},
\qquad
\frac{\partial_{\beta'} \widetilde D}{\widetilde D},
\qquad
1\le |{\beta'}|\le |\beta|.
\]
Since all denominators above are bounded away from zero on $\{|p|\le 3q^0\}$, while the numerators are polynomially bounded in $\langle q\rangle$, it follows that
\[
|F^{-1}\partial_\beta F|
\lesssim
\langle q\rangle^{N_\beta}
\]
for some integer $N_{\beta}>0$. Consequently,
\begin{align}\label{A.12a-new}
|\partial_{\beta} F|
\lesssim
F\,\langle q\rangle^{N_{\beta}}
=
\Big|\frac{\partial q}{\partial\bar u}\Big|
\frac{\mathfrak s\widetilde{\mathfrak B}}{p^0q^0}
\langle q\rangle^{N_{\beta}}.
\end{align}

Next, by the estimate already proved in Step 1 and repeated application of the chain rule,
\begin{align}\label{A.12-new}
\left|\partial_{\beta}\mathbf M_{\mathfrak{c}}^{1/2}(q)\right|
\le C_\beta \langle q\rangle^{N_\beta}\mathbf M_{\mathfrak{c}}^{1/2}(q).
\end{align}
Likewise, every derivative of the cutoff is polynomially bounded:
\begin{align}\label{A.13-new}
|\partial_\beta \chi_A(p,q)|
\le C_\beta \langle q\rangle^{N_\beta}.
\end{align}
Indeed, derivatives of the cutoff depending on $p^0/\mathfrak{c}$ are harmless. For the cutoff depending on $|p|/q^0$, any factor $|p|^{-m}$ arising from differentiating $|p|$ is localized to the support of $\chi^{(k)}(2|p|/(3q^0))$, where $|p|\approx q^0$. Hence $|p|^{-m}\lesssim (q^0)^{-m}$ there, and the resulting terms are again bounded by a polynomial in $\langle q\rangle$.

Applying Leibniz's rule to the definition of $K_{\beta_0}^A$, and using \eqref{A.12a-new}--\eqref{A.13-new}, we obtain
\[
K_{\beta_0}^A(\bar u,p,\omega)
\lesssim
|\omega\cdot\bar u|
\Big|\frac{\partial q}{\partial\bar u}\Big|
\frac{\mathfrak s\widetilde{\mathfrak B}}{p^0q^0}
\langle q\rangle^n
\]
for some integer $n>0$ depending only on $\beta_0$.

Finally, reverting to the original variable $\bar u=p^0q-q^0p$ and using
\[
\Big|\frac{\partial(p^0q-q^0p)}{\partial q}\Big|
\Big|\frac{\partial q}{\partial\bar u}\Big|=1,
\qquad
\frac{\widetilde{\mathfrak B}\,|\omega\cdot\bar u|}{p^0q^0}
=
\mathfrak B(p,q,\omega),
\]
we conclude that
\[
\Big|\frac{\partial(p^0q-q^0p)}{\partial q}\Big|
\left|
K_{\beta_0}^A\big(p^0q-q^0p,p,\omega\big)
\right|
\lesssim
\frac{\mathfrak s\mathfrak B}{p^0q^0}\langle q\rangle^n,
\]
which is exactly \eqref{A.8}. The proof is complete.
\end{proof}
	\begin{lemma}\label{lA.2}
		For any $\ell \geq 9$, we have the following estimate for the collision operator for $|\beta|>0$
		\begin{align}\label{A.10}
			\left|\left\langle\partial_{\beta} \Gamma_{\mathfrak{c}}^A\left(f_1, f_2\right), \partial_{\beta} f_3\right\rangle\right| \lesssim \left\|\partial_{\beta} f_3\right\|_{\infty, \ell}\sum_{\beta_1+\beta_2 \leq \beta}\left\|\partial_{\beta_1} f_1\right\|\left\|\partial_{\beta_2} f_2\right\| .
		\end{align}
		Furthermore, if $\chi(p)$ satisfies $|\partial_{\beta}\chi(p)| \lesssim e^{-\delta_1|p|}$ for some positive constant $\delta_1>0$, then we have the following estimate 
		\begin{align}\label{A.11}
			\left|\langle\partial_{\beta} \Gamma_{\mathfrak{c}}^A\left(f_1, \chi\right),\partial_{\beta} f_3\rangle \right|+\left|\left\langle\partial_{\beta} \Gamma_{\mathfrak{c}}^A\left(\chi, f_1\right),\partial_{\beta} f_3\right\rangle \right| \lesssim \sum_{\beta_1\leq \beta}\left\|\partial_{\beta_1} f_1\right\|_{\nu_{\mathfrak{c}}}\left\|\partial_{\beta} f_3\right\|_{\nu_{\mathfrak{c}}} ,
		\end{align}
		where the constants are independent of $\mathfrak{c}$.
	\end{lemma}
	\begin{proof}
		We notice from \eqref{A.6} and Lemma \ref{lA.1} that
		\begin{align*}
			\left|\partial_{\beta} \Gamma_{\mathfrak{c}}^A\right| \lesssim & \sum_{\beta_1+\beta_2 \leq \beta} \int_{\mathbb{R}^3 \times \mathbb{S}^2} \mathbf{1}_{|p| \leq 3 q^0} \frac{\mathfrak{s} \mathfrak{B}}{p^0 q^0}\langle q\rangle^n \sqrt{\mathbf{M}_{\mathfrak{c}}(q)}\left|\left(\partial_{\beta_1} f_1\right)\left(p^{\prime\prime}\right)\left(\partial_{\beta_2} f_2\right)\left(q^{\prime\prime}\right) \right|d \omega d q \\
			& +\sum_{\beta_1+\beta_2 \leq \beta} \int_{\mathbb{R}^3 \times \mathbb{S}^2} \mathbf{1}_{|p| \leq 3 q^0}  \frac{\mathfrak{s} \mathfrak{B}}{p^0 q^0}\langle q\rangle^n\sqrt{\mathbf{M}_{\mathfrak{c}}(q)} \left|\partial_{\beta_1} f_1(p) \partial_{\beta_2} f_2(q)\right|d \omega d q:=\uppercase\expandafter{\romannumeral1} +\uppercase\expandafter{\romannumeral2} .
		\end{align*}
		We estimate first the second term:
		\begin{align*}
			\left|\langle\uppercase\expandafter{\romannumeral2},\partial_{\beta} f_3\rangle\right|\lesssim\left\|\partial_{\beta} f_3\right\|_{\infty,\ell}\sum_{\beta_1+\beta_2 \leq \beta} \int_{\mathbb{R}^3 \times\mathbb{R}^3\times\mathbb{R}^3\times \mathbb{S}^2}\frac{\left|\partial_{\beta_1} f_1(p) \partial_{\beta_2} f_2(q)\right|}{\langle p \rangle^{\ell-1}\langle q \rangle^{\ell-1}}d\omega d q d p d x,
		\end{align*}
		where we have used $\frac{\mathfrak{s} \mathfrak{B}}{p^0 q^0}\lesssim \frac{\mathfrak{c}\left(q^0\right)^3\left(|q|+|p|\right)}{\mathfrak{c}^4}\lesssim\langle p \rangle\langle q\rangle^4$ on the region $|p| \leq 3 q^0$ and the trivial estimate $\langle q\rangle^{n+5}\sqrt{\mathbf{M}_{\mathfrak{c}}(q)}\lesssim \langle q\rangle^{-\ell}$. Then the estimate \eqref{A.10} follows from the Cauchy-Schwarz inequality.
		
		Similarly, for term  $\uppercase\expandafter{\romannumeral1}$ we have
		\begin{align*}
			\left|\langle\uppercase\expandafter{\romannumeral1},\partial_{\beta} f_3\rangle\right|&\lesssim
			\left\|\partial_{\beta} f_3\right\|_{\infty,\ell}\sum_{\beta_1+\beta_2 \leq \beta} \int_{\mathbb{R}^3 \times\mathbb{R}^3\times\mathbb{R}^3\times \mathbb{S}^2}\frac{\mathfrak{s} \mathfrak{B}}{p^0 q^0}\frac{\left|\partial_{\beta_1} f_1(p^{\prime\prime}) \partial_{\beta_2} f_2(q^{\prime\prime})\right|}{\langle p \rangle^{\ell}\langle q \rangle^{\ell}}d\omega d q d p d x\\
			&\lesssim \left\|\partial_{\beta} f_3\right\|_{\infty,\ell}\sum_{\beta_1+\beta_2 \leq \beta} \prod_{i=1,2}\Bigg\{\int_{\mathbb{R}^3 \times\mathbb{R}^3\times\mathbb{R}^3\times \mathbb{S}^2}\frac{\mathfrak{s} \mathfrak{B}}{p^0 q^0}\frac{\left|\partial_{\beta_i} f_i(p^{\prime\prime})\right|^2}{\langle p \rangle^{\ell}\langle q \rangle^{\ell}}d\omega d q d p d x\Bigg\}^{1/2}.
		\end{align*}
		Above, we used the $\left(p^{\prime\prime}, q^{\prime\prime}\right)$ symmetry to interchange the values of $p^{\prime\prime}$ and $q^{\prime\prime}$ in the integral involving $f_2\left(p^{\prime\prime}\right)$ above. Next, by the pre-post collisional change of variables, which is $d p d q=\frac{p^0q^0}{p^{\prime\prime0}q^{\prime\prime0}}d p^{\prime\prime} d q^{\prime\prime}$, we have
		\begin{align*}
			\int_{\mathbb{R}^3 \times\mathbb{R}^3\times\mathbb{R}^3\times \mathbb{S}^2}\frac{\mathfrak{s} \mathfrak{B}}{p^0 q^0}\frac{\left|\partial_{\beta_i} f_i(p^{\prime\prime})\right|^2}{\langle p \rangle^{\ell}\langle q \rangle^{\ell}}d\omega d q d p d x=\int_{\mathbb{R}^3 \times\mathbb{R}^3\times\mathbb{R}^3\times \mathbb{S}^2}\frac{\mathfrak{s} \mathfrak{B}}{p^0 q^0}\frac{\left|\partial_{\beta_i} f_i(p)\right|^2}{\langle p^{\prime\prime} \rangle^{\ell}\langle q^{\prime\prime} \rangle^{\ell}}d\omega d q d p d x.
		\end{align*}
		Here we used the fact that the kernel of the integral is invariant with respect to the relativistic pre-post collisional change of variables from \cite{Glassey-TTSP-1991}.
		
		Recall that $\langle p\rangle\lesssim\langle p^{\prime\prime}\rangle\langle q^{\prime\prime}\rangle$ and $\langle q\rangle\lesssim\langle p^{\prime\prime}\rangle\langle q^{\prime\prime}\rangle$. Since $\ell\geq9$, $\langle p \rangle^2\langle q \rangle^7\lesssim \langle p^{\prime\prime}\rangle^9\langle q^{\prime\prime}\rangle^9\lesssim\langle p^{\prime\prime} \rangle^{\ell}\langle q^{\prime\prime} \rangle^{\ell}$, then we obtain
		\begin{align*}
			\int_{\mathbb{R}^3 \times\mathbb{R}^3\times\mathbb{R}^3\times \mathbb{S}^2}\frac{\mathfrak{s} \mathfrak{B}}{p^0 q^0}\frac{\left|\partial_{\beta_i} f_i(p)\right|^2}{\langle p^{\prime\prime} \rangle^{\ell}\langle q^{\prime\prime} \rangle^{\ell}}d\omega d q d p d x\lesssim\int_{\mathbb{R}^3 \times\mathbb{R}^3\times\mathbb{R}^3\times \mathbb{S}^2}\frac{\langle p \rangle\langle q \rangle^3}{\langle p \rangle^2\langle q \rangle^7}\left|\partial_{\beta_i} f_i(p)\right|^2d\omega d q d p d x\lesssim\left\|\partial_{\beta_i} f_i\right\|^2.
		\end{align*}
		Thus we complete the proof of \eqref{A.10}.
		
		Now we are ready for the proof the estimate \eqref{A.11}. First we prove that $\frac{\mathfrak{s} \mathfrak{B}}{p^0 q^0}\lesssim \nu_{\mathfrak{c}}(p)\langle q\rangle^4$ on the region $A$. Note that $\mathfrak{s}\leq4p^0q^0$ and $\big[\left(p^0+q^0\right)^2-(\omega \cdot[p+q])^2\big]^2\geq (4\mathfrak{c}^2)^2$ holds generally.
		
		{\it Case 1:} $|p|\leq \mathfrak{c}$. It follows form $p^0\leq2\mathfrak{c}\lesssim q^0$ and the definition of $\mathfrak{B}$ \eqref{01-0} that
		\begin{align*}
			\frac{\mathfrak{s} \mathfrak{B}}{p^0 q^0}\lesssim \frac{\mathfrak{c}\left(q^0\right)^3\left(|q|+|p|\right)}{\mathfrak{c}^4}\lesssim\left(|q|+|p|\right)\langle q\rangle^3\lesssim\nu_{\mathfrak{c}}(p)\langle q\rangle^4.
		\end{align*}
		
		{\it Case 2:} $|p|\geq \mathfrak{c}$ and $|p|\leq 3q^0$. It follows from Lemma \ref{l2.1} that $\nu_{\mathfrak{c}}(p)\simeq \mathfrak{c}$, then there holds
		\begin{align*}
			\frac{\mathfrak{s} \mathfrak{B}}{p^0 q^0}\lesssim \frac{\mathfrak{c}\left(q^0\right)^4}{\mathfrak{c}^4}\lesssim\mathfrak{c}\langle q\rangle^4\lesssim\nu_{\mathfrak{c}}(p)\langle q\rangle^4.
		\end{align*}
		
		Then we will estimate the $\uppercase\expandafter{\romannumeral1}$ and $\uppercase\expandafter{\romannumeral2}$ terms from the top of this proof with $f_2=\chi$ in
		the first case. In particular, as before, we have that
		\begin{align*}
			\left|\langle\uppercase\expandafter{\romannumeral2},\partial_{\beta} f_3\rangle\right|\lesssim&\sum_{\beta_1 \leq \beta} \int_{\mathbb{R}^3 \times\mathbb{R}^3\times\mathbb{R}^3 \times \mathbb{S}^2} \mathbf{1}_{|p| \leq 3 q^0}\nu_{\mathfrak{c}}(p)\langle q\rangle^n\sqrt{\mathbf{M}_{\mathfrak{c}}(q)}e^{-\delta_1|q|} \left|\partial_{\beta_1} f_1(p) \partial_{\beta} f_3(p)\right| d \omega d q d p d x\\
			\lesssim &\sum_{\beta_1\leq \beta}\left\|\partial_{\beta_1} f_1\right\|_{\nu_{\mathfrak{c}}}\left\|\partial_{\beta} f_3\right\|_{\nu_{\mathfrak{c}}}.
		\end{align*}
		We also have
		\begin{align*}
			\left|\langle\uppercase\expandafter{\romannumeral1},\partial_{\beta} f_3\rangle\right|\lesssim&\sum_{\beta_1\leq \beta} \int_{\mathbb{R}^3 \times\mathbb{R}^3\times\mathbb{R}^3 \times \mathbb{S}^2} \mathbf{1}_{|p| \leq 3 q^0} \frac{\mathfrak{s} \mathfrak{B}}{p^0 q^0}\langle q\rangle^n \sqrt{\mathbf{M}_{\mathfrak{c}}(q)}e^{-\delta_1|q^{\prime\prime}|}\left|\left(\partial_{\beta_1} f_1\right)\left(p^{\prime\prime}\right)\partial_{\beta} f_3\left(p\right) \right|d \omega d q d p d x\\
			\lesssim&\sum_{\beta_1\leq \beta}\prod_{i=1,3} \Bigg\{\int_{\mathbb{R}^3 \times\mathbb{R}^3\times\mathbb{R}^3 \times \mathbb{S}^2} \mathbf{1}_{|p| \leq 3 q^0} \frac{\mathfrak{s} \mathfrak{B}}{p^0 q^0}\left|\frac{\left|\partial_{\beta_i} f_i(p)\right|^2}{\langle q\rangle^\ell\langle q^{\prime\prime}\rangle^\ell} \right|d \omega d q d p d x\Bigg\}^{1/2},
		\end{align*}
		where we have used similar reasoning as we did for the previous terms. After the pre-post
		change of variables, as before, we conclude that $\left|\langle\uppercase\expandafter{\romannumeral1},\partial_{\beta} f_3\rangle\right|\lesssim\sum_{\beta_1\leq \beta}\left\|\partial_{\beta_1} f_1\right\|_{\nu_{\mathfrak{c}}}\left\|\partial_{\beta} f_3\right\|_{\nu_{\mathfrak{c}}}$ as desired. We have thus shown that $\left|\langle\partial_{\beta} \Gamma_{\mathfrak{c}}^A\left(f_1, \chi\right),\partial_{\beta} f_3\rangle \right|\lesssim\sum_{\beta_1\leq \beta}\left\|\partial_{\beta_1} f_1\right\|_{\nu_{\mathfrak{c}}}\left\|\partial_{\beta} f_3\right\|_{\nu_{\mathfrak{c}}}$. The last estimate involving the term $\Gamma_{\mathfrak{c}}^A\left(\chi, f_1\right)$ follows in exactly the same way.
	\end{proof}
\subsection{Center of momentum frame}

In this section we estimate the term $\Gamma_{\mathfrak{c}}^{A^c}$. 
Taking $\beta$ momentum derivatives of $\Gamma_{\mathfrak{c}}^{A^c}$, we obtain
\begin{align}\label{A.12}
\left|\partial_{\beta} \Gamma_{\mathfrak{c}}^{A^c}\right|
& \lesssim \sum \int_{\mathbb{R}^3 \times \mathbb{S}^2} d q\, d \omega \,
\mathbf{1}_{|p|\ge \frac32 q^0}
\left|\partial_{\beta_0} \left(v_{\phi}\chi_{A^c}\right)\right|
\sqrt{\mathbf{M}_{\mathfrak{c}}(q)}
\left|\left(\partial_{\beta_1} f_1\right)\left(p^{\prime}\right)\left(\partial_{\beta_2} f_2\right)\left(q^{\prime}\right) \kappa_{\beta_2}^{\beta_1}\right|
\nonumber\\
&\quad
+\sum
\left|\partial_{\beta_1} f_1(p)\right|
\int_{\mathbb{R}^3 \times \mathbb{S}^2} d q\, d \omega \,
\mathbf{1}_{|p|\ge \frac32 q^0}
\left|\partial_{\beta_0} \left(v_{\phi}\chi_{A^c}\right)\right|
\sqrt{\mathbf{M}_{\mathfrak{c}}(q)}
\left|\partial_{\beta_2} f_2(q)\right|.
\end{align}
Here $\kappa_{\beta_2}^{\beta_1}$ denotes the finite sum of products of momentum derivatives of $p'$ and $q'$ arising from repeated applications of the chain rule. As before, the summation is over multi-indices satisfying
\[
\beta_0+\beta_1+\beta_2\le \beta .
\]

We first derive the basic derivative bounds in the \emph{center of momentum} variables.

\begin{lemma}\label{lA.3}
Let $|p|\ge \frac32 q^0$. Then for any multi-index $\beta\neq 0$, there exists an integer $n=n(\beta)\ge 1$ such that
\[
\frac{\left|\partial_{\beta} \left(v_{\phi}\chi_{A^c}\right)\right|}{\left|v_\phi\right|}
+\left|\partial_{\beta} p^{\prime}\right|
+\left|\partial_{\beta} q^{\prime}\right|
\lesssim \langle q\rangle^n .
\]
\end{lemma}

To prove Lemma \ref{lA.3}, we use the following auxiliary estimate.

\begin{lemma}\label{lA.4}
Let $|p|\ge \frac32 q^0$. Then for any multi-index $\beta\neq 0$,
\[
\left|\partial_{\beta} g\right| \lesssim \langle q\rangle^{|\beta|},
\qquad
\left|\partial_{\beta}\left(\frac{1}{g}\right)\right|
\lesssim \frac{\langle q\rangle^{|\beta|}}{g^2},
\]
and
\[
\left|\partial_{\beta} \sqrt{\mathfrak{s}}\right| \lesssim \langle q\rangle^{|\beta|},
\qquad
\left|\partial_{\beta}\left(\frac{1}{\sqrt{\mathfrak{s}}}\right)\right|
\lesssim \frac{\langle q\rangle^{|\beta|}}{g^2}.
\]
\end{lemma}

\begin{proof}
We argue by induction on $|\beta|$. If $|\beta|=1$, then
\begin{align*}
\left|\partial_{p_i} g\right|
&=\left|\frac{1}{g}\left\{\frac{p_i}{p^0}q^0-q_i\right\}\right|
\lesssim \frac{q^0}{g},\\
\left|\partial_{p_i}\left\{\frac{1}{g}\right\}\right|
&=\left|-\frac{\partial_{p_i}g}{g^2}\right|
\lesssim \frac{q^0}{g^3},\\
\left|\partial_{p_i}\sqrt{\mathfrak{s}}\right|
&=\left|\frac{g\,\partial_{p_i}g}{\sqrt{\mathfrak{s}}}\right|
\lesssim \left|\partial_{p_i}g\right|
\lesssim \frac{q^0}{g},\\
\left|\partial_{p_i}\left\{\frac{1}{\sqrt{\mathfrak{s}}}\right\}\right|
&=\left|\frac{g\,\partial_{p_i}g}{\sqrt{\mathfrak{s}}^{\,3}}\right|
\lesssim \frac{q^0}{g^3}.
\end{align*}
Here we use that on $|p|\ge \frac32 q^0$,
\begin{align*}
\mathfrak{c}\lesssim\frac{\mathfrak{c}|p-q|}{\sqrt{p^0 q^0}} \lesssim g \lesssim|p-q| \lesssim \langle p\rangle,
\end{align*}
where $\frac{|p|}{3}\le |p-q|\le 2|p|$ follows from $|p|\ge \frac32 q^0$. In particular, on this region,
\[
q^0\lesssim |p|,\qquad p^0\lesssim |p|,\qquad \sqrt{\mathfrak{s}}\ge g\gtrsim \mathfrak{c},
\]
and therefore
\begin{align}\label{A.13-0}
\frac{q^0}{g}\lesssim \frac{q^0}{\mathfrak{c}}
=\sqrt{1+\frac{|q|^2}{\mathfrak{c}^2}}
\lesssim \langle q\rangle .
\end{align}
Hence the desired bounds follow when $|\beta|=1$.

Assume now that the stated estimates hold for all multi-indices of length at most $|\beta|$. We prove them for derivatives of order $|\beta|+1$.

Since
\[
\partial_{p_i}g=\frac{1}{g}\left\{\frac{p_i}{p^0}q^0-q_i\right\},
\]
and for every multi-index $0\le \tilde\beta\le \beta$,
\[
\left|\partial_{\tilde\beta}\left\{\frac{p_i}{p^0}q^0-q_i\right\}\right|
\lesssim q^0,
\]
Leibniz's rule and the induction hypothesis imply
\begin{align*}
\left|\partial_{\beta}\partial_{p_i}g\right|
&\lesssim \sum_{\beta_1\le \beta}
\left|\partial_{\beta_1}\left\{\frac{1}{g}\right\}\right|
\left|\partial_{\beta-\beta_1}\left\{\frac{p_i}{p^0}q^0-q_i\right\}\right|\\
&\lesssim \frac{q^0}{g}
+\sum_{0\neq \beta_1\le \beta}
\frac{\langle q\rangle^{|\beta_1|}}{g^2}\,q^0
\lesssim \langle q\rangle^{|\beta|+1}.
\end{align*}

For the inverse of $g$, note that by the induction hypothesis and $g\gtrsim \mathfrak{c}$,
\[
\left|\partial_{\tilde\beta}\left\{\frac{1}{g}\right\}\right|
\lesssim \frac{\langle q\rangle^{|\tilde\beta|}}{g},
\qquad 0\le \tilde\beta\le \beta.
\]
Then
\begin{align*}
\left|\partial_{\beta}\partial_{p_i}\left\{\frac{1}{g}\right\}\right|
&=\left|\partial_{\beta}\left\{\frac{\partial_{p_i}g}{g^2}\right\}\right|\\
&\lesssim \sum_{\beta_1+\beta_2+\beta_3=\beta}
\left|\partial_{\beta_1}\left\{\frac{1}{g}\right\}\right|
\left|\partial_{\beta_2}\left\{\frac{1}{g}\right\}\right|
\left|\partial_{\beta_3}\partial_{p_i}g\right|\\
&\lesssim \sum_{\beta_1+\beta_2+\beta_3=\beta}
\frac{\langle q\rangle^{|\beta_1|}}{g}
\frac{\langle q\rangle^{|\beta_2|}}{g}
\langle q\rangle^{|\beta_3|+1}
\lesssim \frac{\langle q\rangle^{|\beta|+1}}{g^2}.
\end{align*}

Similarly, using
\[
\partial_{p_i}\sqrt{\mathfrak{s}}
=\frac{g\,\partial_{p_i}g}{\sqrt{\mathfrak{s}}},
\qquad
\partial_{p_i}\left\{\frac{1}{\sqrt{\mathfrak{s}}}\right\}
=-\frac{g\,\partial_{p_i}g}{\sqrt{\mathfrak{s}}^{\,3}},
\]
together with
\[
\left|\partial_{\tilde\beta}\left\{\frac{1}{\sqrt{\mathfrak{s}}}\right\}\right|
\lesssim \frac{\langle q\rangle^{|\tilde\beta|}}{g},
\qquad
\left|\partial_{\tilde\beta}g\right|
\lesssim \max\{g,\langle q\rangle^{|\tilde\beta|}\},
\qquad 0\le \tilde\beta\le \beta,
\]
we obtain
\begin{align*}
\left|\partial_{\beta}\partial_{p_i}\sqrt{\mathfrak{s}}\right|
&\lesssim
\sum_{\beta_1+\beta_2+\beta_3=\beta}
\left|\partial_{\beta_1}g\right|
\left|\partial_{\beta_2}\left\{\frac{1}{\sqrt{\mathfrak{s}}}\right\}\right|
\left|\partial_{\beta_3}\partial_{p_i}g\right|\lesssim \langle q\rangle^{|\beta|+1},
\\
\left|\partial_{\beta}\partial_{p_i}\left\{\frac{1}{\sqrt{\mathfrak{s}}}\right\}\right|
&\lesssim
\sum_{\beta_1+\beta_2+\beta_3+\beta_4+\beta_5=\beta}
\left|\partial_{\beta_1}\left\{\frac{1}{\sqrt{\mathfrak{s}}}\right\}\right|
\left|\partial_{\beta_2}\left\{\frac{1}{\sqrt{\mathfrak{s}}}\right\}\right|
\left|\partial_{\beta_3}\left\{\frac{1}{\sqrt{\mathfrak{s}}}\right\}\right|
\left|\partial_{\beta_4}g\right|
\left|\partial_{\beta_5}\partial_{p_i}g\right|\\
&\lesssim \frac{\langle q\rangle^{|\beta|+1}}{g^2}.
\end{align*}
This closes the induction and proves the lemma.
\end{proof}

We now use Lemma \ref{lA.4} to prove Lemma \ref{lA.3}.

\begin{proof}[\bf Proof of Lemma \ref{lA.3}]
We work throughout on the region $|p|\ge \frac32 q^0$, where
\begin{align}\label{A.13}
q^0\lesssim |p|,\qquad p^0\lesssim |p|,\qquad |p+q|\approx |p|,\qquad g\gtrsim \mathfrak{c}.
\end{align}

We first estimate $\partial_{\beta}(v_\phi\chi_{A^c})$. Recall from \eqref{04} that
\[
v_{\phi}=\frac{\mathfrak{c}}{4}\frac{g\sqrt{\mathfrak{s}}}{p^0 q^0}.
\]
By Leibniz's rule, Lemma \ref{lA.4}, and the standard bounds for momentum derivatives of $p^0$, we obtain, for some $n=n(\beta)$,
\begin{align*}
\left|\partial_{\beta} v_{\phi}\right|
&\lesssim
\sum_{\beta_1+\beta_2+\beta_3=\beta}
\left|\partial_{\beta_1}g\right|
\left|\partial_{\beta_2}\sqrt{\mathfrak{s}}\right|
\left|\partial_{\beta_3}\left(\frac{\mathfrak{c}}{p^0 q^0}\right)\right|\\
&\lesssim \langle q\rangle^n \frac{\mathfrak{c}}{4}\frac{g\sqrt{\mathfrak{s}}}{p^0 q^0}
\lesssim \langle q\rangle^n |v_\phi|.
\end{align*}
Next, by the explicit form of the smooth cutoff $\chi_{A^c}$ and the same argument used previously for cutoff derivatives, we have
\[
|\partial_\beta \chi_{A^c}(p,q)|\lesssim \langle q\rangle^n,
\qquad |\beta|\ge 1,
\]
with support contained in $\{|p|\ge \frac32 q^0\}$. Therefore, by Leibniz's rule,
\[
\frac{|\partial_\beta(v_\phi\chi_{A^c})|}{|v_\phi|}
\lesssim \langle q\rangle^n .
\]

It remains to estimate $\partial_{\beta}p'$ and $\partial_{\beta}q'$. By \eqref{05}, for $|\beta|>0$,
\begin{align}\label{A.13-1}
\left|\partial_{\beta} p^{\prime}_i\right|
\lesssim
\frac{\delta_{i j}}{2}\mathbf{1}_{\beta=e_j}
+\left|\frac{\partial_{\beta}g}{2}\omega_i\right|
+\sum_{\tilde\beta\le \beta}
\left|
\partial_{\beta-\tilde\beta}\left\{\frac{g}{2}(\gamma_0-1)\right\}
\partial_{\tilde\beta}\left(
(p_i+q_i)\frac{(p+q)\cdot\omega}{|p+q|^2}
\right)
\right|.
\end{align}
The first term is trivial, and the second term is controlled by Lemma \ref{lA.4}. For the geometric factor, using $|p+q|\approx |p|$ on $A^c$,
\begin{align}\label{A.13-2}
\left|
\partial_{\tilde\beta}\left(
(p_i+q_i)\frac{(p+q)\cdot\omega}{|p+q|^2}
\right)
\right|
\lesssim |p+q|^{-|\tilde\beta|}
\lesssim \langle p\rangle^{-|\tilde\beta|}.
\end{align}

It therefore suffices to control derivatives of
\[
H(p,q):=\frac{g}{2}(\gamma_0-1)
=\frac{g}{2}\left(\frac{p^0+q^0}{\sqrt{\mathfrak{s}}}-1\right).
\]
We first note that
\[
|H|
\le \frac{|g|}{2}\left(\frac{p^0+q^0}{\sqrt{\mathfrak{s}}}+1\right)
\lesssim \langle p\rangle.
\]
Hence, if $\beta-\tilde\beta=0$, then
\begin{align}\label{A.13-3}
\left|
H(p,q)\,
\partial_{\beta}\left(
(p_i+q_i)\frac{(p+q)\cdot\omega}{|p+q|^2}
\right)
\right|
\lesssim \langle p\rangle^{1-|\beta|}
\lesssim 1.
\end{align}

It remains to consider the case $|\beta-\tilde\beta|>0$. We claim that for every multi-index $\mu\neq 0$, there exists $n=n(\mu)$ such that
\begin{align}\label{A.13-4}
|\partial_\mu H|\lesssim \langle q\rangle^n.
\end{align}
For $|\mu|=1$, using Lemma \ref{lA.4},
\begin{align*}
\left|\partial_{p_i}H\right|
&\lesssim
\left|\left(\partial_{p_i}g\right)(\gamma_0-1)\right|
+\left|g\,\partial_{p_i}\left\{\frac{1}{\sqrt{\mathfrak{s}}}\right\}(p^0+q^0)\right|
+\left|\frac{g}{\sqrt{\mathfrak{s}}}\partial_{p_i}p^0\right|\\
&\lesssim
\frac{q^0}{g}\left|\frac{p^0+q^0}{\sqrt{\mathfrak{s}}}-1\right|
+\frac{q^0(p^0+q^0)}{g^2}
+1.
\end{align*}
Since $\sqrt{\mathfrak{s}}\ge g$ and 
\[
g^2\gtrsim \frac{\mathfrak{c}^2|p-q|^2}{p^0q^0}\gtrsim \frac{\mathfrak{c}^2|p|^2}{p^0q^0},
\]
we obtain from \eqref{A.13}
\[
\frac{q^0}{g}\left|\frac{p^0+q^0}{\sqrt{\mathfrak{s}}}-1\right|
\lesssim \frac{q^0(p^0+q^0)}{g^2}+\frac{q^0}{g}
\lesssim \frac{(q^0)^2}{\mathfrak{c}^2}+\langle q\rangle
\lesssim \langle q\rangle^2,
\]
and similarly
\[
\frac{q^0(p^0+q^0)}{g^2}
\lesssim \frac{(q^0)^2}{\mathfrak{c}^2}
\lesssim \langle q\rangle^2.
\]
Therefore $|\partial_{p_i}H|\lesssim \langle q\rangle^2$.

For higher derivatives, we use the same induction pattern as in Lemma \ref{lA.4}. In particular, by repeating the previous argument with one extra derivative and using \eqref{A.13-0}, one obtains the schematic bounds
\[
\left|\partial_{\mu}\partial_{p_i}g\right|
\lesssim \frac{\langle q\rangle^{n(\mu)}}{g},
\qquad
\left|\partial_{\mu}\partial_{p_i}\left\{\frac{1}{\sqrt{\mathfrak{s}}}\right\}\right|
\lesssim \frac{\langle q\rangle^{n(\mu)}}{g^3},
\qquad \mu\neq 0.
\]
Applying Leibniz's rule to
\[
H=\frac{g}{2}\left(\frac{p^0+q^0}{\sqrt{\mathfrak{s}}}-1\right),
\]
and using again $p^0+q^0\lesssim |p|$, $g^2\gtrsim \frac{\mathfrak{c}^2|p|^2}{p^0q^0}$, and $\frac{q^0}{\mathfrak{c}}\lesssim \langle q\rangle$, we conclude that
\[
|\partial_\mu H|\lesssim \langle q\rangle^{n(\mu)},
\qquad \mu\neq 0.
\]

Combining \eqref{A.13-1}-\eqref{A.13-4}, we conclude that 
\[
\left|\partial_{\beta} p^{\prime}_i\right|
\lesssim \langle q\rangle^n.
\]
The proof for $\left|\partial_{\beta} q^{\prime}_i\right|$ is identical. This completes the proof of Lemma \ref{lA.3}.
\end{proof}

\begin{lemma}\label{lA.5}
For any $\ell \geq 9$, we have the following estimate for the collision operator for $|\beta|>0$
\begin{align}\label{A.14}
\left|\left\langle\partial_{\beta} \Gamma_{\mathfrak{c}}^{A^c}\left(f_1, f_2\right), \partial_{\beta} f_3\right\rangle\right| \lesssim \left\|\partial_{\beta} f_3\right\|_{\infty, \ell}\sum_{\beta_1+\beta_2 \leq \beta}\left\|\partial_{\beta_1} f_1\right\|\left\|\partial_{\beta_2} f_2\right\| .
\end{align}
Furthermore, if $\chi(p)$ satisfies $|\partial_{\beta}\chi(p)| \lesssim e^{-\delta_1|p|}$ for some positive constant $\delta_1>0$, then we have the following estimate 
\begin{align}\label{A.15}
\left|\langle\partial_{\beta} \Gamma_{\mathfrak{c}}^{A^c}\left(f_1, \chi\right),\partial_{\beta} f_3\rangle \right|+\left|\left\langle\partial_{\beta} \Gamma_{\mathfrak{c}}^{A^c}\left(\chi, f_1\right),\partial_{\beta} f_3\right\rangle \right| \lesssim \sum_{\beta_1\leq \beta}\left\|\partial_{\beta_1} f_1\right\|_{\nu_{\mathfrak{c}}}\left\|\partial_{\beta} f_3\right\|_{\nu_{\mathfrak{c}}} ,
\end{align}
where the constants are independent of $\mathfrak{c}$.
\end{lemma}

\begin{proof}
From \eqref{A.12} and Lemma \ref{lA.3}, we have
\begin{align*}
\left|\partial_{\beta} \Gamma_{\mathfrak{c}}^{A^c}\right|
& \lesssim \sum_{\beta_1+\beta_2 \leq \beta}
\int_{\mathbb{R}^3 \times \mathbb{S}^2}
\mathbf{1}_{|p|\geq\frac{3}{2} q^0}
|v_{\phi}|\langle q\rangle^n\sqrt{\mathbf{M}_{\mathfrak{c}}(q)}
\left|\left(\partial_{\beta_1} f_1\right)\left(p^{\prime}\right)\left(\partial_{\beta_2} f_2\right)\left(q^{\prime}\right)\right|
\,d \omega d q
\\
&\quad
+\sum_{\beta_1+\beta_2 \leq \beta}
\left|\partial_{\beta_1} f_1(p)\right|
\int_{\mathbb{R}^3 \times \mathbb{S}^2}
\mathbf{1}_{|p|\geq\frac{3}{2} q^0}
|v_{\phi}|\langle q\rangle^n\sqrt{\mathbf{M}_{\mathfrak{c}}(q)}
\left|\partial_{\beta_2} f_2(q)\right|
\,d \omega d q
=: \uppercase\expandafter{\romannumeral1}+\uppercase\expandafter{\romannumeral2}.
\end{align*}

We first estimate $\uppercase\expandafter{\romannumeral2}$. Since
\[
\langle q\rangle^n\sqrt{\mathbf{M}_{\mathfrak{c}}(q)}\lesssim \langle q\rangle^{-\ell+1},
\qquad
v_\phi\le \frac{|p-q|}{2}\lesssim \langle p\rangle,
\]
by Lemma \ref{l2.1}, we obtain
\begin{align*}
\left|\langle \uppercase\expandafter{\romannumeral2},\partial_\beta f_3\rangle\right|
&\lesssim
\|\partial_\beta f_3\|_{\infty,\ell}
\sum_{\beta_1+\beta_2\le \beta}
\int_{\mathbb{R}^3\times\mathbb{R}^3\times\mathbb{R}^3\times\mathbb{S}^2}
\frac{|\partial_{\beta_1}f_1(p)\partial_{\beta_2}f_2(q)|}
{\langle p\rangle^{\ell-1}\langle q\rangle^{\ell-1}}
\,d\omega\,dq\,dp\,dx \\
&\lesssim
\|\partial_\beta f_3\|_{\infty,\ell}
\sum_{\beta_1+\beta_2\le \beta}
\|\partial_{\beta_1}f_1\|\,\|\partial_{\beta_2}f_2\|,
\end{align*}
where the last step follows from Cauchy-Schwarz (since $\ell\ge 9$).

For $\uppercase\expandafter{\romannumeral1}$, we similarly have
\begin{align*}
\left|\langle \uppercase\expandafter{\romannumeral1},\partial_\beta f_3\rangle\right|
&\lesssim
\|\partial_\beta f_3\|_{\infty,\ell}
\sum_{\beta_1+\beta_2\le \beta}
\int_{\mathbb{R}^3\times\mathbb{R}^3\times\mathbb{R}^3\times\mathbb{S}^2}
v_\phi
\frac{|\partial_{\beta_1}f_1(p')\partial_{\beta_2}f_2(q')|}
{\langle p\rangle^\ell\langle q\rangle^\ell}
\,d\omega\,dq\,dp\,dx \\
&\lesssim
\|\partial_\beta f_3\|_{\infty,\ell}
\sum_{\beta_1+\beta_2\le \beta}
\prod_{i=1,2}
\Bigg\{
\int_{\mathbb{R}^3\times\mathbb{R}^3\times\mathbb{R}^3\times\mathbb{S}^2}
v_\phi
\frac{|\partial_{\beta_i}f_i(p')|^2}
{\langle p\rangle^\ell\langle q\rangle^\ell}
\,d\omega\,dq\,dp\,dx
\Bigg\}^{1/2},
\end{align*}
where we used the $(p',q')$ symmetry in one factor. By the pre-post collisional change of variables \eqref{08},
\[
\int v_\phi
\frac{|\partial_{\beta_i}f_i(p')|^2}
{\langle p\rangle^\ell\langle q\rangle^\ell}
\,d\omega\,dq\,dp\,dx
=
\int v_\phi
\frac{|\partial_{\beta_i}f_i(p)|^2}
{\langle p'\rangle^\ell\langle q'\rangle^\ell}
\,d\omega\,dq\,dp\,dx.
\]
Using
\[
\frac{1}{\langle p'\rangle^\ell\langle q'\rangle^\ell}
\lesssim
\frac{1}{\langle p\rangle^{\ell/2}\langle q\rangle^{\ell/2}},
\qquad
v_\phi\lesssim \langle p\rangle,
\]
we obtain, since $\ell\ge 9$,
\[
\int v_\phi
\frac{|\partial_{\beta_i}f_i(p)|^2}
{\langle p'\rangle^\ell\langle q'\rangle^\ell}
\,d\omega\,dq\,dp\,dx
\lesssim
\int
\frac{\langle p\rangle}{\langle p\rangle^{\ell/2}\langle q\rangle^{\ell/2}}
|\partial_{\beta_i}f_i(p)|^2
\,d\omega\,dq\,dp\,dx
\lesssim \|\partial_{\beta_i}f_i\|^2.
\]
Combining the estimates for $\uppercase\expandafter{\romannumeral1}$ and $\uppercase\expandafter{\romannumeral2}$ yields \eqref{A.14}.

We now prove \eqref{A.15}. We only consider $\Gamma_{\mathfrak{c}}^{A^c}(f_1,\chi)$, since the term $\Gamma_{\mathfrak{c}}^{A^c}(\chi,f_1)$ is treated in the same way. By the above decomposition with $f_2=\chi$, we have
\begin{align*}
\left|\langle \uppercase\expandafter{\romannumeral2},\partial_\beta f_3\rangle\right|
&\lesssim
\sum_{\beta_1\le \beta}
\int_{\mathbb{R}^3\times\mathbb{R}^3\times\mathbb{R}^3\times\mathbb{S}^2}
\mathbf{1}_{|p|\ge\frac32 q^0}
v_\phi \langle q\rangle^n \sqrt{\mathbf{M}_{\mathfrak{c}}(q)} e^{-\delta_1|q|}
|\partial_{\beta_1}f_1(p)\partial_\beta f_3(p)|
\,d\omega\,dq\,dp\,dx \\
&\lesssim
\sum_{\beta_1\le \beta}
\|\partial_{\beta_1}f_1\|_{\nu_{\mathfrak{c}}}
\|\partial_\beta f_3\|_{\nu_{\mathfrak{c}}},
\end{align*}
where we used $\langle q\rangle^n e^{-\delta_1|q|}\lesssim 1$ and
\[
\int_{\mathbb{R}^3\times\mathbb{S}^2} v_\phi \sqrt{\mathbf{M}_{\mathfrak{c}}(q)}\,d\omega\,dq
\lesssim \nu_{\mathfrak{c}}(p).
\]

For $\uppercase\expandafter{\romannumeral1}$, Cauchy--Schwarz and the pre-post change of variables give
\begin{align*}
\left|\langle \uppercase\expandafter{\romannumeral1},\partial_\beta f_3\rangle\right|
&\lesssim
\sum_{\beta_1\le \beta}
\prod_{i=1,3}
\Bigg\{
\int_{\mathbb{R}^3\times\mathbb{R}^3\times\mathbb{R}^3\times\mathbb{S}^2}
\mathbf{1}_{|p|\ge\frac32 q^0}
v_\phi e^{-\delta_1|q|}e^{-\delta_1|q'|}
|\partial_{\beta_i}f_i(p)|^2
\,d\omega\,dq\,dp\,dx
\Bigg\}^{1/2} \\
&\lesssim
\sum_{\beta_1\le \beta}
\|\partial_{\beta_1}f_1\|_{\nu_{\mathfrak{c}}}
\|\partial_\beta f_3\|_{\nu_{\mathfrak{c}}},
\end{align*}
where we used
\[
\int_{\mathbb{R}^3\times\mathbb{S}^2} v_\phi e^{-\delta_1|q|}\,d\omega\,d q
\lesssim \nu_{\mathfrak{c}}(p).
\]
This proves \eqref{A.15}, and hence completes the proof.
\end{proof}
    \subsection{Proof of Lemma \ref{l2.5}}\label{appendix A.3.}

\begin{proof}
To prove \eqref{A.16}, we estimate the two parts $\Gamma_{\mathfrak{c}}^A$ and $\Gamma_{\mathfrak{c}}^{A^c}$ separately.

For the contribution on $A$, it follows from \eqref{A.6} and Lemma \ref{lA.1} that
\begin{align*}
\left|\partial_{\beta}^{\alpha}\Gamma_{\mathfrak{c}}^A(f_1,f_2)(p)\right|
\lesssim I_1+I_2,
\end{align*}
where
\begin{align*}
I_1
:={}&
\sum_{\substack{\alpha_1+\alpha_2\le \alpha\\ \beta_1+\beta_2\le \beta}}
\int_{\mathbb R^3\times\mathbb S^2}
\mathbf 1_{|p|\le 3q^0}\,
\frac{\mathfrak s\mathfrak B}{p^0q^0}\,
\langle q\rangle^n\sqrt{\mathbf M_{\mathfrak{c}}(q)}
\left|
(\partial_{\beta_1}^{\alpha_1}f_1)(p'')
(\partial_{\beta_2}^{\alpha_2}f_2)(q'')
\right|
\,d\omega\,dq,\\
I_2
:={}&
\sum_{\substack{\alpha_1+\alpha_2\le \alpha\\ \beta_1+\beta_2\le \beta}}
\int_{\mathbb R^3\times\mathbb S^2}
\mathbf 1_{|p|\le 3q^0}\,
\frac{\mathfrak s\mathfrak B}{p^0q^0}\,
\langle q\rangle^n\sqrt{\mathbf M_{\mathfrak{c}}(q)}
\left|
(\partial_{\beta_1}^{\alpha_1}f_1)(p)
(\partial_{\beta_2}^{\alpha_2}f_2)(q)
\right|
\,d\omega\,dq.
\end{align*}
On $A$ we use
\[
\frac{\mathfrak s\mathfrak B}{p^0q^0}\lesssim \nu_{\mathfrak{c}}(p)\langle q\rangle^4.
\]
Hence, by the assumptions on $\partial_{\beta_1}^{\alpha_1}f_1$ and $\partial_{\beta_2}^{\alpha_2}f_2$,
\begin{align*}
I_1
&\lesssim
\nu_{\mathfrak{c}}(p)
\int_{\mathbb R^3\times\mathbb S^2}
\langle q\rangle^{n+4}\sqrt{\mathbf M_{\mathfrak{c}}(q)}
\mathbf M_{\mathfrak{c}}^{\lambda/2}(p'')
\mathbf M_{\mathfrak{c}}^{\lambda/2}(q'')
\,d\omega\,dq.
\end{align*}
By conservation of energy-momentum,
\[
\mathbf M_{\mathfrak{c}}^{\lambda/2}(p'')\mathbf M_{\mathfrak{c}}^{\lambda/2}(q'')
=
\mathbf M_{\mathfrak{c}}^{\lambda/2}(p)\mathbf M_{\mathfrak{c}}^{\lambda/2}(q),
\]
and therefore
\begin{align*}
I_1
&\lesssim
\nu_{\mathfrak{c}}(p)\mathbf M_{\mathfrak{c}}^{\lambda/2}(p)
\int_{\mathbb R^3\times\mathbb S^2}
\langle q\rangle^{n+4}\mathbf M_{\mathfrak{c}}^{(1+\lambda)/2}(q)\,d\omega\,dq
\lesssim
\nu_{\mathfrak{c}}(p)\mathbf M_{\mathfrak{c}}^{\lambda/2}(p).
\end{align*}
The same argument yields
\[
I_2\lesssim \nu_{\mathfrak{c}}(p)\mathbf M_{\mathfrak{c}}^{\lambda/2}(p).
\]
Consequently,
\begin{align}\label{A.18}
\left|\partial_{\beta}^{\alpha}\Gamma_{\mathfrak{c}}^A(f_1,f_2)(p)\right|
\lesssim
\nu_{\mathfrak{c}}(p)\mathbf M_{\mathfrak{c}}^{\lambda/2}(p).
\end{align}

For the contribution on $A^c$, \eqref{A.12} and Lemma \ref{lA.3} imply
\begin{align*}
\left|\partial_{\beta}^{\alpha}\Gamma_{\mathfrak{c}}^{A^c}(f_1,f_2)(p)\right|
\lesssim J_1+J_2,
\end{align*}
where
\begin{align*}
J_1
:={}&
\sum_{\substack{\alpha_1+\alpha_2\le \alpha\\ \beta_1+\beta_2\le \beta}}
\int_{\mathbb R^3\times\mathbb S^2}
\mathbf 1_{|p|\ge \frac32 q^0}\,
|v_\phi|\,
\langle q\rangle^n\sqrt{\mathbf M_{\mathfrak{c}}(q)}
\left|
(\partial_{\beta_1}^{\alpha_1}f_1)(p')
(\partial_{\beta_2}^{\alpha_2}f_2)(q')
\right|
\,d\omega\,dq,\\
J_2
:={}&
\sum_{\substack{\alpha_1+\alpha_2\le \alpha\\ \beta_1+\beta_2\le \beta}}
\left|(\partial_{\beta_1}^{\alpha_1}f_1)(p)\right|
\int_{\mathbb R^3\times\mathbb S^2}
\mathbf 1_{|p|\ge \frac32 q^0}\,
|v_\phi|\,
\langle q\rangle^n\sqrt{\mathbf M_{\mathfrak{c}}(q)}
\left|
(\partial_{\beta_2}^{\alpha_2}f_2)(q)
\right|
\,d\omega\,dq.
\end{align*}
Using again the hypotheses on $f_1$ and $f_2$, together with
\[
\mathbf M_{\mathfrak{c}}^{\lambda/2}(p')\mathbf M_{\mathfrak{c}}^{\lambda/2}(q')
=
\mathbf M_{\mathfrak{c}}^{\lambda/2}(p)\mathbf M_{\mathfrak{c}}^{\lambda/2}(q),
\]
we obtain
\begin{align*}
J_1
&\lesssim
\mathbf M_{\mathfrak{c}}^{\lambda/2}(p)
\int_{\mathbb R^3\times\mathbb S^2}
\mathbf 1_{|p|\ge \frac32 q^0}\,
|v_\phi|\,
\langle q\rangle^n
\mathbf M_{\mathfrak{c}}^{(1+\lambda)/2}(q)
\,d\omega\,dq
\lesssim
\nu_{\mathfrak{c}}(p)\mathbf M_{\mathfrak{c}}^{\lambda/2}(p),
\end{align*}
where we used the exponential decay of $\mathbf M_{\mathfrak{c}}$ to absorb the polynomial weight in $q$. The same estimate holds for $J_2$. Hence,
\begin{align}\label{A.19}
\left|\partial_{\beta}^{\alpha}\Gamma_{\mathfrak{c}}^{A^c}(f_1,f_2)(p)\right|
\lesssim
\nu_{\mathfrak{c}}(p)\mathbf M_{\mathfrak{c}}^{\lambda/2}(p).
\end{align}
Combining \eqref{A.18} and \eqref{A.19}, we obtain \eqref{A.16}.

To prove \eqref{19}, recall that
\[
\mathbf L_{\mathfrak{c}}f
=
-\Gamma_{\mathfrak{c}}\bigl(\sqrt{\mathbf M_{\mathfrak{c}}},f\bigr)
-\Gamma_{\mathfrak{c}}\bigl(f,\sqrt{\mathbf M_{\mathfrak{c}}}\bigr).
\]
Hence
\begin{align}\label{A.20}
[\partial_\beta^\alpha,\mathbf L_{\mathfrak{c}}]f
={}&
-\Big\{
\partial_\beta^\alpha\Gamma_{\mathfrak{c}}\bigl(\sqrt{\mathbf M_{\mathfrak{c}}},f\bigr)
-\Gamma_{\mathfrak{c}}\bigl(\sqrt{\mathbf M_{\mathfrak{c}}},\partial_\beta^\alpha f\bigr)
\Big\} \nonumber\\
&\quad
-\Big\{
\partial_\beta^\alpha\Gamma_{\mathfrak{c}}\bigl(f,\sqrt{\mathbf M_{\mathfrak{c}}}\bigr)
-\Gamma_{\mathfrak{c}}\bigl(\partial_\beta^\alpha f,\sqrt{\mathbf M_{\mathfrak{c}}}\bigr)
\Big\}:=Z_1+Z_2.
\end{align}
By definitions of $\Gamma_{\mathfrak{c}}^A$ and $\Gamma_{\mathfrak{c}}^{A^c}$ in \eqref{A.2}-\eqref{A.3}, together with Lemmas \ref{lA.1}-\ref{lA.3}, we obtain
\begin{align*}
|Z_1|
\lesssim&
\sum_{\substack{\alpha_1+\alpha_2\le \alpha,\ \beta_1+\beta_2\le \beta\\
|\alpha_2|+|\beta_2|<|\alpha|+|\beta|}}\Big\{
\int_{\mathbb R^3\times\mathbb S^2}
\mathbf 1_{|p|\le 3q^0}\,
\frac{\mathfrak s\mathfrak B}{p^0q^0}\,
\langle q\rangle^n\sqrt{\mathbf M_{\mathfrak{c}}(q)}
\left|
(\partial_{\beta_1}^{\alpha_1}\sqrt{\mathbf M_{\mathfrak{c}}})(p'')
(\partial_{\beta_2}^{\alpha_2}f_2)(q'')
\right|
\,d\omega\,dq\\
&\qquad+
\int_{\mathbb R^3\times\mathbb S^2}
\mathbf 1_{|p|\le 3q^0}\,
\frac{\mathfrak s\mathfrak B}{p^0q^0}\,
\langle q\rangle^n\sqrt{\mathbf M_{\mathfrak{c}}(q)}
\left|
(\partial_{\beta_1}^{\alpha_1}\sqrt{\mathbf M_{\mathfrak{c}}})(p)
(\partial_{\beta_2}^{\alpha_2}f_2)(q)
\right|
\,d\omega\,dq\\
&\qquad+
\int_{\mathbb R^3\times\mathbb S^2}
\mathbf 1_{|p|\ge \frac32 q^0}\,
|v_\phi|\,
\langle q\rangle^n\sqrt{\mathbf M_{\mathfrak{c}}(q)}
\left|
(\partial_{\beta_1}^{\alpha_1}\sqrt{\mathbf M_{\mathfrak{c}}})(p')
(\partial_{\beta_2}^{\alpha_2}f_2)(q')
\right|
\,d\omega\,dq\\
&\qquad+
\int_{\mathbb R^3\times\mathbb S^2}
\mathbf 1_{|p|\ge \frac32 q^0}\,
|v_\phi|\,
\langle q\rangle^n\sqrt{\mathbf M_{\mathfrak{c}}(q)}\left|(\partial_{\beta_1}^{\alpha_1}\sqrt{\mathbf M_{\mathfrak{c}}})(p)
(\partial_{\beta_2}^{\alpha_2}f_2)(q)
\right|
\,d\omega\,dq \Big\}.
\end{align*}
By \eqref{e2.35},
\[
\sum_{\widetilde\alpha\le \alpha,\ \widetilde\beta\le \beta}
\bigl|
\partial_{\widetilde\beta}^{\widetilde\alpha}
\sqrt{\mathbf M_{\mathfrak{c}}}(p)
\bigr|
\lesssim
\mathbf M_{\mathfrak{c}}^{\lambda/2}(p),
\]
and by hypothesis,
\[
\sum_{\substack{\widetilde\alpha\le \alpha,\ \widetilde\beta\le \beta\\
|\widetilde\alpha|+|\widetilde\beta|<|\alpha|+|\beta|}}
\bigl|
\partial_{\widetilde\beta}^{\widetilde\alpha}f(p)
\bigr|
\lesssim
\mathbf M_{\mathfrak{c}}^{\lambda/2}(p).
\]
A derivation similar to that in proof of \eqref{A.16} gives
\[
|Z_1|
\lesssim
\nu_{\mathfrak{c}}(p)\mathbf M_{\mathfrak{c}}^{\lambda/2}(p).
\]
By symmetry, the same argument applies to $Z_1$. Thus we obtain
\[
|[\partial_\beta^\alpha,\mathbf L_{\mathfrak{c}}]f|
\lesssim
\nu_{\mathfrak{c}}(p)\mathbf M_{\mathfrak{c}}^{\lambda/2}(p),
\]
which proves \eqref{19}.

Finally, we prove \eqref{20}. By the definition of $\nu_{\mathfrak{c}}(p)$, we decompose
\[
\nu_{\mathfrak{c}}(p)=\nu_{\mathfrak{c}}^A(p)+\nu_{\mathfrak{c}}^{A^c}(p),
\]
where
\begin{align*}
\nu_{\mathfrak{c}}^A(p)
&=
\int_{\mathbb R^3\times\mathbb S^2}
\frac{\mathfrak s\mathfrak B}{p^0q^0}\,
\sqrt{\mathbf M_{\mathfrak{c}}(q)}\chi_A(p, q)
\,d\omega\,d q,\\
\nu_{\mathfrak{c}}^{A^c}(p)
&=
\int_{\mathbb R^3\times\mathbb S^2}
|v_\phi|\,
\sqrt{\mathbf M_{\mathfrak{c}}(q)}\chi_{A^c}(p, q)
\,d\omega\,d q.
\end{align*}
Differentiating under the integral sign and using Lemmas \ref{lA.1} and \ref{lA.3}, we obtain
\begin{align*}
\left|\partial_{\beta}\nu_{\mathfrak{c}}^A(p)\right|
\lesssim &
\int_{\mathbb R^3\times\mathbb S^2}
\mathbf 1_{|p|\le 3q^0}\,
\frac{\mathfrak s\mathfrak B}{p^0q^0}\,
\langle q\rangle^n
\sqrt{\mathbf M_{\mathfrak{c}}(q)}
\,d\omega\,d q\\
=& \int_{\mathbb R^3\times\mathbb S^2}
\mathbf 1_{|p|\le 3q^0}\,
v_\phi\,
\langle q\rangle^n
\sqrt{\mathbf M_{\mathfrak{c}}(q)}
\,d\omega\,d q
\lesssim \nu_{\mathfrak{c}}(p).\\
\left|\partial_{\beta}\nu_{\mathfrak{c}}^{A^c}(p)\right|
\lesssim &
\int_{\mathbb R^3\times\mathbb S^2}
\mathbf 1_{|p|\ge \frac32 q^0}\,
|v_\phi|\,
\langle q\rangle^n
\sqrt{\mathbf M_{\mathfrak{c}}(q)}
\,d\omega\,dq
\lesssim \nu_{\mathfrak{c}}(p).
\end{align*}
Therefore,
\[
\left|\partial_{\beta}\nu_{\mathfrak{c}}(p)\right|
\lesssim \nu_{\mathfrak{c}}(p),
\]
which proves \eqref{20}. The proof of Lemma \ref{l2.5} is complete.
\end{proof}
 \subsection{Proof of Lemma \ref{l2.13}}\label{appendix A.4}	
	\begin{proof}
We first prove \eqref{e2.44}. We begin with $\mathcal{K}_{\mathfrak{c}_1}$.
By definition,
\[
\mathcal{K}_{\mathfrak{c}_1}f
=
\int_{\mathbb R^3\times\mathbb S^2}
v_\phi(p,q)\,\sqrt{J_\mathbf{M}(q)}\,\frac{\mathbf M_{\mathfrak{c}}(p)}{\sqrt{J_\mathbf{M}(p)}}\,f(q)\,d\omega\,d q.
\]
Since $p,q,p',q'$ are independent of $x$, we write
\[
\partial_{x_j}\big(\mathcal{K}_{\mathfrak{c}_1}f\big)
=
\big[\partial_{x_j},\mathcal{K}_{\mathfrak{c}_1}\big]f
+
\mathcal{K}_{\mathfrak{c}_1}(\partial_{x_j}f),
\]
where
\[
\big[\partial_{x_j},\mathcal{K}_{\mathfrak{c}_1}\big]f
=
\int_{\mathbb R^3\times\mathbb S^2}
v_\phi(p,q)\,\sqrt{J_\mathbf{M}(q)}\,
\frac{\partial_{x_j}\mathbf M_{\mathfrak{c}}(p)}{\sqrt{J_\mathbf{M}(p)}}
f(q)\,d\omega\,d q.
\]
Using \eqref{e2.35} together with \eqref{e1.15}, we obtain
\[
\left|
\big[\partial_{x_j},\mathcal{K}_{\mathfrak{c}_1}\big]f
\right|
\lesssim
\int_{\mathbb R^3\times\mathbb S^2}
v_\phi(p,q)\langle p\rangle^2\,J_\mathbf{M}^{\frac{1}{2}}(q)J_\mathbf{M}^{\alpha-\frac{1}{2}}(p)\,|f(q)|\,d\omega\,d q
\lesssim
\int_{\mathbb R^3}\bar{k}_1(p,q)|f(q)|\,d q,
\]
where 
$$
\bar{k}_1(p,q):=\int_{\mathbb S^2}
v_\phi(p,q)\,J_\mathbf{M}^{\frac{\delta}{2}}(p)J_\mathbf{M}^{\frac{\delta}{2}}(q)\,d\omega
$$
with $\delta:=\alpha-\frac{1}{2}>0$. By similar arguments as in \cite[Lemmas 4.3]{Wang-Xiao-JLMS-2026}, one has
	\begin{align*}
	    \bar{k}_1(p, q)\lesssim \hat{k}_1(p, q),
	\end{align*}
	where 
	\begin{align}\label{K.0}
	\hat{k}_1(p, q):=|p-q| e^{-\frac{\delta}{2T_M}|p|} e^{-\frac{\delta}{2T_M}|q|}.
	\end{align}

On the other hand, by \eqref{e2.42}-\eqref{e2.43},
\[
\left|\mathcal{K}_{\mathfrak{c}_1}(\partial_{x_j}f)\right|
\lesssim
\int_{\mathbb R^3} k_1(p,q)\,|\partial_{x_j}f(q)|\,d q.
\]
Hence
\begin{align}\label{K.1}
\left|\partial_{x_j}\big(\mathcal{K}_{\mathfrak{c}_1}f\big)\right|
\lesssim
\int_{\mathbb R^3}\hat k_1(p,q)|f(q)|\,d q
+
\int_{\mathbb R^3} k_1(p,q)\,|\partial_{x_j}f(q)|\,d q.
\end{align}

We next consider $\mathcal{K}_{\mathfrak{c}_2}$. By definition,
\begin{align*}
\mathcal{K}_{\mathfrak{c}_2}f
=
\int_{\mathbb R^3\times\mathbb S^2}
v_\phi(p,q)
\Bigg\{
\mathbf M_{\mathfrak{c}}(p')\frac{\sqrt{J_\mathbf{M}(q')}}{\sqrt{J_\mathbf{M}(p)}}\,f(q')
+
\mathbf M_{\mathfrak{c}}(q')\frac{\sqrt{J_\mathbf{M}(p')}}{\sqrt{J_\mathbf{M}(p)}}\,f(p')
\Bigg\}
\,d\omega\,d q .
\end{align*}
Again $p,q,p',q'$ are independent of $x$, so
\[
\partial_{x_j}\big(\mathcal{K}_{\mathfrak{c}_2}f\big)
=\big[\partial_{x_j},\mathcal{K}_{\mathfrak{c}_2}\big]f
+\mathcal{K}_{\mathfrak{c}_2}(\partial_{x_j}f),
\]
where
$$
\big[\partial_{x_j},\mathcal{K}_{\mathfrak{c}_2}\big]f
=
\int_{\mathbb R^3\times\mathbb S^2}
v_\phi(p,q)
\Big\{
\partial_{x_j}\mathbf M_{\mathfrak{c}}(p')\frac{\sqrt{J_\mathbf{M}(q')}}{\sqrt{J_\mathbf{M}(p)}}\,f(q')
+\partial_{x_j}\mathbf M_{\mathfrak{c}}(q')\frac{\sqrt{J_\mathbf{M}(p')}}{\sqrt{J_\mathbf{M}(p)}}\,f(p')
\Big\}
\,d\omega\,d q .
$$
Using \eqref{e2.35} and \eqref{e1.15} exactly as above, together with $J_\mathbf{M}^{\frac{1}{2}}(p')J_\mathbf{M}^{\frac{1}{2}}(q')=J_\mathbf{M}^{\frac{1}{2}}(p)J_\mathbf{M}^{\frac{1}{2}}(q)$, we obtain
\begin{align}\label{K.2}
\left|
\big[\partial_{x_j},\mathcal{K}_{\mathfrak{c}_2}\big]f
\right| \lesssim &
\int_{\mathbb R^3\times\mathbb S^2}
v_\phi(p,q)J_\mathbf{M}^{\frac{1}{2}}(q)
\Big\{\langle p^\prime \rangle^2 J_\mathbf{M}^{\alpha-\frac{1}{2}}(p^\prime)\,|f(q')|
+\langle q^\prime \rangle^2 J_\mathbf{M}^{\alpha-\frac{1}{2}}(q^\prime)\,|f(p')|\Big\}\,d\omega\,d q \nonumber\\
\lesssim &\int_{\mathbb R^3\times\mathbb S^2}
v_\phi(p,q)J_\mathbf{M}^{\frac{\delta}{2}}(q)
\Big\{ J_\mathbf{M}^{\frac{\delta}{2}}(p^\prime)\,|f(q')|+ J_\mathbf{M}^{\frac{\delta}{2}}(q^\prime)\,|f(p')|\Big\}\,d\omega\,d q.
\end{align}
Similar to \eqref{e2.3}, by the symmetry of the collision kernel for identical particles, the two gain terms in \eqref{K.2} give the same contribution after relabeling the post-collisional variables, and hence can be written in the following unified kernel form
\begin{align*}
    \left|
\big[\partial_{x_j},\mathcal{K}_{\mathfrak{c}_2}\big]f
\right| \lesssim
\int_{\mathbb R^3}\bar k_2(p,q)|f(q)|\,dq,
\end{align*}
where
$$
\bar k_2(p,q)=\frac{\mathfrak{c}}{p^0 q^0} \int_{\mathbb{R}^3} \frac{d q^{\prime}}{q^{\prime 0}} \int_{\mathbb{R}^3} \frac{d p^{\prime}}{p^{\prime 0}} \bar{\mathfrak{s}} \delta^{(4)}\left(p^\mu+p^{\prime \mu}-q^\mu-q^{\prime \mu}\right) J_\mathbf{M}^{\frac{\delta}{2}}(p')J_\mathbf{M}^{\frac{\delta}{2}}(q').
$$
Then, By similar arguments as in \cite[Lemmas 4.3]{Wang-Xiao-JLMS-2026}, one has
	\begin{align}\label{K.3}
	    \bar{k}_2(p, q)\lesssim \hat{k}_2(p, q),
	\end{align}
	where 
	\begin{align*}
	\hat{k}_2(p, q):=\frac{1}{|p-q|} e^{-\frac{\delta}{4T_M}|p-q|}.
	\end{align*}

Moreover, by \eqref{e2.42}-\eqref{e2.43},
\[
\left|\mathcal{K}_{\mathfrak{c}_2}(\partial_{x_j}f)\right|
\lesssim
\int_{\mathbb R^3} k_2(p,q)\,|\partial_{x_j}f(q)|\,dq.
\]
Therefore,
\begin{align*}
\left|\partial_{x_j}\big(\mathcal{K}_{\mathfrak{c}_2}f\big)\right|
\lesssim
\int_{\mathbb R^3}\hat k_2(p,q)|f(q)|\,dq
+
\int_{\mathbb R^3} k_2(p,q)\,|\partial_{x_j}f(q)|\,dq,
\end{align*}
which, together with \eqref{K.1}, yields \eqref{e2.44}.

We now prove \eqref{e2.40-1}. As in the proof of Lemmas \ref{lA.2}--\ref{lA.5}, we decompose
\[
\mathcal{K}_{\mathfrak{c}_1}f
=
\mathcal{K}_{\mathfrak{c}_1}^{A}f+\mathcal{K}_{\mathfrak{c}_1}^{A^c}f,
\]
where
\begin{align*}
\mathcal{K}_{\mathfrak{c}_1}^{A}f
&=
\int_{\mathbb R^3\times\mathbb S^2}
\frac{\mathfrak s\mathfrak B}{p^0q^0}\chi_A(p,q)\,
\sqrt{J_\mathbf{M}(q)}\,\frac{\mathbf M_{\mathfrak{c}}(p)}{\sqrt{J_\mathbf{M}(p)}}\,f(q)\,d\omega\,dq,\\
\mathcal{K}_{\mathfrak{c}_1}^{A^c}f
&=
\int_{\mathbb R^3\times\mathbb S^2}
v_\phi\,\chi_{A^c}(p,q)\,
\sqrt{J_\mathbf{M}(q)}\,\frac{\mathbf M_{\mathfrak{c}}(p)}{\sqrt{J_\mathbf{M}(p)}}\,f(q)\,d\omega\,dq.
\end{align*}

Differentiating $\mathcal{K}_{\mathfrak{c}_1}^{A}f$ with respect to $p_i$, we get
\begin{align*}
\left|\partial_{p_i}\mathcal{K}_{\mathfrak{c}_1}^{A}f\right|
\leq {}&
\int_{\mathbb R^3\times\mathbb S^2}
\mathbf 1_{|p|\leq 3q^0}
\left|
\partial_{p_i}\!\left(
\frac{\mathfrak s\mathfrak B}{p^0q^0}\chi_A(p,q)
\right)
\right|
\sqrt{J_\mathbf{M}(q)}\,\frac{\mathbf M_{\mathfrak{c}}(p)}{\sqrt{J_\mathbf{M}(p)}}\,|f(q)|\,d\omega\,d q\\
&+
\int_{\mathbb R^3\times\mathbb S^2}
\mathbf 1_{|p|\leq 3q^0}
\frac{\mathfrak s\mathfrak B}{p^0q^0}\chi_A(p,q)\,
\sqrt{J_\mathbf{M}(q)}\,
\Big|
\partial_{p_i}\!\Big(\frac{\mathbf M_{\mathfrak{c}}(p)}{\sqrt{J_\mathbf{M}(p)}}\Big)
\Big|
\,|f(q)|\,d\omega\,d q.
\end{align*}
By Lemma \ref{lA.1}, \eqref{e2.35} and \eqref{e.15}, we obtain
\begin{align}\label{K.4}
\left|\partial_{p_i}\mathcal{K}_{\mathfrak{c}_1}^{A}f\right|
\lesssim
\int_{\mathbb R^3\times\mathbb S^2}
\mathbf 1_{|p|\leq 3q^0}
\frac{\mathfrak s\mathfrak B}{p^0q^0}\,
J_\mathbf{M}^{\frac{\delta}{2}}(p)\,J_\mathbf{M}^{\frac{\delta}{2}}(q)\,|f(q)|\,d\omega\,d q.
\end{align}

Similarly,
\begin{align*}
\left|\partial_{p_i}\mathcal{K}_{\mathfrak{c}_1}^{A^c}f\right|
\leq {}&
\int_{\mathbb R^3\times\mathbb S^2}
\mathbf 1_{|p|\geq \frac32 q^0}
\left|
\partial_{p_i}\big(v_\phi\chi_{A^c}(p,q)\big)
\right|
\sqrt{J_\mathbf{M}(q)}\,\frac{\mathbf M_{\mathfrak{c}}(p)}{\sqrt{J_\mathbf{M}(p)}}\,|f(q)|\,d\omega\,d q\\
&+
\int_{\mathbb R^3\times\mathbb S^2}
\mathbf 1_{|p|\geq \frac32 q^0}
v_\phi\chi_{A^c}(p,q)\,
\sqrt{J_\mathbf{M}(q)}\,
\Big|
\partial_{p_i}\!\Big(\frac{\mathbf M_{\mathfrak{c}}(p)}{\sqrt{J_\mathbf{M}(p)}}\Big)
\Big|
\,|f(q)|\,d\omega\,d q.
\end{align*}
By Lemma \ref{lA.3}, \eqref{e2.35} and \eqref{e.15},
\begin{align}\label{K.5}
\left|\partial_{p_i}\mathcal{K}_{\mathfrak{c}_1}^{A^c}f\right|
\lesssim
\int_{\mathbb R^3\times\mathbb S^2}
\mathbf 1_{|p|\geq \frac32 q^0}
v_\phi\,J_\mathbf{M}^{\frac{\delta}{2}}(p)\,J_\mathbf{M}^{\frac{\delta}{2}}(q)\,|f(q)|\,d\omega\,d q.
\end{align}

Combining \eqref{K.4} and \eqref{K.5} by \eqref{07}, we deduce
\[
\left|\partial_{p_i}\mathcal{K}_{\mathfrak{c}_1}f\right|
\lesssim
\int_{\mathbb R^3}\bar k_1(p,q)|f(q)|\,d q
\lesssim
\int_{\mathbb R^3}\hat k_1(p,q)|f(q)|\,d q,
\]
which proves \eqref{e2.40-1}.

Finally, we prove \eqref{e2.40-2}. Decompose
\[
\mathcal{K}_{\mathfrak{c}_2}f
=
\mathcal{K}_{\mathfrak{c}_2}^{A}f+\mathcal{K}_{\mathfrak{c}_2}^{A^c}f,
\]
where
\begin{align*}
\mathcal{K}_{\mathfrak{c}_2}^{A}f
:={}&
\int_{\mathbb R^3\times\mathbb S^2}
\frac{\mathfrak s\mathfrak B}{p^0q^0}\chi_A(p,q)
\Bigg\{
\mathbf M_{\mathfrak{c}}(p'')\frac{\sqrt{J_\mathbf{M}(q'')}}{\sqrt{J_\mathbf{M}(p)}}\,f(q'')
+
\mathbf M_{\mathfrak{c}}(q'')\frac{\sqrt{J_\mathbf{M}(p'')}}{\sqrt{J_\mathbf{M}(p)}}\,f(p'')
\Bigg\}
\,d\omega\,d q,\\
\mathcal{K}_{\mathfrak{c}_2}^{A^c}f
:={}&
\int_{\mathbb R^3\times\mathbb S^2}
v_\phi\chi_{A^c}(p,q)
\Bigg\{
\mathbf M_{\mathfrak{c}}(p')\frac{\sqrt{J_\mathbf{M}(q')}}{\sqrt{J_\mathbf{M}(p)}}\,f(q')
+
\mathbf M_{\mathfrak{c}}(q')\frac{\sqrt{J_\mathbf{M}(p')}}{\sqrt{J_\mathbf{M}(p)}}\,f(p')
\Bigg\}
\,d\omega\,d q .
\end{align*}

We first treat the $A$-region. Differentiating with respect to $p_i$, we split
\[
\partial_{p_i}\mathcal{K}_{\mathfrak{c}_2}^{A}f
=
\mathcal T_{A,1}+\mathcal T_{A,2},
\]
where $\mathcal T_{A,1}$ denotes the sum of all terms in which $\partial_{p_i}$ falls on the kernel, the cutoff, or the Maxwellian coefficients, while $\mathcal T_{A,2}$ denotes the sum of the terms in which $\partial_{p_i}$ falls on $f(q'')$ or $f(p'')$.

For $\mathcal T_{A,1}$, arguing exactly as in the proof of \eqref{K.2}, and using Lemma \ref{lA.1}, we obtain
\begin{align}\label{K.6}
|\mathcal T_{A,1}|
\lesssim &\int_{\mathbb R^3\times\mathbb S^2}
\mathbf 1_{|p|\leq 3q^0}
\frac{\mathfrak s\mathfrak B}{p^0q^0}\,J_\mathbf{M}^{\frac{\delta}{2}}(q)
\Big\{ J_\mathbf{M}^{\frac{\delta}{2}}(p'')\,|f(q'')|+ J_\mathbf{M}^{\frac{\delta}{2}}(q'')\,|f(p'')|\Big\}\,d\omega\,d q .
\end{align}
For $\mathcal T_{A,2}$, the chain rule yields factors of $\partial_{p_i}p''$ and $\partial_{p_i}q''$, which are bounded by $\langle q\rangle^n$ by Lemma \ref{lA.1}. And the polynomial growth $\langle q\rangle^n$ can be absorbed by $J_\mathbf{M}^{\frac{1}{2}}(q)$, then we obtain
\begin{align}\label{K.7}
|\mathcal T_{A,2}|
\lesssim &\int_{\mathbb R^3\times\mathbb S^2}
\mathbf 1_{|p|\leq 3q^0}
\frac{\mathfrak s\mathfrak B}{p^0q^0}\,J_\mathbf{M}^{\delta}(q)
\Big\{ J_\mathbf{M}^{\delta}(p'')\,|\partial_{p_i}f(q'')|+ J_\mathbf{M}^{\delta}(q'')\,|\partial_{p_i}f(p'')|\Big\}\,d\omega\,d q .
\end{align}

We next consider the $A^c$-region. Similarly,
\[
\partial_{p_i}\mathcal{K}_{\mathfrak{c}_2}^{A^c}f
=
\mathcal T_{A^c,1}+\mathcal T_{A^c,2},
\]
where $\mathcal T_{A^c,1}$ contains the terms in which $\partial_{p_i}$ falls on $v_\phi\chi_{A^c}$ or on the Maxwellian coefficients, and $\mathcal T_{A^c,2}$ contains the terms in which $\partial_{p_i}$ falls on $f(q')$ or $f(p')$. By Lemma \ref{lA.3}, we immediately obtain
\begin{align}\label{K.8}
|\mathcal T_{A^c,1}|
\lesssim &\int_{\mathbb R^3\times\mathbb S^2}
\mathbf 1_{|p|\geq \frac32 q^0}
v_\phi(p,q)\,J_\mathbf{M}^{\frac{\delta}{2}}(q)
\Big\{ J_\mathbf{M}^{\frac{\delta}{2}}(p')\,|f(q')|+ J_\mathbf{M}^{\frac{\delta}{2}}(q')\,|f(p')|\Big\}\,d\omega\,d q 
\end{align}
and
\begin{align}\label{K.9}
|\mathcal T_{A^c,2}|
\lesssim \int_{\mathbb R^3\times\mathbb S^2}
\mathbf 1_{|p|\geq \frac32 q^0}
v_\phi(p,q)\,J_\mathbf{M}^{\delta}(q)
\Big\{ J_\mathbf{M}^{\delta}(p')\,|\partial_{p_i}f(q')|+ J_\mathbf{M}^{\delta}(q')\,|\partial_{p_i}f(p')|\Big\}\,d\omega\,d q .
\end{align}

Combining \eqref{K.6}-\eqref{K.9} by \eqref{07}, we deduce
\begin{align*}
\left|\partial_{p_i}\mathcal{K}_{\mathfrak{c}_2}^{A^c}f\right|
\lesssim &\int_{\mathbb R^3\times\mathbb S^2}
v_\phi(p,q)\,J_\mathbf{M}^{\frac{\delta}{2}}(q)
\Big\{ J_\mathbf{M}^{\frac{\delta}{2}}(p')\,|f(q')|+ J_\mathbf{M}^{\frac{\delta}{2}}(q')\,|f(p')|\Big\}\,d\omega\,d q \\
&+\int_{\mathbb R^3\times\mathbb S^2}
v_\phi(p,q)\,J_\mathbf{M}^{\delta}(q)
\Big\{ J_\mathbf{M}^{\delta}(p')\,|\partial_{p_i}f(q')|+ J_\mathbf{M}^{\delta}(q')\,|\partial_{p_i}f(p')|\Big\}\,d\omega\,d q .
\end{align*}
By the same transformation as in \eqref{K.2}, we thus obtain
\[
\left|\partial_{p_i}\mathcal{K}_{\mathfrak{c}_2}^{A^c}f\right|
\lesssim
\int_{\mathbb R^3}\hat k_2(p,q)|f(q)|\,dq
+
\int_{\mathbb R^3}k_2(p,q)\,|\partial_{p_i}f(q)|\,dq,
\]
which proves \eqref{e2.40-2}. Therefore we complete the proof of Lemma \ref{l2.13}.
\end{proof}

    \section{Derivation and Symmetrization of the Equations for the Macroscopic Quantities}
    \subsection{The positive definiteness of the matrix $A_0$}\label{Appendix0}
    Since the system is considered in the physical state space where $n_\mathfrak{e} > 0$ and $h(n_{\mathfrak{e}}) > 0$, we have $n_\mathfrak{e}h(n_{\mathfrak{e}}) > 0$ as in \eqref{e.1}.
	From \cite[Proposition 3.4]{Speck-CMP-2011}, the Gibbs relation \eqref{e.18} is equivalent to
	\begin{align*}
	\left.\frac{\partial e_{\mathfrak{e}}}{\partial n_{\mathfrak{e}}}\right|_{S_\mathfrak{e}}=\frac{e_{\mathfrak{e}}+P_{\mathfrak{e}}}{n_{\mathfrak{e}}},\left.\quad \frac{\partial e_{\mathfrak{e}}}{\partial S_\mathfrak{e}}\right|_{n_{\mathfrak{e}}}=n_{\mathfrak{e}} T_\mathfrak{e},
	\end{align*}
	then one has
	\begin{align*}
	\left.n_{\mathfrak{e}} \frac{\partial P_{\mathfrak{e}}}{\partial n_{\mathfrak{e}}}\right|_{S_\mathfrak{e}}=\left.\left.n_{\mathfrak{e}} \frac{\partial P_{\mathfrak{e}}}{\partial e_{\mathfrak{e}}}\right|_{S_\mathfrak{e}} \cdot \frac{\partial e_{\mathfrak{e}}}{\partial n_{\mathfrak{e}}}\right|_{S_\mathfrak{e}}=\frac{a^2}{\mathfrak{c}^2}\left(e_{\mathfrak{e}}+P_{\mathfrak{e}}\right),
	\end{align*}
	where $a^2:=\left.\mathfrak{c}^2 \frac{\partial P_{\mathfrak{e}}}{\partial e_{\mathfrak{e}}}\right|_{S_\mathfrak{e}}$ is the square of sound speed. Since $S_\mathfrak{e}$ is a constant and $h(n_{\mathfrak{e}})=\mathfrak{c}^2\frac{K_3(\gamma)}{K_2(\gamma)}$, we then have $n_{\mathfrak{e}}
    h^{\prime}=P_{\mathfrak{e}}^{\prime}\left(n_{\mathfrak{e}}\right)=a^2\frac{K_3(\gamma)}{K_2(\gamma)} $. Using the fact that $a \in\left(0, \frac{\mathfrak{c}}{\sqrt{3}}\right)$ (see \cite{Speck-CMP-2011}), one can show that
    \begin{align}\label{e.17}
    h\left(n_{\mathfrak{e}}\right)>n_{\mathfrak{e}}
    h^{\prime}.
    \end{align}
    
	For $\mathfrak{c}\gg 1$, it holds from \eqref{2.5} that
	\begin{align*}
	\left(\frac{K_3(\gamma)}{K_2(\gamma)}\right)^2-\frac{5}{\gamma} \frac{K_3(\gamma)}{K_2(\gamma)}-1<\left(\frac{K_3(\gamma)}{K_2(\gamma)}\right)^2-\frac{5}{\gamma} \frac{K_3(\gamma)}{K_2(\gamma)}+\frac{1}{\gamma^2}-1<0,
	\end{align*}
	which yields that
	\begin{gather*}
	     \operatorname{det}\left(A_0\right)_{1 \times1}=h^{\prime}=\frac{T_\mathfrak{e}}{n_{\mathfrak{e}}}\frac{\left(\frac{K_3(\gamma)}{K_2(\gamma)}\right)^2-\frac{5}{\gamma} \frac{K_3(\gamma)}{K_2(\gamma)}-1}{\left(\frac{K_3(\gamma)}{K_2(\gamma)}\right)^2-\frac{5}{\gamma} \frac{K_3(\gamma)}{K_2(\gamma)}+\frac{1}{\gamma^2}-1}>0, \\
	    \operatorname{det}\left(A_0\right)_{2 \times 2}=\frac{n_{\mathfrak{e}}h^{\prime}h}{\mathfrak{c}^2}-\frac{n_{\mathfrak{e}}h^{\prime}h}{\mathfrak{c}^2}\frac{u_{\mathfrak{e},1}^2}{(u_{\mathfrak{e}}^0)^2}-\frac{n_{\mathfrak{e}}^2(h^{\prime})^2}{(u_{\mathfrak{e}}^0)^2}\frac{u_{\mathfrak{e},1}^2}{(u_{\mathfrak{e}}^0)^2}>\frac{n_{\mathfrak{e}}h^{\prime}h}{\mathfrak{c}^2}\left(1-\frac{2u_{\mathfrak{e},1}^2}{(u_{\mathfrak{e}}^0)^2}\right)>0,\\
	    \operatorname{det}\left(A_0\right)_{3 \times 3}=\frac{n_{\mathfrak{e}}h}{\mathfrak{c}^2}\left[\frac{n_{\mathfrak{e}}h^{\prime}h}{\mathfrak{c}^2}-\left(\frac{n_{\mathfrak{e}}h^{\prime}h}{\mathfrak{c}^2}+\frac{n_{\mathfrak{e}}^2(h^{\prime})^2}{(u_{\mathfrak{e}}^0)^2}\right)\left(\frac{u_{\mathfrak{e},1}^2+u_{\mathfrak{e},2}^2}{(u_{\mathfrak{e}}^0)^2}\right)\right]>\frac{n_{\mathfrak{e}}^2h^{\prime}h^2}{\mathfrak{c}^4}\left(1-\frac{2u_{\mathfrak{e},1}^2+2u_{\mathfrak{e},2}^2}{(u_{\mathfrak{e}}^0)^2}\right)>0,\\
	    \det (A_0)_{4\times 4} = \left(\frac{n_\mathfrak{e} h}{\mathfrak{c}^2}\right)^2 \left[ \frac{n_\mathfrak{e} h' h}{\mathfrak{c}^2} - \left( \frac{n_\mathfrak{e} h' h}{\mathfrak{c}^2} + \frac{n_\mathfrak{e}^2 (h')^2}{(u_\mathfrak{e}^0)^2} \right) \frac{u_{\mathfrak{e},1}^2 + u_{\mathfrak{e},2}^2 + u_{\mathfrak{e},3}^2}{(u_\mathfrak{e}^0)^2} \right] \\
        > \frac{n_\mathfrak{e}^3 h^3 h'}{\mathfrak{c}^6} \left( 1 - \frac{2(u_{\mathfrak{e},1}^2 + u_{\mathfrak{e},2}^2 + u_{\mathfrak{e},3}^2)}{(u_\mathfrak{e}^0)^2} \right) > 0,
	\end{gather*}
	where we have used \eqref{e.17} and $(u_\mathfrak{e}^0)^2>2u_{\mathfrak{e}}^2$ for sufficiently large $\mathfrak{c}$. Therefore, $A_0$ is a positive definite matrix.
	\subsection{Derivation of the equations for
	$\left(a_{n+1}, b_{n+1}, c_{n+1}\right)$}\label{AppendixC}
     Using \eqref{e1.14} and Lemma \ref{l5.2}, by tedious calculations, one has
     {\small
	\begin{align*}
		\int_{\mathbb{R}^3} F_{n+1}^{\mathfrak{c}} d p
		&=\int_{\mathbb{R}^3}\Big\{a_{n+1}+b_{n+1} \cdot p+c_{n+1} \frac{p^0}{\mathfrak{c}}\Big\} \mathbf{M}_{\mathfrak{c}} d p \\
		& =\frac{n_{\mathfrak{e}} u_{\mathfrak{e}}^0}{\mathfrak{c}} a_{n+1}+\frac{e_{\mathfrak{e}}+P_{\mathfrak{e}}}{\mathfrak{c}^3} u_{\mathfrak{e}}^0\left(u_{\mathfrak{e}} \cdot b_{n+1}\right)+\frac{e_{\mathfrak{e}}\big(u_{\mathfrak{e}}^0\big)^2+P_{\mathfrak{e}}|u_{\mathfrak{e}}|^2}{\mathfrak{c}^4} c_{n+1},\\
		\int_{\mathbb{R}^3} \frac{p_j p}{p^0} F_{n+1}^{\mathfrak{c}} d p
		&=\int_{\mathbb{R}^3} \frac{p_j p}{p^0}\Big\{a_{n+1}+b_{n+1} \cdot p+c_{n+1} \frac{p^0}{\mathfrak{c}}\Big\} \mathbf{M}_{\mathfrak{c}} d p+\int_{\mathbb{R}^3} \frac{p_j p}{p^0} \sqrt{\mathbf{M}_{\mathfrak{c}}}\left\{\mathbf{I}-\mathbf{P}_{\mathfrak{c}}\right\}\Big(\frac{F_{n+1}^{\mathfrak{c}}}{\sqrt{\mathbf{M}_{\mathfrak{c}}}}\Big) d p \\
		& =\frac{n_{\mathfrak{e}}}{\mathfrak{c} \gamma K_2(\gamma)}\left(6 K_3(\gamma)+\gamma K_2(\gamma)\right) u_{\mathfrak{e},j} u_{\mathfrak{e}}\Big\{\left(u_{\mathfrak{e}} \cdot b_{n+1}\right)+\frac{u_{\mathfrak{e}}^0}{\mathfrak{c}} c_{n+1}\Big\} \\
		& \quad+\frac{e_{\mathfrak{e}}+P_{\mathfrak{e}}}{\mathfrak{c}^3} u_{\mathfrak{e},j} u_{\mathfrak{e}} a_{n+1}+\mathbf{e}_j a_{n+1} \frac{P_{\mathfrak{e}}}{\mathfrak{c}}+\frac{\mathfrak{c} n_{\mathfrak{e}} K_3(\gamma)}{\gamma K_2(\gamma)}\left(u_{\mathfrak{e}} b_{n+1, j}+u_{\mathfrak{e},j} b_{n+1}\right) \\
		&\quad +\mathbf{e}_j \frac{\mathfrak{c} n_{\mathfrak{e}} K_3(\gamma)}{\gamma K_2(\gamma)}\Big\{\left(u_{\mathfrak{e}} \cdot b_{n+1}\right)+\frac{u_{\mathfrak{e}}^0}{\mathfrak{c}} c_{n+1}\Big\}+\int_{\mathbb{R}^3} \frac{p_j p}{p^0} \sqrt{\mathbf{M}_{\mathfrak{c}}}\left\{\mathbf{I}-\mathbf{P}_{\mathfrak{c}}\right\}\Big(\frac{F_{n+1}^{\mathfrak{c}}}{\sqrt{\mathbf{M}_{\mathfrak{c}}}}\Big) d p,\\
		\int_{\mathbb{R}^3} \hat{p}_j F_{n+1}^{\mathfrak{c}} d p
		&=\int_{\mathbb{R}^3} \hat{p}_j\Big\{a_{n+1}+b_{n+1} \cdot p+c_{n+1} \frac{p^0}{\mathfrak{c}}\Big\} \mathbf{M}_{\mathfrak{c}} d p+\int_{\mathbb{R}^3} \hat{p}_j \sqrt{\mathbf{M}_{\mathfrak{c}}}\left\{\mathbf{I}-\mathbf{P}_{\mathfrak{c}}\right\}\Big(\frac{F_{n+1}^{\mathfrak{c}}}{\sqrt{\mathbf{M}_{\mathfrak{c}}}}\Big) d p \\
		& =n_{\mathfrak{e}} u_{\mathfrak{e},j} a_{n+1}+\frac{e_{\mathfrak{e}}+P_{\mathfrak{e}}}{\mathfrak{c}^2} u_{\mathfrak{e},j}\left(u_{\mathfrak{e}} \cdot b_{n+1}\right)+P_{\mathfrak{e}} b_{n+1, j}+\frac{e_{\mathfrak{e}}+P_{\mathfrak{e}}}{\mathfrak{c}^3} u_{\mathfrak{e}}^0 u_{\mathfrak{e},j} c_{n+1} \\
		& \quad+\int_{\mathbb{R}^3} \hat{p}_j \sqrt{\mathbf{M}_{\mathfrak{c}}}\left\{\mathbf{I}-\mathbf{P}_{\mathfrak{c}}\right\}\Big(\frac{F_{n+1}^{\mathfrak{c}}}{\sqrt{\mathbf{M}_{\mathfrak{c}}}}\Big) d p, \\
		\int_{\mathbb{R}^3} p_j F_{n+1}^{\mathfrak{c}} d p
		&=\int_{\mathbb{R}^3} p_j\Big\{a_{n+1}+b_{n+1} \cdot p+c_{n+1} \frac{p^0}{\mathfrak{c}}\Big\} \mathbf{M}_{\mathfrak{c}} d p \\
		& =\frac{n_{\mathfrak{e}}}{\mathfrak{c} \gamma K_2(\gamma)}\Big\{\left(6 K_3(\gamma)+\gamma K_2(\gamma)\right) u_{\mathfrak{e}}^0 u_{\mathfrak{e},j}\left(u_{\mathfrak{e}} \cdot b_{n+1}\right)+\mathfrak{c}^2 K_3(\gamma) u_{\mathfrak{e}}^0 b_{n+1, j}\Big\} \\
		& \quad+\frac{n_{\mathfrak{e}}}{\mathfrak{c}^2 \gamma K_2(\gamma)}\Big\{\left(5 K_3(\gamma)+\gamma K_2(\gamma)\right)\big(u_{\mathfrak{e}}^0\big)^2+K_3(\gamma)|u_{\mathfrak{e}}|^2\Big\} u_{\mathfrak{e},j} c_{n+1} \\
		&\quad +\frac{e_{\mathfrak{e}}+P_{\mathfrak{e}}}{\mathfrak{c}^3} u_{\mathfrak{e}}^0 u_{\mathfrak{e},j} a_{n+1}, \\
		\int_{\mathbb{R}^3} p^0 F_{n+1}^{\mathfrak{c}} d p
		&=\int_{\mathbb{R}^3} p^0\Big\{a_{n+1}+b_{n+1} \cdot p+c_{n+1} \frac{p^0}{\mathfrak{c}}\Big\} \mathbf{M}_{\mathfrak{c}} d p \\
		& =\frac{n_{\mathfrak{e}}}{\mathfrak{c} \gamma K_2(\gamma)}\Big\{\left(5 K_3(\gamma)+\gamma K_2(\gamma)\right)\big(u_{\mathfrak{e}}^0\big)^2+K_3(\gamma)|u_{\mathfrak{e}}|^2\Big\}\left(u_{\mathfrak{e}} \cdot b_{n+1}\right) \\
		& \quad+\frac{n_{\mathfrak{e}}}{\mathfrak{c}^2 \gamma K_2(\gamma)}\Big\{\left(3 K_3(\gamma)+\gamma K_2(\gamma)\right)\big(u_{\mathfrak{e}}^0\big)^2+3 K_3(\gamma)|u_{\mathfrak{e}}|^2\Big\} u_{\mathfrak{e}}^0 c_{n+1} \\
		&\quad+\frac{e_{\mathfrak{e}}\big(u_{\mathfrak{e}}^0\big)^2+P_{\mathfrak{e}}|u_{\mathfrak{e}}|^2}{\mathfrak{c}^3} a_{n+1},
	\end{align*}}
	where $\mathbf{e}_j(j=1,2,3)$ are the unit base vectors in $\mathbb{R}^3$.
	
	Next, we shall derive the equations for $\left(a_{n+1}, b_{n+1}, c_{n+1}\right)$.
	Integrating $\eqref{N.2}_1$ with respect to $p$, we have
	{\small
	\begin{align}\label{N.4}
		& \partial_t\Big(\frac{n_{\mathfrak{e}} u_{\mathfrak{e}}^0}{\mathfrak{c}} a_{n+1}+\frac{e_{\mathfrak{e}}+P_{\mathfrak{e}}}{\mathfrak{c}^3} u_{\mathfrak{e}}^0\left(u_{\mathfrak{e}} \cdot b_{n+1}\right)+\frac{e_{\mathfrak{e}}\big(u_{\mathfrak{e}}^0\big)^2+P_{\mathfrak{e}}|u_{\mathfrak{e}}|^2}{\mathfrak{c}^4} c_{n+1}\Big) \nonumber\\
		& \quad+\operatorname{div}\Big(n_{\mathfrak{e}} u_{\mathfrak{e}} a_{n+1}+\frac{e_{\mathfrak{e}}+P_{\mathfrak{e}}}{\mathfrak{c}^2} u_{\mathfrak{e}}\left(u_{\mathfrak{e}} \cdot b_{n+1}\right)+P_{\mathfrak{e}} b_{n+1}+\frac{e_{\mathfrak{e}}+P_{\mathfrak{e}}}{\mathfrak{c}^3} u_{\mathfrak{e}}^0 u_{\mathfrak{e}} c_{n+1}\Big)\nonumber\\
		&\quad+\operatorname{div} \int_{\mathbb{R}^3} \hat{p} \sqrt{\mathbf{M}_{\mathfrak{c}}}\left\{\mathbf{I}-\mathbf{P}_{\mathfrak{c}}\right\}\Big(\frac{F_{n+1}^{\mathfrak{c}}}{\sqrt{\mathbf{M}_{\mathfrak{c}}}}\Big) d p=0 .
	\end{align}}
	Note that
	{\small
	\begin{align}
		-\int_{\mathbb{R}^3} p_j\big(E_{n+1}^{\mathfrak{c}}+\frac{p}{p^0} \times B_{n+1}^{\mathfrak{c}}\big) \cdot \nabla_p F_0^{\mathfrak{c}} d p\nonumber
		&=\int_{\mathbb{R}^3} E_{n+1, j}^{\mathfrak{c}} F_0^{\mathfrak{c}} d p+\int_{\mathbb{R}^3}\Big(\frac{p}{p^0} \times B_{n+1}^{\mathfrak{c}}\Big)_j F_0^{\mathfrak{c}} d p \nonumber\\
		& =\frac{n_{\mathfrak{e}} u_{\mathfrak{e}}^0}{\mathfrak{c}} E_{n+1, j}^{\mathfrak{c}}+\left(\frac{n_{\mathfrak{e}} u_{\mathfrak{e}}}{\mathfrak{c}} \times B_{n+1}^{\mathfrak{c}}\right)_j \text {, } \nonumber
	\end{align}}
	and
	{\small
	\begin{align}
		&-\int_{\mathbb{R}^3} p_j\big(E_0^{\mathfrak{c}}+\frac{p}{p^0} \times B_0^{\mathfrak{c}}\big) \cdot \nabla_p F_{n+1}^{\mathfrak{c}} d p\nonumber\\
		=&\int_{\mathbb{R}^3} E_{0, j}^{\mathfrak{c}} F_{n+1}^{\mathfrak{c}} d p+\int_{\mathbb{R}^3}\Big(\frac{p}{p^0} \times B_0^{\mathfrak{c}}\Big)_j F_{n+1}^{\mathfrak{c}} d p \nonumber\\
		=&E_{0, j}^{\mathfrak{c}}\Big(\frac{n_{\mathfrak{e}} u_{\mathfrak{e}}^0}{\mathfrak{c}} a_{n+1}+\frac{e_{\mathfrak{e}}+P_{\mathfrak{e}}}{\mathfrak{c}^3} u_{\mathfrak{e}}^0\left(u_{\mathfrak{e}} \cdot b_{n+1}\right)+\frac{e_{\mathfrak{e}} \big(u_{\mathfrak{e}}^0\big)^2+P_{\mathfrak{e}}|u_{\mathfrak{e}}|^2}{\mathfrak{c}^4} c_{n+1}\Big)\nonumber \\
		&+ \left\{\Big(\frac{n_{\mathfrak{e}} u_{\mathfrak{e}} a_{n+1}}{\mathfrak{c}}+\frac{e_{\mathfrak{e}}+P_{\mathfrak{e}}}{\mathfrak{c}^3} u_{\mathfrak{e}}\left(u_{\mathfrak{e}} \cdot b_{n+1}\right)+\frac{P_{\mathfrak{e}} b_{n+1}}{\mathfrak{c}}+\frac{e_{\mathfrak{e}}+P_{\mathfrak{e}}}{\mathfrak{c}^4} u_{\mathfrak{e}}^0 u_{\mathfrak{e}} c_{n+1}\Big) \times B_0^{\mathfrak{c}}\right\}_j \nonumber\\
		&+ \int_{\mathbb{R}^3}\Big(\frac{p}{p^0} \times B_0^{\mathfrak{c}}\Big)_j \sqrt{\mathbf{M}_{\mathfrak{c}}}\left\{\mathbf{I}-\mathbf{P}_{\mathfrak{c}}\right\}\Big(\frac{F_{n+1}^{\mathfrak{c}}}{\sqrt{\mathbf{M}_{\mathfrak{c}}}}\Big) d p , \nonumber
	\end{align}}
	multiplying $\eqref{N.2}_1$ by $p_j$ and integrating over $\mathbb{R}^3$, one gets
	{\small
	\begin{align}\label{N.5}
		\partial_t&\left\{\frac{e_{\mathfrak{e}}+P_{\mathfrak{e}}}{\mathfrak{c}^3} u_{\mathfrak{e}}^0 u_{\mathfrak{e},j} a_{n+1}+\frac{n_{\mathfrak{e}}}{\mathfrak{c} \gamma K_2(\gamma)}\Big\{\left(6 K_3(\gamma)+\gamma K_2(\gamma)\right) u_{\mathfrak{e}}^0 u_{\mathfrak{e},j}\left(u_{\mathfrak{e}} \cdot b_{n+1}\right)+\mathfrak{c}^2 K_3(\gamma) u_{\mathfrak{e}}^0 b_{n+1, j}\Big\}\right. \nonumber\\
		& \left.+\frac{n_{\mathfrak{e}}}{\mathfrak{c}^2 \gamma K_2(\gamma)}\Big\{\left(5 K_3(\gamma)+\gamma K_2(\gamma)\right)\big(u_{\mathfrak{e}}^0\big)^2+K_3(\gamma)|u_{\mathfrak{e}}|^2\Big\} u_{\mathfrak{e},j} c_{n+1}\right\} \nonumber\\
		& +\operatorname{div}\Big\{\frac{e_{\mathfrak{e}}+P_{\mathfrak{e}}}{\mathfrak{c}^2} u_{\mathfrak{e},j} u_{\mathfrak{e}} a_{n+1}+\frac{n_{\mathfrak{e}}}{\gamma K_2(\gamma)}\left(6 K_3(\gamma)+\gamma K_2(\gamma)\right) u_{\mathfrak{e},j} u_{\mathfrak{e}}\Big(\left(u_{\mathfrak{e}} \cdot b_{n+1}\right)+\frac{u_{\mathfrak{e}}^0}{\mathfrak{c}} c_{n+1}\Big)\Big\} \nonumber\\
		& +\partial_{x_j}\left(P_{\mathfrak{e}} a_{n+1}\right)+\operatorname{div}\Big\{\frac{\mathfrak{c}^2 n_{\mathfrak{e}} K_3(\gamma)}{\gamma K_2(\gamma)}\left(u_{\mathfrak{e}} b_{n+1, j}+u_{\mathfrak{e},j} b_{n+1}\right)\Big\} \nonumber\\
		& +\partial_{x_j}\Big\{\frac{\mathfrak{c}^2 n_{\mathfrak{e}} K_3(\gamma)}{\gamma K_2(\gamma)}\Big(\left(u_{\mathfrak{e}} \cdot b_{n+1}\right)+\frac{u_{\mathfrak{e}}^0}{\mathfrak{c}} c_{n+1}\Big)\Big\}+\operatorname{div}\Big\{ \int_{\mathbb{R}^3} p_j \hat{p} \sqrt{\mathbf{M}_{\mathfrak{c}}}\left\{\mathbf{I}-\mathbf{P}_{\mathfrak{c}}\right\}\Big(\frac{F_{n+1}^{\mathfrak{c}}}{\sqrt{\mathbf{M}_{\mathfrak{c}}}}\Big) d p\Big\}\nonumber\\
		&+E_{0, j}^{\mathfrak{c}}\Big(\frac{n_{\mathfrak{e}} u_{\mathfrak{e}}^0}{\mathfrak{c}} a_{n+1}+\frac{e_{\mathfrak{e}}+P_{\mathfrak{e}}}{\mathfrak{c}^3} u_{\mathfrak{e}}^0\left(u_{\mathfrak{e}} \cdot b_{n+1}\right)+\frac{e_{\mathfrak{e}}\big(u_{\mathfrak{e}}^0\big)^2+P_{\mathfrak{e}}|u_{\mathfrak{e}}|^2}{\mathfrak{c}^4} c_{n+1}\Big)\nonumber\\
		&+\left\{\Big(\frac{n_{\mathfrak{e}} u_{\mathfrak{e}} a_{n+1}}{\mathfrak{c}}+\frac{e_{\mathfrak{e}}+P_{\mathfrak{e}}}{\mathfrak{c}^3} u_{\mathfrak{e}}\left(u_{\mathfrak{e}} \cdot b_{n+1}\right)+\frac{P_{\mathfrak{e}} b_{n+1}}{\mathfrak{c}}+\frac{e_{\mathfrak{e}}+P_{\mathfrak{e}}}{\mathfrak{c}^4}u_{\mathfrak{e}}^0 u_{\mathfrak{e}}\left(c_{n+1}\right)\Big)\times B_0^{\mathfrak{c}}\right\}_j\nonumber\\
		&+\int_{\mathbb{R}^3}\Big(\frac{p}{p^0} \times B_0^{\mathfrak{c}}\Big)_j \sqrt{\mathbf{M}_{\mathfrak{c}}}\left\{\mathbf{I}-\mathbf{P}_{\mathfrak{c}}\right\}\Big(\frac{F_{n+1}^{\mathfrak{c}}}{\sqrt{\mathbf{M}_{\mathfrak{c}}}}\Big) d p+\frac{n_{\mathfrak{e}} u_{\mathfrak{e}}^0}{\mathfrak{c}} E_{n+1, j}^{\mathfrak{c}}+\left(\frac{n_{\mathfrak{e}} u_{\mathfrak{e}}}{\mathfrak{c}} \times B_{n+1}^{\mathfrak{c}}\right)_j \nonumber\\
		&+\sum_{\substack{k+l=n+1 \\
				k, l \geq 1}} E^{\mathfrak{c}}_{k, j}\Big(\frac{n_{\mathfrak{e}} u_{\mathfrak{e}}^0}{\mathfrak{c}} a_l+\frac{e_{\mathfrak{e}}+P_{\mathfrak{e}}}{\mathfrak{c}^3} u_{\mathfrak{e}}^0\left(u_{\mathfrak{e}} \cdot b_l\right)+\frac{e_{\mathfrak{e}}\big(u_{\mathfrak{e}}^0\big)^2+P_{\mathfrak{e}}|u_{\mathfrak{e}}|^2}{\mathfrak{c}^4} c_l\Big) \nonumber\\
		&+\sum_{\substack{k+l=n+1 \\
				k, l \geq 1}}\left\{\Big(\frac{n_{\mathfrak{e}} u_{\mathfrak{e}} a_l}{\mathfrak{c}}+\frac{e_{\mathfrak{e}}+P_{\mathfrak{e}}}{\mathfrak{c}^3} u_{\mathfrak{e}}\left(u_{\mathfrak{e}} \cdot b_l\right)+\frac{P_{\mathfrak{e}} b_l}{\mathfrak{c}}+\frac{e_{\mathfrak{e}}+P_{\mathfrak{e}}}{\mathfrak{c}^4} u_{\mathfrak{e}}^0 u_{\mathfrak{e}} c_l\Big) \times B_k^{\mathfrak{c}}\right\}_j\nonumber\\
		&+\sum_{\substack{k+l=n+1 \\
				k, l \geq 1}} \int_{\mathbb{R}^3}\Big(\frac{p}{p^0} \times B_k^{\mathfrak{c}}\Big)_j \sqrt{\mathbf{M}_{\mathfrak{c}}}\left\{\mathbf{I}-\mathbf{P}_{\mathfrak{c}}\right\}\Big(\frac{F_l^{\mathfrak{c}}}{\sqrt{\mathbf{M}_{\mathfrak{c}}}}\Big) d p=0
	\end{align}	}
	for $j=1,2,3$ with $b_{n+1}=\left(b_{n+1,1}, b_{n+1,2}, b_{n+1,3}\right)^t$.
	
	Note that
	{\small
	\begin{align}
		-\int_{\mathbb{R}^3} p^0\big(E_{n+1}^{\mathfrak{c}}+\frac{p}{p^0} \times B_{n+1}^{\mathfrak{c}}\big) \cdot \nabla_p F_0^{\mathfrak{c}} d p
		&=\int_{\mathbb{R}^3} E_{n+1}^{\mathfrak{c}} \cdot \frac{p}{p^0} F_0^{\mathfrak{c}} d p
		=E_{n+1}^{\mathfrak{c}} \cdot \frac{n_{\mathfrak{e}} u_{\mathfrak{e}}}{\mathfrak{c}}, \nonumber
	\end{align}}
	and
	{\small
	\begin{align}
		&-\int_{\mathbb{R}^3} p^0\big(E_0^{\mathfrak{c}}+\frac{p}{p^0} \times B_0^{\mathfrak{c}}\big) \cdot \nabla_p F_{n+1}^{\mathfrak{c}} d p \nonumber\\
		=&\int_{\mathbb{R}^3} E_0^{\mathfrak{c}} \cdot \frac{p}{p^0} F_{n+1} d p \nonumber\\
		=&\frac{n_{\mathfrak{e}} a_{n+1}}{\mathfrak{c}} u_{\mathfrak{e}} \cdot E_0^{\mathfrak{c}}+\frac{e_{\mathfrak{e}}+P_{\mathfrak{e}}}{\mathfrak{c}^3}\left(u_{\mathfrak{e}} \cdot b_{n+1}\right)\left(E_0^{\mathfrak{c}} \cdot u_{\mathfrak{e}}\right)+\frac{P_{\mathfrak{e}} b_{n+1}}{\mathfrak{c}} \cdot E_0^{\mathfrak{c}}+ \nonumber\\
		& \frac{e_{\mathfrak{e}}+P_{\mathfrak{e}}}{\mathfrak{c}^4} u_{\mathfrak{e}}^0 c_{n+1} u_{\mathfrak{e}} \cdot E_0^{\mathfrak{c}}+\int_{\mathbb{R}^3} \frac{p}{p^0} \cdot E_0^{\mathfrak{c}} \sqrt{\mathbf{M}_{\mathfrak{c}}}\left\{\mathbf{I}-\mathbf{P}_{\mathfrak{c}}\right\}\Big(\frac{F_{n+1}^{\mathfrak{c}}}{\sqrt{\mathbf{M}_{\mathfrak{c}}}}\Big) d p ,\nonumber
	\end{align}}
	multiplying $\eqref{N.2}_1$ by $\frac{p^0}{\mathfrak{c}}$ and integrating over $\mathbb{R}^3$, one obtains that
	{\small
	\begin{align}\label{N.6}
		\partial_t&\Big\{  \frac{e_{\mathfrak{e}}\big(u_{\mathfrak{e}}^0\big)^2+P_{\mathfrak{e}}|u_{\mathfrak{e}}|^2}{\mathfrak{c}^4} a_{n+1}+\frac{n_{\mathfrak{e}}}{\mathfrak{c}^2 \gamma K_2(\gamma)}\Big\{\left(5 K_3(\gamma)+\gamma K_2(\gamma)\right)\big(u_{\mathfrak{e}}^0\big)^2+K_3(\gamma)|u_{\mathfrak{e}}|^2\Big\}\left(u_{\mathfrak{e}} \cdot b_{n+1}\right) \nonumber\\
		& +\frac{n_{\mathfrak{e}}}{\mathfrak{c}^3 \gamma K_2(\gamma)}\Big\{\left(3 K_3(\gamma)+\gamma K_2(\gamma)\right)\big(u_{\mathfrak{e}}^0\big)^2+3 K_3(\gamma)|u_{\mathfrak{e}}|^2\Big\} u_{\mathfrak{e}}^0 c_{n+1}\Big\} \nonumber\\
		& +\operatorname{div}\left\{\frac{e_{\mathfrak{e}}+P_{\mathfrak{e}}}{\mathfrak{c}^3} u_{\mathfrak{e}}^0 u_{\mathfrak{e}} a_{n+1}+\frac{n_{\mathfrak{e}}}{\mathfrak{c} \gamma K_2(\gamma)}\Big\{\left(6 K_3(\gamma)+\gamma K_2(\gamma)\right) u_{\mathfrak{e}}^0 u_{\mathfrak{e}}\left(u_{\mathfrak{e}} \cdot b_{n+1}\right)+\mathfrak{c}^2 K_3(\gamma) u_{\mathfrak{e}}^0 b_{n+1}\Big\}\right. \nonumber\\
		& \left.+\frac{n_{\mathfrak{e}}}{\mathfrak{c}^2 \gamma K_2(\gamma)}\Big\{\left(5 K_3(\gamma)+\gamma K_2(\gamma)\right)\big(u_{\mathfrak{e}}^0\big)^2+K_3(\gamma)|u_{\mathfrak{e}}|^2\Big\} u_{\mathfrak{e}} c_{n+1}\right\}+\frac{n_{\mathfrak{e}} u_{\mathfrak{e}}}{\mathfrak{c}^2} \cdot E_{n+1}^{\mathfrak{c}}\nonumber\\
		&+\frac{n_{\mathfrak{e}} a_{n+1}}{\mathfrak{c}^2} u_{\mathfrak{e}} \cdot E_0^{\mathfrak{c}} +\frac{e_{\mathfrak{e}}+P_{\mathfrak{e}}}{\mathfrak{c}^4}\left(u_{\mathfrak{e}} \cdot b_{n+1}\right)\left(E_0^{\mathfrak{c}} \cdot u_{\mathfrak{e}}\right)+\frac{P_{\mathfrak{e}} b_{n+1}}{\mathfrak{c}^2} \cdot E_0^{\mathfrak{c}}+\frac{e_{\mathfrak{e}}+P_{\mathfrak{e}}}{\mathfrak{c}^5} u_{\mathfrak{e}}^0 c_{n+1}\left(u_{\mathfrak{e}} \cdot E_0^{\mathfrak{c}}\right) \nonumber\\
		&+\int_{\mathbb{R}^3} \frac{p \cdot E_0^{\mathfrak{c}}}{\mathfrak{c} p^0} \sqrt{\mathbf{M}_{\mathfrak{c}}}\left\{\mathbf{I}-\mathbf{P}_{\mathfrak{c}}\right\}\Big(\frac{F_{n+1}^{\mathfrak{c}}}{\sqrt{\mathbf{M}_{\mathfrak{c}}}}\Big) d p+\sum_{\substack{k+l=n+1 \\
				k, l \geq 1}}\left\{\frac{n_{\mathfrak{e}} a_l}{\mathfrak{c}^2}\left(u_{\mathfrak{e}} \cdot E_k^{\mathfrak{c}}\right)+\frac{e_{\mathfrak{e}}+P_{\mathfrak{e}}}{\mathfrak{c}^4}\left(u_{\mathfrak{e}} \cdot b_l\right)\left(E_k^{\mathfrak{c}} \cdot u_{\mathfrak{e}}\right)\right. \nonumber\\
		&\left.+\frac{P_{\mathfrak{e}}}{\mathfrak{c}^2} E_k^{\mathfrak{c}} \cdot b_l+\frac{e_{\mathfrak{e}}+P_{\mathfrak{e}}}{\mathfrak{c}^5} u_{\mathfrak{e}}^0 c_l\left(u_{\mathfrak{e}} \cdot E_k^{\mathfrak{c}}\right)+\int_{\mathbb{R}^3} \frac{p \cdot E_k^{\mathfrak{c}}}{\mathfrak{c} p^0} \sqrt{\mathbf{M}_{\mathfrak{c}}}\left\{\mathbf{I}-\mathbf{P}_{\mathfrak{c}}\right\}\Big(\frac{F_l^{\mathfrak{c}}}{\sqrt{\mathbf{M}_{\mathfrak{c}}}}\Big) d p\right\}=0 . 
	\end{align}}
	Moreover, $\eqref{N.2}_2-\eqref{N.2}_5$ can be rewritten as 
	{\small
	\begin{align}\label{N.7}
		\begin{cases}
			\partial_t E_{n+1}^{\mathfrak{c}}-\mathfrak{c} \nabla_x \times B_{n+1}^{\mathfrak{c}}=4 \pi n_{\mathfrak{e}} u_{\mathfrak{e}} a_{n+1}+4 \pi \frac{e_{\mathfrak{e}}+P_{\mathfrak{e}}}{\mathfrak{c}^2} u_{\mathfrak{e}}\left(u_{\mathfrak{e}} \cdot b_{n+1}\right)+4 \pi P_{\mathfrak{e}} b_{n+1}+4 \pi \frac{e_{\mathfrak{e}}+P_{\mathfrak{e}}}{\mathfrak{c}^3} u_{\mathfrak{e}}^0 u_{\mathfrak{e}} c_{n+1}\\
			\qquad\qquad\qquad\qquad\qquad\quad+4 \pi\int_{\mathbb{R}^3} \hat{p}\sqrt{\mathbf{M}_{\mathfrak{c}}}\left\{\mathbf{I}-\mathbf{P}_{\mathfrak{c}}\right\}\Big(\frac{F_{n+1}^{\mathfrak{c}}}{\sqrt{\mathbf{M}_{\mathfrak{c}}}}\Big) d p ,\\
			\partial_t B_{n+1}^{\mathfrak{c}}+\mathfrak{c} \nabla_x \times E_{n+1}^{\mathfrak{c}}=0 ,\\
			\operatorname{div} E_{n+1}^{\mathfrak{c}}=-4 \pi \frac{n_{\mathfrak{e}} u_{\mathfrak{e}}^0}{\mathfrak{c}} a_{n+1}-4 \pi \frac{e_{\mathfrak{e}}+P_{\mathfrak{e}}}{\mathfrak{c}^3} u_{\mathfrak{e}}^0\left(u_{\mathfrak{e}} \cdot b_{n+1}\right)-4 \pi \frac{e_{\mathfrak{e}}\big(u_{\mathfrak{e}}^0\big)^2+P_{\mathfrak{e}} | u_{\mathfrak{e}}|^2}{\mathfrak{c}^4} c_{n+1}, \\
			\operatorname{div} B_{n+1}^{\mathfrak{c}}=0.
		\end{cases}
	\end{align}}
	
	We simplify equations \eqref{N.4},\eqref{N.5} and \eqref{N.6} as follows:
	{\small
		\begin{align}\label{N.8}
			& \frac{n_{\mathfrak{e}} u_{\mathfrak{e}}^0}{\mathfrak{c}}\Big\{\partial_t a_{n+1}+\frac{h}{\mathfrak{c}^2}\left(u_{\mathfrak{e}} \cdot \partial_t b_{n+1}\right)+\frac{h}{\mathfrak{c}^3} u_{\mathfrak{e}}^0 \partial_t c_{n+1}\Big\}-\frac{P_{\mathfrak{e}}}{\mathfrak{c}^2} \partial_t c_{n+1}+\frac{n_{\mathfrak{e}} u_{\mathfrak{e}}^0}{\mathfrak{c}}\left\{\frac{1}{\mathfrak{c}^2} \partial_t(h u_{\mathfrak{e}}) \cdot b_{n+1}\right. \nonumber\\
			& \left.+\frac{1}{\mathfrak{c}^3} \partial_t\left(h u_{\mathfrak{e}}^0\right) c_{n+1}\right\}-\frac{c_{n+1}}{\mathfrak{c}^2} \partial_{t} P_{\mathfrak{e}}+n_{\mathfrak{e}} u_{\mathfrak{e}} \cdot\Big\{\nabla_x a_{n+1}+\frac{h}{\mathfrak{c}^2} \nabla_x b_{n+1} \cdot u_{\mathfrak{e}}+\frac{h}{\mathfrak{c}^3} u_{\mathfrak{e}}^0 \nabla_x c_{n+1}\Big\} \nonumber\\
			& +P_{\mathfrak{e}} \operatorname{div} b_{n+1}+n_{\mathfrak{e}} u_{\mathfrak{e}} \cdot\Big\{\nabla_x\Big(\frac{h}{\mathfrak{c}^2} u_{\mathfrak{e}}\Big) \cdot b_{n+1}+\nabla_x\Big(\frac{h}{\mathfrak{c}^3} u_{\mathfrak{e}}^0\Big) c_{n+1}\Big\}+b_{n+1} \cdot \nabla_x P_{\mathfrak{e}} \nonumber\\
			&+\operatorname{div} \int_{\mathbb{R}^3} \hat{p} \sqrt{\mathbf{M}_{\mathfrak{c}}}\left\{\mathbf{I}-\mathbf{P}_{\mathfrak{c}}\right\}\Big(\frac{F_{n+1}^{\mathfrak{c}}}{\sqrt{\mathbf{M}_{\mathfrak{c}}}}\Big) d p=0,
	\end{align}}
	
	{\small
	\begin{align}\label{N.9}
		&\frac{n_{\mathfrak{e}} u_{\mathfrak{e}}^0}{\mathfrak{c}}\left\{\frac{h}{\mathfrak{c}^2} u_{\mathfrak{e},j} \partial_{t} a_{n+1}+\frac{1}{\gamma K_2(\gamma)}\Big\{\left(6 K_3(\gamma)+\gamma K_2(\gamma)\right) u_{\mathfrak{e},j}\Big(u_{\mathfrak{e}} \cdot \partial_t b_{n+1}+\frac{u_{\mathfrak{e}}^0}{\mathfrak{c}} \partial_t c_{n+1}\Big)+\mathfrak{c}^2 K_3(\gamma) \partial_t b_{n+1, j}\Big\}\right\} \nonumber\\
		& -\frac{n_{\mathfrak{e}} K_3(\gamma)}{\gamma K_2(\gamma)} u_{\mathfrak{e},j} \partial_t c_{n+1}-\partial_t\Big(\frac{n_{\mathfrak{e}} K_3(\gamma)}{\gamma K_2(\gamma)} u_{\mathfrak{e},j}\Big) c_{n+1}  +\frac{n_{\mathfrak{e}} u_{\mathfrak{e}}^0}{\mathfrak{c}}\left\{\partial_t\Big(\frac{h}{\mathfrak{c}^2} u_{\mathfrak{e},j}\Big) a_{n+1}\right.\nonumber\\
		&\left.+\partial_t\Big(\frac{6 K_3(\gamma)+\gamma K_2(\gamma)}{\gamma K_2(\gamma)} u_{\mathfrak{e},j} u_{\mathfrak{e}}\Big) \cdot b_{n+1}+\partial_t\Big(\frac{6 K_3(\gamma)+\gamma K_2(\gamma)}{\mathfrak{c} \gamma K_2(\gamma)} u_{\mathfrak{e},j} u_{\mathfrak{e}}^0\Big) c_{n+1}+\partial_t\Big(\frac{\mathfrak{c}^2 K_3(\gamma)}{\gamma K_2(\gamma)}\Big) b_{n+1, j}\right\} \nonumber\\
		& +n_{\mathfrak{e}} u_{\mathfrak{e}} \cdot\left\{\frac{h}{\mathfrak{c}^2} u_{\mathfrak{e},j} \nabla_x a_{n+1}+\frac{6 K_3(\gamma)+\gamma K_2(\gamma)}{\gamma K_2(\gamma)} u_{\mathfrak{e},j}\Big\{\left(\nabla_x b_{n+1} \cdot u_{\mathfrak{e}}\right)+\frac{u_{\mathfrak{e}}^0}{\mathfrak{c}} \nabla_x c_{n+1}\Big\}+\frac{\mathfrak{c}^2 K_3(\gamma)}{\gamma K_2(\gamma)} \nabla_x b_{n+1, j}\right\}\nonumber\\
		&+P_{\mathfrak{e}} \partial_{x_j} a_{n+1}+\partial_{x_j} P_{\mathfrak{e}} a_{n+1}  +n_{\mathfrak{e}} u_{\mathfrak{e}} \cdot\left\{\nabla_x\Big(\frac{h}{\mathfrak{c}^2} u_{\mathfrak{e},j}\Big) a_{n+1}+\nabla_x\Big(\frac{6 K_3(\gamma)+\gamma K_2(\gamma)}{\gamma K_2(\gamma)} u_{\mathfrak{e},j} u_{\mathfrak{e}}\Big) \cdot b_{n+1}\right.\nonumber\\
		&+\left.\nabla_x\Big(\frac{6 K_3(\gamma)+\gamma K_2(\gamma)}{\mathfrak{c} \gamma K_2(\gamma)} u_{\mathfrak{e},j} u_{\mathfrak{e}}^0\Big) c_{n+1}+\nabla_x\Big(\frac{\mathfrak{c}^2 K_3(\gamma)}{\gamma K_2(\gamma)}\Big) b_{n+1, j}\right\} 
		+\frac{\mathfrak{c}^2 n_{\mathfrak{e}} K_3(\gamma)}{\gamma K_2(\gamma)} u_{\mathfrak{e},j} \operatorname{div} b_{n+1}\nonumber\\
		&+\nabla_x\Big(\frac{\mathfrak{c}^2 n_{\mathfrak{e}} K_3(\gamma)}{\gamma K_2(\gamma)} u_{\mathfrak{e},j}\Big) \cdot b_{n+1}  +\partial_{x_j}\Big(\frac{\mathfrak{c}^2 n_{\mathfrak{e}} K_3(\gamma)}{\gamma K_2(\gamma)} u_{\mathfrak{e}}\Big) \cdot b_{n+1}+\partial_{x_j}\Big(\frac{\mathfrak{c} n_{\mathfrak{e}} K_3(\gamma)}{\gamma K_2(\gamma)} u_{\mathfrak{e}}^0\Big) c_{n+1}\nonumber\\
		&+\frac{\mathfrak{c}^2 n_{\mathfrak{e}} K_3(\gamma)}{\gamma K_2(\gamma)}\Big\{u_{\mathfrak{e}} \cdot \partial_{x_j} b_{n+1}+\frac{u_{\mathfrak{e}}^0}{\mathfrak{c}} \partial_{x_j} c_{n+1}\Big\} \nonumber \\
		&+E_{0, j}\Big\{\frac{n_{\mathfrak{e}} u_{\mathfrak{e}}^0}{\mathfrak{c}} a_{n+1}+\frac{e_{\mathfrak{e}}+P_{\mathfrak{e}}}{\mathfrak{c}^3} u_{\mathfrak{e}}^0\left(u_{\mathfrak{e}} \cdot b_{n+1}\right)+\frac{e_{\mathfrak{e}}\big(u_{\mathfrak{e}}^0\big)^2+P_{\mathfrak{e}}|u_{\mathfrak{e}}|^2}{\mathfrak{c}^4} c_{n+1}\Big\} \nonumber\\
		& +\Big\{\Big(\frac{n_{\mathfrak{e}} u_{\mathfrak{e}}}{\mathfrak{c}} a_{n+1}+\frac{e_{\mathfrak{e}}+P_{\mathfrak{e}}}{\mathfrak{c}^3} u_{\mathfrak{e}}\left(u_{\mathfrak{e}} \cdot b_{n+1}\right)+\frac{P_{\mathfrak{e}} b_{n+1}}{\mathfrak{c}}+\frac{e_{\mathfrak{e}}+P_{\mathfrak{e}}}{\mathfrak{c}^4} u_{\mathfrak{e}}^0 u_{\mathfrak{e}} c_{n+1}\Big) \times B_0^{\mathfrak{c}}\Big\}_j \nonumber\\
		& +\frac{n_{\mathfrak{e}} u_{\mathfrak{e}}^0}{\mathfrak{c}} E_{n+1, j}^{\mathfrak{c}}+\left(\frac{n_{\mathfrak{e}} u_{\mathfrak{e}}}{\mathfrak{c}} \times B_{n+1}^{\mathfrak{c}}\right)_j \nonumber\\
		& +\sum_{\substack{k+l=n+1 \\
				k, l \geq 1}} E_{k, j}\Big\{\frac{n_{\mathfrak{e}} u_{\mathfrak{e}}^0}{\mathfrak{c}} a_l+\frac{e_{\mathfrak{e}}+P_{\mathfrak{e}}}{\mathfrak{c}^3} u_{\mathfrak{e}}^0\left(u_{\mathfrak{e}} \cdot b_l\right)+\frac{e_{\mathfrak{e}}\big(u_{\mathfrak{e}}^0\big)^2+P_{\mathfrak{e}}|u_{\mathfrak{e}}|^2}{\mathfrak{c}^4} c_l\Big\}\nonumber \\
		& +\sum_{\substack{k+l=n+1 \\
				k, l \geq 1}}\Big\{\Big(\frac{n_{\mathfrak{e}} u_{\mathfrak{e}}}{\mathfrak{c}} a_l+\frac{e_{\mathfrak{e}}+P_{\mathfrak{e}}}{\mathfrak{c}^3} u_{\mathfrak{e}}\left(u_{\mathfrak{e}} \cdot b_l\right)+\frac{P_{\mathfrak{e}}}{\mathfrak{c}} b_l+\frac{e_{\mathfrak{e}}+P_{\mathfrak{e}}}{\mathfrak{c}^4} u_{\mathfrak{e}}^0 u_{\mathfrak{e}} c_l\Big) \times B_k^{\mathfrak{c}}\Big\}_j \nonumber\\
		& +\operatorname{div} \left\{\int_{\mathbb{R}^3} p_j \hat{p} \sqrt{\mathbf{M}_{\mathfrak{c}}}\left\{\mathbf{I}-\mathbf{P}_{\mathfrak{c}}\right\}\Big(\frac{F_{n+1}^{\mathfrak{c}}}{\sqrt{\mathbf{M}_{\mathfrak{c}}}}\Big) d p\right\}+\int_{\mathbb{R}^3}\Big(\frac{p}{p^0} \times B_0^{\mathfrak{c}}\Big)_j \sqrt{\mathbf{M}_{\mathfrak{c}}}\left\{\mathbf{I}-\mathbf{P}_{\mathfrak{c}}\right\}\Big(\frac{F_{n+1}^{\mathfrak{c}}}{\sqrt{\mathbf{M}_{\mathfrak{c}}}}\Big) d p \nonumber\\
		& +\sum_{\substack{k+l=n+1 \\
				k, l \geq 1}} \int_{\mathbb{R}^3}\Big(\frac{p}{p^0} \times B_k^{\mathfrak{c}}\Big)_j \sqrt{\mathbf{M}_{\mathfrak{c}}}\left\{\mathbf{I}-\mathbf{P}_{\mathfrak{c}}\right\}\Big(\frac{F_l^{\mathfrak{c}}}{\sqrt{\mathbf{M}_{\mathfrak{c}}}}\Big) d p=0 ,
	\end{align}}
	and
	{\small
	\begin{align}\label{N.10}
		& \frac{n_{\mathfrak{e}} u_{\mathfrak{e}}^0}{\mathfrak{c}}\Big\{\frac{h u_{\mathfrak{e}}^0}{\mathfrak{c}^3} \partial_t a_{n+1}+\frac{u_{\mathfrak{e}}^0}{\mathfrak{c} \gamma K_2(\gamma)}\left(6 K_3(\gamma)+\gamma K_2(\gamma)\right)\Big(u_{\mathfrak{e}} \cdot \partial_t b_{n+1}+\frac{u_{\mathfrak{e}}^0}{\mathfrak{c}} \partial_t c_{n+1}\Big)\Big\} + \frac{n_{\mathfrak{e}} u_{\mathfrak{e}}^0}{\mathfrak{c}}\Big\{\partial_t\Big(\frac{h u_{\mathfrak{e}}^0}{\mathfrak{c}^3}\Big) a_{n+1}\nonumber\\
		&+\partial_t\Big(\frac{u_{\mathfrak{e}}^0}{\mathfrak{c} \gamma K_2(\gamma)}\left(6 K_3(\gamma)+\gamma K_2(\gamma)\right) u_{\mathfrak{e}}\Big) \cdot b_{n+1}+\partial_t\Big(\frac{u_{\mathfrak{e}}^0}{\mathfrak{c}^2 \gamma K_2(\gamma)}\left(6 K_3(\gamma)+\gamma K_2(\gamma)\right) u_{\mathfrak{e}}^0\Big) c_{n+1}\Big\} \nonumber\\
		& -\frac{P_{\mathfrak{e}}}{\mathfrak{c}^2} \partial_t a_{n+1}-\partial_t\Big(\frac{P_{\mathfrak{e}}}{\mathfrak{c}^2}\Big) a_{n+1}-\frac{n_{\mathfrak{e}} K_3(\gamma)}{ \gamma K_2(\gamma)} u_{\mathfrak{e}} \cdot \partial_t b_{n+1}-\partial_t\Big(\frac{n_{\mathfrak{e}} K_3(\gamma) u_{\mathfrak{e}}}{\gamma K_2(\gamma)}\Big) \cdot b_{n+1}-\frac{3 n_{\mathfrak{e}} K_3(\gamma) u_{\mathfrak{e}}^0}{\mathfrak{c} \gamma K_2(\gamma)} \partial_t c_{n+1} \nonumber\\
		& -\partial_t\Big(\frac{3 n_{\mathfrak{e}} K_3(\gamma) u_{\mathfrak{e}}^0}{\mathfrak{c} \gamma K_2(\gamma)}\Big)c_{n+1}+n_{\mathfrak{e}} u_{\mathfrak{e}} \cdot\Big\{\frac{h u_{\mathfrak{e}}^0}{\mathfrak{c}^3} \nabla_x a_{n+1}+\frac{u_{\mathfrak{e}}^0}{\mathfrak{c} \gamma K_2(\gamma)}\left(6 K_3(\gamma)+\gamma K_2(\gamma)\right) \nonumber\\
		&\Big\{u_{\mathfrak{e}} \cdot \nabla_x b_{n+1}+\frac{u_{\mathfrak{e}}^0}{\mathfrak{c}} \nabla_x c_{n+1}\Big\}\Big\}+n_{\mathfrak{e}} u_{\mathfrak{e}} \cdot\Big\{\nabla_x\Big(\frac{h u_{\mathfrak{e}}^0}{\mathfrak{c}^3}\Big) a_{n+1}+\nabla_x\Big\{\frac{u_{\mathfrak{e}}^0}{\mathfrak{c} \gamma K_2(\gamma)}\left(6 K_3(\gamma)+\gamma K_2(\gamma)\right) u_{\mathfrak{e}}\Big\} \cdot b_{n+1}\nonumber\\
		&+\nabla_x\Big\{\frac{u_{\mathfrak{e}}^0}{\mathfrak{c}^2 \gamma K_2(\gamma)}\left(6 K_3(\gamma)+\gamma K_2(\gamma)\right) u_{\mathfrak{e}}^0\Big\} c_{n+1}\Big\} \nonumber\\
		& +\frac{\mathfrak{c} K_3(\gamma) n_{\mathfrak{e}} u_{\mathfrak{e}}^0}{\gamma K_2(\gamma)} \operatorname{div} b_{n+1}+\nabla_x\Big(\frac{\mathfrak{c} K_3(\gamma) n_{\mathfrak{e}} u_{\mathfrak{e}}^0}{\gamma K_2(\gamma)}\Big) \cdot b_{n+1}-\frac{K_3(\gamma)}{\gamma K_2(\gamma)} n_{\mathfrak{e}} u_{\mathfrak{e}} \cdot \nabla_x c_{n+1}- \operatorname{div}\Big(\frac{K_3(\gamma)}{\gamma K_2(\gamma)}n_{\mathfrak{e}} u_{\mathfrak{e}}\Big) c_{n+1} \nonumber\\
		& +\frac{n_{\mathfrak{e}} u_{\mathfrak{e}}}{\mathfrak{c}^2} \cdot E_{n+1}^{\mathfrak{c}}+\frac{n_{\mathfrak{e}} a_{n+1}}{\mathfrak{c}^2} u_{\mathfrak{e}} \cdot E_0^{\mathfrak{c}}+\frac{e_{\mathfrak{e}}+P_{\mathfrak{e}}}{\mathfrak{c}^4}\left(u_{\mathfrak{e}} \cdot b_{n+1}\right)\left(E_0^{\mathfrak{c}} \cdot u_{\mathfrak{e}}\right)+\frac{P_{\mathfrak{e}} E_0^{\mathfrak{c}} \cdot b_{n+1}}{\mathfrak{c}^2}\nonumber\\
		&+\frac{e_{\mathfrak{e}}+P_{\mathfrak{e}}}{\mathfrak{c}^5} u_{\mathfrak{e}}^0 c_{n+1}\left(u_{\mathfrak{e}} \cdot E_0^{\mathfrak{c}}\right)+\int_{\mathbb{R}^3} \frac{p}{\mathfrak{c} p^0} \cdot E_0^{\mathfrak{c}} \sqrt{\mathbf{M}_{\mathfrak{c}}}\left\{\mathbf{I}-\mathbf{P}_{\mathfrak{c}}\right\}\Big(\frac{F_{n+1}^{\mathfrak{c}}}{\sqrt{\mathbf{M}_{\mathfrak{c}}}}\Big) d p \nonumber\\
		& +\sum_{\substack{k+l=n+1 \\
				k, l \geq 1}}\Big\{\frac{n_{\mathfrak{e}} a_l}{\mathfrak{c}^2} u_{\mathfrak{e}} \cdot E_k^{\mathfrak{c}}+\frac{e_{\mathfrak{e}}+P_{\mathfrak{e}}}{\mathfrak{c}^4}\left(u_{\mathfrak{e}} \cdot b_l\right)\left(E_k^{\mathfrak{c}} \cdot u_{\mathfrak{e}}\right)+\frac{P_{\mathfrak{e}}}{\mathfrak{c}^2} E_k^{\mathfrak{c}} \cdot b_l+\frac{e_{\mathfrak{e}}+P_{\mathfrak{e}}}{\mathfrak{c}^5} u_{\mathfrak{e}}^0 c_l\left(u_{\mathfrak{e}} \cdot E_k^{\mathfrak{c}}\right) \nonumber\\
		& +\int_{\mathbb{R}^3} \frac{p \cdot E_k^{\mathfrak{c}}}{\mathfrak{c} p^0} \sqrt{\mathbf{M}_{\mathfrak{c}}}\left\{\mathbf{I}-\mathbf{P}_{\mathfrak{c}}\right\}\Big(\frac{F_l^{\mathfrak{c}}}{\sqrt{\mathbf{M}_{\mathfrak{c}}}}\Big) d p\Big\}=0.
	\end{align}}
	Despite their cumbersome form, equations \eqref{N.8}-\eqref{N.10} can be reduced to a symmetric hyperbolic system.
   \section{Representation Formulae for the Electromagnetic Fields and Their Derivatives}
\subsection{Proof of Lemma \ref{lW.2}}\label{appendix C.1}
Recall \eqref{W0.1}. Up to the harmless multiplicative constant \(4\pi\), which can be absorbed into the generic constant \(C\) in all subsequent estimates, we derive the following explicit representation formulas for the electromagnetic fields. Using the Poisson formula for the wave equation in \(\mathbb{R}^3\), the solutions to \eqref{W0.1} can be written as follows:
\begin{align}
E_{R,i}^{\varepsilon,\mathfrak{c}}(t,x)
=&\frac{1}{\mathfrak{c}^2}\partial_t\Big(\frac{1}{t}\int_{|x-y|=\mathfrak{c}t}E_{R,i}^{\varepsilon,\mathfrak{c}}(0,y)\,d S_y\Big)
+\frac{1}{\mathfrak{c}^2 t}\int_{|x-y|=\mathfrak{c}t}\partial_t E_{R,i}^{\varepsilon,\mathfrak{c}}(0,y)\,d S_y \nonumber\\
&+\frac{1}{\mathfrak{c}^2}\int_{|x-y|<\mathfrak{c}t}\frac{d y}{|x-y|}\int_{\mathbb{R}^3}
\Big(\mathfrak{c}^2\partial_{y_i}F_R^{\varepsilon,\mathfrak{c}}+\hat{p}_i\partial_t F_R^{\varepsilon,\mathfrak{c}}\Big)
\Big(t-\frac{|x-y|}{\mathfrak{c}},y,p\Big)\,d p, \label{G.1}\\
B_{R,i}^{\varepsilon,\mathfrak{c}}(t,x)
=&\frac{1}{\mathfrak{c}^2}\partial_t\Big(\frac{1}{t}\int_{|x-y|=\mathfrak{c}t}B_{R,i}^{\varepsilon,\mathfrak{c}}(0,y)\,d S_y\Big)
+\frac{1}{\mathfrak{c}^2 t}\int_{|x-y|=\mathfrak{c}t}\partial_tB_{R,i}^{\varepsilon,\mathfrak{c}}(0,y)\,d S_y \nonumber\\
&-\frac{1}{\mathfrak{c}^2}\int_{|x-y|<\mathfrak{c}t}\frac{d y}{|x-y|}
\Big(\nabla_y\times\int_{\mathbb{R}^3}\mathfrak{c}\hat{p}\,F_R^{\varepsilon,\mathfrak{c}}
\Big(t-\frac{|x-y|}{\mathfrak{c}},y,p\Big)\,d p\Big)_i. \label{G.2}
\end{align}

Following \cite{Glassey-ARMA-1986}, we define
\begin{equation}\label{G.3}
\left\{
\begin{aligned}
S&=\partial_t+\hat{p}\cdot\nabla_y,\\
T&=-\frac{\omega}{\mathfrak{c}}\partial_t+\nabla_y,
\end{aligned}
\right.
\end{equation}
and $\omega=\frac{y-x}{|x-y|}$. Moreover,
\begin{equation}\label{G.4}
\left\{
\begin{aligned}
\partial_t&=\Big(1+\frac{\hat{p}}{\mathfrak{c}}\cdot\omega\Big)^{-1}(S-\hat{p}\cdot T),\\
\nabla_y&=T+\frac{\omega}{\mathfrak{c}}\Big(1+\frac{\hat{p}}{\mathfrak{c}}\cdot\omega\Big)^{-1}(S-\hat{p}\cdot T).
\end{aligned}
\right.
\end{equation}
By \eqref{G.4}, the last term on the right-hand side of \eqref{G.1} can be written as
\begin{align}\label{G.5}
&\frac{1}{\mathfrak{c}^2}\iint_{|x-y|<\mathfrak{c}t}\frac{dydp}{|x-y|}
\frac{\mathfrak{c}\omega_i+\hat{p}_i}{1+\frac{\hat{p}}{\mathfrak{c}}\cdot\omega}
\big(SF_R^{\varepsilon,\mathfrak{c}}\big)\Big(t-\frac{|x-y|}{\mathfrak{c}},y,p\Big)\nonumber\\
&+\frac{1}{\mathfrak{c}^2}\iint_{|x-y|<\mathfrak{c}t}\frac{dydp}{|x-y|}
\Big(\mathfrak{c}^2T_iF_R^{\varepsilon,\mathfrak{c}}
-\frac{\mathfrak{c}\omega_i+\hat{p}_i}{1+\frac{\hat{p}}{\mathfrak{c}}\cdot\omega}\,\hat{p}\cdot TF_R^{\varepsilon,\mathfrak{c}}\Big)
\Big(t-\frac{|x-y|}{\mathfrak{c}},y,p\Big).
\end{align}
Noting that
\begin{align*}
\nabla_y\Big[f\Big(t-\frac{|x-y|}{\mathfrak{c}},y,p\Big)\Big]
&=\big(\nabla_y f\big)\Big(t-\frac{|x-y|}{\mathfrak{c}},y,p\Big)
-\frac{\omega}{\mathfrak{c}}\big(\partial_tf\big)\Big(t-\frac{|x-y|}{\mathfrak{c}},y,p\Big)\\
&=(Tf)\Big(t-\frac{|x-y|}{\mathfrak{c}},y,p\Big),
\end{align*}
we may use the divergence theorem to deduce that  
\begin{align}\label{G.6}
\eqref{G.5}_2=&-\iint_{|x-y|<\mathfrak{c}t}
\left\{
\partial_{y_i}\Big(\frac{1}{|x-y|}\Big)
-\nabla_y\cdot\Big(
\frac{\frac{\hat{p}}{\mathfrak{c}}}{|x-y|}
\frac{\omega_i+\frac{\hat{p}_i}{\mathfrak{c}}}{1+\frac{\hat{p}}{\mathfrak{c}}\cdot\omega}
\Big)
\right\}
F_R^{\varepsilon,\mathfrak{c}}\Big(t-\frac{|x-y|}{\mathfrak{c}},y,p\Big)\,dydp \nonumber\\
&+\iint_{|x-y|=\mathfrak{c}t}\frac{dS_y\,dp}{|x-y|}
\left(
\omega_i-\frac{\omega_i+\frac{\hat{p}_i}{\mathfrak{c}}}{1+\frac{\hat{p}}{\mathfrak{c}}\cdot\omega}
\frac{\hat{p}}{\mathfrak{c}}\cdot\omega
\right)
F_R^{\varepsilon,\mathfrak{c}}(0,y,p).
\end{align}

To compute the derivatives appearing in \eqref{G.6}, note first that
\[
\frac{\partial |x-y|}{\partial y_j}=\omega_j,
\qquad
\frac{\partial \omega_i}{\partial y_j}
=\frac{\delta_{ij}-\omega_i\omega_j}{|x-y|}.
\]
Now we compute
\begin{align*}
&\nabla_y\cdot
\left(
\frac{\big(\omega_i+\frac{\hat{p}_i}{\mathfrak{c}}\big)\frac{\hat{p}}{\mathfrak{c}}}
{|x-y|\big(1+\frac{\hat{p}}{\mathfrak{c}}\cdot\omega\big)}
\right)
=\sum_{j=1}^3
\frac{\partial}{\partial y_j}
\left(
\frac{\big(\omega_i+\frac{\hat{p}_i}{\mathfrak{c}}\big)\frac{\hat{p}_j}{\mathfrak{c}}}
{|x-y|\big(1+\frac{\hat{p}}{\mathfrak{c}}\cdot\omega\big)}
\right)\\
&=|x-y|^{-2}\Big(1+\frac{\hat{p}}{\mathfrak{c}}\cdot\omega\Big)^{-2}
\left\{
\Big(1-\frac{|\hat{p}|^2}{\mathfrak{c}^2}\Big)\frac{\hat{p}_i}{\mathfrak{c}}
-2\frac{\hat{p}}{\mathfrak{c}}\cdot\omega\,\omega_i
-\omega_i\Big(\frac{\hat{p}}{\mathfrak{c}}\cdot\omega\Big)^2
-\omega_i\frac{|\hat{p}|^2}{\mathfrak{c}^2}
\right\}.
\end{align*}
Hence,
\begin{align}\label{G.7}
\frac{\partial}{\partial y_i}\Big(\frac{1}{|x-y|}\Big)
-\nabla_y\cdot
\left(
\frac{\big(\omega_i+\frac{\hat{p}_i}{\mathfrak{c}}\big)\frac{\hat{p}}{\mathfrak{c}}}
{|x-y|\big(1+\frac{\hat{p}}{\mathfrak{c}}\cdot\omega\big)}
\right)
=
-|x-y|^{-2}\Big(1+\frac{\hat{p}}{\mathfrak{c}}\cdot\omega\Big)^{-2}
\Big(1-\frac{|\hat{p}|^2}{\mathfrak{c}^2}\Big)
\Big(\omega_i+\frac{\hat{p}_i}{\mathfrak{c}}\Big).
\end{align}
Substituting \eqref{G.7} into \eqref{G.6}, we obtain
\begin{align}\label{G.8}
\eqref{G.5}_2=&
\iint_{|x-y|<\mathfrak{c}t}\frac{d y d p}{|x-y|^2}
\frac{\Big(1-\frac{|\hat{p}|^2}{\mathfrak{c}^2}\Big)\Big(\omega_i+\frac{\hat{p}_i}{\mathfrak{c}}\Big)}
{\Big(1+\frac{\hat{p}}{\mathfrak{c}}\cdot\omega\Big)^2}
F_R^{\varepsilon,\mathfrak{c}}\Big(t-\frac{|x-y|}{\mathfrak{c}},y,p\Big)\nonumber\\
&+\iint_{|x-y|=\mathfrak{c}t}\frac{d S_y\,d p}{|x-y|}
\left(
\omega_i-\frac{\omega_i+\frac{\hat{p}_i}{\mathfrak{c}}}{1+\frac{\hat{p}}{\mathfrak{c}}\cdot\omega}
\frac{\hat{p}}{\mathfrak{c}}\cdot\omega
\right)
F_R^{\varepsilon,\mathfrak{c}}(0,y,p).
\end{align}

Substituting \eqref{G.5} and \eqref{G.8} into \eqref{G.1}, we immediately obtain
\begin{align}\label{G.10}
E_{R,i}^{\varepsilon,\mathfrak{c}}(t,x)
=&\frac{1}{\mathfrak{c}^2}\partial_t\Big(\frac{1}{t}\int_{|x-y|=\mathfrak{c}t}E_{R,i}^{\varepsilon,\mathfrak{c}}(0,y)\,dS_y\Big)
+\frac{1}{\mathfrak{c}^2 t}\int_{|x-y|=\mathfrak{c}t}\partial_tE_{R,i}^{\varepsilon,\mathfrak{c}}(0,y)\,dS_y \nonumber\\
&+\iint_{|x-y|=\mathfrak{c}t}\frac{dS_y\,dp}{|x-y|}
\left(
\omega_i-\frac{\omega_i+\frac{\hat{p}_i}{\mathfrak{c}}}{1+\frac{\hat{p}}{\mathfrak{c}}\cdot\omega}
\frac{\hat{p}}{\mathfrak{c}}\cdot\omega
\right)
F_R^{\varepsilon,\mathfrak{c}}(0,y,p)\nonumber\\
&+\iint_{|x-y|<\mathfrak{c}t}\frac{dydp}{|x-y|^2}
\frac{\Big(1-\frac{|\hat{p}|^2}{\mathfrak{c}^2}\Big)\Big(\omega_i+\frac{\hat{p}_i}{\mathfrak{c}}\Big)}
{\Big(1+\frac{\hat{p}}{\mathfrak{c}}\cdot\omega\Big)^2}
F_R^{\varepsilon,\mathfrak{c}}\Big(t-\frac{|x-y|}{\mathfrak{c}},y,p\Big)\nonumber\\
&+\frac{1}{\mathfrak{c}^2}\iint_{|x-y|<\mathfrak{c}t}\frac{dydp}{|x-y|}
\frac{\mathfrak{c}\omega_i+\hat{p}_i}{1+\frac{\hat{p}}{\mathfrak{c}}\cdot\omega}
\big(SF_R^{\varepsilon,\mathfrak{c}}\big)\Big(t-\frac{|x-y|}{\mathfrak{c}},y,p\Big).
\end{align}
Then \eqref{W.2-1} follows directly from \eqref{G.10}.

The representation for $B_{R,i}^{\varepsilon,\mathfrak{c}}$ can be derived in the same manner as that for $E_{R,i}^{\varepsilon,\mathfrak{c}}$, and we omit the details. This completes the proof of Lemma \ref{lW.2}.
$\Box$

\subsection{Proof of Lemma \ref{lW.6}}\label{Appendix C.2}
Applying $\partial_j:=\partial_{x_j}$ to \eqref{W.2-1}, for any $x\in\mathbb{R}^3$, we obtain
\begin{align}\label{G.11}
\partial_j E_{R,i}^{\varepsilon,\mathfrak{c}}(t,x)
=\partial_j\Big(\mathbf{I}_{E_{R,i}^{\varepsilon,\mathfrak{c}}}(t,x)\Big)
+\partial_j\Big(E_{R_T,i}^{\varepsilon,\mathfrak{c}}(t,x)\Big)
+\partial_j\Big(E_{R_S,i}^{\varepsilon,\mathfrak{c}}(t,x)\Big),
\end{align}
where
\begin{align}\label{G.12}
&\partial_j\Big(E_{R_T,i}^{\varepsilon,\mathfrak{c}}(t,x)\Big)
+\partial_j\Big(E_{R_S,i}^{\varepsilon,\mathfrak{c}}(t,x)\Big)\nonumber\\
=&\frac{1}{\mathfrak{c}^2}\iint_{|x-y|\leq \mathfrak{c}t}\frac{d y d p}{|x-y|^2}
\frac{\big(\omega_i+\frac{\hat{p}_i}{\mathfrak{c}}\big)(\mathfrak{c}^2-|\hat{p}|^2)}
{\big(1+\frac{\hat{p}}{\mathfrak{c}}\cdot\omega\big)^2}
\big(\partial_{y_j}F_R^{\varepsilon,\mathfrak{c}}\big)\Big(t-\frac{|x-y|}{\mathfrak{c}},y,p\Big)\nonumber\\
&+\frac{1}{\mathfrak{c}^2}\iint_{|x-y|\leq \mathfrak{c}t}\frac{dydp}{|x-y|}
\frac{\mathfrak{c}\omega_i+\hat{p}_i}{1+\frac{\hat{p}}{\mathfrak{c}}\cdot\omega}
\Big(\partial_{y_j}(SF_R^{\varepsilon,\mathfrak{c}})\Big)\Big(t-\frac{|x-y|}{\mathfrak{c}},y,p\Big).
\end{align}
This follows by changing variables $z=x-y$ in the integrals defining $E_{R_T,i}^{\varepsilon,\mathfrak{c}}$ and $E_{R_S,i}^{\varepsilon,\mathfrak{c}}$.

Now we perform the same procedure as before: substituting \eqref{G.4} for $\partial_j$ and integrating by parts the two terms involving $T_k$, we obtain
{\small
\begin{align}\label{G.13}
&\frac{1}{\mathfrak{c}^2}\iint_{|x-y|\leq \mathfrak{c}t}
\frac{\big(\omega_i+\frac{\hat{p}_i}{\mathfrak{c}}\big)(\mathfrak{c}^2-|\hat{p}|^2)}
{\big(1+\frac{\hat{p}}{\mathfrak{c}}\cdot\omega\big)^2}
\big(\partial_{y_j}F_R^{\varepsilon,\mathfrak{c}}\big)\Big(t-\frac{|x-y|}{\mathfrak{c}},y,p\Big)\frac{dydp}{|x-y|^2}\nonumber\\
&+\frac{1}{\mathfrak{c}^2}\iint_{|x-y|\leq \mathfrak{c}t}
\frac{\mathfrak{c}\omega_i+\hat{p}_i}{1+\frac{\hat{p}}{\mathfrak{c}}\cdot\omega}
\Big(\partial_{y_j}(SF_R^{\varepsilon,\mathfrak{c}})\Big)\Big(t-\frac{|x-y|}{\mathfrak{c}},y,p\Big)\frac{dydp}{|x-y|}\nonumber\\
=&\frac{1}{\mathfrak{c}^2}\iint_{|x-y|\leq \mathfrak{c}t}
\frac{\big(\omega_i+\frac{\hat{p}_i}{\mathfrak{c}}\big)(\mathfrak{c}^2-|\hat{p}|^2)}
{\big(1+\frac{\hat{p}}{\mathfrak{c}}\cdot\omega\big)^2}
\left\{
\Big(\delta_{jk}-\frac{\omega_j\frac{\hat{p}_k}{\mathfrak{c}}}{1+\frac{\hat{p}}{\mathfrak{c}}\cdot\omega}\Big)T_k
+\frac{\frac{\omega_j}{\mathfrak{c}}}{1+\frac{\hat{p}}{\mathfrak{c}}\cdot\omega}S
\right\}
F_R^{\varepsilon,\mathfrak{c}}\Big(t-\frac{|x-y|}{\mathfrak{c}},y,p\Big)\frac{dydp}{|x-y|^2}\nonumber\\
&+\frac{1}{\mathfrak{c}^2}\iint_{|x-y|\leq \mathfrak{c}t}
\frac{\mathfrak{c}\omega_i+\hat{p}_i}{1+\frac{\hat{p}}{\mathfrak{c}}\cdot\omega}
\left\{
\Big(\delta_{jk}-\frac{\omega_j\frac{\hat{p}_k}{\mathfrak{c}}}{1+\frac{\hat{p}}{\mathfrak{c}}\cdot\omega}\Big)T_k
+\frac{\frac{\omega_j}{\mathfrak{c}}}{1+\frac{\hat{p}}{\mathfrak{c}}\cdot\omega}S
\right\}
(SF_R^{\varepsilon,\mathfrak{c}})\Big(t-\frac{|x-y|}{\mathfrak{c}},y,p\Big)\frac{dydp}{|x-y|}\nonumber\\
=&-\frac{1}{\mathfrak{c}^2}\iint_{|x-y|\leq \mathfrak{c}t}
\partial_{y_k}\left\{
\frac{\big(\omega_i+\frac{\hat{p}_i}{\mathfrak{c}}\big)(\mathfrak{c}^2-|\hat{p}|^2)}
{\big(1+\frac{\hat{p}}{\mathfrak{c}}\cdot\omega\big)^2}
\Big(\delta_{jk}-\frac{\omega_j\frac{\hat{p}_k}{\mathfrak{c}}}{1+\frac{\hat{p}}{\mathfrak{c}}\cdot\omega}\Big)\frac{1}{|x-y|^2}
\right\}
F_R^{\varepsilon,\mathfrak{c}}\Big(t-\frac{|x-y|}{\mathfrak{c}},y,p\Big)\,dydp\nonumber\\
&+\frac{1}{\mathfrak{c}^2}\iint_{|x-y|=\mathfrak{c}t}
\frac{\big(\omega_i+\frac{\hat{p}_i}{\mathfrak{c}}\big)(\mathfrak{c}^2-|\hat{p}|^2)}
{\big(1+\frac{\hat{p}}{\mathfrak{c}}\cdot\omega\big)^2}
\Big(\omega_j-\frac{\omega_j\frac{\hat{p}}{\mathfrak{c}}\cdot\omega}{1+\frac{\hat{p}}{\mathfrak{c}}\cdot\omega}\Big)
F_R^{\varepsilon,\mathfrak{c}}(0,y,p)\frac{dS_y\,dp}{|x-y|^2}\nonumber\\
&+\frac{1}{\mathfrak{c}^2}\iint_{|x-y|\leq \mathfrak{c}t}
\frac{\big(\omega_i+\frac{\hat{p}_i}{\mathfrak{c}}\big)(\mathfrak{c}^2-|\hat{p}|^2)\frac{\omega_j}{\mathfrak{c}}}
{\big(1+\frac{\hat{p}}{\mathfrak{c}}\cdot\omega\big)^3}
SF_R^{\varepsilon,\mathfrak{c}}\Big(t-\frac{|x-y|}{\mathfrak{c}},y,p\Big)\frac{dydp}{|x-y|^2}\nonumber\\
&-\frac{1}{\mathfrak{c}^2}\iint_{|x-y|\leq \mathfrak{c}t}
\partial_{y_k}\left\{
\frac{\mathfrak{c}\omega_i+\hat{p}_i}{1+\frac{\hat{p}}{\mathfrak{c}}\cdot\omega}
\Big(\delta_{jk}-\frac{\omega_j\frac{\hat{p}_k}{\mathfrak{c}}}{1+\frac{\hat{p}}{\mathfrak{c}}\cdot\omega}\Big)\frac{1}{|x-y|}
\right\}
(SF_R^{\varepsilon,\mathfrak{c}})\Big(t-\frac{|x-y|}{\mathfrak{c}},y,p\Big)\,dydp\nonumber\\
&+\frac{1}{\mathfrak{c}^2}\iint_{|x-y|=\mathfrak{c}t}
\frac{\mathfrak{c}\omega_i+\hat{p}_i}{1+\frac{\hat{p}}{\mathfrak{c}}\cdot\omega}
\Big(\omega_j-\frac{\omega_j\frac{\hat{p}}{\mathfrak{c}}\cdot\omega}{1+\frac{\hat{p}}{\mathfrak{c}}\cdot\omega}\Big)
(SF_R^{\varepsilon,\mathfrak{c}})(0,y,p)\frac{dS_y\,dp}{|x-y|}\nonumber\\
&+\frac{1}{\mathfrak{c}^2}\iint_{|x-y|\leq \mathfrak{c}t}
\frac{\big(\omega_i+\frac{\hat{p}_i}{\mathfrak{c}}\big)\omega_j}
{\big(1+\frac{\hat{p}}{\mathfrak{c}}\cdot\omega\big)^2}
(S^2F_R^{\varepsilon,\mathfrak{c}})\Big(t-\frac{|x-y|}{\mathfrak{c}},y,p\Big)\frac{dydp}{|x-y|}.
\end{align}
}
The last term is precisely the one appearing in \eqref{W.11}, and it yields \eqref{W.11-5}.

Since
\begin{align}\label{G.14}
&-\partial_{y_k}\left\{
\frac{\big(\omega_i+\frac{\hat{p}_i}{\mathfrak{c}}\big)(\mathfrak{c}^2-|\hat{p}|^2)}
{\mathfrak{c}^2\big(1+\frac{\hat{p}}{\mathfrak{c}}\cdot\omega\big)^2}
\Big(\delta_{j k}-\frac{\omega_j\frac{\hat{p}_k}{\mathfrak{c}}}{1+\frac{\hat{p}}{\mathfrak{c}}\cdot\omega}\Big)\frac{1}{|x-y|^2}
\right\}\nonumber\\
=&\frac{a_A(\omega,\hat{p})}{|x-y|^3}
\equiv
\frac{
3\big(\omega_i+\frac{\hat{p}_i}{\mathfrak{c}}\big)
\left[\omega_j\Big(1-\frac{|\hat{p}|^2}{\mathfrak{c}^2}\Big)
+\frac{\hat{p}_j}{\mathfrak{c}}\Big(1+\frac{\hat{p}}{\mathfrak{c}}\cdot\omega\Big)\right]
-\Big(1+\frac{\hat{p}}{\mathfrak{c}}\cdot\omega\Big)^2\delta_{i j}
}
{\big(1+\frac{|p|^2}{\mathfrak{c}^2}\big)\big(1+\frac{\hat{p}}{\mathfrak{c}}\cdot\omega\big)^4|x-y|^3},
\end{align}
we obtain the fourth term in \eqref{W.11} and \eqref{W.11-3}. Moreover,
\begin{align*}
&-\partial_{y_k}\left\{
\frac{\omega_i+\frac{\hat{p}_i}{\mathfrak{c}}}{1+\frac{\hat{p}}{\mathfrak{c}}\cdot\omega}
\Big(\delta_{jk}-\frac{\omega_j\frac{\hat{p}_k}{\mathfrak{c}}}{1+\frac{\hat{p}}{\mathfrak{c}}\cdot\omega}\Big)\frac{1}{|x-y|}
\right\}\\
=&-\frac{
\big(\omega_i+\frac{\hat{p}_i}{\mathfrak{c}}\big)
\left[2\omega_j\Big(\frac{|\hat{p}|^2}{\mathfrak{c}^2}-1\Big)
-2\frac{\hat{p}_j}{\mathfrak{c}}\Big(1+\frac{\hat{p}}{\mathfrak{c}}\cdot\omega\Big)\right]
+\Big(1+\frac{\hat{p}}{\mathfrak{c}}\cdot\omega\Big)^2\delta_{ij}
}
{\big(1+\frac{\hat{p}}{\mathfrak{c}}\cdot\omega\big)^3|x-y|^2}.
\end{align*}
If we add the kernel
\[
\frac{\big(\omega_i+\frac{\hat{p}_i}{\mathfrak{c}}\big)\Big(1-\frac{|\hat{p}|^2}{\mathfrak{c}^2}\Big)\omega_j}
{\big(1+\frac{\hat{p}}{\mathfrak{c}}\cdot\omega\big)^3|x-y|^2},
\]
then we obtain
\begin{align}\label{G.15}
\frac{b_A(\omega,\hat{p})}{|x-y|^2}
=
-\frac{
\big(\omega_i+\frac{\hat{p}_i}{\mathfrak{c}}\big)
\left[3\omega_j\Big(\frac{|\hat{p}|^2}{\mathfrak{c}^2}-1\Big)
-2\frac{\hat{p}_j}{\mathfrak{c}}\Big(1+\frac{\hat{p}}{\mathfrak{c}}\cdot\omega\Big)\right]
+\Big(1+\frac{\hat{p}}{\mathfrak{c}}\cdot\omega\Big)^2\delta_{ij}
}
{\big(1+\frac{\hat{p}}{\mathfrak{c}}\cdot\omega\big)^3|x-y|^2},
\end{align}
which is the fifth term in \eqref{W.11}, and clearly \eqref{W.11-4} holds.

Arranging the remaining terms, define
\begin{align*}
\mathbf{I}_{E_{R,ij}^{\varepsilon,\mathfrak{c}}}(t,x)
:=&\partial_j\Big(\mathbf{I}_{E_{R,i}^{\varepsilon,\mathfrak{c}}}(t,x)\Big)+\frac{1}{\mathfrak{c}^2}\iint_{|x-y|=\mathfrak{c}t}
\frac{\mathfrak{c}\omega_i+\hat{p}_i}{1+\frac{\hat{p}}{\mathfrak{c}}\cdot\omega}
\Big(\omega_j-\frac{\omega_j\frac{\hat{p}}{\mathfrak{c}}\cdot\omega}{1+\frac{\hat{p}}{\mathfrak{c}}\cdot\omega}\Big)
(SF_R^{\varepsilon,\mathfrak{c}})(0,y,p)\frac{dS_y\,dp}{|x-y|}\\
&+\frac{1}{\mathfrak{c}^2}\iint_{|x-y|=\mathfrak{c}t}
\frac{\big(\omega_i+\frac{\hat{p}_i}{\mathfrak{c}}\big)(\mathfrak{c}^2-|\hat{p}|^2)}
{\big(1+\frac{\hat{p}}{\mathfrak{c}}\cdot\omega\big)^2}
\Big(\omega_j-\frac{\omega_j\frac{\hat{p}}{\mathfrak{c}}\cdot\omega}{1+\frac{\hat{p}}{\mathfrak{c}}\cdot\omega}\Big)
F_R^{\varepsilon,\mathfrak{c}}(0,y,p)\frac{dS_y\,dp}{|x-y|^2},
\end{align*}
which yields \eqref{W.11-1}. Arranging the remaining terms, we obtain \eqref{W.11}.

The representation for $\partial_jB_{R,i}^{\varepsilon,\mathfrak{c}}$ can be derived in the same manner as that for $\partial_jE_{R,i}^{\varepsilon,\mathfrak{c}}$, and we omit the details.

With the help of \eqref{W.5}, the estimate \eqref{W.12} is immediate. It remains to establish \eqref{W.13}. For this purpose, using the formula for $a_A(\omega,\hat{p})$ in \eqref{W.11-3}, we write
\[
a_A(\omega,\hat{p})=I+II+III,
\]
where
\begin{align*}
I
=&\frac{3\big(\omega_i+\frac{\hat{p}_i}{\mathfrak{c}}\big)\omega_j\big(1-\frac{|\hat{p}|^2}{\mathfrak{c}^2}\big)}
{\big(1+\frac{|p|^2}{\mathfrak{c}^2}\big)\big(1+\frac{\hat{p}}{\mathfrak{c}}\cdot\omega\big)^4}
=\frac{3\omega_j\big(\omega_i+\frac{\hat{p}_i}{\mathfrak{c}}\big)}
{\big(1+\frac{|p|^2}{\mathfrak{c}^2}\big)^2\big(1+\frac{\hat{p}}{\mathfrak{c}}\cdot\omega\big)^4}\\
=&\frac{3\omega_j\big(\omega_i+\frac{\hat{p}_i}{\mathfrak{c}}\big)}
{\Big(\sqrt{1+\frac{|p|^2}{\mathfrak{c}^2}}+\frac{p}{\mathfrak{c}}\cdot\omega\Big)^4}
=-\frac{\partial}{\partial p_i}
\left\{
\frac{(\mathfrak{c}\omega_j+\hat{p}_j)-\hat{p}_j}
{\Big(\sqrt{1+\frac{|p|^2}{\mathfrak{c}^2}}+\frac{p}{\mathfrak{c}}\cdot\omega\Big)^3}
\right\}\\
=&\frac{\partial}{\partial p_i}
\left\{
\frac{\partial}{\partial p_j}
\left[
\frac{\mathfrak{c}^2}{2}
\Big(\sqrt{1+\frac{|p|^2}{\mathfrak{c}^2}}+\frac{p}{\mathfrak{c}}\cdot\omega\Big)^{-2}
\right]
+\hat{p}_j\Big(\sqrt{1+\frac{|p|^2}{\mathfrak{c}^2}}+\frac{p}{\mathfrak{c}}\cdot\omega\Big)^{-3}
\right\},\\
II
=&\frac{3\big(\omega_i+\frac{\hat{p}_i}{\mathfrak{c}}\big)\frac{\hat{p}_j}{\mathfrak{c}}}
{\big(1+\frac{|p|^2}{\mathfrak{c}^2}\big)\big(1+\frac{\hat{p}}{\mathfrak{c}}\cdot\omega\big)^3}
=\frac{3\big(\omega_i+\frac{\hat{p}_i}{\mathfrak{c}}\big)\frac{p_j}{\mathfrak{c}}}
{\Big(\sqrt{1+\frac{|p|^2}{\mathfrak{c}^2}}+\frac{p}{\mathfrak{c}}\cdot\omega\Big)^3}\\
=&-\frac{3}{2}p_j\frac{\partial}{\partial p_i}
\Big(\sqrt{1+\frac{|p|^2}{\mathfrak{c}^2}}+\frac{p}{\mathfrak{c}}\cdot\omega\Big)^{-2},\\
III
=&-\frac{\delta_{ij}}
{\big(1+\frac{|p|^2}{\mathfrak{c}^2}\big)\big(1+\frac{\hat{p}}{\mathfrak{c}}\cdot\omega\big)^2}
=-\frac{\delta_{ij}}
{\Big(\sqrt{1+\frac{|p|^2}{\mathfrak{c}^2}}+\frac{p}{\mathfrak{c}}\cdot\omega\Big)^2}.
\end{align*}

A direct computation gives
\begin{align}\label{G.16}
\int_{|\omega|=1}\Big(\sqrt{1+\frac{|p|^2}{\mathfrak{c}^2}}+\frac{p}{\mathfrak{c}}\cdot\omega\Big)^{-2}d\omega
=&\Big(1+\frac{|p|^2}{\mathfrak{c}^2}\Big)^{-1}
\int_{|\omega|=1}\frac{d\omega}{\big(1+\frac{\hat{p}}{\mathfrak{c}}\cdot\omega\big)^2}\nonumber\\
=&2\pi\Big(1+\frac{|p|^2}{\mathfrak{c}^2}\Big)^{-1}
\int_0^\pi\frac{\sin\phi\,d\phi}{\big(1+\frac{|\hat{p}|}{\mathfrak{c}}\cos\phi\big)^2}
=4\pi,
\end{align}
and similarly
\begin{align}\label{G.17}
\int_{|\omega|=1}\Big(\sqrt{1+\frac{|p|^2}{\mathfrak{c}^2}}+\frac{p}{\mathfrak{c}}\cdot\omega\Big)^{-3}\hat{p}_j\,d\omega
=4\pi p_j.
\end{align}
Thus,
\[
\int I\,d\omega
=\frac{\mathfrak{c}^2}{2}\frac{\partial^2}{\partial p_i\partial p_j}(4\pi)
+\frac{\partial}{\partial p_i}(4\pi p_j)
=4\pi\delta_{ij},
\qquad
\int II\,d\omega
=-\frac{3}{2}p_j\frac{\partial}{\partial p_i}(4\pi)=0,
\]
and
\[
\int III\,d\omega=-4\pi\delta_{ij}.
\]
This proves the first equality in \eqref{W.13}.

We now prove the second equality in \eqref{W.13} for the kernel $a_B(\omega,\hat{p})$. Write
\[
a_B(\omega,\hat{p})=I'+II'+III',
\]
where
\[
I'=\frac{3\omega_j\big(\omega_3\frac{\hat{p}_2}{\mathfrak{c}}-\omega_2\frac{\hat{p}_3}{\mathfrak{c}}\big)}
{\big(1+\frac{|p|^2}{\mathfrak{c}^2}\big)^2\big(1+\frac{\hat{p}}{\mathfrak{c}}\cdot\omega\big)^4},
\quad
II'=\frac{3\frac{\hat{p}_j}{\mathfrak{c}}\big(\omega_3\frac{\hat{p}_2}{\mathfrak{c}}-\omega_2\frac{\hat{p}_3}{\mathfrak{c}}\big)}
{\big(1+\frac{|p|^2}{\mathfrak{c}^2}\big)\big(1+\frac{\hat{p}}{\mathfrak{c}}\cdot\omega\big)^3},
\quad
III'=\frac{\frac{\hat{p}_3}{\mathfrak{c}}\delta_{2j}-\frac{\hat{p}_2}{\mathfrak{c}}\delta_{3j}}
{\big(1+\frac{|p|^2}{\mathfrak{c}^2}\big)\big(1+\frac{\hat{p}}{\mathfrak{c}}\cdot\omega\big)^2}.
\]
From \eqref{G.16}, we already know that
\[
\int III'\,d\omega
=4\pi\left(\frac{\hat{p}_3}{\mathfrak{c}}\delta_{2j}-\frac{\hat{p}_2}{\mathfrak{c}}\delta_{3j}\right).
\]

Next, consider $II'$. We have
\begin{align*}
\frac{1}{3}\int II'\,d\omega
=&\frac{\hat{p}_j}{\mathfrak{c}}p_2
\int\frac{\omega_3+\frac{\hat{p}_3}{\mathfrak{c}}}
{\Big(\sqrt{1+\frac{|p|^2}{\mathfrak{c}^2}}+\frac{p}{\mathfrak{c}}\cdot\omega\Big)^3}\,d\omega
-\frac{\hat{p}_j}{\mathfrak{c}}p_3
\int\frac{\omega_2+\frac{\hat{p}_2}{\mathfrak{c}}}
{\Big(\sqrt{1+\frac{|p|^2}{\mathfrak{c}^2}}+\frac{p}{\mathfrak{c}}\cdot\omega\Big)^3}\,d\omega\\
=&-\frac{\hat{p}_jp_2}{2}\frac{\partial}{\partial p_3}
\int\frac{d\omega}{\Big(\sqrt{1+\frac{|p|^2}{\mathfrak{c}^2}}+\frac{p}{\mathfrak{c}}\cdot\omega\Big)^2}
+\frac{\hat{p}_jp_3}{2}\frac{\partial}{\partial p_2}
\int\frac{d\omega}{\Big(\sqrt{1+\frac{|p|^2}{\mathfrak{c}^2}}+\frac{p}{\mathfrak{c}}\cdot\omega\Big)^2}.
\end{align*}
The integral above is the constant $4\pi$ by \eqref{G.16}. Hence $\int II'\,d\omega=0$.

Finally, the integral of $I'$ can be written as
\[
\int I'\,d\omega
=3\int
\frac{\big(\omega_3\frac{\hat{p}_2}{\mathfrak{c}}-\omega_2\frac{\hat{p}_3}{\mathfrak{c}}\big)\big(\omega_j+\frac{\hat{p}_j}{\mathfrak{c}}\big)}
{\big(1+\frac{|p|^2}{\mathfrak{c}^2}\big)^2\big(1+\frac{\hat{p}}{\mathfrak{c}}\cdot\omega\big)^4}\,d\omega
-3\frac{\hat{p}_j}{\mathfrak{c}}
\int
\frac{\omega_3\frac{\hat{p}_2}{\mathfrak{c}}-\omega_2\frac{\hat{p}_3}{\mathfrak{c}}}
{\big(1+\frac{|p|^2}{\mathfrak{c}^2}\big)^2\big(1+\frac{\hat{p}}{\mathfrak{c}}\cdot\omega\big)^4}\,d\omega.
\]
The first integral equals
\begin{align*}
&-\int\frac{\omega_3p_2-\omega_2p_3}{\sqrt{1+\frac{|p|^2}{\mathfrak{c}^2}}}
\frac{\partial}{\partial p_j}
\Big(\sqrt{1+\frac{|p|^2}{\mathfrak{c}^2}}+\frac{p}{\mathfrak{c}}\cdot\omega\Big)^{-3}\,d\omega\\
=&-\frac{\partial}{\partial p_j}
\left[
\frac{1}{\sqrt{1+\frac{|p|^2}{\mathfrak{c}^2}}}
\int
\frac{\omega_3\hat{p}_2-\omega_2\hat{p}_3}
{\big(1+\frac{|p|^2}{\mathfrak{c}^2}\big)\big(1+\frac{\hat{p}}{\mathfrak{c}}\cdot\omega\big)^3}
\,d\omega
\right]\\
&-\frac{\frac{p_j}{\mathfrak{c}^2}}{\big(1+\frac{|p|^2}{\mathfrak{c}^2}\big)^{3/2}}
\int
\frac{\omega_3\hat{p}_2-\omega_2\hat{p}_3}
{\big(1+\frac{|p|^2}{\mathfrak{c}^2}\big)\big(1+\frac{\hat{p}}{\mathfrak{c}}\cdot\omega\big)^3}
\,d\omega\\
&+\frac{1}{\sqrt{1+\frac{|p|^2}{\mathfrak{c}^2}}}
\int
\frac{\omega_3\delta_{2j}-\omega_2\delta_{3j}}
{\Big(\sqrt{1+\frac{|p|^2}{\mathfrak{c}^2}}+\frac{p}{\mathfrak{c}}\cdot\omega\Big)^3}
\,d\omega.
\end{align*}
The integrals in the first two terms have the same structure as that in $\int II'\,d\omega$; after factoring out the $p$-dependent coefficients, they reduce to the same type of angular integrals, and hence vanish because $\int II'\,d\omega=0$. Therefore,
\begin{align*}
\int I'\,d\omega
=&\frac{1}{\sqrt{1+\frac{|p|^2}{\mathfrak{c}^2}}}
\int
\frac{\omega_3\delta_{2j}-\omega_2\delta_{3j}}
{\Big(\sqrt{1+\frac{|p|^2}{\mathfrak{c}^2}}+\frac{p}{\mathfrak{c}}\cdot\omega\Big)^3}
\,d\omega-3\frac{\hat{p}_j}{\mathfrak{c}}
\int
\frac{\omega_3\frac{\hat{p}_2}{\mathfrak{c}}-\omega_2\frac{\hat{p}_3}{\mathfrak{c}}}
{\big(1+\frac{|p|^2}{\mathfrak{c}^2}\big)^2\big(1+\frac{\hat{p}}{\mathfrak{c}}\cdot\omega\big)^4}
\,d\omega\\
=&\frac{1}{\sqrt{1+\frac{|p|^2}{\mathfrak{c}^2}}}
\int
\frac{
\big(\omega_3+\frac{\hat{p}_3}{\mathfrak{c}}\big)\delta_{2j}
-\big(\omega_2+\frac{\hat{p}_2}{\mathfrak{c}}\big)\delta_{3j}
+\big(\frac{\hat{p}_2}{\mathfrak{c}}\delta_{3j}-\frac{\hat{p}_3}{\mathfrak{c}}\delta_{2j}\big)
}
{\Big(\sqrt{1+\frac{|p|^2}{\mathfrak{c}^2}}+\frac{p}{\mathfrak{c}}\cdot\omega\Big)^3}
\,d\omega\\
&-3\frac{\hat{p}_j}{\mathfrak{c}}
\int
\frac{
\big(\omega_3+\frac{\hat{p}_3}{\mathfrak{c}}\big)\frac{\hat{p}_2}{\mathfrak{c}}
-\big(\omega_2+\frac{\hat{p}_2}{\mathfrak{c}}\big)\frac{\hat{p}_3}{\mathfrak{c}}
}
{\big(1+\frac{|p|^2}{\mathfrak{c}^2}\big)^2\big(1+\frac{\hat{p}}{\mathfrak{c}}\cdot\omega\big)^4}
\,d\omega.
\end{align*}
By the same argument as above, we obtain
\[
\int I'\,d\omega
=4\pi\left(\frac{\hat{p}_2}{\mathfrak{c}}\delta_{3j}-\frac{\hat{p}_3}{\mathfrak{c}}\delta_{2j}\right),
\]
which exactly cancels $\int III'\,d\omega$. This completes the proof of \eqref{W.13}.

Therefore, the proof of Lemma \ref{lW.6} is completed. 
$\Box$

\

    
	\noindent{\bf Acknowledgments.}
    Yong Wang's research is partially supported by the National Key Research and Development Program of China, grant 2021YFA1000800; the National Natural Science Foundation of China, grants 12288201 and 12421001; and the CAS Project for Young Scientists in Basic Research, grant YSBR-031. 
    
\

    \noindent{\bf Conflict of Interest:} The authors declare that they have no conflict of interest.


\begin{thebibliography}{99}

\bibitem{Caflisch-CPAM-1980}
Caflisch, R.~E.: The fluid dynamic limit of the nonlinear Boltzmann equation.
{\it Comm. Pure Appl. Math.} {\bf 33}, 651--666, 1980.

\bibitem{Calogero-JMP-2004}
Calogero, S.: The Newtonian limit of the relativistic Boltzmann equation.
{\it J. Math. Phys.} {\bf 45}, 4042--4052, 2004.

\bibitem{Cao-Ouyang-Wang-Xiao-JFA-2026}
Cao, C., Ouyang, J., Wang, Y., Xiao, C.: The global well-posedness and Newtonian limit for the relativistic Boltzmann equation in a periodic box.
{\it J. Funct. Anal.} {\bf 290}, Art. 111251, 2026.

\bibitem{Cercignani-2002}
Cercignani, C., Kremer, G.~M.: {\it The Relativistic Boltzmann Equation: Theory and Applications}.
Birkh\"{a}user, Basel 2002.

\bibitem{Chapman-CMP-2021}
Chapman, J., Jang, J.~W., Strain, R.~M.: On the determinant problem for the relativistic Boltzmann equation.
{\it Comm. Math. Phys.} {\bf 384}, 1913--1943, 2021.

\bibitem{Yvonne-OMM-2009}
Choquet-Bruhat, Y.: {\it General Relativity and the Einstein Equations}.
Oxford University Press, Oxford 2009.

\bibitem{Dafermos-2016}
Dafermos, C.~M.: {\it Hyperbolic conservation laws in continuum physics}.
Springer, Berlin 2016.

\bibitem{Groot-RKT-1980}
de Groot, S.~R., van Leeuwen, W.~A., van Weert, Ch.~G.: {\it Relativistic Kinetic Theory: Principles and Applications}.
North-Holland, Amsterdam 1980.



\bibitem{Evans-1998}
Evans, L.~C.: {\it Partial differential equations}.
Amer. Math. Soc., Providence, RI 1998.

\bibitem{Friedrichs-CPAM-1954}
Friedrichs, K.~O.: Symmetric hyperbolic linear differential equations.
{\it Comm. Pure Appl. Math.} {\bf 7}, 345--392, 1954.

\bibitem{Glassey-ARMA-1986}
Glassey, R.~T., Strauss, W.~A.: Singularity formation in a collisionless plasma could occur only at high velocities.
{\it Arch. Ration. Mech. Anal.} {\bf 92}, 59--90, 1986.

\bibitem{Glassey-TTSP-1991}
Glassey, R.~T., Strauss, W.~A.: On the derivatives of the collision map of relativistic particles.
{\it Transp. Theory Stat. Phys.} {\bf 20}, 55--68, 1991.

\bibitem{Glassey-PRIMS-1993}
Glassey, R.~T., Strauss, W.~A.: Asymptotic stability of the relativistic Maxwellian.
{\it Publ. Res. Inst. Math. Sci.} {\bf 29}, 301--347, 1993.

\bibitem{Guo-CMP-1998}
Guo, Y.: Smooth irrotational flows in the large to the Euler-Poisson system in $\mathbb{R}^{3+1}$.
{\it Comm. Math. Phys.} {\bf 195}, 249--265, 1998.


\bibitem{Guo-AOM-2016}
Guo, Y., Ionescu, A. D., Pausader, B.: Global solutions of the Euler-Maxwell two-fluid system in 3D.
{\it Ann. of Math. (2)} {\bf 183}, 377--498, 2016.

\bibitem{Guo-KRM-2009}
Guo, Y., Jang, J., Jiang, N.: Local Hilbert expansion for the Boltzmann equation.
{\it Kinet. Relat. Models} {\bf 2}, 205--214, 2009.

\bibitem{Guo-CPAM-2010}
Guo, Y., Jang, J., Jiang, N.: Acoustic limit for the Boltzmann equation in optimal scaling.
{\it Comm. Pure Appl. Math.} {\bf 63}, 337--361, 2010.

\bibitem{Guo-CMP-2010}
Guo, Y., Jang, J.: Global Hilbert expansion for the Vlasov-Poisson-Boltzmann system.
{\it Comm. Math. Phys.} {\bf 299}, 469--501, 2010.

\bibitem{Guo-CMP-2011}
Guo, Y., Pausader, B.: Global smooth ion dynamics in the Euler-Poisson system.
{\it Comm. Math. Phys.} {\bf 303}, 89--125, 2011.

\bibitem{Guo-CMP-2012}
Guo, Y., Strain, R.~M.: Momentum regularity and stability of the relativistic Vlasov-Maxwell-Boltzmann system.
{\it Comm. Math. Phys.} {\bf 310}, 649--673, 2012.

\bibitem{Guo-CMP-2021}
Guo, Y., Xiao, Q.: Global Hilbert expansion for the relativistic Vlasov-Maxwell-Boltzmann system.
{\it Comm. Math. Phys.} {\bf 384}, 341--401, 2021.

\bibitem{Jiang-Lei-Zhao-JFA-2024}
Jiang, N., Lei, Y., Zhao, H.: On the Vlasov--Poisson--Boltzmann limit of the Vlasov--Maxwell--Boltzmann system.
{\it J. Funct. Anal.} {\bf 287}, Art. 110529, 2024.

\bibitem{Jiang-CMP-2025}
Jiang, N., Luo, Y.-L., Tang, S.: Grad-Caflisch pointwise decay estimates revisited.
{\it Comm. Math. Phys.} {\bf 406}, Art. 225, 73 pp., 2025.

\bibitem{Majda-AMS-1984}
Majda, A.~J.: {\it Compressible fluid flow and systems of conservation laws in several space variables}.
Springer, New York 1984.

\bibitem{Nishida-CMP-1978}
Nishida, T.: Fluid dynamical limit of the nonlinear Boltzmann equation to the level of the compressible Euler equation.
{\it Comm. Math. Phys.} {\bf 61}, 119--148, 1978.

\bibitem{Olver-1997}
Olver, F.~W.~J.: {\it Asymptotics and special functions}.
A K Peters, Wellesley, MA 1997.
Reprint of the 1974 original.

\bibitem{Peng-CAMSB-2007}
Peng, Y.-J., Wang, S.: Convergence of compressible Euler-Maxwell equations to compressible EP equations.
{\it Chinese Ann. Math. Ser. B} {\bf 28}, 583--602, 2007.

\bibitem{Speck-CMP-2011}
Speck, J.~R., Strain, R.~M.: Hilbert expansion from the Boltzmann equation to relativistic fluids.
{\it Comm. Math. Phys.} {\bf 304}, 229--280, 2011.

\bibitem{Strain-CMP-2010}
Strain, R.~M.: Asymptotic stability of the relativistic Boltzmann equation for the soft potentials.
{\it Comm. Math. Phys.} {\bf 300}, 529--597, 2010.

\bibitem{Strain-JMA-2010}
Strain, R.~M.: Global Newtonian limit for the relativistic Boltzmann equation near vacuum.
{\it SIAM J. Math. Anal.} {\bf 42}, 1568--1601, 2010.

\bibitem{Strain-KRM-2011}
Strain, R.~M.: Coordinates in the relativistic Boltzmann theory.
{\it Kinet. Relat. Models} {\bf 4}, 345--359, 2011.

\bibitem{Wang-Xiao-JLMS-2026}
Wang, Y. and Xiao, C.: Hydrodynamic limit and Newtonian limit from the
relativistic Boltzmann equation to the classical Euler equations.
{\it J. London Math. Soc.} {\bf 113}, 2026.

\bibitem{Watson-1980}
Watson, G.~N.: {\it A Treatise on the Theory of Bessel Functions}.
2nd ed., Cambridge University Press, Cambridge 1944.

\end{thebibliography}
\end{document}